\documentclass[12pt,a4paper]{report}
\usepackage{graphicx}
\usepackage[italian]{babel}
\usepackage{ae,aecompl}
\usepackage[T1]{fontenc}
\usepackage[utf8]{inputenc}
\usepackage{times}
\usepackage{epsfig}
\usepackage{fancyhdr}
\usepackage{amssymb}
\usepackage{amsmath}
\usepackage{amsfonts}
\usepackage{amsthm}
\usepackage{amscd}
\usepackage{hyperref}
\usepackage[active]{srcltx}

\hypersetup{
  pdftitle={Tesi},
  pdfauthor={Alberto Giulio Setti, Daniele Valtorta},
  pdfsubject={Potenziali di Evans su varietà paraboliche},
  pdfkeywords={Setti 7i evans potential potentials parabolic manifold tesi matematica}
  pdfpagelayout=SinglePage,
  pdfpagemode=UseOutlines}

\newcommand{\N}{\mathbb{N}}
\newcommand{\R}{\mathbb{R}}
\newcommand{\Rp}{\mbox{\footnotesize{$\mathbb{ R}$}}}

\newcommand{\C}{\mathbb{C}}

\newcommand{\norm}[1]{\left\|#1\right\|}
\newcommand{\normop}[1]{\left\|#1\right\|^*}

\newcommand{\ps}[2]{\left\langle#1\middle\vert#2\right\rangle}
\newcommand{\psb}[2]{\left\langle#1\middle\vert#2\right\rangle_B}

\newcommand{\gr}[1]{\left\{#1\right\}}
\newcommand{\ton}[1]{\left(#1\right)}
\newcommand{\abs}[1]{\left|#1\right|}

\newcommand{\Tc}{Tchebycheff}
\newcommand{\A}{\mathcal{A}}
\newcommand{\B}{\mathcal{B}}
\newcommand{\U}{\mathcal{U}}
\newcommand{\M}{\mathcal{M}}
\newcommand{\F}{\mathcal{F}}

\newcommand{\Ro}{\mathbb{M}(R)}
\newcommand{\Roo}{\mathbb{M}_0(R)}
\newcommand{\Rod}{\mathbb{M}_{\Delta}(R)}

\newtheorem{prop}{Proposizione}[chapter]
\newtheorem{teo}[prop]{Teorema}
\newtheorem{deph}[prop]{Definizione}
\newtheorem{lemma}[prop]{Lemma}
\newtheorem{oss}[prop]{Osservazione}

\fancypagestyle{plain}{
  \fancyhead{}
  
}

\begin{document}
\begin{titlepage}

\begin{center}
UNIVERSIT\`A DEGLI STUDI DELL'INSUBRIA
\end{center}

\begin{center}

\ 
Dipartimento di scienze MM. FF. NN. Como \\
Anno accademico 2008 - 2009\\
Sessione di laurea del 30 Settembre 2009\\

\ \newline
{\Huge 

Tesi di Matematica: \\
}
\ \\

\Huge{\textbf{Potenziali di Evans\\ su varietà paraboliche}}

\ \newline
\newline
\includegraphics[width=.5\textwidth]{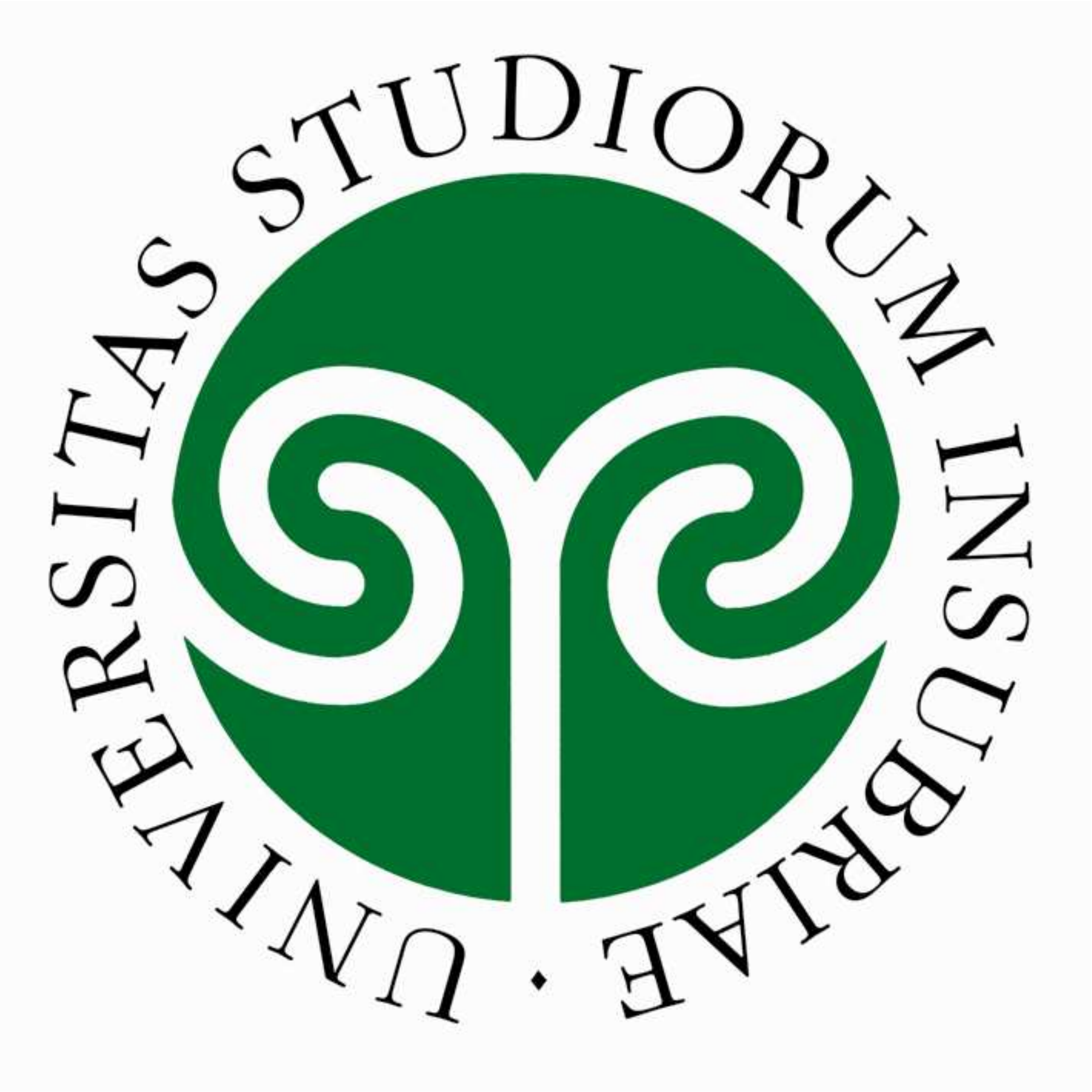}
\end{center}

\begin{flushright}
\ \ \ \ \ \ \ \ \ \ \ \ \ \ \ 
\newline
Autore:\\

Daniele Valtorta, matricola 617528\\
e-mail: \href{mailto:danielevaltorta@gmail.com}{danielevaltorta@gmail.com}
\\

Relatore:\\

Alberto Giulio Setti\\

Co-relatore:\\

Stefano Pigola\\

\end{flushright}
\end{titlepage}
\setcounter{page}{2} 
\section{Riassunto}
La tesi si occupa di una particolare caratterizzazione delle varietà paraboliche, in particolare una varietà Riemanniana $R$ è parabolica se e solo se ammette potenziali di Evans, funzioni armoniche proprie definite fuori da un compatto.\\
Il principale riferimento bibliografico è il libro di Sario e Nakai \cite{5}, libro nel quale gli autori dimostrano l'esistenza di potenziali di Evans su superfici riemanniane. La generalizzazione al caso di varietà di dimensione qualsiasi presenta alcune differenze tecniche non irrilevanti, ad esempio leggermente diverse sono le tecniche utilizzate nel paragrafo \ref{subsec_energy} relativo al principio dell'energia. Si ringrazia il professor Wolfhard Hansen (University of Bielefeld) per i suggerimenti gentilmente forniti riguardo alla teoria del potenziale, che gioca un ruolo fondamentale in queste dimostrazioni.\\
La tesi si articola in una parte introduttiva in cui vengono ripresi alcuni concetti di geometria Riemanniana e di teoria delle funzioni armoniche. Speciale attenzione è stata posta sulla soluzione del problema di Dirichlet su varietà, in particolare è stata trovata una dimostrazione geometrica della solubilità del problema di Dirichlet per domini particolari, dimostrazione riportata nella sezione \ref{sec_dir}. In seguito vengono definite e studiate l'algebra di Royden $\Ro$ e la compattificazione di Royden $R^*$ di una varietà Riemanniana $R$, strumento essenziale per la dimostrazione di esistenza dei potenziali di Evans. Essendo $R^*$ uno spazio topologico non primo numerabile, un capitolo della tesi (il capitolo \ref{chap_filtri}) è dedicato allo studio di filtri e ultrafiltri, concetti utili a descrivere questa topologia e soprattutto a fornire esempi espliciti di elementi in $R^*$ non banali (vedi paragrafo \ref{subsec_charroy}).\\
Importanti strumenti tecnici introdotti nella tesi sono le formule di Green e il principio di Dirichlet, relazioni che sono state dimostrate con ipotesi poco restrittive sulla regolarità degli insiemi e delle funzioni in gioco (vedi paragrafo \ref{sec_dirpri2}). Anche il principio del massimo ha un ruolo essenziale, sia nelle sue forme ``standard'', sia nelle sue versioni più sofisticate presentate nel paragrafo \ref{sec_max2}.\\
Nel capitolo \ref{chap_parab} viene introdotta la capacità di un insieme, e questo concetto è utilizzato per la caratterizzazione delle varietà paraboliche, seguendo il percorso tracciato nella monografia \cite{20}. Sono introdotti anche il diametro transfinito e la costante di \Tc, e viene dimostrato che il diametro transfinito di ogni compatto contenuto nel bordo irregolare di una varietà parabolica è infinito. Questa dimostrazione è basata fortemente sulla teoria del potenziale, teoria che tratta delle proprietà di misure di Borel regolari, della loro energia e delle proprietà dei relativi potenziali di Green. La non finitezza del diametro transfinito è essenziale per la costruzione di funzioni armoniche che tendono a infinito sul bordo irregolare di una varietà iperbolica, e grazie al legame tra queste varietà e le varietà paraboliche esposto nel paragrafo \ref{subset_irreg}, otteniamo la dimostrazione (costruttiva) dell'esistenza dei potenziali di Evans, riportata nel teorema \ref{teo_evans}.\\
In tutte queste dimostrazioni la linearità è essenziale, in particolare si sfrutta il fatto che la somma di due funzioni armoniche è armonica. Un possibile sviluppo futuro di questa tesi è cercare di caratterizzare le varietà $p$-paraboliche con l'esistenza di potenziali di Evans $p$-armonici, che richiederebbe di adattare le tecniche utilizzate in questo lavoro al caso non lineare.
\tableofcontents

\chapter{Richiami di matematica}
\section{Geometria Riemanniana}
Lo scopo di questo capitolo è passare brevemente in rassegna alcune definizioni e risultati che riguardano le varietà Riemanniane. Nel farlo considereremo solo funzioni lisce e varietà $C^{\infty}$.\\
Verranno dati per scontati i concetti di varietà differenziale e di metrica riemanniana, come referenza su questi argomenti consigliamo \cite{7}, \cite{21} e \cite{27}. Ricordiamo che una varietà riemanniana $M$ di dimensione $m$ può essere dotata di una metrica (cioè una funzione $g:T(M)\times T(M)\to \R$ simmetrica definita positiva e liscia), e indicheremo con $g_{ij}$ la matrice definita da $g_{ij}(p)=g(\frac{\partial}{\partial x^i},\frac{\partial}{\partial x^j})$, con $\abs g$ il suo determinante e con $g^{ij}$ la sua inversa.
\subsection{Gradiente e laplaciano di una funzione}
Data una varietà Riemanniana $(M,g)$, è possibile associare a ogni funzione differenziabile $f:M\to \R$ il suo differenziale $df:T(M)\to T(\R)$, che in carte locali assume la forma
\begin{gather*}
 df|_x (v)= \left.\frac{\partial \tilde f}{\partial x^i}\right|_x v^i
\end{gather*}
dove $\tilde f$ è la rappresentazione locale della funzione $f$ rispetto a una qualsiasi carta, e $v^i$ le componenti del vettore $v$ nella stessa carta. Il gradiente della funzione $f$ è il duale del suo differenziale, nel senso che $\nabla f$ è l'unico elemento di $T(M)$ tale che per ogni campo vettoriale $V\in T(M)$:
\begin{gather*}
 g_{ij}(\nabla f)^i v^j=\ps{\nabla f}{ V} = D_V(f)\equiv V(f) = \sum_{i=1}^m v^i\frac{\partial f}{\partial x^i}
\end{gather*}
da questa relazione risulta che il gradiente in coordinate locali assume la forma
\begin{gather}\label{eq_grad}
 (\nabla f)^i=g^{ij}\frac{\partial f}{\partial x^j}
\end{gather}
Oltre al gradiente, per una funzione reale possiamo definire anche l'operatore di Laplace-Beltrami (o laplaciano) $\Delta: C^{m+2}(M,\R)\to C^{m}(M,\R)$. La forma locale di questo operatore (che è l'unico aspetto di interesse in questa tesi) è la seguente:
\begin{gather}\label{eq_lap}
 \Delta f = \frac{1}{\sqrt g} \frac{\partial}{\partial x^i} \ton{g^{ij}\sqrt g \frac{\partial f}{\partial x^j}}
\end{gather}
dove si sottointende la somma degli indici ripetuti (convenzione di Einstein).\\
\`E utile definire anche l'operatore divergenza. Dato un campo vettoriale $V$ liscio, la sua divergenza è una funzione reale e in coordinate locali questa funzione è descritta da:
\begin{gather*}
 div(V)=\frac{1}{\sqrt g} \frac{\partial}{\partial x^i} \ton{\sqrt g V^i}
\end{gather*}
dalla definizione osserviamo subito che
\begin{gather*}
 \Delta f = div(\nabla f)
\end{gather*}
e per un campo vettoriale $V$ e una funzione $f$ qualsiasi vale che:
\begin{gather*}
 div(f V)=\ps{\nabla f}{ V} + f\Delta V
\end{gather*}
A questo punto ha senso parlare di \textbf{funzioni armoniche}, che per definizione sono le funzioni per le quali $\Delta f=0$. \\
Per approfondimenti sul laplaciano consigliamo il libro \cite{21} (nell'esercizio 10 cap. 2.8 pag 57 si trova la maggior parte delle informazioni necessarie per questo lavoro).
\subsection{Punti critici, valori critici e teorema di Sard}\label{sec_sard}
In questa sezione riportiamo brevemente le definizioni di punto critico, valore critico e il teorema di Sard. Per approfondimenti e chiarimenti rimandiamo a \cite{18}. Tutte le definizioni e i teoremi sono riportati con la generalità sufficiente agli scopi di questo lavoro, in modo da non appesantire la tesi con dettagli eccessivi, quindi in particolare tutte le varietà differenziali saranno varietà lisce.\\
La prima domanda che ci poniamo è se le funzioni $f:M\to \R$ possono essere utili per definire in qualche senso una sottovarietà di $M$.
\begin{deph}
Data una varietà differenziale $M$ e una funzione $f\in C^1(M,\R)$, si dice che $p\in M$ è un \textbf{punto critico} se $df|_p=0$ (o equivalentemente $\nabla f |_p=0$). Si dice invece che $x\in \R$ è un \textbf{valore critico} per $f$ se $f^{-1}(x)$ contiene almeno un punto critico.\\
Al contrario si dice che $p\in M$ è un \textbf{valore regolare} per $f$ se non è un punto critico, cioè se $df|_p\neq 0$, mentre $x\in \R$ è un \textbf{valore regolare} per $f$ se $f^{-1}(p)$ non contiene punti critici, cioè contiene solo punti regolari \footnote{nel caso banale $f^{-1}(x)=\emptyset$, si dice che $x$ è un valore regolare per $f$}.\\
Indichiamo $\mathcal{C}(f)\subset \R$ l'insieme dei valori critici di $f:M\to \R$.
\end{deph} 
\begin{deph}
 Data una funzione $f:M\to \R$, per ogni valore $c\in \R$, definiamo \textbf{insieme di sottolivello} di $f$ rispetto a $c$ l'insieme $f^{-1}(-\infty,c]$, \textbf{insieme di sopralivello} di $f$ rispetto a $c$ l'insieme $f^{-1}[c,\infty)$ e \textbf{insieme di livello} di $f$ rispetto a $c$ l'insieme $f^{-1}(c)$.\\
Osserviamo subito che il bordo degli insiemi di sotto/sopra livello è contenuto nel relativo insieme di livello, cioè $\partial (f^{-1}(-\infty,c])\subset f^{-1}(c)$, e che se $c$ è un valore regolare per $f$, questa inclusione si trasforma in uguaglianza, cioè $\partial (f^{-1}(-\infty,c])= f^{-1}(c)$.
\end{deph}
Ha senso chiedersi quando gli insiemi di sottolivello e di livello sono sottovarietà di $M$.
\begin{prop}
 Se $c\in \R$ è un valore regolare di $f\in C^{r}(M,\R)$ \footnote{$r\geq 1$}, allora $f^{-1}(c)$ è una sottovarietà $m-1$ dimensionale di $M$ di regolarità almeno $C^r$, e quindi $f^{-1}(-\infty,c]$ è una sottovarietà con bordo.
\end{prop}
Questa proposizione corrisponde al teorema 3.2 cap. 1.3 pag 22 di \cite{18}, e rimandiamo a questo libro per la dimostrazione.\\
Ora possiamo chiederci ``quanti'' siano i valori regolari di una funzione indipendentemente dalla funzione data. La risposta è contenuta nel \textit{teorema di Sard}. Prima di enunciarlo definiamo gli insiemi di misura nulla su una varietà. Consigliamo il paragrafo 3.1 pag 68 di \cite{18} per approfondimenti.
\begin{deph}
 Un insieme $A\subset \R^m$ si dice essere di \textbf{misura nulla} se e solo se per ogni $\epsilon>0$ esiste un'insieme al più numerabile di cubi $C_n$ in $\R^m$ \footnote{cioè di insiemi della forma $\prod_{i=1}^m [a_i,b_i]$ che hanno misura $\lambda(\prod_{i=1}^m [a_i,b_i])=\prod_{i=1}^m (b_i-a_i)$} tali che
\begin{enumerate}
 \item $A\subset \bigcup_n C_n$
 \item $\sum_{n=1}^{\infty} \lambda (C_n)<\epsilon$
\end{enumerate}
Ricordiamo che questa definizione è equivalente alla richiesta che la \textbf{misura di Lebesgue} di $A$ sia nulla.
\end{deph}
Valgono le seguenti proprietà
\begin{prop}
 Per gli insiemi di misura nulla vale che
\begin{enumerate}
 \item L'unione numerabile di insiemi di misura nulla ha misura nulla
 \item Gli insiemi aperti non vuoti non hanno misura nulla, e il complementare di insiemi di misura nulla è denso
 \item I sottoinsiemi di insiemi di misura nulla hanno misura nulla
 \item L'immagine attraverso una funzione localmente Lipschitziana di un insieme di misura nulla ha misura nulla
\end{enumerate}
\end{prop}
Grazie alla proprietà (1) ha senso definire
\begin{deph}
 Un insieme $A\subset M$ si dice essere di \textbf{misura nulla} se per ogni carta locale $(U,\phi)$ di $M$, l'insieme $\phi(A)$ ha misura nulla. Osserviamo che per paracompattezza di $M$, esiste sempre un atlante al più numerabile $\{(U_n,\phi_n)\}$ di $M$ \footnote{anzi se la dimensione di $M$ è $m$, esiste sempre un atlante formato da al più $m+1$ carte, vedi problema 2.8 pag 21 di \cite{19}}, quindi dato che $A=\cup_n (A\cap U_n)$, ha senso la definizione di insieme di misura nulla e non dipende dall'atlante scelto.
\end{deph}
Il teorema di Sard (o meglio una sua versione non molto generale) garantisce che
\begin{teo}[Teorema di Sard]
 Data $M$ varietà differenziale $m-$dimensionale, sia $r\geq \max\{1,m-1\}$. Allora se $f\in C^r(M,\R)$, l'insieme dei valori critici $\mathcal C (f)$ ha misura nulla in $\R$, quindi il suo complementare è denso.
\end{teo}

\subsection{Coordinate polari geodetiche e modelli}\label{sec_polar}
Un sistema di coordinate che utilizzeremo spesso sulle varietà sono le \textit{coordinate polari geodetiche}. In questa sezione ci occupiamo di dare una breve carrellata sulle coordinate polari e sulle varietà modello, ovvero varietà sfericamente simmetriche. Come referenze per questa sezione consigliamo il capitolo 3 di \cite{7} e il capitolo 3 di \cite{20}. In tutta la sezione lavoreremo con una varietà riemannana $R$ di dimensione $m$.\\
Fissato un punto $p\in R$, esiste sempre un intorno normale $U$ di $p$, un intorno cioè dove la mappa esponenziale $\exp:TR\to R$ è un diffeomorfismo. Visto che $TR$ è isomorfo a $\R^m$, possiamo definire su $TR\setminus \{0\}$ le classiche coordinate polari $m$-dimensionali \footnote{ogni punto può essere individuato dalla distanza dall'origine e da una coordinata su $S^{m-1}$ che rappresenta la direzione del punto} e grazie alla mappa esponenziale possiamo portare queste coordinate sulla varietà. Sull'aperto $U\setminus \{p\}$ chiamiamo le coordinate definite da
\begin{gather*}
 \phi(q)=\{r(\exp^{-1}(q)),\theta(\exp^{-1}(q))\}=(r,\theta)
\end{gather*}
\textbf{coordinate polari geodetiche}. In questo sistema di coordinate la metrica assume la forma:
\begin{gather}\label{eq_gpol}
 g=\begin{bmatrix} 1 & 0 \\ 0 & A(q) \end{bmatrix}
\end{gather}
oppure utilizzando una notazione più famigliare alla geometria riemanniana:
\begin{gather*}
 ds^2 = dr^2 + A(r,\theta)_{ij} d\theta ^i d \theta ^j
\end{gather*}
dove $A(q)$ è una matrice definita positiva \footnote{vedi equazione 3.1 di \cite{20} oppure pagina 136 di \cite{21}}. In coordinate polari inoltre l'operatore laplaciano assume la forma:
\begin{gather}\label{eq_lappol}
 \Delta (f)  = \frac{1}{\sqrt{\abs{A}}} \frac{\partial}{\partial r} \ton{\sqrt{\abs A} \frac{\partial f}{\partial r}} + \Delta_S (f) = \frac{\partial^2 f}{\partial r^2} + \frac{1}{2}\frac{\partial \log(\abs{A})}{\partial r} \frac{\partial f}{\partial r} + \Delta_S (f)
\end{gather}
Dove $\Delta_S$ indica il laplaciano sulla sottovarietà $r=$costante \footnote{vedi equazione 3.4 di \cite{20}}.\\
Si nota quindi che per le funzioni ``radiali'', ossia quelle funzioni che dipendono solo da $r$, il laplaciano assume una forma abbastanza semplificata. Grazie alla semplice forma di $g$ in coordinate polari anche il gradiente di una funzione radiale è molto semplice da calcolare, infatti grazie alla \ref{eq_grad} notiamo che per le funzioni radiali $f$:
\begin{gather}\label{eq_gradpol}
 (\nabla f)^1 = \frac{\partial f}{\partial r}
\end{gather}
mentre tutte le altre componenti del gradiente sono nulle.\\
Queste coordinate risultano particolarmente adatte per descrivere le \textit{varietà modello}. Queste varietà sono varietà sfericamente simmetriche rispetto a rotazioni attorno a un punto fisso.
\begin{deph}
 Definiamo $R$ una varietà con polo $o$, se la mappa esponenziale $\exp|_o$ è un diffeomorfismo globale tra $\R^m$ e $R$. Inoltre se la metrica di questa varietà rispetto alle coordinate polari geodetiche in $o$ è della forma
\begin{gather*}
 ds^2=dr^2 + \sigma (r)^2 d\theta ^2
\end{gather*}
dove $d\theta^2$ è la metrica euclidea standard di $S^{m-1}$, definiamo questa una \textbf{varietà modello}
\end{deph}
La definizione trova la sua giustificazione in questa osservazione. Se consideriamo una rotazione $\rho:S^{m-1}\to S^{m-1}$, possiamo definire una funzione \begin{gather*}
\psi_{\rho}:R\to R \ \ \psi_{\rho}(r,\theta)=\psi_{\rho}(r,\rho(\theta))
\end{gather*}
questa funzione (che possiamo considerare a tutti gli effetti una rotazione su $R$), ha la proprietà di tenere la metrica invariata, cioè
\begin{gather*}
 \psi_{\rho}^*(g)=g
\end{gather*}
dove $\psi_{\rho}^*$ è l'operatore di push-forward, definito da:
\begin{gather*}
\psi_{\rho}^*(g) (v_1,v_2)=g(d \psi_{\rho} (v_1), d \psi_{\rho} (v_2))
\end{gather*}
\subsection{Integrazione su varietà e formula di Green}\label{sec_int}
In questa sezione riportiamo brevemente la definizione di integrazione su varietà riemanniane e la formula di Green per domini regolari. Rimandiamo al testo \cite{9} per approfondimenti.\\
Ai nostri scopi interessa solo ricordare brevemente la definizione di integrale di una funzione su una varietà e di una 1-forma su un bordo. Data una funzione reale $f:R\to \R$, per definire il suo integrale utilizziamo l'integrale di Lebesgue su $\R^m$. Non introduciamo la teoria dell'operatore duale di Hodge, ma ne utilizziamo comunque il simbolo ($\ast$) definendolo quando necessario.
\begin{deph}
 Sia $f:R\to \R$ una funzione continua a supporto compatto $supp(f)\Subset U$ con $(U,\phi)$ carta locale. Allora definiamo
\begin{gather*}
 \int_R f dV = \int_{\phi(U)} \tilde f(x) \sqrt{\abs g} dx^1\cdots dx^m=\int_{\phi(U)} f\phi^{-1}(x) \sqrt{\abs g} dx^1\cdots dx^m
\end{gather*}
Osserviamo che la presenza di $\sqrt{\abs g}$ garantisce che la definizione non dipenda dalla carta locale scelta.
Data una funzione $f:R\to \R$ senza ulteriori condizioni sul supporto, definiamo:
\begin{gather*}
 \int_R f dV = \sum_{n=1}^{\infty} \int_{R} f\lambda_n dV
\end{gather*}
dove $\{\lambda_n\}$ è una partizione dell'unità di $R$ subordinata a un ricoprimento di carte locali qualsiasi.
\end{deph}
\begin{deph}\label{deph_hodge}
Date due funzioni $f,h:R\to \R$, e un dominio regolare $\Omega\subset R$ \footnote{con regolare si intende un dominio aperto relativamente compatto con $\partial \Omega$ bordo liscio (almeno a pezzi)}, definiamo
\begin{gather*}
 \int_{\partial \Omega} h\ast df =\int_{\phi(\partial \Omega)} h(x)g^{im}(x)\frac{\partial f}{\partial x^i}(x) \sqrt{\abs g} dx^1\cdots dx^{m-1}
\end{gather*}
 quando $f$ ha supporto contenuto in una carta $(U,\phi)$ regolare per $\partial \Omega$, cioè una carta in cui $\phi(\partial\Omega)\subset \{x^m=0\}$. Se $f$ non ha queste caratteristiche, l'integrale si ottiene come sopra grazie alle partizioni dell'unità.
\end{deph}
Osserviamo che la notazione $\ast df$ indica il \textit{duale di Hodge} della forma $df$. In questa tesi però questo operatore verrà usato solo in casi simili a quello appena descritto, quindi per brevità tralasciamo la sua definizione e lo utilizziamo solo per comodità di scrittura.\\
Grazie a quanto appena descritto possiamo enunciare la prima identità di Green per varietà Riemanniane:
\begin{prop}[Prima identità di Green]\label{prop_green1}
 Dato un dominio regolare $\Omega\subset R$, e due funzioni lisce $f,h:R\to \R$, si ha che:
\begin{gather*}
 \int_{\Omega} \ps{\nabla f}{\nabla h} dV +\int_{\Omega} f \Delta h dV =\int_{\partial \Omega} f \ast dh
\end{gather*}
\end{prop}
In seguito rilasseremo le ipotesi sulla regolarità di $f$ in questa identità. La dimostrazione di questa identità è una facile conseguenza del teorema di Stokes (consigliamo come referenza il capitolo 7 di \cite{9}).

\subsection{Coordinate di Fermi e applicazioni}
In questo paragrafo descriveremo un sistema di coordinate locali intorno a una sottovarietà regolare particolarmente utile in quanto per molti aspetti molto simile alle coordinate euclidee di $\R^n$. Il riferimento principale in questa sezione è il libro \cite{27}, in particolare il II capitolo, dal quale riporteremo alcuni risultati senza dimostrazione.\\
Data una varietà riemanniana $R$ di dimensione $m$ e una sua sottovarietà regolare $S$ di dimensione $s$, cerchiamo un sistema di coordinate locali $(x_1,\cdots,x_m)$ in un intorno $U$ di $x_0\in S$ tali che ogni punto $p\in U\cap S$ abbia $x^{m-s+1}(p)=\cdots=x^m(p)=0$ e tali che la metrica in queste coordinate assuma la forma:
\begin{gather*}
 g(p)=\begin{bmatrix} A(p) & 0 \\ 0 & B(p) \end{bmatrix}
\end{gather*}
dove $A$ è una matrice $(m-s)\times(m-s)$ e $B$ è una matrice $s\times s$. In realtà per gli scopi di questo lavoro siamo interessati solo al caso $s=m-1$, quindi per comodità di notazione tratteremo solo questo caso, anche se i risultati di questo paragrafo possono facilmente essere estesi a sottovarietà di dimensioni qualsiasi.\\
Per comodità, indicheremo con $p, \ q$ i punti sulla varietà riemanniana, con $$(x^1,\cdots,x^{m-1})\equiv \vec x$$ le prime $m-1$ funzioni coordinate sulla varietà e con $y$ l'ultima coordinata.
\begin{prop}\label{prop_coord_fermi}
 Data una sottovarietà $S$ di codimensione $1$ in $R$, per ogni punto $p_0\in S$, esiste un suo intorno $U$ e delle coordinate locali tali che:
\begin{enumerate}
 \item $y(U\cap S)=0$
 \item la metrica assume la forma
\begin{gather*}
 g(p)=\begin{bmatrix} A(p) & 0 \\ 0 & 1 \end{bmatrix}
\end{gather*}
dove $A$ è una matrice $(m-1)\times (m-1)$. Inoltre è possibile scegliere le coordinate in modo che su $S$ le funzioni $x^1\cdots x^{m-1}$ siano coordinate qualsiasi relative alla sottovarietà $S$. 
\end{enumerate}
\begin{proof}
Il sistema di coordinate cercato prende il nome di \textit{coordinate di Fermi} per la sottovarietà $S$. L'esistenza di queste coordinate è dimostrata in \cite{27}. In particolare a pagina 17 si trova la definizione di coordinate di Fermi, e grazie al lemma 2.3 di pagina 18 e al corollario 2.14 di pagina 31 \footnote{il \textit{lemma di Gauss generalizzato}}, si dimostra la proprietà (2).
\end{proof}
\end{prop}
Introduciamo ora gli intorni geodeticamente convessi su una varietà riemanniana.
\begin{deph}
 Un insieme aperto $U\subset R$ si dice \textbf{geodeticamente convesso} se e solo se per ogni coppia di punti $(p,q)\subset U\times U$ la geodetica che minimizza la distanza tra questi due punti esiste unica e la sua traccia è contenuta in $U$.
\end{deph}
L'esempio più classico di intorno geodeticamente convesso sono le bolle in $\R^n$. L'esistenza di intorni geodeticamente convessi è garantita da questa proposizione, di cui non riportiamo la dimostrazione.
\begin{prop}\label{prop_geoconv}
 Dato un aperto $W\in R$ e $x\in W$, esiste sempre un intorno aperto geodeticamente convesso $U$ tale che $x\in U\subset W$.
\begin{proof}
 La dimostrazione di questa proposizione può essere trovata su \cite{29}, proposizione 10.5.4 pagina 334.
\end{proof}
\end{prop}
Ora siamo pronti per dimostrare un lemma molto tecnico che servirà in seguito per dimostrare una proprietà dei potenziali di Green.
\begin{prop}\label{prop_proj1}
 Dato $K$ insieme aperto relativamente compatto in $R$ con bordo liscio, definiamo per ogni $q\in R$:
\begin{gather*}
 \pi(q) \ \ \ t.c. \ \ \ d(\pi(q),q)=\inf\{d(q,p) \ t.c. \ p\in K\}
\end{gather*}
Dimostriamo che per ogni $p_0\in \partial K$, esiste un intorno $V$ di $p_0$ dove $\pi$ è una funzione continua ben definita, e per ogni $\epsilon>0$, esiste un intorno $U_{\epsilon}$ di $p_0$ per il quale per ogni $p\in K$ e per ogni $q\in U_{\epsilon}$:
\begin{gather*}
 d(p,q)\leq (1+\epsilon)d(p,\pi(q))
\end{gather*}
\begin{proof}
 Consideriamo un intorno $V_1$ di $p_0$ dove siano definite le coordinate di Fermi per la sottovarietà regolare $\partial K$, e in particolare consideriamo un sistema di coordinate per cui la metrica nel punto $p_0$ assume la forma euclidea standard, cioè
\begin{gather*}
 g_{ij}(p_0)=\delta_{ij}
\end{gather*}
Definiamo $K_1=K\cap \overline V_1$ e $K_2\equiv K\cap V_1^C$. L'insieme definito da
\begin{gather*}
 A=\{p\in R \ t.c. \ d(p,K_1)<d(p,K_2)
\end{gather*}
è un'insieme aperto per la continuità della funzione distanza, non vuoto perché $p_0\in A$. Se $p\in A$, sicuramente il punto (o i punti) $\pi(p)$ sono da cercare sono nell'insieme $K_1$. Consideriamo un cilindro $V$ rispetto all'ultima coordinata \footnote{quindi un insieme della forma $B\times (-\epsilon,\epsilon)$ con $B $ aperto contenente la proiezione $p_0$} contenuto nell'aperto $A\cap V_1$. Data la particolare forma delle coordinate di Fermi, è facile dimostrare che se $q\in V$, allora il punto $\pi(q)$ è il punto di coordinate
\begin{gather*}
 y(\pi(q))=0, \ \ x^i(\pi(q))=x^i(q) \ \ \forall i=1\cdots m-1
\end{gather*}
Infatti supponiamo per assurdo che $y(\pi(q))\neq 0$ \footnote{dove l'indice $1$ può essere sostituito da qualsiasi altro indice}, e sia $\gamma$ la geodetica che unisce $\gamma(0)=\pi(q)$ e $\gamma(1)=q$. La curva:
\begin{gather*}
 \tilde\gamma(t)=(x^1(q),\cdots,x^{m-1}(q), y(\gamma(t)))
\end{gather*}
è una curva che unisce i punti $(x^1(q),\cdots,m^{m-1}(q),0)$ con $q$, e la sua lunghezza è
\begin{gather*}
 L(\tilde\gamma)=\int_{0}^1 g_{ij} \dot{\tilde\gamma}^i(t),\dot{\tilde\gamma}^j(t) dt=\int_0^1 (\dot \gamma^m (t))^2 dt< \int_{0}^1 g_{ij} \dot\gamma^i(t),\dot\gamma^j(t) dt=L(\gamma)
\end{gather*}
Per dimostrare la seconda parte della proposizione, ragioniamo sulle prime $m-1$ righe e colonne della matrice $g_{ij}$, cioè sulla matrice $A_{ij}$ per $i,j=1\cdots m-1$. Possiamo sempre scrivere che per ogni punto di $V$:
\begin{gather*}
 A_{ij}(p)=A_{ij}(p_0)+\alpha_{ij}(p)=\delta_{ij}+\alpha_{ij}(p)
\end{gather*}
Per continuità della metrica, per ogni $\epsilon>0$, esiste un intorno $U'_{\epsilon}(p_0)$ contenuto in $V$ per il quale il modulo di tutti gli autovalori di $\alpha_{ij}(p)$ è minore di $\epsilon$. Sia $U''_{\epsilon}(p_0)$ un cilindro rispetto all'ultima coordinata della varietà contenuto in $U'_{\epsilon}$, e $U_{\epsilon}$ un intorno di $p_0$ geodeticamente convesso contenuto in $U''_{\epsilon}$. Siano $p\in U_{\epsilon}\cap K$ e $q\in U_{\epsilon}$. Se $q\in K$, non c'è niente da dimostrare in quanto $q=\pi(q)$, quindi consideriamo solo il caso $q\not \in K$.\\
Sia $\gamma$ una geodetica normalizzata (cioè con $g_{ij}(\gamma(t))\dot \gamma^i (t)\dot \gamma^j(t) =1$) che congiunge $p$ a $q$. Questa geodetica è necessariamente contenuta in $U_{\epsilon}$. Se chiamiamo $l=d(p,q)$, sappiamo che
\begin{gather*}
 l=\int_{0}^l \sqrt{g_{ij}|_{\gamma(t)}\dot \gamma^i(t)\dot \gamma^j(t)} dt=\int_{0}^l \sqrt{A_{ij}|_{\gamma(t)}\dot \gamma^i(t)\dot \gamma^j(t) + \abs{\dot \gamma^{m}(t)}^2}dt 
\end{gather*}
Sia $\eta(t)$ una curva che unisce i punti $p$ e $\pi(q)$ descritta da:
\begin{gather*}
 \eta(t)=\ton{\gamma^1(t),\cdots,\gamma^{m-1}(t), \frac{y(p)}{y(p)-y(q)}(\gamma^m(t)-y(q))}
\end{gather*}
questa curva è contenuta nell'insieme $U''_{\epsilon}$ e anche se non è necessariamente una geodetica, vale che:
\begin{gather*}
 d(\pi(q),p)\leq \int_0^l \sqrt{g_{ij}|_{\eta(t)}\dot \eta^i(t)\dot \eta^j(t)}dt=\\
=\int_0^l \sqrt{A_{ij}|_{\eta(t)}\dot \gamma^i(t)\dot \gamma^j(t) + \kappa^2 \abs{\dot \gamma(t)}^2}dt\leq  \int_0^l \sqrt{A_{ij}|_{\eta(t)}\dot \gamma^i(t)\dot \gamma^j(t) + \abs{\dot \gamma(t)}^2}dt
\end{gather*}
dove $\kappa=\abs{\frac{y(p)}{y(p)-y(q)}}<1$ grazie al fatto che $p\in K$, quindi $y(p)\leq 0$, e $q\not \in K$, quindi $y(q)>0$.\\
Consideriamo inoltre che, poch\`e $\gamma$ \`e normalizzata:
\begin{gather*}
 A_{ij}|_{\gamma(t)}\dot \gamma^i(t)\dot\gamma^j(t)\leq 1
\end{gather*}
Se indichiamo $$\norm{\dot\gamma(t)}_{m-1}=\sum_{i,j\leq m-1}\delta_{ij}\dot\gamma^i(t)\dot\gamma^j(t)$$ ricaviamo che per ogni $t$:
\begin{gather*}
 \delta_{ij}\dot\gamma^i(t)\dot\gamma^i(t)+\alpha|_{\gamma(t)}\dot\gamma^i(t)\dot\gamma^j(t)\leq 1\\ \ \\
\norm{\dot\gamma(t)}_{m-1}\leq 1-\alpha|_{\gamma(t)}\dot\gamma^i(t)\dot\gamma^j(t)\leq 1+\epsilon \norm{\dot\gamma(t)}_{m-1}\\ \ \\
\norm{\dot\gamma(t)}_{m-1}\leq \frac{1}{1-\epsilon}
\end{gather*}
dove abbiamo utilizzato il fatto che il modulo di tutti gli autovalori di tutte le matrici $\alpha_{ij}$ in $U'_{\epsilon}$ è minore di $\epsilon$. Osserviamo che questa stima è indipendente dalla scelta di $p$ e $q$.\\
Grazie alla definizione delle matrici $\alpha_{ij}$ abbiamo anche che:
\begin{gather*}
 A_{ij}|_{\eta(t)}\dot \gamma^i(t)\dot \gamma^j(t)=[A_{ij}|_{\gamma(t)}-\alpha_{ij}|_{\gamma(t)}+\alpha_{ij}|_{\eta(t)}]\dot \gamma^i(t)\dot \gamma^j(t)
\end{gather*}
quindi:
\begin{gather*}
 [A_{ij}|_{\eta(t)}-A_{ij}|_{\gamma(t)}]\dot \gamma^i(t)\dot \gamma^j(t)=[-\alpha_{ij}|_{\gamma(t)}+\alpha_{ij}|_{\eta(t)}]\dot \gamma^i(t)\dot \gamma^j(t) \\ \ \\
 \abs{[A_{ij}|_{\eta(t)}-A_{ij}|_{\gamma(t)}]\dot \gamma^i(t)\dot \gamma^j(t)}\leq 2 \epsilon \norm{\dot\gamma(t)}_{m-1}\leq \frac{2\epsilon}{1-\epsilon}
\end{gather*}
Grazie a quest'ultima disuguaglianza e ricordando che $\sqrt{a+b}\leq \sqrt{a}+\sqrt{\abs b}$ per ogni coppia di numeri reali tali che $a\geq 0, \ a+b\geq 0$, possiamo concludere che:
\begin{gather*}
 d(\pi(q),p)\leq \int_0^l \sqrt{A_{ij}|_{\eta(t)}\dot \gamma^i(t)\dot \gamma^j(t) + \abs{\dot \gamma(t)}^2}dt =\\
=\int_0^l\sqrt{A_{ij}|_{\gamma(t)}\dot \gamma^i(t)\dot \gamma^j(t) + \abs{\dot \gamma(t)}^2 + [A_{ij}|_{\eta(t)}-A_{ij}|_{\gamma(t)}]\dot \gamma^i(t)\dot \gamma^j(t)}dt\leq\\
\leq \int_0^l\sqrt{A_{ij}|_{\gamma(t)}\dot \gamma^i(t)\dot \gamma^j(t) + \abs{\dot \gamma(t)}^2}dt + \int_{0}^l \sqrt{\frac{2\epsilon}{1-\epsilon}}dt \leq l\ton{1+\sqrt{\frac{2\epsilon}{1-\epsilon}}}
\end{gather*}
ricordiamo che $l=d(p,q)$. La tesi segue dal fatto che
\begin{gather*}
 f(\epsilon)=\sqrt{\frac{2\epsilon}{1-\epsilon}}
\end{gather*}
è una funzione continua per $0\leq\epsilon<1$ e $f(0)=0$. Osserviamo che $\epsilon$ non è stata scelta in funzione di $p$ e $q$, ma solo in funzione dell'intorno $U_{\epsilon}$, quindi la stima vale per ogni coppia di punti $p$ e $q$ in $U_{\epsilon}(p_0)$.
\end{proof}
\end{prop}
Osserviamo che questa proposizione vale anche in un altro caso, nel caso in cui $K$ sia una sottovarietà di $R$ con bordo regolare. La dimostrazione appena conclusa non tiene conto della possibilità che $K$ abbia un bordo regolare, ma può essere adattata facilmente. Basta considerare che se $\pi(q)\in \partial K$ (bordo inteso come bordo della sottovarietà regolare), allora la geodetica che unisce $p$ e $\pi(q)$ è di certo più corta della geodetica che unisce $p$ alla proiezione di $q$ sul piano contenente la parte di bordo considerata. Riportiamo quindi la proposizione lasciando i dettagli della dimostrazione al lettore.
\begin{prop}\label{prop_proj2}
  Data $K$ sottovarietà regolare di $R$ possibilmente con bordo liscio, definiamo per ogni $q\in R$:
\begin{gather*}
 \pi(q) \ \ \ t.c. \ \ \ d(\pi(q),q)=\inf\{d(q,p) \ t.c. \ p\in K\}
\end{gather*}
Per ogni $p_0\in K$, esiste un intorno $V$ di $p_0$ dove $\pi$ è una funzione continua ben definita, e per ogni $\epsilon>0$, esiste un intorno $U_{\epsilon}$ di $p_0$ per il quale per ogni $p\in K$ e per ogni $q\in U_{\epsilon}$:
\begin{gather*}
 d(p,q)\leq (1+\epsilon)d(p,\pi(q))
\end{gather*}
\end{prop}

\subsection{Esaustioni regolari}
In questa sezione ci occupiamo di dimostrare l'esistenza di \textbf{esaustioni regolari} per una qualsiasi varietà differenziale liscia connessa $M$ \footnote{in realtà come si può dedurre dalle dimostrazioni, un'esaustione esiste per ogni varietà $M$, e per ogni varietà esiste un'esaustione regolare quanto lo è la varietà}. Per prima cosa dimostriamo l'esistenza di esaustioni non regolari, poi sfruttiamo la densità delle funzioni lisce nelle funzioni continue su una varietà per ottenere l'esistenza di esaustioni regolari. Come risulta evidente dalle definizioni, ha senso procurarsi un'esaustione di $M$ solo nel caso $M$ non compatta.\\
Oltre che ad esaustioni regolari, avremo bisogno anche di esaustioni che si comportano bene in relazione ad una sottovarietà di $R$. In particolare, vogliamo dimostrare che per ogni sottovarietà regolare con bordo $S\subset R$ ($dim(S)=dim(R)$), esiste un'esausione regolare $K_n$ tale che $C=\partial K_n\cap S$ è ancora una sottovarietà regolare di codimensione $1$ con bordo.\\
Ovviamente cominciamo con le definizioni del caso.
\begin{deph}
 Data una varietà differenziale non compatta $M$, diciamo che $\{K_n\}$ è un'esaustione di $M$ se $K_n$ sono insiemi aperti connessi relativamente compatti tali che
\begin{enumerate}
 \item $\overline{K_n}\Subset K_{n+1}$
 \item $M=\bigcup_n K_n$
\end{enumerate}
Spesso nella definzione si confondono $K_n$ e $\overline{K_n}$.\\
Diciamo che questa esaustione è \textbf{regolare} se i bordi $\partial K_n$ sono $C^{\infty}$, o equivalentemente è una sottovarietà di $M$.
\end{deph}
Dimostriamo per prima cosa l'esistenza di un'esaustione.
\begin{prop}
 Ogni varietà differenziale connessa $M$ ammette un'esaustione $K_n$.
\begin{proof}
 Per ogni punto $p\in M$ consideriamo un intorno aperto connesso relativamente compatto $U(p)$ \footnote{che esiste poiché $M$ è localmente euclidea}. Al variare di $p\in M$, $U(p)$ è un ricoprimento aperto di $M$, quindi poiché $M$ è II numerabile, esiste un sottoricoprimento numerabile $U(p_n)\equiv U_n$. Costruiamo per induzione l'esaustione $K_n$.\\
Per prima cosa rinumeriamo gli insiemi $U_n$ in modo che per ogni $N$, $S_N\equiv \bigcup_{n=1}^N U_n$ sia un insieme connesso. Questo è sicuramente vero se $N=1$. Supponiamo che sia vero per $N$. Dato che $M=\cup_n U_n$, se non esistesse nessun $\bar n\geq N$ tale che
\begin{gather*}
 S_N\cap U_{\bar n} \neq \emptyset
\end{gather*}
allora $M$ sarebbe disconnessa perché $S_N\bigcap \ton{\bigcup_{n=N+1}^{\infty} U_n} =\emptyset$. Rinominando $\bar n =N+1$, si ottiene che $S_{N+1}$ è l'unione di insiemi connessi a intersezione non vuota, quindi è connesso.\\
Consideriamo ora $K_1=U_1$, e costruiamo per induzione una successione a valori interi strettamente crescente $k(n)$ tale che
\begin{enumerate}
 \item $k(1)=1$
 \item per ogni $n$, $\bigcup_{j=1}^{k(n)}\overline{U_j} \Subset\bigcup_{j=1}^{k(n+1)}U_j $
\end{enumerate}
Se definiamo $K_n\equiv \bigcup_{j=1}^{k(n)}U_j$, il gioco è fatto.\\
Questo è possibile perché $\overline{K_n} =\bigcup_{j=1}^{k(n)}\overline{U_j}$ è un insieme compatto (unione finita di compatti), ricoperto dagli aperti $U_j$, quindi esiste un sottoricoprimento finito di $U_j$ che ricopre $K_n$. Se chiamiamo $k(n+1)$ l'indice massimo di questo sottoricoprimento, la successione $k(n)$ ha le proprietà desiderate.
\end{proof}
\end{prop}
Grazie all'esistenza di questa esaustione, possiamo costruire una funzione di esaustione $f:M\to \R$ in questo modo:
\begin{prop}
 Per ogni varietà $M$, esiste una funzione continua positiva propria $f:M\to \R$ \footnote{ricordiamo che una funzione è detta \textit{propria} se la controimmagine di un qualsiasi compatto è compatta}.
\begin{proof}
 Definiamo $f(\overline{K_1})=1$, e $f(\partial K_n)=n$. Per ogni insieme $\overline{K_{n+1}}\setminus K_n$ ($n=1,2,\cdots$) estendiamo $f$ a una funzione continua $n\leq f \leq n+1$ grazie al lemma di Urysohn \footnote{vedi teorema 4.1 pag 146 di \cite{1}. Ricordiamo che ogni spazio metrizzabile (e quindi ogni varietà o ogni suo sottoinsieme) è normale}. In questo modo otteniamo una funzione continua $f:M\to \R$ tale che $\overline{K_n}\Subset f^{-1}(-\infty, n]\Subset K_{n+1}$, il che garantisce che $f$ sia propria.
\end{proof}
\end{prop}
Osserviamo che ci sono altri modi di definire una funzione continua positiva propria su $M$. Ad esempio se $M$ è riemanniana geodeticamente completa (e per ogni varietà è possibile definire una metrica $g$ che la renda completa), la funzione distanza da un punto soddisfa queste caratteristiche.\\
Trovata questa funzione, siamo in grado di dimostrare che
\begin{prop}\label{prop_exreg}
 Ogni varietà liscia $M$ ammette un'esaustione regolare $C_n$.
\begin{proof}
 Data un'esaustione $K_n$ e definita $f$ come nella proposizione precedente, sia $h\in C^{\infty} (M,\R)$ tale che $\norm{h-f}_{\infty}<1/4$ \footnote{questo è possibile per densità delle funzioni lisce nello spazio delle funzioni continue con la norma uniforme, vedi teorema 2.6 pag. 49 di \cite{18}, oppure è possibile dimostrare questo fatto adattando la dimostrazione della proposizione \ref{prop_dens} a pagina \pageref{prop_dens}}. Per ogni $n$ consideriamo l'insieme $\Delta_n=(n+1/4,n+3/4)$ aperto in $\R$. Grazie al teorema di Sard, sappiamo che esiste un numero reale $x_n\in \Delta_n$ che sia valore regolare di $h$. Grazie alle proposizioni di questa sezione, questo garantisce che l'insieme $A_n\equiv h^{-1}(-\infty,x_n)$ sia un aperto di $M$ con bordo $\partial A_n=h^{-1}(x_n)$ sottovarietà regolare di $M$.\\
Per definizione di $h$, si verifica facilmente che $K_n\subset A_n\subset K_{n+1}$, quindi $A_n$ è relativamente compatto, però non è detto che $A_n$ sia connesso. Per questo motivo definiamo $C_n$ la componente connessa di $A_n$ che contiene $K_n$. Poiché $\partial C_n$ è una componente connessa di $\partial A_n$, mantiene le sue proprietà di regolarità. Infine, dato che per ogni $n$ $K_n\subset C_n\subset K_{n+1}$, vale che
\begin{enumerate}
 \item $\overline C_n \Subset C_{n+1}$
 \item $M=\bigcup_n C_n$
\end{enumerate}
\end{proof}
\end{prop}
Oltre che ad esaustioni regolari, avremo bisogno anche di esaustioni che si comportano bene in relazione ad una sottovarietà di $R$. In particolare, vogliamo dimostrare che per ogni sottovarietà regolare con bordo $S\subset R$ ($dim(S)=dim(R)$), esiste un'esausione regolare $K_n$ tale che $C=\partial K_n\cap S$ è ancora una sottovarietà regolare di codimensione $1$ con bordo.\\
A questo scopo, assumiamo che $S$ sia chiusa ma non compatta (se $S$ compatta, $S\subset K_n^{\circ}$ definitivamente). Dato che $C^{\circ}=\partial K_n \cap S^{\circ}$ è di certo una sottovarietà regolare, resta da dimostrare solo che anche il bordo $\partial C=\partial K_n\cap \partial S$ è regolare. Riportiamo ora un risultato riguardo alla regolarità dell'intersezione tra due sottovarietà regolari.
\begin{deph}
 Data $R$ varietà riemanniana e $S,M$ sue sottovarietà regolari, diciamo che $M$ è \textbf{trasversale} rispetto a $S$ se e solo se per ogni $x\in S\cap M$:
\begin{gather*}
 T_x(M)+T_x(S)=T_x(R)
\end{gather*}
\end{deph}
\begin{prop}
 Data $R$ varietà riemanniana e date $S,M$ sue sottovarietà regolari senza bordo, se $S$ ed $M$ sono trasversali, allora $S\cap M$ è ancora una sottovarietà regolare di $R$ con:
\begin{gather*}
 codim(S\cap M)=codim(S)+codim(M)
\end{gather*}
dove $codim(S)$ indica la codimensione della sottovarietà, cioè $$codim(S)\equiv dim(R)-dim(S)$$
\begin{proof}
 La dimostrazione di questo risultato ed alcuni approfondimenti sulla teoria della trasversalità possono essere trovati su \cite{31}. In particolare questo teorema è il teorema pag. 30 del cap. 1.
\end{proof}
\end{prop}
Con questo risultato possiamo dimostrare che:
\begin{prop}\label{prop_exreg2}
 Data una sottovarietà regolare $S\subset R$ della stessa dimensione di $R$, chiusa e con bordo regolare $\partial S$, esiste un'esaustione regolare $K_n$ tale che $\partial K_n\cap S$ è una sottovarietà regolare di codimensione $1$ con bordo.
\begin{proof}
 Come affermato sopra, è sufficiente dimostrare che $\partial K_n\cap \partial S$ sia una sottovarietà regolare di codimensione $2$ (è il bordo della sottovarietà $\partial K_m\cap S$).\\
Consideriamo una funzione di esaustione liscia positiva $h$ come nella proposizione \ref{prop_exreg}. La funzione $\bar h \equiv h|_{\partial S}$ è anch'essa una funzione liscia propria e positiva. Sia $\{A_n,\phi_n\}$ una successione di aperti relativamente compatti di $R$ su cui siano definite le coordinate di Fermi $\phi_n$ relative a $\partial S$ tali che
\begin{gather*}
 \partial S\subset \bigcup_n A_n
\end{gather*}
Chiamiamo le funzioni coordinate di Fermi $(x_1,\cdots,x_{m-1},y)$ e definiamo delle funzioni $f_n$ in modo che $f_n:A_n\to \R$, $f_n|_{\partial S}=\bar h$ e che rispetto alle carte $\phi_n$ si abbia:
\begin{gather*}
\tilde {f_n}(x,y)\equiv f_n\circ \phi_n^{-1}(x,y)=f_n(x,0)=\bar h (x) 
\end{gather*}
Cioè definiamo $f_n$ in modo che siano costanti sulle geodetiche perpendicolari a $\partial S$.\\
Consideriamo ora il ricoprimento di aperti di $R$ dato da $\{A_n,(\partial S)^C\}$, e sia $\{\lambda_n,\lambda\}$ una partizione dell'unità subordinata a questo ricoprimento. Definiamo:
\begin{gather*}
 F(p)\equiv \sum_{n} f_n(p)\lambda_n(p) + h(p)\lambda(p)
\end{gather*}
Osserviamo subito che $F\in C^{\infty}(R,\R^+)$. Inoltre vale che per ogni punto $p\in \partial S$, $\nabla(F)|_p\in T_p(S)$. Infatti, sia $p\in A_n\cap S$, indichiamo le sue coordinate rispetto alla carta $\phi_n$ come $p=(\bar x, 0)$ \footnote{$\bar x$ rappresenta il vettore delle prime $m-1$ coordinate} :
\begin{gather*}
 \left.\frac{\partial F}{\partial y}\right|_p = \frac{\partial}{\partial t} F(\bar x,t)=\frac{\partial}{\partial t} \ton{\sum_n f_n(\bar x, t)\lambda_n(\bar x, t)}=\\
=\frac{\partial}{\partial t}\ton{\sum_n h(\bar x)\lambda_n(\bar x, t)}=h(\bar x) \frac{\partial}{\partial t} \sum_n \lambda_n (\bar x, t)=0
\end{gather*}
dove abbiamo sfruttato il fatto che in un intorno di $\partial S$, la funzione $\lambda=0$, quindi $\sum_n \lambda_n (\bar x, t)=1$. Questo garantische che l'ultima componente del differenziale sia nulla, e data la particolare forma della metrica, anche l'ultima componente di $\nabla(F)|_p$ è nulla. Dato che in questa carta l'insieme $\partial S$ è il piano avente ultima coordinata nulla, questo dimostra che $\nabla(F)|_p\in T_p(\partial S)$.\\
Consideriamo ora un valore regolare $c\in \R$ di $F$ \footnote{grazie al teorema di Sard, i valori regolari sono densi in $\R$}. Allora l'insieme $F^{-1}(c)$ è una sottovarietà regolare $F_c$ di $R$. Ricordiamo che per queste sottovarietà:
\begin{gather*}
 T_p(F_c)=\{\nabla(F)|_p\}^{\perp}
\end{gather*}
Dato che $\nabla(F)|_p\in T_p(\partial S)$:
\begin{gather*}
 T_p(F_c)+ T_p(\partial S)\supset \{\nabla(F)|_p\}^{\perp} + \{\nabla(F)|_p\} =T_p(R)
\end{gather*}
quindi le sottovarietà $F_c$ e $\partial S$ sono trasversali, e la loro intersezione è regolare.\\
Se verifichiamo anche che la funzione $F$ è una funzione propria, allora possiamo scegliere una qualsiasi successione $c_n\nearrow \infty$ di valori regolari di $F$ e definire
\begin{gather*}
 K_n \equiv F^{-1}[0, c_n]
\end{gather*}
per ottenere la tesi.\\
Per dimostrare che $F$ è propria, dimostriamo che per ogni $x\in R^+$, esiste un compatto $C_x$ tale che
\begin{gather*}
 p\not\in C_x \ \Rightarrow \ F(p)>x
\end{gather*}
Sia $K=h^{-1}[0,x]$, insieme compatto per il fatto che $h$ è propria. Per locale finitezza delle partizioni dell'unità, esiste un numero finito di indici $n$ per il quale i supporti di $\lambda_n$ intersecano l'insieme $K$. Indichiamo con $I_x$ l'insieme finito di questi indici e consideriamo l'insieme compatto:
\begin{gather*}
 C_x=K\bigcup_{n\in I_x} \overline {A_n}
\end{gather*}
Supponiamo che $p\not\in C_x$. Allora $h(p)>x$, e anche $f_n(p)=h(\pi(p))>x$, dove $\pi(p)$ indica la proiezione di $p$ su $\partial S$, cioè il punto che nelle coordinate di Fermi è caratterizzato da:
\begin{gather*}
 \bar x(\pi(p))=\bar x(p), \ \ y(\pi(p))=0
\end{gather*}
Infatti, $p \not\in C_x \ \Rightarrow \ \pi(p)\not \in K$. Supponiamo per assurdo il contrario. Dato che deve esistere $\bar n$ per cui $p\in A_{\bar n}$, e dato che $\pi(p)\in A_{\bar n}$ per costruzione delle coordinate di Fermi, allora se $\pi(p)\in K$ necessariamente $\bar n\in I_x$, assurdo poiché abbiamo assunto $p\not \in C_x$.\\
Ricordando la definizione di $F$, è immediato verificare che $F(p)>x$.
\end{proof}
\end{prop}
Ripetendo questa costruzione per una successione di sottovarietà $\{S_n\}$ con bordi disgiunti possiamo ottenere:
\begin{prop}\label{prop_exreg3}
 Sia $\{S_n\}$ una successione di sottovarietà di $R$ della sua stessa dimensione con bordi $\partial S_n$ lisci disgiunti tra loro. Allora esiste un'esaustione regolare $K_m$ tale che per ogni $n$, $m$ l'insieme $\partial K_m \cap S_n$ è una sottovarietà di $R$ $m-1$ dimensionale con bordo liscio.
\begin{proof}
 L'idea della dimostrazione è la stessa della dimostrazione precedente. Scegliamo per ogni sottovarietà $\partial S_n$ un insieme numerabile di suoi intorni $\{I_n^m\}_{m=1}^\infty$ che ammettono coordinate di Fermi, restringendoli in modo che questi intorni siano disgiunti dagli altri bordi $\partial S_k$ (con $k\neq n$). Fissiamo inoltre una partizione dell'unità subordinata al ricoprimento di aperti
$$\left\{\{I_n^m\}_{n,m=1}^\infty ; \ \bigcap_n \partial S_n^C\right\}$$
Ripetendo una costruzione del tutto analoga alla precedente otteniamo la tesi.
\end{proof}
\end{prop}

\subsection{Insiemi chiusi e funzioni lisce}
Lo scopo di questo paragrafo è ottenere l'osservazione \ref{oss_an}, che sarà utile nel seguito del lavoro.
\begin{lemma}\label{lemma_T4liscio}
 Dato un insieme chiuso $C\subset \R^n$, esiste una funzione $f:\R^n\to \R$ liscia positiva tale che $C=f^{-1}(0)$.
\begin{proof}
 Grazie al fatto che $\R^n$ è uno spazio metrico, è facile trovare una funzione $g:\R^n\to\R$ continua per cui $C=g^{-1}(0)$, basta considerare infatti la funzione
\begin{gather*}
 g(x)=d(x,C)\equiv \inf_{y\in C} d(x,y)
\end{gather*}
ha le caratteristiche cercate. Consideriamo ora gli insiemi
\begin{gather*}
 C_{\epsilon}\equiv \{x\in \R^n \ t.c. \ d(x,C)\leq \epsilon\}
\end{gather*}
data la continuità della funzione $d(\cdot,C)$, tutti questi insiemi sono insiemi chiusi, quindi esistono funzioni $g_{\epsilon}:\R^n\to \R$ tali che $g_{\epsilon}^{-1}(0)=C_{\epsilon}$. Notiamo che per ogni $y\in C$, $B_{\epsilon}(y)\subset C_{\epsilon}$, quindi in particolare $g(B_{\epsilon(y)})=0$ per ogni $y\in C$, e che
\begin{gather*}
 C=\bigcap_{\epsilon>0} C_{\epsilon}=\bigcap_{n=1}^{\infty} C_{1/n} 
\end{gather*}
Grazie al lemma \ref{lemma_conv_0}, esiste un nucleo di convoluzione $\Theta_{\epsilon}$ con supporto compatto contenuto in $B_{\epsilon}(0)$, e grazie al lemma \ref{lemma_conv_1}, la funzione
\begin{gather*}
 g'_{\epsilon}\equiv g_{\epsilon}\ast \Theta_{\epsilon}
\end{gather*}
è una funzione liscia da $\R^n$ a $\R$ per la quale $g'_{\epsilon}(y)=0$ per ogni $y\in C$ e $g'_{\epsilon}(x)>0$ per ogni $x\in \R^n$ tale che $d(x,C)>\epsilon$.\\
Definiamo la funzione
\begin{gather*}
 f(x)\equiv \sum_{n=1}^{\infty} \frac{g'_{1/n}(x)}{A_n} \frac{1}{2^n}
\end{gather*}
dove $A_n=\max_{i=1}^n \norm{D^i g'_{1/n}}_{\infty, B_n(0)}$. Grazie a questa definizione è facile verificare che la serie che definisce $f(x)$ converge localmente uniformemente in $\R^n$ assieme alla serie delle derivate di qualunque ordine, quindi la funzione $f$ è una funzione liscia. Ovviamente $f$ è positiva, e dato che per ogni $n$, $g'_{1/n}(y)=0$ per ogni $y\in C$, anche $f(y)=0$ per ogni $y\in C$. Se $y\not\in C$, grazie al fatto che $C$ è chiuso, esiste $\bar n$ sufficientemente grande per cui $y\not \in C_{1/\bar n}$, quindi dalle considerazioni precedenti abbiamo che $g'_{1/\bar n}(y)>0$, quindi anche $f(y)>0$. Da cui la tesi.
\end{proof}
\end{lemma}
\begin{oss}
 Utilizzando le partizioni dell'unità, è possibile estendere questo lemma a una varietà differenziale $M$ qualsiasi \footnote{se vogliamo $f$ liscia, la varietà $M$ deve avere una struttura $C^{\infty}$}.
\begin{proof}
 Sia $M$ una varietà differenziale qualsiasi e sia $\lambda_n$ una partizione dell'unità di $M$ subordinata a un ricoprimento di carte locali $(U_n,\phi_n)$. Sia $f_n$ la funzione descritta nel lemma precedente relativa all'insieme chiuso in $\R^m$ $\phi_n(C\cap~ supp(\lambda_n))$. Allora la funzione
\begin{gather*}
 f(p)\equiv \sum_n f_n\circ \phi_n(p)\cdot \lambda_n(p)
\end{gather*}
ha le caratteristiche cercate. Infatti dato che la somma è localmente finita, $f$ è una funzione liscia. La positività di $f$ è ovvia conseguenza della positività delle funzioni $f_n$ e $\lambda_n$; se $p\in C$, allora $f_n(\phi_n(p))\cdot \lambda_n(p)=0$, mentre se $p\not \in C$, esiste $\bar n$ tale che $\lambda_{\bar n}(p)>0$ \footnote{quindi in particolare $p\in supp(\lambda_{\bar n})$} e quindi $f_n(\phi_n(p))\cdot\lambda_{\bar n}(p)>0$.
\end{proof}
\end{oss}
\begin{oss}\label{oss_an}
Grazie alla precedente osservazione e al teorema di Sard, possiamo concludere che per ogni insieme $C$ chiuso in una varietà $R$, esiste una successione di insiemi aperti con bordo liscio $A_n\subset A_{n-1}$ in $R$ tali che $C=\cap_n A_n$. Sia infatti $f$ una funzione con le caratteristiche appena descritte. Consideriamo una successione di valori regolari $a_n$ di $f$ tali che $a_n\searrow 0$, allora è facile verificare che gli insiemi
\begin{gather*}
 A_n=f^{-1}[0,a_n)
\end{gather*}
soddisfano le proprietà richieste. Osserviamo inoltre che se $C$ è compatto, gli insiemi possono essere scelti relativamente compatti (lasciamo al lettore i facili dettagli di questa dimostrazione).
\end{oss}
\begin{oss}\label{oss_an2}
 Sia $\Omega$ un aperto di $R$. Allora esiste una successione di aperti $A_n$ con bordo liscio tali che:
\begin{gather*}
  A_n\subset \Omega \ \ \ A_n\subset A_{n+1} \ \ \ \bigcup_n A_n=\Omega
\end{gather*}
Se $\Omega$ è relativamente compatto, anche gli insiemi $A_n$ lo sono.
\begin{proof}
 Sia $f$ una funzione liscia positiva tale che $f^{-1}(0)=\Omega^C$, e sia $\epsilon_n\searrow 0$ una successione di valori regolari per la funzione $f$. È facile verificare che gli insiemi
\begin{gather*}
 A_n\equiv f^{-1}(\epsilon_n,\infty)
\end{gather*}
hanno le caratteristiche desiderate.
\end{proof}
\end{oss}
\section{Funzioni assolutamente continue}\label{sec_ac}
In questa sezione riportiamo alcuni risultati sulle funzioni assolutamente continue senza dimostrazione. Come referenza consigliamo il paragrafo 5.4 di \cite{10}.
\begin{deph}\label{deph_ac}
 Una funzione $f:U\to \R$ (dove $U\subset \R$ è un aperto) è detta assolutamente continua se per ogni $\epsilon>0$, esiste $\delta>0$ tale che per ogni scelta di $\{x_i,x_i'\}\subset U$ \footnote{in questa definizione è equivalente chiedere che $i$ vari su un insieme finito o numerabile di indici} tali che $x_i'>x_i>x_{i-1}'$ e tali che $\sum_i (x_i'-x_i)<\delta$, allora:
\begin{gather*}
 \sum_i \abs{f(x_i')-f(x_i)} < \epsilon
\end{gather*}
\end{deph}
\begin{prop}\label{prop_ac}
 Se una funzione $f$ è assolutamente continua, allora la sua derivata esiste finita quasi ovunque ed è Lebesgue integrabile sui compatti. Inoltre $f$ è l'integrale della sua derivata se e solo se è assolutamente continua. Tutte le funzioni di Lipschitz sono assolutamente continue.
\end{prop}
\begin{prop}\label{prop_ac_times}
 Date due funzioni assolutamente continue limitate $f$ e $g$, il loro prodotto è assolutamente continuo.
\begin{proof}
 La dimostrazione segue facilmente dalla disuguaglianza:
\begin{gather*}
 \abs{f(x)g(x)-f(y)g(y)}\leq \abs{f(x)\cdot (g(x)-g(y))}+\abs{g(y)\cdot (f(x)-f(y))}\leq\\
\leq \max\{\norm f _{\infty},\norm g _{\infty}\}\cdot(\abs{f(x)-f(y)}+\abs{g(x)-g(y)})
\end{gather*}
\end{proof}
\end{prop}

\section[Spazi L2 di forme]{Spazi $\mathcal{L}^2$ di forme}\label{sec_L2}
In teoria dell'intergrazione sono famosi gli spazi $L^p$, spazi di funzioni misurabili il cui modulo elevato alla $p$ è integrabile sullo spazio. Per un'introduzione su questa teoria consigliamo \cite{12} (capitolo 3). I risultati che ci interessano sono comunque riportati in questa breve rassegna:
\begin{deph}
 Dato uno spazio di misura $(X,\mu)$ con $\mu$ misura positiva, definiamo per $1\leq p < \infty$ $L^p(X,\mu)$ come lo spazio delle funzioni $f:X\to \R$ (o anche $f:X\to \C$) misurabili \footnote{rispetto alla misura $\mu$ e all'algebra di Borel su $\R$} tali che
\begin{gather*}
 \norm f _p ^p \equiv \int_X \abs f ^p d\mu <\infty
\end{gather*}
Se introduciamo la relazione di equivalenza $f\sim g$ se $f=g$ quasi ovunque, lo spazio quoziente $L^p/\sim$ con la norma $\norm \cdot _p$, $L^p$ è uno spazio di Banach, e per $p=2$ uno spazio di Hilbert.
\end{deph}
\begin{prop}\label{prop_lp_1}
 Se $f_n\to f$ rispetto a una qualsiasi norma $\norm \cdot _p$, allora esiste una sottosuccesione $f_{n_k}$ che converge puntualmente quasi ovunque a $f$. Questo implica anche che se $f_n$ converge uniformemente a $f$ e $f_n$ converge in norma $p$ ad $h$, allora $f=h$ quasi ovunque. 
\end{prop}
\begin{prop}\label{prop_Ldens}
Data una varietà Riemanniana $R$, gli spazi $L^p(R)$ sono separabili per $1\leq p<\infty$.
\begin{proof}
 Grazie al teorema 3.14 pag 68 di \cite{12}, sappiamo che l'insieme delle funzioni continue a supporto compatto in $R$ (che indicheremo $C_C(R)$ è denso nello spazio $L^p(R)$. Sia $K_n$ un'esaustione di $R$, allora:
$$
C_C(R)=\bigcup_n C(K_n)
$$
dato che $C(K_n)$ dotato della norma del sup è separabile, e dato che la norma del sup è più forte della norma $L^p$ su insiemi compatti, allora $C(K_n)$ è separabile anche nella norma $L^p$. Siano $D_n$ insiemi numerabili densi in $C(K_n)$, e sia $D=\cup_n D_n$. Allora $D$ è numerabile e denso in $L^p(R)$, infatti data $f\in L^p(R)$ e dato $\epsilon>0$, esiste $g\in C_C(R)$ tale che $\norm{g-f}_p<\epsilon$. Sia $K_n$ tale che $supp(g)\subset K_n$, allora esiste $h\in D_n\subset D$ tale che $\norm{g-h}_p<\epsilon$, quindi:
\begin{gather*}
 \norm{f-h}_p\leq \norm{f-g}_p+\norm{g-h}_p\leq 2 \epsilon
\end{gather*}
data l'arbitrarietà di $\epsilon$, si ottiene la tesi.
\end{proof}

\end{prop}

Una facile generalizzazione di questi spazi sono gli spazi $L^p$ per forme su varietà Riemanniane. La teoria di questi spazi è del tutto analoga agli spazi $L^p$ di funzioni, quindi anche in questo caso non riporteremo tutte le dimostrazioni. Un riferimento su questo argomento può essere \cite{13}, capitolo 7 \footnote{in realtà questo libro si occupa solo di superfici riemanniane, ma in questo frangente non c'è nessuna sostanziale differenza con una varietà di dimensione generica}.
\begin{deph}
 Data una varietà Riemanniana $(R,g)$, definiamo lo spazio $\mathcal{L}^2(R)$ come lo spazio delle 1-forme a quadrato integrabile su $R$. Sia cioè $\alpha\in T^*(R)$ una 1-forma (in coordinate locali $\alpha=\sum_{i=1}^m\alpha_i(x)dx^i$), $\alpha\in \mathcal{L}^2(R)$ se e solo se $\alpha_i$ sono tutte misurabili e
\begin{equation*}
 \int_R \abs{\alpha}^2 dV<\infty
\end{equation*}
\end{deph}
Ricordiamo che l'integrale su una varietà è definito tramite partizioni dell'unità. Sia $\{\lambda_n\}$ una partizione dell'unità di $R$ subordinata a un ricoprimento di aperti coordinati. Per definizione
\begin{gather*}
 \int_R f dV=\sum_{n=1}^{\infty} \int_{supp(\lambda_n)} \lambda_n\cdot f \ dV=\sum_{n=1}^{\infty} \int_{\phi_n(supp(\lambda_n))} \tilde\lambda_n(x) \tilde f(x) \ \sqrt{\abs g} dx^1\dots dx^m
\end{gather*}
Ricoriamo anche che $\abs {\tilde\alpha(x)} ^2 = g_{ij}(x) \tilde\alpha^i(x)\tilde\alpha^j(x)$, dove le funzioni con la tilde indicano le rappresentazioni in coordinate delle relative funzioni.\\
Anche per questo spazio valgono gli analoghi delle proposizioni \ref{prop_lp_1} e \ref{prop_Ldens}:
\begin{prop}\label{prop_L2}
 Lo spazio $\mathcal{L}^2(R)$ è uno spazio di Hilbert, e se $\alpha_n\to \alpha$ in norma, allora esiste una sottosuccesione $\alpha_{n_k}$ che converge puntualmente quasi ovunque ad $\alpha$. Inoltre $\mathcal L^2(R)$ è separabile.
\end{prop}
Anche se in questa tesi non toccheremo l'argomento, osserviamo che $$\abs \alpha ^2 dV =\alpha \ast \alpha$$ dove $\ast$ indica l'operatore duale di Hodge.\\
Osserviamo che l'integrale di Dirichlet di una funzione
\begin{gather*}
 \int_R\abs{\nabla f}^2dV=\int_R\abs{df}^2dV
\end{gather*}
è la norma nello spazio $\mathcal{L}^2(R)$ della 1-forma $df$.

\section{Derivazione sotto al segno d'integrale}
In questa sezione ripotiamo due lemmi che consentono sotto certe ipotesi di scambiare derivata e integrale. Il secondo lemma è una generalizzazione del primo, anche se più difficile da dimostrare.
\begin{lemma}[Derivazione sotto al segno di integrale]\label{diid}
 Sia $U\subset \R$ un aperto e $(\Omega,\mu)$ uno spazio di misura. Supponiamo che $f:U\times \Omega \to \R$ abbia le proprietà:
\begin{enumerate}
 \item $f(x,\omega)$ è Lebesgue-integrabile per ogni $x\in U$
 \item quasi ovunque rispetto a $\mu$ la funzione $\frac{\partial f(x,\omega)}{\partial x}$ esiste per ogni $x\in U$
 \item esiste una funzione $g\in L^1(\Omega,\mu)$ tale che $\abs{\frac{\partial f(x,\omega)}{\partial x}}\leq g(\omega)$ per ogni $x\in U$
\end{enumerate}
Allora vale che:
\begin{gather*}
 \frac{d}{d x} \int_{\Omega}f(x,\omega) d\mu = \int_{\Omega}\frac{\partial}{\partial x}f(x,\omega) d\mu
\end{gather*}
\begin{proof}
 La dimostrazione è una semplice applicazione del teorema di convergenza dominata (vedi teorema 1.34 pag. 26 di \cite{12}).
\end{proof}
\end{lemma}
Prima di enunciare il secondo lemma, riportiamo una proposizione che sarà utile per la sua dimostrazione
\begin{prop}\label{prop_inted}
 Se $f\in L^1(X,\mu)$, dove $(X,\mu)$ è uno spazio di misura positiva qualsiasi, allora per ogni $\epsilon>0$ esiste $\delta>0$ tale che
\begin{gather*}
 \mu(E)<\delta \Rightarrow \int_{E} \abs f d\mu <\epsilon
\end{gather*}
dove $E$ è un qualsiasi sottoinsieme misurabile di $(X,\mu)$.
\begin{proof}
 Osserviamo che questo è il testo dell'esercizio 12 pag 32 di \cite{12}. Supponiamo per assurdo che esista $\epsilon>0$ tale che per ogni $n\in \N$ esiste $E_n$ con $\mu(E_n)<1/n$ tale che $\int_{E_n} \abs f \geq \epsilon$. Consideriamo la successione $\abs f \chi_{E_n}$. Questa successione converge in misura a $0$, quindi per il teorema di convergenza dominata $\int_{E_n}\abs f d\mu=\int_{X} f\chi_{E_n}d\mu\to0$, assurdo.
\end{proof}
\end{prop}
\begin{lemma}[Derivazione sotto al segno di integrale]\label{diid+}
Sia $U\subset \R$ un aperto e $(\Omega,\mu)$ uno spazio di misura. Supponiamo che $f:U\times \Omega \to \R$ abbia le proprietà:
\begin{enumerate}
 \item $f(x,\omega)$ è misurabile su $U\times \Omega$, e per quasi ogni $x\in U$ è integrabile su $\Omega$
 \item quasi ovunque rispetto a $\mu$ la funzione $f(x,\omega)$ è assolutamente continua in $x$
 \item \begin{gather*}
 \int_U dx\int_{\Omega} d\mu\abs{\frac{\partial f}{\partial x}(x,\omega)}<\infty
\end{gather*}
\end{enumerate}
Allora $F(x)\equiv\int_{\Omega} f(x,\omega)d\mu$ è assolutamente continua in $x$ e per quasi ogni $x$:
\begin{gather*}
 \frac{d}{d x} \int_{\Omega}f(x,\omega) d\mu = \int_{\Omega}\frac{\partial}{\partial x}f(x,\omega) d\mu
\end{gather*}
\begin{proof}
Per prima cosa dimostriamo che $F(x)$ è assolutamente continua. Grazie alla proposizione precedente (la \ref{prop_inted}), per ogni $\epsilon>0$ esiste $\delta>0$ per cui se $\lambda(E)<\delta$ \footnote{$\lambda$ indica la misura di Lebesgue su $\R$}, allora
\begin{gather*}
\int_{E} dx\int_{\Omega} d\mu\abs{\frac{\partial f}{\partial x}(x,\omega)}<\epsilon 
\end{gather*}
Siano $\{a^i,b^i\}$ tali che $b_i>a_i>b_{i-1}$ e $\sum_i (b_i-a_i)<\delta$, e chiamiamo $E=\cup_i (a^i,b^i)$. Allora
\begin{gather*}
 \sum_i \abs{F(b^i)-F(a^i)}=\sum_i\abs{\int_{\Omega} f(b^i,\omega)-f(a^i,\omega)d\mu}=\\
=\sum_i\abs{\int_{\Omega} \int_{a^i}^{b^i}\frac{\partial f}{\partial x}(x,\omega)\ dx\ d\mu}\leq \int_{E} dx\int_{\Omega} d\mu\abs{\frac{\partial f}{\partial x}(x,\omega)}<\epsilon 
\end{gather*}
quindi $F$ è assolutamente continua in $x$. Questo significa che esiste quasi ovunque $F'(x)\equiv \partial F(x) /\partial x$ e che $F(b)-F(a)=\int_{a}^b F'(t)dt$  ogni volta che $[a,b]\subset U$. Ma
\begin{gather*}
 F(b)-F(a)=\int_{\Omega} f(b,\omega)-f(a,\omega)\ d\mu =\int_{\Omega}\int_a^b \frac{\partial f}{\partial x}(x,\omega) \ dx\ d\mu=\\
=\int_a^b \int_{\Omega} \frac{\partial f}{\partial x}(x,\omega) \ d\mu\ dx
\end{gather*}
dove l'ultimo passaggio è giustificato dal teorema di Fubini (vedi teorema 7.8 pag 140 di \cite{12}). Ma allora si ha che per ogni $[a,b]\subset U$:
\begin{gather*}
 \int_a^b \ton{F'(x)-\int_{\Omega} \frac{\partial f}{\partial x}(x,\omega) \ d\mu} dx =0
\end{gather*}
quindi grazie a una proprietà nota degli integrali \footnote{vedi teorema 1.39 (b) pag 29 di \cite{12}}
\begin{gather*}
F'(x)=\int_{\Omega} \frac{\partial f}{\partial x}(x,\omega) \ d\mu 
\end{gather*}
dove l'uguaglianza è intesa quasi ovunque in $x$ in ogni componente connessa dell'insieme $U$, quindi in tutto $U$.
\end{proof}
\end{lemma}
\section{Convoluzioni}\label{sec_conv}
In questa sezione introdurremo la convoluzione tra funzioni reali (non nella forma più generale possibile, ma nella forma utile ai nostri scopi), e esploreremo alcune tecniche di regolarizzazione di funzioni tramite convoluzione. Per prima cosa dimostriamo l'esistenza dei \textit{nuclei di convoluzione}.
\begin{lemma}\label{lemma_conv_0}
 Per ogni $\alpha>0$ esiste una funzione positiva $\Theta_{\alpha}\in C^{\infty}(\R^m,\R)$ con $supp(\Theta_{\alpha})\Subset \overline{B_{\alpha}(0)}$ e $\int_{\Rp^m}\Theta_{\alpha}(x)dx=1$. Chiamiamo queste funzioni \textbf{nuclei di convoluzione}.
\begin{proof}
 \`E sufficiente trovare una funzione con le caratteristiche descritte per $\alpha=1$. Infatti è facile verificare che la funzione $$\Theta_{\alpha}(x)\equiv \frac{1}{\alpha^m}\Theta\left(\frac{x}{\alpha}\right)$$ verifica tutte le richieste.\\
Per trovare la funzione $\Theta_1$ basta considerare la funzione
\begin{gather*}
 \tilde\Theta_1(x)\equiv\begin{cases} exp\left(\frac{1}{\norm x ^2 -1} \right) & se \norm x \leq 1 \\
               0 & se \norm x \geq 1
              \end{cases}\\
\Theta_1(x)\equiv\frac{\tilde \Theta_1(x)}{\int_{\Rp^m}\tilde\Theta_1(x)dx}
\end{gather*}
oppure è possibile sfruttare l'esistenza delle partizioni dell'unità. Trovata una funzione a supporto compatto $0\leq\lambda(x)\leq 1$ (non identicamente nulla) con le tecniche descritte qui sopra la si può traslare e scalare in modo da ottenere la funzione desiderata.
\end{proof}

\end{lemma}

Passiamo a definire la convoluzione, operazione che si rivelerà molto utile per approssimare funzioni abbastanza generiche con una funzioni lisce. La dimostrazione di queste (e altre) proprietà può essere trovata su \cite{12} al \textsection 7.13, oppure su \cite{14} al \textsection 8.2.
\begin{prop}\label{prop_mate_young}
 \textbf{Disuguaglianza di Young:} Date $f \in L_1(\R^m)$ e $g \in L_p(\R^m)$ ($1\leq p \leq \infty$), si definisce:
\begin{equation*}
 f*g(x)=\int_{\Rp} \ f(x-y)g(y) dy
\end{equation*}
Questa definizione ha senso solo quasi ovunque e vale che:
\begin{equation}
 \norm{f*g}_p\leq \norm f _1 \cdot \norm g _p
\end{equation}
\end{prop}
\begin{prop}
 Date $f \in L_p(\R^m)$ e $g \in L_q(\R^m)$ con $\frac{1}{p}+\frac{1}{q}=1, \ 1\leq p \leq \infty$, ha senso definire $\forall x$:
\begin{displaymath}
 f*g(x)=\int_{\Rp} \ f(x-y)g(y) dy
\end{displaymath}
e vale che $f*g$ è uniformemente continua su $\R^m$ e:
\begin{displaymath}
 \norm{f*g}_{\infty} \leq \norm f _p \cdot \norm g _q
\end{displaymath}
\end{prop}
Osserviamo subito che $\ast$ è un'operazione commutativa, cioè $f\ast g =g\ast f$. La convoluzione può essere utilizzata per regolarizzare una funzione, cioè per trovare una funzione liscia vicina a piacere alla funzione data. I dettagli sono nei lemmi seguenti.\\
Per comodità di notazione, ricordiamo brevemente la definizione di multiindice. Un vettore di numeri interi non negativi $\vec k=(k_1,\cdots k_m)$ è detto multiindice. La sua lunghezza è per definizione $\abs {\vec k}=\sum_{i=1}^m k_i$, inoltre con il simbolo $D^{\vec k} f$ intendiamo:
\begin{gather*}
 D^{\vec k} f \equiv \frac{\partial^{\abs {\vec k}} f}{\partial x_1^{k_1}\cdots \partial x_m^{k_m}}
\end{gather*}

\begin{lemma}\label{lemma_conv_1}
 Sia $\Theta\in C^r(\R^m,\R)$ con $supp(\Theta)\Subset \overline{B_{\alpha}(0)}$ e $\int \Theta(x)dx=1$ con $0\leq r\leq \infty$, sia $f\in L^1(\R^m)$. Allora si ha che:
\begin{enumerate}
 \item $(\Theta \ast f) \in C^r(\R^m,\R)$
 \item Per ogni $\abs {\vec k}\leq r$, $D^{\vec k}(\Theta \ast f)\vert_x=(D^{\vec k} (\Theta) \ast f)\vert_x$
\end{enumerate}
Inoltre se $f\in C^s(\R^m,\R)$ ($0\leq s\leq \infty$) con $supp(f)=C$, allora si ha anche che:
\begin{enumerate}
 \item $(\Theta \ast f) \in C^s(\R^m,\R)$
 \item Per ogni $\abs{\vec k}\leq s$, $D^{\vec k}(\Theta \ast f)\vert_x=(\Theta \ast(D^{\vec k} f))\vert_x$
 \item $supp(\Theta\ast f)\subset C+\alpha=\{x\in \R^m \ t.c. \ d(x,C)\leq \alpha\}$
\end{enumerate}
\begin{proof}
 La dimostrazione è una semplice applicazione del lemma \ref{diid}.
\end{proof}
\end{lemma}
Nel seguito avremo bisogno di una versione più raffinata di questo lemma, in particolare con richieste meno stringenti sulla regolarità di $f$. A questo scopo dimostriamo che:
\begin{lemma}\label{lemma_conv_1b}
Sia $f\in C(\R^m)$ con $supp(f)=K\Subset \prod_{i=1}^m(a_i,b_i)\equiv U$. Sia inoltre $f(\bar x^1,\cdots, x^i,\cdots,\bar x^m)$ assolutamente continua rispetto a $x^i$ quasi ovunque rispetto a $(\bar x^1,\cdots,\bar x^{i-1},\bar x^{i+1},\cdots,\bar x^m)$, con $\partial f /\partial x^i\in L^1(\R^m)$. Allora se $\Theta$ è un nucleo di convoluzione:
\begin{gather*}
 \frac{\partial}{\partial x^i}(\Theta\ast f)=\Theta \ast \frac{\partial f}{\partial x^i}
\end{gather*}
\begin{proof}
Questo lemma è conseguenza del lemma \ref{diid+}, ma per completezza riportiamo anche una dimostrazione più ``elementare''.\\ 
Per comodità di notazione indicheremo
\begin{gather*}
 \frac{\partial}{\partial x^i}\equiv \partial_i
\end{gather*}
Osserviamo che quasi ovunque rispetto alle $\bar x$, $f(\bar x ,x^i)=\int_{-\infty} ^{x^i} \partial_i f(\bar x,t)dt$. Consideriamo $f'_n$ successione di funzioni continue a supporto compatto $f'_n:\R^m\to \R$ tali che $\norm{f'_n-\partial_i f}_{L^1}\to 0$ \footnote{questa successione esiste per densità di $C_C(\R^m)$ in $L^1(\R^m)$, vedi teorema 3.14 pag 68 di \cite{12}}. Definiamo
\begin{gather*}
 f_n(\bar x, x^i)=\int_{-\infty}^{x^i} f'_n(\bar x,t)dt
\end{gather*}
in questo modo la successione $f_n$ converge in norma $L_1$ a $f$, infatti:
\begin{gather*}
 \norm{f_n-f}_{L^1}\leq\int_{\Rp^m} \int_{\infty}^{x^i} \abs{f'_n(\bar x,x^i)-\partial_if(\bar x,x^i)}dt\  d\bar x dx^i\leq\\
\leq (b_i-a_i) \norm{f'_n-\partial_i f}_{L^1}\to 0
\end{gather*}
Grazie ai lemmi precedenti osserviamo che
\begin{gather*}
 \partial_i (\Theta\ast f_n) = (\partial_i \Theta)\ast f_n \to (\partial_i \Theta)\ast f=\partial_i(\Theta \ast f)
\end{gather*}
e anche
\begin{gather*}
 \partial_i (\Theta\ast f_n) = \Theta\ast (\partial_i f_n)=\Theta\ast f'_n= \to \Theta\ast (\partial_i f)
\end{gather*}
dove la convergenza nelle ultime due uguaglianze è intesa nel senso di $L^1$. Per unicità del limite, $\partial_i(\Theta \ast f)=\Theta\ast (\partial_i f)$ nel senso di $L^1$, quindi quasi ovunque. Ma le funzioni $\Theta\ast (\partial_i f)$ e $\partial_i(\Theta \ast f)$ sono funzioni lisce, quindi uguaglianza quasi ovunque implica uguaglianza ovunque.
\end{proof}
\end{lemma}
Ora analizziamo le proprietà di approssimazione della convoluzione. A questo scopo utilizzeremo la notazione
\begin{gather*}
 \norm{f}_{\infty}\equiv \max_{x\in \R^m} \abs{f(x)}
\end{gather*}
\begin{lemma}\label{lemma_conv_2}
 Sia $f\in C^r(\R^m,\R)$ con $supp(f)=K$ compatto e sia $\vec k$ un multiindice di lunghezza $\abs{\vec k}\leq r$. Dato $\epsilon>0$, esiste $\alpha>0$ tale che per ogni $\Theta$ con supporto in $\overline{B_{\alpha}(0)}$ e $\int \Theta(x)dx=1$ si ha che $\norm{D^{\vec k}(\Theta \ast f) -D^{\vec k}(f)}_{\infty}<\epsilon$.
\begin{proof}
Per uniforme continuità di $D^{\vec k}(f)$ su $K$, vale che dato $\epsilon>0$, esiste $\alpha>0$ tale che $\norm{x-y}\leq\alpha \Rightarrow \abs{D^{\vec k} f|_{y}-D^{\vec k}f|_{x}}<\epsilon$. Dato che:
\begin{gather*}
 D^{\vec k} (\Theta \ast f)|_x-D^{\vec k}(f)|_x=\int_{B_{\alpha}(0)}{\Theta(y)} (D^{\vec k} (f)|_{x-y}-D^{\vec k}(f)|_x) dy
\end{gather*}
si ottiene:
\begin{gather*}
 \abs{D^{\vec k} (\Theta \ast f)|_x-D^{\vec k}(f)|_x}\leq\\
\leq \int_{B_{\alpha}(0)}\abs{{\Theta(y)}}\abs{ (D^{\vec k} (f)|_{x-y}-D^{\vec k}(f)|_x)} dy\leq \epsilon \int_{B_{\alpha}(0)}\abs{{\Theta(y)}}=\epsilon
\end{gather*}
\end{proof}
\end{lemma}
Questo lemma garantisce che ogni funzione a supporto compatto può essere approssimata uniformemente a ogni suo ordine di regolarità con una funzione liscia.
\section{Duali di spazi di Banach}
In questa sezione ricordiamo brevemente alcuni risultati riguardo agli spazi duali degli spazi di Banach, in particolare la definizione di topologia debole-* e il teorema di Banach-Alaoglu. Per approfondimenti rimandiamo ai capitoli 3 e 4 di \cite{4}.
\begin{deph}
 Dato uno spazio di Banach (reale) $(X,\norm \cdot)$, si definisce $X^*$ il suo spazio duale, cioè:
\begin{gather*}
 X^*=\{\phi:X\to \R \ \ t.c. \ \phi \ lineare \ e \ continuo\}
\end{gather*}
\end{deph}
Possiamo rendere questo spazio uno spazio di Banach con la norma:
\begin{gather*}
\normop \phi \equiv \sup_{x\neq 0} \frac{\abs{\phi(x)}}{\norm x} =\sup_{\norm x \leq 1} \abs{\phi(x)}=\sup_{\norm x = 1} \abs{\phi(x)}
\end{gather*}
Notiamo che per un funzionale lineare, essere continuo, essere continuo nel punto 0 ed essere limitato (cioè avere $\normop{\phi}<\infty$) sono proprietà equivalenti. Sullo spazio $X^*$ però è possibile definire anche una topologia vettoriale più debole della topologia indotta da questa norma, la topologia debole-*.
\begin{deph}
 Sullo spazio $X^*$ definiamo $\tau^*$ la topologia debole-*, una topologia invariante per traslazioni tale che una base di intorni del punto $0$ è costituita dagli insiemi
\begin{gather*}
 V(0,\epsilon,x_1,\cdots,x_n)=\{\phi \ t.c. \ \abs{\phi(x_i)}<\epsilon \ \forall \ 1\leq i\leq n\}
\end{gather*}
al variare di $\epsilon$ e di $x_1,\cdots,x_n\in X$ \footnote{al variare degli elementi e del numero degli elementi, che può essere un numero \textit{finito} qualsiasi}. Questa topologia rende $X^*$ uno spazio vettoriale topologico.
\end{deph}
Ricordiamo che per gli spazi di Hilbert reali esiste un'isometria lineare tra lo spazio $H$ e il suo duale $H^*$, isometria data dal teorema di rappresentazione di Riesz (vedi ad esempio teorema 3.4 di \cite{35}).\\
Per la topologia debole-* vale un teorema molto famoso, il teorema di \textit{Banach Alaoglu} (vedi teorema 3.15 pag 68 di \cite{4}). In questa rassegna riportiamo una versione del teorema adatta ai nostri scopi
\begin{teo}[Teorema di Banach-Alaoglu]\label{teo_BA}
 Nello spazio $(X^*,\tau^*)$, l'insieme $$B^*=\{\normop \phi \leq 1\}$$ è compatto. Inoltre se $X$ è separabile, allora $B^*$ è anche sequenzialmente compatto.
\begin{proof}
La compattezza di $B^*$ è dimostrata nel teorema 3.15 pag 68 di \cite{4}, per quanto riguarda la compattezza per successioni, il teorema 3.17 a pagina 70 di \cite{4} garantisce che se $X$ è separabile, allora gli insiemi compatti di $X^*$ sono metrizzabili, quindi anche sequenzialmente compatti.
\end{proof}
\end{teo}

\section{Funzioni armoniche}
In questa sezione riportiamo alcuni risultati riguardo alle funzioni armoniche su varietà riemanniane. In tutta la sezione $R$ sarà una varietà Riemanniana liscia senza bordo di dimensione $m$.
\subsection{Principio del massimo}\label{sec_max}
In questo paragrafo riportiamo alcuni risultati sugli operatori ellittici, in particolare il principio del massimo.
\begin{deph}
Un operatore differenziale lineare del secondo ordine $$D:C^{2}(\Omega,\R)\to C(\Omega,\R)$$ dove $\Omega$ è un aperto in $\R^m$ è detto \textbf{ellittico} se ha la forma
\begin{gather*}
 D(f)=a^{ij}(x) \partial_i \partial_j f +b^i(x)\partial_i f + c(x) f
\end{gather*}
dove le funzioni $a$ e $b$ sono lisce in $\Omega$ e la matrice $a^{ij}$ è definita positiva (quindi ha tutti gli autovalori strettamente maggiori di $0$).\\
L'operatore $D$ è detto \textbf{strettamente ellittico} se esiste un numero positivo $\lambda>0$ tale che per ogni vettore $v$ e per ogni $x\in \Omega$
\begin{gather*}
 a^{ij}(x) v_i v_j\geq \lambda \sum_i \abs{v_i}^2
\end{gather*}
o equivalentemente se l'autovalore minimo della matrice $a^{ij}$ è limitato dal basso sull'insieme $\Omega$.\\
L'operatore $D$ è detto \textbf{uniformemente ellittico} su $\Omega$ se il rapporto tra l'autovalore massimo e l'autovalore minimo di $a^{ij}(x)$ è limitato (indipendentemente da $x$).\\
La definizione di operatore ellittico può essere estesa facilmente a operatori su varietà Riemanniane chiedendo semplicemente che ogni loro rappresentazione locale abbia le caratteristiche descritte sopra.
\end{deph}
La teoria di questi operatori viene sviluppata in maniera esaustiva su \cite{23}, testo dal quale estraiamo solo i risultati che serviranno in seguito per questa tesi.\\
Una proprietà interessante di questi operatori è che se i coefficienti $a^{ij}, b^i, c$ sono funzioni lisce e una funzione $C^2$ soddisfa $Du=f$ con $f\in C^\infty(\Omega)$, allora automaticamente la funzione $u\in C^{\infty}(\Omega)$. Anzi si può dimostrare che questo continua a valere anche se $Du=f$ solo nel senso delle distrubuzioni.\\
Un'altra proprietà che sfrutteremo molto in questo lavoro è il principio del massimo, che permette di controllare il valore di una funzione con i suoi valori al bordo dell'insieme di definizione.
\begin{prop}[Principio del massimo]\label{prop_max1}
 Sia $D$ un operatore ellittico in un dominio relativamente compatto $\Omega$. Se la funzione $u:\Omega\to \R$ soddisfa:
\begin{enumerate}
 \item $u\in C^{2}(\Omega)\cap C^0(\overline{\Omega})$
 \item $Du=0$
\end{enumerate}
allora il massimo e il minimo di $u$ su $\overline{\Omega}$ sono raggiunti sul bordo $\partial \Omega$. Cioè:
\begin{gather*}
 \sup_{x\in \Omega} u(x)=\sup_{x\in \partial\Omega} u(x) \ \ \
 \inf_{x\in \Omega} u(x)=\inf_{x\in \partial\Omega} u(x) 
\end{gather*}
Inoltre se $u$ non è continua su $\overline{\Omega}$, la conclusione può essere sostituita da:
\begin{gather*}
  \sup_{x\in \Omega} u(x)=\limsup_{x\to \partial \Omega} u(x) \ \ \
 \inf_{x\in \Omega} u(x)=\liminf_{x\to \partial\Omega} u(x) 
\end{gather*}
dove con $\limsup_{x\to \partial \Omega} u(x)$ intendiamo il limite di $\sup_{x\in K_n^C} u(x)$ quando $K_n$ è un'esaustione di $\Omega$.
\end{prop}
Rimandiamo al teorema 3.1 pagina 31 di \cite{23} per la dimostrazione di questo teorema.\\
Vale un principio simile anche se $\Omega$ non è relativamente compatto, infatti:
\begin{prop}\label{prop_max+}
Sia $D$ un operatore ellittico in un dominio $\Omega$, e sia $u$ una funzione $u\in C^{2}(\Omega)\cap C^0(\overline{\Omega})$. Allora
\begin{gather*}
  \sup_{x\in \Omega} u(x)\leq\limsup_{x\to \infty} u(x) \ \ \
 \inf_{x\in \Omega} u(x)\geq\liminf_{x\to \infty} u(x) 
\end{gather*}
dove con $\limsup_{x\to \infty} u(x)$ intendiamo il limite di $\sup_{x\in K_n^C} u(x)$ quando $K_n$ è un'esaustione di $\Omega$.
\begin{proof}
 Supponiamo per assurdo che $\sup_{x\in \Omega} u(x)>\limsup_{x\to \infty} u(x)\equiv L$, e consideriamo l'insieme $A=u^{-1}(L,\infty)$, aperto non vuoto per ipotesi. Dato che $\limsup_{x\to \infty} u(x)= L$, $A$ è relativamente compatto, e applicando il precedente principio a questo insieme si ottiene che $u\equiv L$ in $A$, assurdo.
\end{proof}

\end{prop}

Una forma leggermente più forte del principio del massimo è la seguente:
\begin{prop}[Principio del massimo]\label{prop_max2}
 Sia $D$ un operatore uniformemente ellittico con $c=0$. Se una funzione $u$ soddisfa $Du=0$ assume il suo massimo in un punto interno a $\Omega$, allora è constante
\end{prop}
Rimandiamo al teorema 3.5 pagina 34 di \cite{23} per la dimostrazione di questo teorema.\\
Il principio del massimo può essere espresso anche in forme più sofisticate, alcune delle quali verranno esposte nel seguito della tesi.\\
Osserviamo subito che questi principi possono essere applicati nel caso dell'operatore laplaciano su varietà Riemanniane. Infatti la rappresentazione locale del laplaciano (vedi \ref{eq_lap}) dimostra che questo è un operatore ellittico \footnote{ricordiamo che per definizione di metrica la matrice $g^{ij}$ è definita postiva sulla varietà $R$}, quindi vale il corollario:
\begin{prop}\label{prop_maxq}
 Sia $u$ una funzione armonica in $\Omega$ dominio relativamente compatto in $R$ e continua fino al bordo. Allora il massimo e il minimo della funzione sono assunti sul bordo $\partial \Omega$.
\end{prop}
Inoltre è facile verificare che su ogni insieme compatto $K\Subset R$ l'operatore $\Delta$ è uniformemente ellittico, quindi:
\begin{prop}\label{prop_max3}
 Sia $u$ una funzione armonica in $\Omega$ dominio relativamente compatto in $R$. Se $u$ assume il suo massimo in un punto interno a $\Omega$, allora $u$ è costante su $\Omega$.
\end{prop}
Ovviamente qualunque funzione armonica su tutta la varietà $R$ è armonica su ogni dominio relativamente compatto, quindi per una funzione di questo tipo che non sia costante vale che:
\begin{gather*}
 \inf_{x\in \partial \Omega} u(x) < u(p)< \sup_{x\in \partial \Omega} u(x)
\end{gather*}
per ogni punto $p\in \Omega$.
\subsection{Stime sul gradiente}
Uno strumento importantissimo nello studio delle funzioni armoniche sono le \textit{stime sul gradiente}, cioé stime sul modulo del gradiente di una funzione armonica positiva. Riportiamo solo il risultato, la cui dimostrazione può essere trovata su \cite{26} \footnote{teorema 3.1 pagina 17}
\begin{prop}\label{prop_gradest}
 Sia $R$ una varietà Riemanniana completa con $dim(R)\equiv m\geq 2$, e sia $B_{2r}(x_0)$ la bolla geodetica di raggio $2r$ centrata in $x_0$. Supponiamo che $u$ sia una funzione armonica positiva su $B_{2r}(x_0)$, e sia $Ric(M)\geq -(m-1)K$ su $B_{2r}$ \footnote{essendo la curvatura di Ricci una funzione continua su $R$, su ogni insieme compatto assume un minimo finito, quindi la costante $K$ si può sempre trovare} dove $K\geq 0$ è una costante. Allora:
\begin{gather*}
 \frac{\abs{\nabla u}}{u}\leq C_m \frac{1+r\sqrt K}{r}
\end{gather*}
sull'insieme $B_r(x_0)$, dove $C_m$ è una costante che dipende solo dalla dimensione $m$ della varietà.
\end{prop}
Osserviamo che sul testo \cite{26}, il teorema è enunciato in una forma diversa, si richiede infatti che il limite inferiore sulla curvatura valga su tutta la varietà $R$. Dalla dimostrazione però è evidente che questa ipotesi può essere rilassata.
\subsection{Disuguaglianza, funzione e principio di Harnack}
Un'altra proprietà che vale per le funzioni armoniche su varietà è la disuguaglianza di Harnack:
\begin{prop}[Disuguaglianza di Harnack]\label{prop_harnack}
 Dato un dominio $\Omega$ e un insieme $K\Subset \Omega$, esiste una costante $\Lambda$ \footnote{chiamata \textbf{costante di Harnack}} che dipende solo da $\Omega$ e $K$ tale che per ogni funzione armonica positiva $u$ su $\Omega$ vale che:
\begin{gather*}
 \sup_{x\in K} u(x) \leq \Lambda \inf_{x\in K} u(x)
\end{gather*}
\end{prop}
Questa disuguaglianza è una conseguenza del corollario 8.21 pagina 189 di \cite{23}.\\
Grazie a questa disuguaglianza siamo in grado di provare il \textit{principio di Harnack}, che riguarda successioni di funzioni armoniche positive su varietà
\begin{prop}[Principio di Harnack]\label{prop_harnackpri}
 Sia $u_m$ una successione crescente di funzioni armoniche positive su $\Omega$ dominio in $R$. Allora o $u_m$ diverge localmente uniformemente, o converge localmente uniformemente in $\Omega$.\\
Se la successione di funzioni $u_m$ è uniformemente limitata (non necessariamente positiva o crescente), allora esiste una sua sottosuccessione che converge localmente uniformemente su $\Omega$.
\end{prop}
Questo principio di può trovare su \cite{24} pagina 49 \footnote{in realtà il testo \cite{24} tratta il caso del laplaciano standard in $\R^m$, ma per molte dimostrazioni le tecniche usate si basano su proprietà (come la disuguaglianza di Harnack) che valgono anche nel caso di laplaciano su varietà, quindi possono essere facilmente estese}.\\
Questo principio è molto significativo anche perché vale che
\begin{prop}\label{prop_conv_lu}
 Se una successione di funzioni armoniche $u_n$ converge localmente uniformemente a una funzione $u$, allora $u$ è armonica.
\end{prop}
Oltre alla disuguaglianza di Harnack e alla relativa costante, possiamo definire una \textbf{funzione di Harnack} in questo modo:
\begin{deph}
 Su una varietà Riemanniana $R$ dato un aperto $\Omega$, possiamo definire una funzione $k:\Omega\times \Omega \to [1,\infty)$ nel seguente modo:
\begin{gather*}
 k(x,y)\equiv \sup\{c \ t.c. \ c^{-1}u(x)\leq u(y)\leq cu(x) \ \ \forall \ u\in HP(\Omega)\}
\end{gather*}
dove $HP(\Omega)$ è l'insieme delle funzioni armoniche positive su $\Omega$.\\
Chiamiamo questa funzione \textbf{funzione di Harnack}.
\end{deph}
Grazie all'esistenza della costante di Harnack $\Lambda(K,\Omega)$, sappiamo che la funzione $k$ è ben definita su $\Omega$, infatti dati due punti $x,y\in \Omega$, se consideriamo un compatto $K$ che li contiene, abbiamo che:
\begin{gather*}
 k(x,y)\leq \Lambda(K)
\end{gather*}
però vale anche un'altra importante proprietà di questa funzione:
\begin{prop}\label{prop_harnackfun}
 Per ogni $\Omega$ dominio in $R$, e per ogni $x\in \Omega$, vale che:
\begin{gather*}
 \lim_{y\to x, \ y\in \Omega} k(x,y)=1
\end{gather*}
\begin{proof}
 Questa dimostrazione è un'applicazione delle stime sul gradiente \ref{prop_gradest}. Consideriamo una bolla geodetica $B_{2r}(x)\subset \Omega$, e sia $K$ la costante descritta nella proposizione \ref{prop_gradest}. Sia $y\in B_r(x)$ tale che $d(x,y)\equiv d$, sia $\gamma:[0,d]\to R$ la geodetica che unisce $x$ e $y$ e per una qualsiasi funzione $u$ armonica positiva su $\Omega$ definiamo la funzione
\begin{gather*}
 \phi:[0,d]\to R \ \ \ \ \ \phi(t)\equiv \log\circ u\circ \gamma (t)=\log(u(\gamma(t)))
\end{gather*}
Sappiamo che
\begin{gather*}
 \phi(d)-\phi(0)=\log\ton{\frac{u(y)}{u(x)}}
\end{gather*}
e grazie alle stime sul gradiente possiamo osservare che:
\begin{gather*}
 \abs{\phi(d)-\phi(0)}=\abs{\int_0^d \frac{d\phi}{dt} (s) ds }\leq \int_0^d \abs{\frac{d\phi}{dt}(s)} ds=\int_0^d \abs{\ps{\nabla \log(u)}{\frac{d\gamma}{dt}}} ds=\\
= \int_{0}^d \frac{\abs{\nabla u}}{ u} ds \leq C_m \frac{1+r\sqrt K}{r}d
\end{gather*}
Questo significa che per qualunque funzione armonica positiva $u$:
\begin{gather*}
 \frac{u(y)}{u(x)}\leq \exp\ton{C_m \frac{1+r\sqrt K}{r}d(x,y)}
\end{gather*}
e anche:
\begin{gather*}
 \frac{u(x)}{u(y)}\leq \exp\ton{C_m \frac{1+r\sqrt K}{r}d(x,y)}
\end{gather*}
quindi questa quantità tende a $1$ se $d(x,y)\to 0$.
\end{proof}

\end{prop}

\subsection{Funzioni di Green}\label{subsec_green1}
\begin{deph}\label{deph_fgreen}
 Una funzione di Green per un insieme aperto $\Omega\subset R$ rispetto a un punto $p$ è una funzione $G\in C^\infty (\Omega\times\Omega\setminus D; \R)$, dove  $D=\{(p,p)\ t.c. \ p\in \Omega\}$ è la diagonale di $\Omega \times \Omega$, tale che:
\begin{enumerate}
 \item $G$ è strettamente positiva su $\Omega$
 \item $G$ è simmetrica, cioè $G(p,q)=G(q,p)$
 \item Fissato $p\in \Omega$, la funzione $G(p,q)$ è armonica rispetto a $q$ sull'insieme $\Omega\setminus \{p\}$ e superarmonica su tutto $\Omega$.
 \item $G$ soddisfa la condizione di Dirichlet al bordo, cioè per ogni $p\in \Omega$, $G(p,q)=0$ per ogni $q\in \partial \Omega$ \footnote{il valore di $G$ sul bordo è inteso come il limite per $q_n\to q$ dove $q_n\in \Omega$}
 \item $G$ è soluzione fondamentale dell'operatore $\Delta$, cioé per ogni funzione liscia $f$ a supporto compatto in $\Omega$:
\begin{gather*}
 \Delta_x \int_{\Omega} G(x,y)f(y) dy =\int_{\Omega} G(x,y)\Delta_y(f)(y) dy =-f(x)
\end{gather*}
questo significa che nel senso delle distribuzioni $\Delta_y G(x,y)=-\delta_x$
\item Il flusso di $G_{\Omega}(\ast,p)$ attraverso il bordo di un'insieme regolare $K\Subset \Omega$ con $p\not \in \partial K$ vale:
\begin{gather*}
 \int_{\partial K}\ast d G(\cdot,p)=\begin{cases}
                                   -1 & se \ p\in K\\
0 & se \ p\not \in K
                                  \end{cases}
\end{gather*}

\item La funzione $G$ ha un comportamento asintotico della forma:
\begin{gather*}
 G(x,y)\sim C(m) \begin{cases}
                  -log(d(x,y)) & m=2\\
		  d(x,y)^{m-2} & m\geq 3
                 \end{cases}
\end{gather*}
quando $d(x,y)\to 0$. La costante $C(m)$ dipende solo dalla dimensione della varietà e può essere determinata sfruttando la condizione (5).
\end{enumerate}
\end{deph}
Segnaliamo che per domini $\Omega$ con bordo liscio esiste unica la funzione di Green associata a questo dominio. Un riferimento per questa proposizione è \cite{22} a pagina 165.
\begin{prop}\label{prop_fg1}
 Dato un insieme aperto relativamente compatto $\Omega$ con bordo liscio, esiste unica la funzione di Green $G_{\Omega}$
\end{prop}
Osserviamo che grazie al principio del massimo possiamo dimostrare che
\begin{prop}\label{prop_maxG}
 Sia $\Omega$ un dominio relativamente compatto dal bordo liscio in $R$ e sia $G(\cdot,p)$ la funzione di Green relativa a $\Omega$. Se $K$ è un dominio relativamente compatto tale che
\begin{gather*}
 p\in K\Subset \Omega
\end{gather*}
Allora la funzione $G(\cdot,p)$ rispetto all'insieme $\overline{\Omega}\setminus K$ assume il suo massimo su $\partial K$. Inoltre se $p\in K\Subset K'\Subset \Omega$, $G$ assume il suo massimo rispetto a $\overline{K'}\setminus K$ su $\partial K$.
\begin{proof}
 La dimostrazione segue dal principio del massimo applicato all'insieme $\Omega\setminus K$. Essendo $G(\cdot,p)|_{\partial \Omega}=0$ per definizione, ed essendo $\Delta$ un operatore uniformemente ellittico su $\overline{\Omega}\setminus K$ \footnote{grazie alla compattezza di questo insieme}, vale il principio \ref{prop_max2}, e quindi 
\begin{gather}\label{eq_maxG1}
G(\cdot,p)|_{\Omega \setminus \overline K} < \max_{x\in \partial K} G(x,p) 
\end{gather}
Questo dimostra anche che sull'insieme compatto $\partial K'\Subset \Omega \setminus \overline{K}$ il massimo è strettamente minore del massimo di $G$ su $\partial K$.
\end{proof}
\end{prop}
Se $R$ non è compatta, ha senso chiedersi se esiste una funzione con proprietà simili a quelle descritte definita su tutta la varietà. La risposta a questa domanda è legata alla parabolicità della varietà $R$ che introdurremo in seguito, e si trova nella sezione \ref{sec_green2}

\subsection{Singolarità di funzioni armoniche}
Grazie al principio del massimo e alle funzioni di Green siamo in grado di caratterizzare le singolarità delle funzioni armoniche positive.
\begin{prop}\label{prop_sing1}
 Dato un dominio aperto relativamente compatto $\Omega\subset R$ con bordo liscio e una funzione $$v\in H(\Omega \setminus \{x_0\})\cap C^0(\overline{\Omega}\setminus \{x_0\})$$ positiva \footnote{in realtà è sufficiente che sia limitata dal basso}, detta $G_0$ la funzione di Green $G_{\Omega}(\cdot,x_0)$, se
\begin{gather*}
 \lim_{d(x,x_0)\to 0} \frac{v(x)}{G_0(x)}=0
\end{gather*}
allora la funzione $v$ è estendibile a una funzione armonica su tutta $\Omega$.
\begin{proof}
 Sia $\Omega'$ un dominio aperto relativamente compatto con bordo liscio tale che
\begin{gather*}
 x_0\in \Omega' \subset \overline{\Omega'}\Subset \Omega
\end{gather*}
Sia $[v]$ la soluzione del problema di Dirichlet su $\Omega'$ con valore al bordo $v|_{\partial \Omega'}$. L'obiettivo della dimostrazione è mostrare che $v=[v]$ su $\Omega_0\equiv\Omega'\setminus \{x_0\}$. A questo scopo chiamiamo $\delta=v-[v]$ la funzione differenza. Evidentemente
\begin{gather*}
 \delta \in H(\Omega_0)\cap C^0({\overline \Omega' \setminus \{x_0\}}) \ \ \delta|_{\partial \Omega'}=0
\end{gather*}
Consideriamo la funzione $\phi(x)=\frac{\delta(x)}{G_0(x)}$ definita su $\Omega_0$. Questa funzione è identicamente nulla sul bordo di $\Omega'$ e per $x$ che tende a $x_0$ il suo limite vale $0$ \footnote{grazie al fatto che $G_0$ tende a infinito in $x_0$, $[v]$ è limitata e questa proprietà vale per il rapporto $v/G_0$}. Inoltre questa funzione soddisfa:
\begin{gather*}
\nabla\ton{\frac{\delta}{G_0}}=\frac{\nabla \delta}{G_0} - \frac{\delta \nabla G_0}{G_0^2}\\
 \Delta(\phi)=div(\nabla(\phi))=div\ton{\frac{\nabla \delta}{G_0} - \frac{\delta \nabla G_0}{G_0^2}}=\\
=\frac{\Delta \delta}{G_0} -2\frac{\ps{\nabla \delta}{\nabla G_0}}{G_0^2} -\delta\frac{\Delta G_0}{G_0^2} +2 \frac{\delta \norm{\nabla G_0}^2}{G_0^3}= -2\frac{\ps{\nabla \delta}{\nabla G_0}}{G_0^2}  +2 \frac{\delta \norm{\nabla G_0}^2}{G_0^3}=\\
-\frac{2}{G_0}\ton{\ps{\frac{\nabla \delta}{G_0}}{\nabla G_0}-\ps{\frac{\delta \nabla G_0}{G_0^2}}{\nabla G_0}}=-2 \ps{\nabla \phi}{\frac{\nabla G_0}{G_0}}
\end{gather*}
quindi:
\begin{gather*}
  \Delta(\phi)+2 \ps{\nabla \phi}{\frac{\nabla G_0}{G_0}}=0
\end{gather*}
Dato che l'operatore
\begin{gather*}
 D(\cdot)=\Delta (\cdot) + 2 \ps{\nabla \cdot}{\frac{\nabla G_0}{G_0}}
\end{gather*}
è un'operatore ellittico su $\Omega_0$ (i coefficienti sono funzioni lisce e $\Delta$ è ellittico), allora grazie al principio del massimo \ref{prop_max1} otteniamo che $\phi\equiv 0$ su tutto l'insieme $\Omega_0$, da cui la tesi.
\end{proof}
\end{prop}
Grazie a questa proposizione siamo in grado di dimostrare il corollario:
\begin{prop}\label{prop_sing2}
 Dato un dominio aperto relativamente compatto $\Omega\subset R$ con bordo liscio e una funzione $$v\in H(\Omega \setminus \{x_0\})\cap C^0(\overline{\Omega}\setminus \{x_0\})$$ positiva \footnote{in realtà è sufficiente che sia limitata dal basso}, detta $G_0$ la funzione di Green $G_{\Omega}(\cdot,x_0)$, se
\begin{gather*}
 \lim_{d(x,x_0)\to 0} \frac{v(x)}{G_0(x)}=c
\end{gather*}
allora esiste una funzione armonica $\delta$ su $\Omega$ continua fino al bordo tale che
\begin{gather*}
 v(x)=\delta(x) +cG_0(x)
\end{gather*}
\begin{proof}
 Basta applicare la proposizione precedente alla funzione $$v(x)-cG_0(x)$$
\end{proof}

\end{prop}

\subsection{Principio di Dirichlet}
In questo paragrafo ci occupiamo di introdurre il principio di Dirichlet nella sua forma più standard. Questo principio afferma che le funzioni armoniche sono le funzioni che hanno integrale di Dirichlet minimo in una certa famiglia di funzioni. È possibile rilassare le ipotesi di regolarità su queste funzioni, come riportato nella sezione \ref{sec_dirpri2}.\\
Data una funzione su una varietà $f:R\to \R$ e un insieme misurabile $\Omega$, possiamo definire il suo \textit{integrale di Dirichlet} come:
\begin{gather*}
 D_{\Omega}(f)\equiv\int_{\Omega} \abs{\nabla f}^2 dV
\end{gather*}
e in maniera simile definiamo anche
\begin{gather*}
 D_{\Omega}(f,h)\equiv\int_{\Omega} \ps{\nabla f}{\nabla g} dV
\end{gather*}
Il principio di Dirichlet \footnote{o meglio una versione del principio, in seguito dimostreremo versioni con ipotesi meno restrittive sulla regolarità delle funzioni in gioco} afferma che:
\begin{prop}[Principio di Dirichlet]\label{prop_D1}
 Dato un dominio regolare \footnote{un insieme aperto relativamente compatto con bordo liscio a tratti} $\Omega$ e una funzione continua $h:\partial \Omega \to \R$, per ogni funzione $f:\overline \Omega\to \R$ liscia tale che $f=h$ su $\partial \Omega$:
\begin{gather*}
 D_{\Omega}(f) = D_{\Omega} (u) + D_{\Omega} (f-u)
\end{gather*}
dove $u$ è l'unica soluzione del problema di Dirichlet che ha $h$ come valore al bordo \footnote{cioè $u$ è l'unica funzione continua in $\overline \Omega$, tale che $\Delta u=0$ in $\Omega$ e $u=h$ su $\partial \Omega$}. In particolare $u$ ha integrale di Dirichlet \textbf{minimo} tra tutte le funzioni $f$.
\end{prop}
La dimostrazione di questo principio si può trovare su \cite{16}, sezione 7.1 pagina 173.
\subsection{Funzioni super e subarmoniche}
Oltre alle funzioni armoniche, è possibile definire altre due categorie di funzioni legate all'armonicità: le funzioni superarmoniche e le funzioni subarmoniche. Prima di definire queste funzioni, ricordiamo alcune definizioni preliminari.
\begin{deph}
 Una funzione $f:R\to \R\cup\{+\infty\}$ si dice \textbf{semicontinua inferiormente} se vale una delle seguenti proprietà equivalenti:
\begin{enumerate}
 \item per ogni $x\in R$, $f(x)\leq \liminf_{y\to x} f(y)$
 \item per ogni $a\in R$, $\{x\in \R \ t.c. \ f(x)>a\}$ è un'insieme aperto
 \item per ogni $a\in R$, $\{x\in \R \ t.c. \ f(x)\leq a \}$ è un'insieme chiuso
\end{enumerate}
$f$ si dice \textbf{semicontinua superiormente} se e solo se $-f$ è semicontinua inferiormente, o equivalentemente se e solo se $f:R\to \R\cup\{-\infty\}$ e:
\begin{enumerate}
 \item per ogni $x\in R$, $f(x)\geq \limsup_{y\to x} f(y)$
 \item per ogni $a\in R$, $\{x\in \R \ t.c. \ f(x)<a\}$ è un'insieme aperto
 \item per ogni $a\in R$, $\{x\in \R \ t.c. \ f(x)\geq a \}$ è un'insieme chiuso
\end{enumerate}
\end{deph}
Ricordiamo che:
\begin{prop}
Una successione crescente di funzioni $f_n$ continue converge a una funzione semicontinua inferiormente. Una successione decrescente di funzioni continue converge a una funzione semicontinua superiormente. Il minimo e il massimo tra due (o tra un numero finito) di funzioni semicontinue inferiormente o superiormente è ancora una funzione semicontinua inferiormente o superiormente.
\end{prop}
Il seguente lemma sarà utile per confrontare funzioni sub e superarmoniche. La sua formulazione può essere data con ipotesi meno restrittive, ma per gli scopi della tesi è sufficiente assumere di lavorare su spazi metrici. Questo lemma è tratto del lemma 4.3 pag 171 di \cite{1}.
\begin{lemma}\label{lemma_cfrsemi}
 Sia $X$ uno spazio metrico, e siano $G:X\to \R$ semicontinua superiormente, $g:X\to \R$ semicontinua inferiormente. Se $G(x)<g(x)$ per ogni $x\in X$, allora esiste una funzione continua $\phi:X\to \R$ tale che per ogni $x\in X$:
\begin{gather*}
 G(x)< \phi(x)< g(x)
\end{gather*}
\begin{proof}
Costruiamo la funzione $\phi$ grazie alle partizioni dell'unità. Per ogni numero razionale positivo $r>0$, sia
\begin{gather*}
  U_r=\{x\in X \ t.c. \ G(x)<r\}\cap\{x\in X \ t.c. \ g(x)>r\}
\end{gather*}
per semicontinuità tutti questi insiemi sono aperti, inoltre evidentemente formano un ricoprimento dell'insieme $X$. Dato che $X$ è uno spazio metrico, ogni ricoprimento aperto ammette una partizione dell'unità subordinata a tale ricoprimento. Siano $\lambda_r$ le funzioni di questa partizione, definiamo:
\begin{gather*}
 \phi(x)=\sum_r r\lambda_r
\end{gather*}
Dato che la somma è localmente finita, la funzione $\phi:X\to \R$ è continua su $X$. Inoltre, dato che per ogni $r$ $supp(\lambda_r)\subset U_r$, si ha che per ogni $x\in X$:
\begin{gather*}
 G(x)=G(x)\sum_r \lambda_r < \phi(x) < g(x)\sum_n \lambda_r =g(x)
\end{gather*}
\end{proof}
\end{lemma}
Passiamo ora a trattare le funzioni sub e superarmoniche.\\
Intuitivamente, una funzione superarmonica è una funzione che confrontata con una funzione armonica è maggiore di questa funzione. Se la funzione $v$ è continua, si dice che è superarmonica se e solo se per ogni compatto con bordo liscio $K$, la funzione armonica determinata da $v|_{\partial K}$ è minore della funzione $v$ su tutto $K$. È però possibile rilassare l'ipotesi sulla continuità di $v$. Comunque vogliamo che abbia senso confrontare la funzione $v$ con una funzione armonica $f$ su $K$ tale che $f|_{\partial K}\leq v_{\partial K}$. Se vogliamo che $f\leq v$ su tutto $K$, è necessario chiedere che la funzione $v$ sia semicontinua inferiormente. Quindi definiamo:
\begin{deph}
 Una funzione $v:R\to \R\cup \{+\infty\}$ si dice \textbf{superarmonica} se e solo se è una funzione semicontinua inferiormente e se per ogni insieme compatto $K\Subset R$ con bordo liscio e per ogni funzione $f$ armonica in $K^{\circ}$, continua su $K$ e tale che $f|_{\partial K}\leq v|_{\partial K}$, allora $f\leq v$ su tutto l'insieme $K$.\\
 Una funzione $v:R\to \R\cup \{-\infty\}$ si dice \textbf{subarmonica} se e solo se $-v$ è superarmonica, quindi se e solo se è una funzione semicontinua superiormente e se per ogni insieme compatto $K\Subset R$ con bordo liscio e per ogni funzione $f$ armonica in $K^{\circ}$, continua su $K$ e tale che $f|_{\partial K}\geq v|_{\partial K}$, allora $f\geq v$ su tutto l'insieme $K$.
\end{deph}
Osserviamo subito che la richiesta che $K$ abbia bordo liscio può essere rilassata.
\begin{prop}\label{prop_noliscio}
 Sia $\Omega$ un dominio aperto relativamente compatto in $R$. Se $v$ è superarmonica su $\overline \Omega$ e se $f\in H(\Omega)\cap(C(\overline \Omega))$ è tale che:
$$
v|_{\partial \Omega}\geq f|_{\partial \Omega}
$$
allora $v\geq f$ su tutto l'insieme $\Omega$.
\begin{proof}
 Sia $\epsilon>0$. Data la semicontinuità inferiore di $v-f$, l'insieme $$\Omega_{\epsilon}\equiv (v-f)^{-1}(-\epsilon,\infty)=\{x \ t.c. \ v(x)>f(x)-\epsilon\}$$ è aperto, e per ipotesi contiene $\partial \Omega$. Grazie all'osservazione \ref{oss_an2}, esiste una successione di aperti relativamente compatti con bordo liscio $K_n$ tali che
\begin{gather*}
 K_n\Subset \Omega \ \ \ K_n\Subset K_{n+1} \ \ \ \bigcup_n K_n=\Omega
\end{gather*}
È facile verificare che $\Omega_{\epsilon}^C\subset K_n$ definitivamente. Infatti $\Omega_{\epsilon}^C$ è un compatto di $R$ ed è ricoperto dagli aperti $\{K_n\}$.\\
Dato che $f\in H(\overline K_n)$ e dato che definitivamente in $n$, $$v|_{\partial K_n}<f|_{\partial K_n}-\epsilon$$ poiché il bordo di $K_n$ è liscio, per la definizione di superarmonicità si ha che definitivamente in $n$:
\begin{gather*}
 v|_{K_n}\geq f|_{K_n}-\epsilon \ \ \Rightarrow \ \ v|_{\Omega}\geq f|_{\Omega}-\epsilon
\end{gather*}
per l'arbitrarietà di $\epsilon$ si ottiene la tesi.
\end{proof}

\end{prop}
Le funzioni sub e superarmoniche si possono confrontare tra loro, in particolare si ha che:
\begin{prop}\label{prop_cfrsubsup}
 Sia $K$ un compatto con bordo liscio, e siano $u$ subarmonica su $K$ e $v$ superarmonica su $K$, allora se $u|_{\partial K}\leq v|_{\partial K}$ la disuguaglianza vale su tutto l'insieme $K$.
\begin{proof}
 Grazie al lemma \ref{lemma_cfrsemi}, sappiamo che per ogni $\epsilon>0$ esiste una funzione continua $\phi:\partial K\to \R$ tale che
\begin{gather*}
 u|_{\partial K} < \phi_{\epsilon} < v|_{\partial K} + \epsilon
\end{gather*}
Sia $\Phi_{\epsilon}$ la soluzione del problema di Dirichlet su $K$ con condizioni al bordo $\phi_{\epsilon}$. Allora per definizione di superarmonicità, sappiamo che:
\begin{gather*}
 u\leq \Phi_{\epsilon} \leq v+\epsilon
\end{gather*}
su tutto l'insime $K$. Data l'arbitrarietà di $\epsilon$, otteniamo la tesi.
\end{proof}
\end{prop}
\begin{oss}
 Con una tecnica analoga a quella utilizzata per la dimostrazione della proposizione \ref{prop_noliscio}, si può dimostrare che \ref{prop_cfrsubsup} vale anche se si toglie l'ipotesi di liscezza del bordo di $K$.
\end{oss}
Osserviamo che condizione necessaria e sufficiente affinché una funzione $v$ sia armonica su $R$ è che sia contemporaneamente sub e superarmonica.\\
Una proprietà elementare delle funzioni sub e superarmoniche è che:
\begin{prop}\label{sub_max}
 Il minimo in una famiglia finita di funzioni superarmoniche è superarmonico, e il massimo in una famiglia finita di funzioni subarmoniche è subarmonico.
\end{prop}
Lo spazio delle funzioni sub e superarmoniche è un cono in uno spazio vettoriale, più precisamente:
\begin{prop}\label{prop_cono}
 Combinazioni lineari a coefficienti positivi di funzioni superarmoniche sono superarmoniche.
\begin{proof}
 Se $f$ è superarmonica, è ovvio che per ogni $t\geq 0$, anche $tf$ è una funzione superarmonica. Resta da dimostrare che la somma mantiene la proprietà di superarmonicità. È facile dimostrare che la somma di una funzione armonica e una superarmonica è superarmonica, e analogamente la somma di una funzione armonica e una subarmonica è subarmonica. Consideriamo ora due funzioni $v_1$ e $v_2$ entrambe superarmoniche. Sia $K$ un compatto con bordo liscio in $R$ e sia $u\in H(K^{\circ})\cap C(K)$ tale che:
\begin{gather*}
 u|_{\partial K}\leq (v_1+v_2)|_{\partial K}=v_1 |_{\partial K} + v_2 |_{\partial K}
\end{gather*}
dato che la funzione $u-v_1$ è subarmonica mentre $v_2$ è superarmonica, grazie alla proposizione \ref{prop_cfrsubsup} sappiamo che $u-v_1\leq v_2$ su tutto l'insieme $K$, da cui la tesi. 
\end{proof}
\end{prop}
Evidentemente vale un'affermazione analoga per le funzioni subarmoniche.\\
Riportiamo ora due proposizioni che saranno utili in seguito per dimostrare la superarmonicità dei potenziali di Green.
\begin{prop}\label{prop_sup1}
 Data una funzione $f:R\times R \to \R$ continua per cui per ogni $y\in R$, $f(\cdot,y)$ è superarmonica in $R$, e data una misura di Borel positiva a supporto compatto $K$ con $\mu(K)=\mu(R)<\infty$, la funzione
\begin{gather*}
 F(x)\equiv \int_{R} f(x,y)d\mu(y)
\end{gather*}
è una funzione superarmonica.
\begin{proof}
 Grazie al teorema di convergenza dominata, è facile dimostrare che la funzione $F$ è continua su $R$, quindi in particolare semicontinua inferiormente.\\
Consideriamo ora una successione con indice $n$ di partizioni di $K$ costituite da insiemi $E_k^{(n)}$ tali che il diametro di ogni $E_k^{(n)}$ sia minore di $1/n$, e per ogni $n$ e $k$ sia $y_k^{(n)}$ un punto qualsiasi dell'insieme $E_k^{(n)}$. La successione di funzioni
\begin{gather*}
 F^n(x)=\int_{K} f(x,y_k^{(n)})\chi(E_k^{(n)})(y)d\mu(y)=\sum_k f(x,y_k^{(n)})\mu(E_k^{(n)})
\end{gather*}
è una successione di funzioni superarmoniche grazie alla proposizione \ref{prop_cono}, inoltre converge localmente uniformemente alla funzione $F(x)$. Sia infatti $C$ un qualsiasi insieme compatto in $R$. L'insieme $C\times K$ è compatto in $R\times R$, e quindi la funzione $f(x,y)$ è uniformemente continua su questo insieme. Questo significa che per ogni $\epsilon>0$, esiste $\delta>0$ tale che uniformemente in $x$ si ha:
\begin{gather*}
d(y_1,y_2)<\delta \ \Rightarrow \ \abs{f(x,y_1)-f(x,y_2)}<\epsilon
\end{gather*}
Allora per ogni $\epsilon>0$, se scegliamo $n$ in modo che $1/n\leq \delta$, otteniamo:
\begin{gather*}
 \abs{F(x)-F^n(x)}\leq \int_K \abs{f(x,y)-\sum_k f(x,y_k^{(n)})\chi(E_k^{(n)})(y)}d\mu(y)\leq \mu(K) \epsilon
\end{gather*}
Ora osserviamo che una successione di funzioni superarmoniche che converge localmente uniformemente ha limite superarmonico.\\
Infatti sia $C$ un compatto con bordo liscio in $R$ e $u$ una funzione armonica sulla parte interna di $C$ continua fino al bordo tale che
\begin{gather*}
 u|_{\partial C}(x)\leq F|_{\partial C}(x)
\end{gather*}
Allora per ogni $\epsilon>0$, esiste un $N$ tale che per ogni $n\geq N$:
\begin{gather*}
 F^n|_{\partial C}\geq F|_{\partial C}-\epsilon\geq u|_{\partial C} -\epsilon
\end{gather*}
Per superarmonicità di $F_n$, vale che su tutto $C$:
\begin{gather*}
 F^n(x)\geq u(x)-\epsilon
\end{gather*}
Passando al limite su $n$ e grazie all'arbitrarietà di $\epsilon$, otteniamo la tesi.
\end{proof}
\end{prop}
\begin{prop}\label{prop_sup2}
 Una successione crescente di funzioni continue superarmoniche ha limite superarmonico.
\begin{proof}
 Sia $f_n(x)$ una successione crescente di funzioni continue superarmoniche, e sia $f(x)$ il suo limite \footnote{automaticamente $f$ è una funzione semicontinua inferiormente}. Per dimostrare la superarmonicità di $f$, consideriamo $K$ un compatto con bordo liscio in $R$ e una funzione $u$ armonica sull'interno di $K$ continua fino al bordo con $u|_{\partial K}\leq f|_{\partial K}$. Consideriamo la successione di funzioni
\begin{gather*}
 g_n(x)=\min\{f_n|_{\partial K},u|_{\partial K}\}
\end{gather*}
questa successione è una successione di funzioni continue, crescenti, e ha come limite la funzione continua $u|_{\partial K}$. Grazie al teorema di Dini (riportato dopo questa dimostrazione, teorema \ref{teo_dini}) la convergenza è uniforme, quindi per ogni $\epsilon>0$, esiste $n$ tale che
\begin{gather*}
 g_n|_{\partial K}\geq u|_{\partial K}-\epsilon \ \Rightarrow \ f_n|_{\partial K} \geq u|_{\partial K}-\epsilon
\end{gather*}
data la superarmonicità di $f_n$, si ha che
\begin{gather*}
 f_n(x)\geq u(x)-\epsilon
\end{gather*}
per ogni $x\in K$. Per monotonia si ottiene che per ogni $x\in K$:
\begin{gather*}
 f(x)\geq u(x)-\epsilon
\end{gather*}
data l'arbitrarietà di $\epsilon$, si ottiene la tesi.
\end{proof}

\end{prop}
Riportiamo ora il teorema di Dini. Per ulteriori approfondimenti su questo teorema rimandiamo al testo \cite{28} (il teorema di Dini è il teorema 1.3 pag 381).
\begin{teo}[Teorema di Dini]\label{teo_dini}
Sia $X$ uno spazio topologico compatto, e sia $f_n$ una successione crescente di funzioni $f_n:X\to \R$ continue. Se $f_n$ converge puntualmente a una funzione $f:X\to \R$ continua, allora la convergenza è uniforme.
\begin{proof}
Fissato $\epsilon>0$,  consideriamo $X_n$ gli insiemi aperti
\begin{gather*}
 X_n\equiv \{x\in X \ t.c. \ f(x)-f_n(x)<\epsilon\}
\end{gather*}
Dato che $f_n$ è una successione crescente, $X_n\subset X_{n+1}$, e vista la convergenza puntuale di $f$, sappiamo che $X=\cup_n X_n$. Per compattezza di $X$, esiste un indice $\bar n$ per cui $X=X_{\bar n}$, cioè per ogni $x\in X$ e per ogni $n\geq \bar n$ si ha che
\begin{gather*}
 f_n(x)>f(x)-\epsilon
\end{gather*}
quindi la convergenza è uniforme.
\end{proof}

\end{teo}

Concludiamo questo paragrafo con una proposizione che sarà spesso utilizzata nel seguito.
\begin{prop}\label{prop_superH}
Dato $K$ compatto con bordo liscio in $R$, se $f$ è armonica in $R\setminus K$, costante sull'insieme $K$, minore o uguale alla costante fuori da $K$ e continua su $R$, allora è superarmonica.
\begin{proof}
Senza perdita di generalità, supponiamo che $f\equiv 1$ sull'insieme $K$.\\
Osserviamo che se esiste $v$ superarmonica su $R$ tale che $v|_{R\setminus K}=f|_{R\setminus K}$, questa proposizione è conseguenza del fatto che il minimo tra funzioni superarmoniche è superarmonico. Non è necessario però richiedere che esista una tale funzione.\\
Sia $C$ un insieme compatto con bordo liscio in $R$, e $u$ una funzione armonica sulla parte interna di $C$ continua fino al bordo la cui restrizione al bordo sia minore o uguale alla restrizione di $f$. Se $C\subset R\setminus K$ o $C\subset K$ non c'è niente da dimostrare. Negli altri casi, sia $C_1\equiv C\cap K$ e $C_2\equiv C\cap K^C$. Grazie al principio del massimo, sappiamo che $u\leq 1$ su tutto l'insieme $C$, quindi in particolare su $C_1$. Sempre con il principio del massimo (la forma descritta in \ref{prop_maxq}), confrontando $u$ e $f$ sull'insieme $\overline {C_2}$ abbiamo la tesi. 
\end{proof}

\end{prop}

\section{Algebre di Banach e caratteri}
Questa sezione è dedicata a una breve rassegna sulle algebre di Banach e alcuni risultati che saranno utili nello svolgimento della tesi. Per approfondimenti sull'argomento consigliamo il testo \cite{4}, in particolare il capitolo 10.\\
Iniziamo con il ricordare alcune definizioni di base.
\begin{deph}
 Un'\textbf{algebra associativa} è uno spazio vettoriale su un campo $\mathbb{K}$ (che d'ora in avanti noi assumeremo sempre essere $\R$) con un'operazione di moltiplicazione associativa e distributiva rispetto alla somma. In simboli, uno spazio vettoriale $V$ è un'algebra associativa se è definita una funzione $\cdot:V\times V\to V$ tale che:
\begin{enumerate}
 \item $(x\cdot y)\cdot z=x\cdot(y\cdot z)$
 \item $(x+y)\cdot z=x\cdot z +y\cdot z$
 \item $z\cdot (x+y)=z\cdot z +z\cdot y$
 \item $a(x\cdot y)=(ax)\cdot y=x\cdot (ay)$
\end{enumerate}
Per ogni $x,y,z\in V$ e per ogni $a\in \mathbb{K}$.
\end{deph}
Nel seguito il simbolo di moltiplicazione $\cdot$ sarà sottointeso quando questo non causerà confusione.\\
Un'algebra di Banach è semplicemente un'algebra associativa dotata di una norma compatibile con le operazioni di somma e moltiplicazione e che sia completa rispetto a questa norma.
\begin{deph}
 Un'\textbf{algebra di Banach} è un'algebra associativa su cui è definita un'operazione $\norm{\cdot}:V\to \R^+$ tale che $\norm{x\cdot y}\leq \norm x \norm y$ e che lo spazio normato $(V,\norm \cdot)$ sia completo (o di Banach).
\end{deph}
Si dice che l'algebra di Banach $A$ sia dotata di unità se esiste un elemento $e$ tale che per ogni $x\in A$, $ex=xe=x$ e anche $\norm e =1$. \`E facile dimostrare che se esiste, l'elemento $e$ è unico. D'ora in avanti ci occuperemo solo di algebre di Banach con unità sul campo dei numeri reali. 
\paragraph{Elementi invertibili}
Nelle algebre di Banach con unità (che indicheremo con $A$ in tutta la sezione), ha senso parlare di ``elementi invertibili''. Si dice invertibile un elemento $x$ se esiste $x^{-1}$ tale che $xx^{-1}=x^{-1}x=e$. Grazie all'associatività è facile dimostrare che gli elementi invertibili sono chiusi rispetto alla moltiplicazione, e ovviamente $e$ è un elemento invertibile, il cui inverso è sè stesso. \`E facile dimostrare che una categoria particolare di elementi di $A$ è sempre invertibile:
\begin{prop}\label{prop_algb_inv}
 Se $\norm x <1$, allora $(e-x)$ è invertibile in $A$.
\begin{proof}
 Lo scopo di questa dimostrazione è provare che
\begin{gather*}
 (e-x)^{-1}=\sum_{i=0}^{\infty} x^i
\end{gather*}
dove per convenzione $x^0=e$. Per prima cosa, grazie al fatto che $\norm x < 1$, la serie converge totalmente, quindi converge essendo $A$ uno spazio di Banach. Inoltre:
\begin{gather*}
 \ton{\sum_{i=0}^{\infty} x^i}\cdot (e-x)=\ton{\lim_{n\to \infty} \sum_{i=0}^n x^i}\cdot (e-x)=\lim_{n\to \infty} \left[\ton{\sum_{i=0}^n x^i}\cdot (e-x)\right]=\\
=\lim_{n\to \infty} \left[\sum_{i=0}^n x^i -\sum_{i=1}^{n+1} x^i)   \right]=e-\lim_{n\to \infty} x^{n+1}=e
\end{gather*}
e un ragionamento del tutto analogo vale per $(e-x)\cdot\ton{\sum_{i=0}^{\infty} x^i}$.
\end{proof}
\end{prop}
\paragraph{Caratteri}
Introduciamo ora una classe particolare di funzionali lineari sulle algebre di Banach, i caratteri, o funzionali moltiplicativi. Anche in questa sezione riporteremo solo i risultati che interessano agli scopi della tesi, e per ulteriori approfondimenti rimandiamo a \cite{4}.
\begin{deph}
 Un funzionale lineare $\phi:A\to \R$ si dice moltiplicativo se conserva la moltiplicazione, cioè se $\forall x,y \in A$
\begin{gather*}
 \phi(x\cdot y)=\phi(x)\phi(y)
\end{gather*}
\end{deph}
Le proprietà che risultano evidenti dalla definizione sono che:
\begin{enumerate}
\item $\phi(e)=1$
\item per ogni elemento invertibile $x$, $\phi(x)\neq 0$
\item $\phi(x^{-1})=\ton{\phi(x)}^{-1}$
\end{enumerate}
Grazie all'ultima proprietà e alla proposizione \ref{prop_algb_inv} possiamo dimostrare che qualunque carattere è necessariamente continuo e ha norma 1 \footnote{ricordiamo che la norma di un funzionale può essere definita da $\norm \phi \equiv \sup_{\norm x \leq 1}(\abs{\phi(x)})$, un funzionale è continuo se e solo se ha norma finita}.
\begin{prop}
 Ogni carattere $\phi$ ha norma 1.
\begin{proof}
Dal fatto che $\phi(e)=1$ si ricava facilmente che $\norm \phi \geq 1$. Consideriamo ora un qualunque $x$ con norma minore o uguale a 1, e un qualunque numero reale $\lambda$ con $\abs \lambda > 1$. Dalla proposizione \ref{prop_algb_inv} segue che $e-\lambda^{-1}x$ è un elemento invertibile, quindi
\begin{gather*}
 \phi(e-\lambda^{-1}x)\neq 0 \Longleftrightarrow 1-\frac{\phi(x)}{\lambda}\neq 0 \Longleftrightarrow \phi(x)\neq \lambda
\end{gather*}
Questa considerazione è valida per qualunque numero $\abs{\lambda} > 1$, e quindi si ricava che $\sup_{\norm x \leq 1}\{\abs{\phi(x)}\}\leq 1$, da cui la tesi.
\end{proof}
\end{prop}
\`E utile osservare che:
\begin{prop}\label{prop_charK}
 L'insieme $\mathcal{C}$ dei funzionali lineari \textbf{moltiplicativi} è compatto rispetto alla topologia debole-*.
\begin{proof}
 Grazie al teorema \ref{teo_BA}, è sufficiente dimostare che $\mathcal C$ è chiuso nella topologia debole-*\footnote{ricordiamo che questo insieme è limitato in norma.}.
A questo scopo, consideriamo $\phi\in \overline{\mathcal C}$, e verifichiamo se $\phi (xy)=\phi(x)\phi(y)$. Dato che $\phi\in \overline {\mathcal C}$, per ogni $\epsilon>0$, esiste $\tilde \phi \in \mathcal C$ tale che:
\begin{gather*}
 \tilde \phi \in V(\phi,\epsilon,x,y,xy)\equiv \{\psi \in A^* \ t.c. \ \abs{\psi(t)-\phi(t)}<\epsilon \ \ t=x,y,xy\}
\end{gather*}
quindi per ogni $\epsilon>0$:
\begin{gather*}
 \abs{\phi(xy)-\phi(x)\phi(y)}\leq\\
\leq\abs{\phi(xy)-\tilde \phi(xy)}+\abs{\tilde\phi(xy)-\tilde\phi(x)\tilde\phi(y)}+\abs{\tilde\phi(x)\tilde\phi(y)-\phi(x)\phi(y)}\leq\\
\leq 2\epsilon + \abs{\tilde \phi(x)\tilde\phi(y)-\tilde\phi(x)\phi(y)}+\abs{\tilde\phi(x)\phi(y)-\phi(x)\phi(y)}\leq\\
\leq 2\epsilon+(\abs{\tilde\phi(x)}+\abs{\phi(y)})\epsilon \leq 2\epsilon + (\abs{\phi(x)}+\abs{\phi(y)}+\epsilon)\epsilon
\end{gather*}
Data l'arbitrarietà di $\epsilon$, si ottiene la tesi.
\end{proof}

\end{prop}

\section{Problema di Dirichlet}\label{sec_dir}
Lo scopo di questa sezione è illustrare alcuni risultati per risolvere il problema di Dirichlet su domini limitati su varietà riemanniane. Daremo la definizione di \textit{dominio regolare} per il problema di Dirichlet e utilizzeremo il metodo di Perron e le barriere per caratterizzare questi domini, in seguito dimostreremo che alcuni domini particolari sulle varietà sono regolari, in particolare i domini limitati con bordo liscio e i domini della forma $\Omega\setminus K$, dove $\Omega$ è un dominio limitato con bordo liscio e $K$ una sottovarietà di codimensione $1$ con bordo liscio contenuta in $\Omega$. Il lemma \ref{lemma_fiko} (che citeremo senza dimostrazione) sarà lo strumento principale per i nostri scopi, in quanto lega la solubilità del problema di Dirichlet per un operatore ellittico abbastanza generico al problema di Dirichlet del più noto e studiato laplaciano standard in $\R^n$.\\
In tutta la sezione $\Omega$ indicherà un dominio (cioè un insieme aperto connesso relativamente compatto) in $R$ varietà riemanniana o in $\R^n$, mentre $L:C^2(R,\R)\to C(R,\R)$ indicherà un operatore differenziale del II ordine strettamente ellittico (il laplaciano su varietà riemanniane soddisfa queste ipotesi come illustrato nella sezione \ref{sec_max}). Inoltre considereremo solo funzioni sub e superarmoniche \textit{continue} fino alla chiusura dell'insieme di definizione.\\
Prima di cominciare diamo una definizione preliminare.
\begin{deph}
 Un dominio $\Omega\subset R$ ha bordo $C^k$ con $0\leq k\leq \infty$ se e solo se per ogni punto $p\in \partial \Omega$, esiste un intorno $V(p)$ e una funzione $f\in C^k(V,\R)$ tale che:
\begin{gather*}
 \Omega \cap V = f^{-1}(\infty,0) \ \ \ \partial \Omega \cap V =f^{-1}(0)
\end{gather*}
Inoltre $0$ deve essere un valore regolare di $f$.
\end{deph}
Osserviamo che non è sufficiente per un dominio in $R$ avere il bordo costutuito da sottovarietà regolari per essere considerato di bordo liscio. Ad esempio consideriamo l'insieme $A=B(0,2)\setminus \{(x,0) \ -1\leq x \leq 1\}\subset \R^2$. Il bordo di questo dominio è costituito dal bordo della bolla e dal segmento $\{(x,0) \ -1\leq x \leq 1\}$, e sebbene per ogni punto del segmento $A=\{(x,0) \ -1< x < 1\}$ esiste un intorno $V$ e una funzione $f:V\to \R$ tale che $A\cap V=f^{-1}(0)$, non è possibile fare in modo che $\Omega\cap V =f^{-1}(\infty,0)$.
\subsection{Metodo di Perron}
Per prima cosa definiamo il problema di Dirichlet.
\begin{deph}
 Dati $\Omega\in R$ e $L$ con le caratteristiche descritte appena sopra, $\phi:\partial \Omega \to \R$ funzione continua, diciamo che $\Phi:\overline{\Omega}\to \R$ è soluzione del problema di Dirichlet se:
\begin{gather}\label{eq_dir}
 \Phi\in C^2(\Omega)\cap C(\overline \Omega) \ \ \ e \ \ \ L(\Phi)=0 \ \ su \ \ \Omega \ \ \ e \ \ \ \Phi|_{\partial \Omega}=\phi
\end{gather}
\end{deph}
La solubilità del problema di Dirichlet è fortemente legata alla regolarità del dominio $\Omega$. Ad esempio, se $\Omega=B(0,1)\setminus \{0\}\subset \R^n$ \footnote{$B(0,1)$ indica la bolla di raggio $1$ in $\R^n$}, non esistono soluzioni del problema di Dirichlet ``classico'' (cioè con $L=\Delta$) uguali a $0$ sul bordo di $B(0,1)$ e diverse da $0$ ma limitate nell'origine (vedi esercizi 15 e 16 cap 3 pag 57 di \cite{24}, e per un altro esempio vedi es 25 cap 2 pag 54 dello stesso libro).
\begin{deph}
 Un dominio $\Omega\in R$ si dice essere \textbf{regolare} rispetto al problema di Dirichlet se il più generico di questi problemi ha un'unica soluzione.
\end{deph}
Il problema di Dirichlet per domini in $\R^n$ rispetto all'operatore laplaciano standard è un argomento molto studiato in matematica e con molti risultati nella letteratura. Il seguente lemma (che riportiamo senza dimostrazione) lega la solubilità del problema di Dirichlet per un operatore ellittico abbastanza generico a quella del laplaciano standard, il che semplifica molto il problema.
\begin{lemma}\label{lemma_fiko}
 Dato un dominio $\Omega\in \R^n$ con $\Omega\subset \overline \Omega \Subset \Omega_1$ e dato $L$ operatore ellittico della forma:
\begin{gather*}
 L=a^{ij}(x)\partial_i\partial_j+b^i(x)\partial_i
\end{gather*}
con coefficienti $a^{ij}$ e $b^i$ localmente lipschitziani in $\Omega_1$ e per il quale esiste una funzione di Green su $\Omega_1$, allora il generico problema di Dirichlet su $\Omega$ è risolubile se e solo se lo è anche il generico problema di Dirichlet legato all'operatore laplaciano.
\begin{proof}
 Questo lemma è il corollario al teorema 36.3 in \cite{30} (ultimo risultato presente nell'articolo).
\end{proof}
\end{lemma}
Se consideriamo $\Omega\subset R$ dominio contenuto in una carta locale, il laplaciano sulla varietà $R$ soddisfa tutte le ipotesi del teorema. Questo implica in particolare che:
\begin{prop}
 Per ogni dominio $\Omega\subset R$ e per ogni $p\in \Omega$, esiste un intorno aperto $V(p)\subset \overline{V(p)}\Subset \Omega$ tale che il generico problema di Dirichlet rispetto al laplaciano su varietà è risolubile su $V(p)$.
\begin{proof}
 La dimostrazione è un semplice corollario del teorema precedente. Per il laplaciano standard, il problema di Dirichlet su ogni bolla $B(x_0,r)\subset \R^n$ è risolubile (vedi ad esempio teorema 1.17 pag 13 di \cite{24}), quindi se consideriamo un qualsiasi insieme $A$ aperto intorno di $p$ la cui chiusura è contenuta in $\Omega\cap U$ dove $(U,\phi)$ è una carta locale di $R$ e tale che $\phi(A)=B(x_0,r)$, il problema di Dirichlet relativo al laplaciano su varietà è risolubile su $A$.
\end{proof}
\end{prop}
Per caratterizzare alcuni dei domini su cui è sempre possibile risolvere il problema di Dirichlet, utilizzeremo il \textbf{metodo di Perron}.
\begin{deph}\label{deph_perron}
 Dati $\Omega$, $L$, $\phi$ come sopra, indichiamo con $P[\phi]$ la funzione $P[\phi]:\overline \Omega \to \R$ definita da
\begin{gather*}
 P[\phi](x)=\sup_{u\in S_\phi} u(x)
\end{gather*}
dove $S_\phi$ è l'insieme delle funzioni subarmoniche $u:\overline \Omega\to \R$ continue su $\overline \Omega$ tali che
\begin{gather*}
 u|_{\partial \Omega}\leq \phi
\end{gather*}
Osserviamo che per la compattezza di $\partial \Omega$, la funzione $\phi$ ha minimo $m$ e massimo $M$ finiti, e dato che le tutte le funzioni costanti soddisfano $\Delta_R (c)=0$, l'insieme $S_{\phi}$ non è vuoto in quanto contiene la funzione costante uguale a $m$, e tutte le funzioni in $S_{\phi}$, in quanto funzioni subarmoniche sono limitate da $M$, quindi $P[\phi](x)\leq M$ per ogni $x\in \overline\Omega$.
\end{deph}
Il metodo di Perron è illustrato ad esempio nel paragrafo 11.3 pag 226 di \cite{24} o nel paragrafo 2.8 pag 23 di \cite{23} per risolvere il problema di Dirichlet per il laplaciano standard su domini qualsiasi. I cardini essenziali di questo metodo però sono il principio del massimo e la possibilità di risolvere il problema di Dirichlet sulla bolla. Come abbiamo visto, questi principi valgono anche per il laplaciano su varietà, non è quindi difficile immaginare che anche per questo problema il metodo di Perron fornisca una soluzione adeguata.\\
Verifichiamo ora sotto quali condizioni $P[\phi]$ risolve \ref{eq_dir}. Per prima cosa dimostriamo che indipendentemente dal dominio $\Omega$, $P[\phi]$ è una funzione armonica.
\begin{prop}
 $P[\phi]$ è armonica sull'insieme $\Omega$.
\begin{proof}
Per dimostrare che la funzione $P[\phi]$ è armonica, dimostriamo che per ogni punto $p\in \Omega$, esiste un intorno $V=V(p)$ su cui la funzione è armonica.\\ Consideriamo a questo scopo $V(p)$ un aperto con chiusura contenuta in $\Omega\cap U$, dove $U$ è un intorno coordinato di $R$ tale che la rappresentazione in carte locali di $V(p)$ sia una bolla. Grazie a \ref{lemma_fiko}, sappiamo che è possibile risolvere il problema di Dirichlet relativo a $\Delta_R$ in carte locali su $V(p)$. Per ogni funzione continua $v:\Omega\to \R$, possiamo definire il suo \textit{sollevamento armonico} come la funzione data da:
\begin{gather*}
 \bar v (x)=\begin{cases}
             v(x) & se \ x\in V^C\\
	    D[v](x) & se \ x\in \overline V
            \end{cases}
\end{gather*}
dove $D[v]$ è la soluzione del problema di Dirichlet relativo a $\Delta_R$ su $V(p)$ con $v|_{\partial V}$ come condizione al bordo. Osserviamo che questa funzione è continua su $\Omega$ e se $v$ è subarmonica, allora $\bar v$ mantiene questa proprietà in quanto massimo di due funzioni subarmoniche.\\
Sia ora $u_k$ una successione di funzioni in $S_\phi$ tali che $u_k(p)\to P[\phi](p)$. Consideriamo $\bar u_k$ la successione dei sollevamenti armonici di queste funzioni relativamente all'intorno $V(p)$. Dato che $\bar u_k(p)\geq u_k(p)$ e che $\bar u_k \in S_\phi$, si ha che $\bar u_k (p) \to P[\phi](p)$. Data l'uniforme limitatezza delle funzioni $\bar u_k$, per il principio di Harnack \ref{prop_harnackpri} esiste una sottosuccessione che per comodità continueremo a indicare con lo stesso indice che converge localmente uniformemente su $V(p)$ a una funzione armonica. Sia $u=\lim_k \bar u_k$, vogliamo dimostrare che $u|_{V(p)}=P[\phi]|_{V(p)}$.\\
È facile osservare che $u\leq P[\phi]$ su tutto $\Omega$, infatti ogni funzione $\bar u_k\in S_\phi$. Inoltre come osservato in precedenza $u(p)=P[\phi](p)$. Supponiamo per assurdo che esista $q\in V(p)$ tale che $u(q)< P[\phi](q)$. Allora per definizione per ogni $\epsilon>0$, esiste una funzione $\tilde w \in S_\phi$ tale che
\begin{gather*}
 u(q)<\tilde w(q)\leq P[\phi](q) \ \ \ \tilde w(q)>P[\phi](q)-\epsilon
\end{gather*}
Se definiamo $w_k$ come il sollevamento armonico della funzione $\max\{\tilde w, u_k\}$ rispetto all'insieme $V(p)$, otteniamo come prima una successione di funzioni armoniche limitate che, a patto di passare a una sottosuccessione, converge localmente uniformemente a una funzione $w$ armonica su $V(p)$. Poiché tutte le funzioni $w_k\in S_\phi$, sappiamo che $w\leq P[\phi]$, e quindi in particolare $w(p)\leq P[\phi](p)$. Inoltre per costruzione $u_k\leq w_k$, quindi passando al limite otteniamo che $u\leq w$ sull'insieme $V(p)$. Ma dato che $P[\phi](p)=u(p)\leq w(p)\leq P[\phi](p)$, abbiamo che $w(p)=u(p)$, cioè la funzione armonica $u-w$ assume il suo massimo in $V(p)$ in un punto interno all'insieme, quindi per il principio del massimo \ref{prop_max2} $u-w=0$. Dato che per costruzione $w(q)>P[\phi](q)-\epsilon$, abbiamo che:
\begin{gather*}
 u(q)>P[\phi](q)-\epsilon
\end{gather*}
e l'assurdo segue dall'arbitrarietà di $\epsilon$.
\end{proof}
\end{prop}
Resta da verificare se la funzione $P[\phi]$ è continua su $\overline \Omega$ e se $P[\phi]|_{\partial \Omega}=\phi$. Prima di dare condizioni per verificare queste proprietà, sottolineamo che il problema di Dirichlet è risolvibile se e solo se valgono queste proprietà per $P[\phi]$.
\begin{prop}
 Il problema di Dirichlet \ref{eq_dir} è risolubile se e solo se $P[\phi]$ risulta essere continua su $\overline \Omega$ e $P[\phi]|_{\partial \Omega}=\phi$.
\begin{proof}
 Se le condizioni su $P[\phi]$ sono verificate, automaticamente $P[\phi]$ è l'unica soluzione del problema di Dirichlet. Al contrario, supponiamo che esista $\Phi$ soluzione del problema. Allora $\Phi\in S_{\phi}$, anzi $\Phi=P[\phi]$. Infatti la funzione $\Phi$ è armonica, quindi maggiora tutte le funzioni in $S_\phi$. 
\end{proof}
\end{prop}
Come accennato in precedenza, non sempre il problema di Dirichlet è risolubile. Data l'ultima equivalenza questo implica che sebbene $P[\phi]$ sia sempre una funzione armonica, non sempre soddisfa tutte le condizioni del problema di Dirichlet. Un modo per verificare quando il problema è risolubile è il \textit{criterio delle barriere}.
\begin{deph}
 Dato un punto $p\in \partial \Omega$, diciamo che una funzione $S:\overline \Omega \to \R$ è una \textbf{barriera} per il punto $p$ se $S$ è una funzione superarmonica in $\Omega$, continua in $\overline \Omega$ tale che
\begin{gather*}
 S(p)=0 \ \ \ S|_{\overline \Omega \setminus \{p\}}>0
\end{gather*}
 Se la funzione $\beta$ ha queste caratteristiche ma è definita solamente su $\overline \Omega \cap V(p)$, dove $V(p)$ è un intorno qualsiasi del punto $p$, diciamo che $\beta$ è una \textbf{barriera locale}.
\end{deph}
Per prima cosa osserviamo che l'esistenza di una barriera ``globale'' è equivalente all'esistenza di una barriera locale, infatti:
\begin{prop}
 Se esiste una barriera locale $\beta$ per un punto $p\in \partial \Omega$, allora esiste anche una barriera globale per $p$.
\begin{proof}
La dimostrazione è relativamente facile. Sia $V(p)$ un intorno aperto tale che $\beta$ sia definita su $\overline \Omega \cap V(p)$. Sia $W\subset\overline W\subset V$ un secondo intorno di $p$. Per continuità $\beta$ assume minimo $m$ strettamente positivo sull'insieme $W^C$. Consideriamo la funzione $S$ definita da:
\begin{gather*}
 S(q)=\begin{cases}
       \min\{\beta(q),m\} & se \ q\in \overline \Omega \cap W \\
      m & se \ q\in \overline \Omega \cap W^C
      \end{cases}
\end{gather*}
dato che $S$ è il minimo tra funzioni superarmoniche, è ancora una funzione superarmonica. È facile dimostrare che anche le proprietà richieste a una barriera sono verificate.
\end{proof}
\end{prop}
L'utilità del concetto di barriera è contenuta nella seguente proposizione:
\begin{prop}
 Sia $p\in \partial \Omega$. Se esiste una barriera per $p$ rispetto a $\Omega$, allora
\begin{gather*}
 P[\phi](p)=\lim_{x\to p} P[\phi](x)=\phi(p)
\end{gather*}
\begin{proof}
 Poiché $\phi$ è continua sull'insieme $\partial \Omega$, per ogni $\epsilon>0$ esiste un intorno $U_{\epsilon}(p)$ tale che per ogni $x\in U_\epsilon \cap \partial \Omega$:
\begin{gather*}
 \phi(p)-\epsilon<\phi(x)<\phi(p)+\epsilon
\end{gather*}
Dato che $S$ è una funzione strettamente positiva su $\partial \Omega \setminus U_{\epsilon}$, esiste una costante $c>0$ tale che
\begin{gather*}
 \phi(p)-\epsilon-c S(x)<\phi(x)< \phi(p)+\epsilon+c S(x)
\end{gather*}
per ogni $x \in \partial \Omega \setminus U_{\epsilon}$. Data la positività di $S$, questa relazione vale su tutto l'insieme $\partial \Omega$. Dato che $S$ è superarmonica, la funzione $\phi(p)-\epsilon-c S(x)\in S_\phi$, quindi
\begin{gather}\label{eq_a}
 \phi(p)-\epsilon-c S(x)\in S_\phi\leq P[\phi](x)
\end{gather}
ora consideriamo una funzione $u\in S_\phi$. Vale che:
\begin{gather*}
 u|_{\partial \Omega}\leq \phi <\phi(p)+\epsilon+c S(x) \ \ \Rightarrow \ \ [u-cS]|_{\partial \Omega} \leq \phi(p)+\epsilon
\end{gather*}
Data la subarmonicità della funzione continua $[u-cS]$, questa relazione è valida su tutto l'insieme $\overline \Omega$, e data l'arbitrarietà di $u\in S_\phi$ otteniamo che:
\begin{gather}\label{eq_b}
 P[\phi](x)\leq \phi(p)+\epsilon +cS(x)
\end{gather}
Passando al limite per $x\to p$ nelle relazioni \ref{eq_a} e \ref{eq_b}, otteniamo che per ogni $\epsilon>0$:
\begin{gather*}
 \phi(p)-\epsilon\leq\liminf_{x\to p}P[\phi](x)\leq \limsup_{x\to p}P[\phi](x)\leq \phi(p)+\epsilon
\end{gather*}
data l'arbitrarietà di $\epsilon>0$, si ottiene la tesi.
\end{proof}
\end{prop}
Il criterio delle barriere è molto utile perchè l'esistenza delle barriere è condizione necessaria e sufficiente per risolvere il problema di Dirichlet, infatti:
\begin{prop}
 Il generico problema di Dirichlet su $\Omega \subset R$ è risolubile se e solo se per ogni $p\in \partial \Omega$ esiste una barriera relativa a $p$.
\begin{proof}
 Grazie alla proposizione precedente, se ogni punto del bordo di $\Omega$ ha una barriera, $P[\phi]$ soddisfa tutte le condizioni del problema di Dirichlet, quindi esiste unica la soluzione di tale problema. Supponiamo al contrario che il problema sia sempre risolubile. Allora fissato $p\in \partial \Omega$, sia $\phi_p$ una funzione continua sul bordo tale che $\phi_p(p)=0$, e $\phi_p(x)>0$ per ogni $x\neq p$ (ad esempio, $\phi_p(x)=d(x,p)$). Sia $\Phi_p$ la soluzione del relativo problema di Dirichlet. Questa funzione è una barriera per $p$, infatti è una funzione armonica (quindi anche superarmonica) in $\Omega$, continua in $\overline \Omega$ e uguale a $\phi_p$ su $\partial \Omega$, e grazie al principio del massimo, $\Phi_p(x)>0$ per ogni $x\in \Omega$.
\end{proof}
\end{prop}
Come corollario immediato di questo teorema, dimostriamo che
\begin{prop}\label{prop_dircap}
 Se $\Omega$ e $\Omega'$ sono domini regolari rispetto a $\Delta_R$, allora anche la loro intersezione $\Omega''=\Omega\cap \Omega'$ è un dominio regolare.
\begin{proof}
 La regolarità dei domini è equivalente all'esistenza di barriere per ogni punto del bordo. Dato che
\begin{gather*}
 \partial (\Omega\cap \Omega')\subset \partial \Omega \cup \partial \Omega'
\end{gather*}
ogni punto del bordo di $\Omega''$ appartiene ad almeno uno dei due bordi. Senza perdita di generalità, consideriamo $p\in (\partial\Omega''\cap \partial \Omega)$. Per regolarità di $\Omega$, esiste una barriera $B$ per $p$ rispetto a $\Omega$. La funzione $B$ è superarmonica in $\Omega$, quindi automaticamente anche in $\Omega''$, strettamente positiva su $\overline \Omega\setminus\{p\}$, quindi anche su $\overline {\Omega''}\setminus\{p\}$, e ovviamente $B|_{\overline{\Omega''}}$ è una funzione continua. Questo dimostra che ogni punto di $\partial {\Omega''}$ possiede una barriera, $\Omega''$ è regolare.
\end{proof}
\end{prop}
Riassumendo, un dominio $\Omega\subset R$ è regolare rispetto al laplaciano sulla varietà $R$ se e solo se per ogni punto del bordo di $\Omega$ esiste una barriera, barriera che come abbiamo visto è un concetto locale. Se il dominio $\Omega$ è contenuto in una carta locale, allora grazie al lemma \ref{lemma_fiko}, possiamo utilizzare tutti i risultati noti per il laplaciano standard e concludere che ogni insieme regolare per il laplaciano standard è regolare anche per il laplaciano su varietà (e viceversa). Utilizzando assieme queste due tecniche, dimostreremo che i domini con bordo liscio in $R$ sono regolari, ma anche i domini della forma $\Omega=A\setminus K$, dove $A$ è un dominio con bordo liscio, e $K$ una sottovarietà con bordo regolare di $R$ di codimensione $1$ contenuta in $A$.
\subsection{Domini con bordo liscio}
Dalla letteratura sul laplaciano standard in $\R^n$, sappiamo che una condizione geometrica semplice per garantire la regolarità di un dominio è la \textit{condizione della bolla esterna}.
\begin{prop}[Condizione della bolla esterna]
 Dato un dominio $\Omega\subset \R^n$, $\Omega$ è regolare per il laplaciano standard se per ogni punto di $\Omega$ esiste una bolla $B(x(p),\epsilon(p))$ tale che
\begin{gather*}
 \overline {B(x,\epsilon)}\cap \overline \Omega=\{p\}
\end{gather*}
cioè una bolla esterna al dominio $\Omega$ la cui chiusura interseca il bordo di $\Omega$ solo nel punto $p$.
\begin{proof}
 Grazie all'omogeneità del dominio $\R^n$ e dell'operatore $\Delta$, e soprattutto grazie al grande dettaglio con cui questo problema è stato studiato, la trattazione dell'argomento è particolarmente semplice. Per dimostrare questa proposizione, dimostriamo che per ogni punto $p$ che soddisfa la condizione della bolla esterna, esiste una barriera $B_p$. Ad esempio possiamo considerare la funzione:
\begin{gather*}
 B_p(y)=\begin{cases}
                \log(\norm{x(p)-y})-\log(\epsilon(p)) & se \ n=2\\
		-\frac{1}{\norm{x(p)-y}^{n-2}}+\frac{1}{\epsilon(p)} & se \ n\geq 3
               \end{cases}
\end{gather*}
È facile dimostrare che questa funzione soddisfa tutte le proprietà richieste a una barriera per $\Omega$ (vedi ad esempio teorema 11.13 pag 229 di \cite{24}).
\end{proof}
\end{prop}
Come corollario di questa condizione, è facile dimostrare che
\begin{prop}
 Qualunque dominio $\Omega\subset \R^n$ con bordo $C^2$ è regolare per il problema di Dirichlet.
\begin{proof}
 È sufficiente dimostrare che ogni punto del bordo di questi domini soddisfa la condizione della bolla esterna. Per i dettagli vedi corollario 11.13 di \cite{24}.
\end{proof}
\end{prop}
La condizione della bolla esterna implica in particolare che ogni anello in $\R^n$, cioè ogni insieme della forma
\begin{gather*}
 A(x_0,r_1,r_2)\equiv B(x_0,r_2)\setminus \overline {B(x_0,r_1)}=\{x\in \R^n \ t.c. \ r_1<\norm{x-x_0}<r_2\}
\end{gather*}
è un dominio regolare per il laplaciano standard, e grazie al lemma \ref{lemma_fiko}, ogni anello in un insieme coordinato \footnote{cioè ogni insieme $A\subset \overline A\Subset U\subset R$ tale che $(U,\psi)$ sia un intorno coordinato e $\psi(A)$ sia un anello in $\psi(U)$} è regolare rispetto al laplaciano su $R$.\\
Questo ci permette di dimostrare che
\begin{prop}\label{prop_dir_li}
 Ogni dominio $\Omega\subset R$ con bordo liscio \footnote{è sufficiente con bordo $C^2$} è un dominio regolare per $\Delta_R$.
\begin{proof}
 Dimostriamo che ogni punto del bordo possiede una barriera locale. Sia a questo scopo $p\in \partial \Omega$ e $(U,\psi)$ un intorno coordinato di $p$ tale che
\begin{gather*}
 U\cap \Omega =\psi^{-1}\{(x_1,\cdots,x_{m-1},0)\}
\end{gather*}
cioè sia $\psi$ una carta che manda $U\cap \Omega$ in un semispazio chiuso di $\R^n$. È facile immaginare (vedi disegno \ref{fig_1}) che esiste un anello $A(x_0,r_1,r_2)$ in $\psi(U)$ tale che $\overline{B(x_0,r_1)}\cap \psi(\Omega)=\{\psi(p)\}$.
\begin{figure}[ht!]
\begin{center}
     \includegraphics[width=60mm]{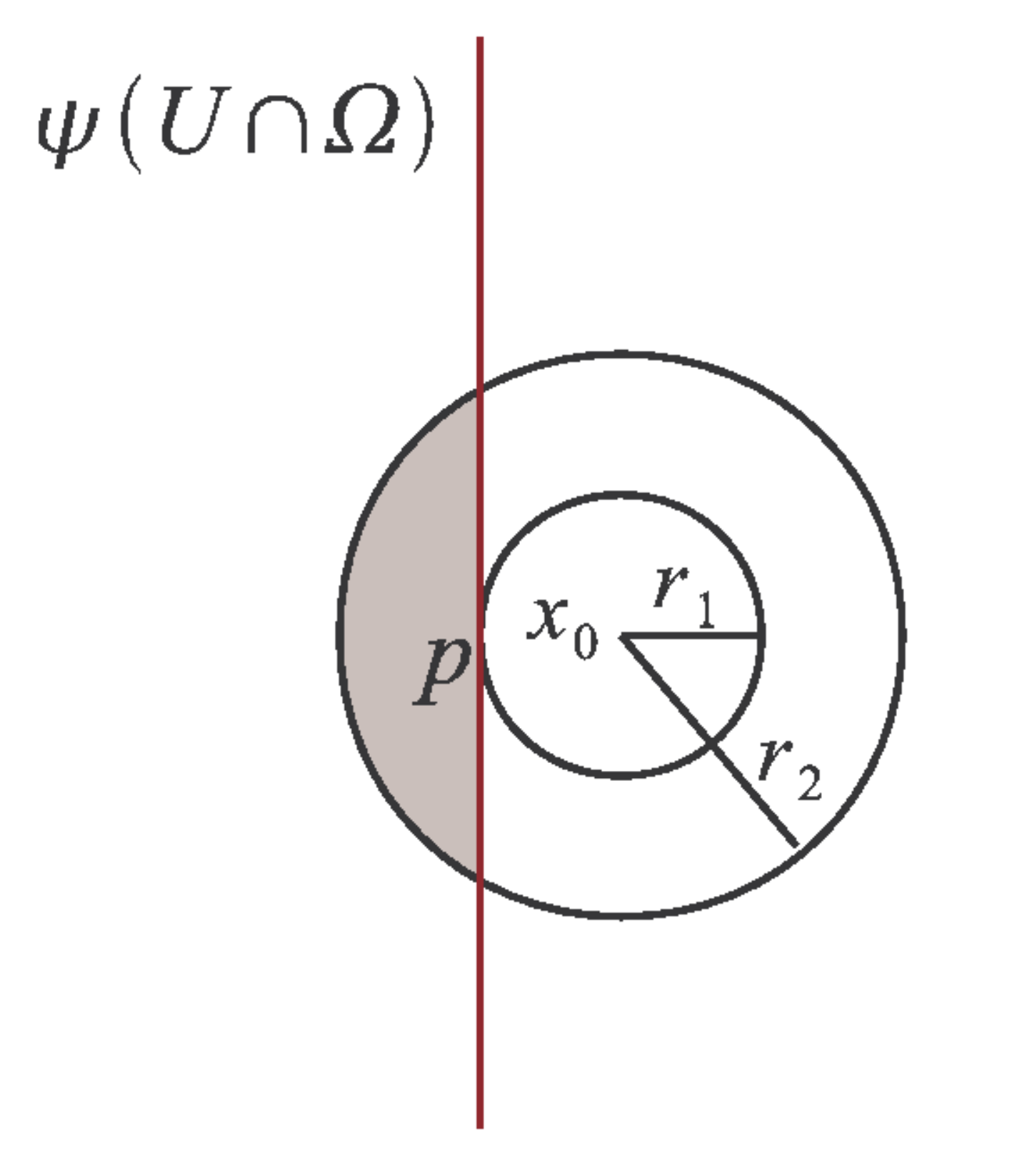}
\caption{L'area in grigio è l'area di definizione della barriera locale $\beta$}\label{fig_1}
\end{center}
\end{figure}
Visto che questo insieme è regolare per il laplaciano $\Delta_R$, esiste una funzione $\beta$ armonica su $A$, continua su $\overline A$, tale che $\beta|_{\partial B(x_0,r_1)}=0$, $\beta|_{\partial B(x_0,r_2)}=1$. Grazie al principio del massimo, è facile dimostrare che questa funzione è una barriera locale per il punto $p$.
\end{proof}
\end{prop}
\subsection{Altri domini regolari}\label{sec_altri_dir}
Lo scopo di questa sezione è dimostrare che ogni dominio della forma $R\supset\Omega= A\setminus K$ dove $A$ è un dominio con bordo liscio e $K$ una sottovarietà regolare di codimensione $1$ con bordo regolare contenuta in $A$, è un dominio regolare per $\Delta_R$. Cercheremo di tenere la dimostrazione al livello più elementare e ``geometrico'' possibile. Il lemma \ref{lemma_fiko} sarà essenziale per ottenere questo risultato.\\
Cominciamo con alcuni risultati preliminari sul laplaciano standard in $\R^n$. Nel seguito chiameremo armoniche le funzioni armoniche nel senso di $\R^n$ e considereremo quindi domini in $\R^n$ ($n\geq 2$), indicando con $B(\bar x,r)$ la bolla euclidea aperta di centro $\bar x$ e raggio $r$, e con $D(\bar x,r)$ la bolla aperta $n-1$ dimensionale centrata in $\bar x$ di raggio $r$, cioè l'insieme
\begin{gather*}
 D(\bar x,r)=\{x\in \R^n \ t.c. \ \norm{\bar x-x}<r \ e \ x_m=\bar x_m\}=B(\bar x,r)\cap \{x_m=\bar x_m\}
\end{gather*}
\begin{prop}
 Sull'insieme $\Omega\equiv B(\bar x, r_2 )\setminus \overline{D(\bar x, r_1)}$ con $r_2>r_1$, esiste una funzione $f:\Omega\to \R$ tale che
\begin{gather*}
 f\in H(\Omega) \cap C(\overline \Omega)=H(\Omega)\cap C(\overline {B(\bar x,r)}) \ \ \ f|_{\partial B(\bar x,r)}=0, \ \ \ f|_{\overline{D(\bar x,r)}}=1
\end{gather*}
\begin{proof}
 Vista l'invarianza per riscalamento di $\R^n$ e dell'operatore laplaciano standard, possiamo assumere senza perdita di generalità $r\equiv r_1 < r_2=1$. In tutta la dimostrazione, $B=B(\bar x,1)$, $D=D(\bar x,r)$.\\
Grazie all'osservazione \ref{oss_an}, esiste una successione di aperti relativamente compatti con bordo liscio $A_n$ tali che $A_n\subset A_{n-1}$ e $\overline D=\cap_n A_n$. Definitivamente in $n$, vale che $\overline{A_n}\subset B(0,1)$. Assumiamo senza perdita di generalità che $\overline{A_1}\subset B$. Per ogni $n$, l'insieme $\Omega_n=B\setminus \overline{A_n}$ è un insieme regolare per il laplaciano (ha bordo liscio), quindi per ogni $n$ esiste una funzione:
\begin{gather*}
 f_n\in C(\overline{B},\R) \ \ \ f_n\in H(\Omega_n) \ \ \ f_n|_{\partial B}=0 \ \ \ f_n|_{\overline{A_n}}=1
\end{gather*}
Grazie al principio del massimo e al fatto che $A_n\subset A_{n-1}$, è facile verificare che questa successione è una successione decrescente, e che per ogni $n$, $0\leq f_n\leq 1$. Grazie al principio di Harnack (vedi \ref{prop_harnackpri}), la successione $f_n$ converge localmente uniformemente su $\Omega$ a una funzione $f$ armonica su $\Omega$. Dimostriamo che questa funzione $f$ è continua e che soddisfa le richieste sui valori al bordo.\\
Per prima cosa consideriamo $\bar y\in D$. Per definizione di $D$, esiste un raggio $r(\bar y)$ sufficientemente piccolo in modo che $D(\bar y,r(\bar y))\subset D$ e $B(\bar y,r(\bar y))\subset B$. Consideriamo la semisfera:
\begin{gather*}
 B^+(\bar y,r(\bar y))=B(\bar y,r(\bar y))\cap \{x\in \R^n \ t.c. \ x_m\geq 0\}
\end{gather*}
Grazie alla proposizione \ref{prop_dircap}, questo insieme è un'insieme regolare per il laplaciano. Consideriamo una qualsiasi funzione continua $h:\partial B^+\to \R$ tale che $h|_{\partial B^+}=0$, $h(\bar y)=1$ e $0\leq h \leq 1$, e sia $u$ la soluzione del relativo problema di Dirichlet (cioè $u|_{\partial B^+}=h$). Grazie al principio del massimo, è facile concludere che per ogni funzione $f_n$, $f_n\geq u$ su $B^+\cap A_n$, quindi anche su tutta $B^+$. Passando al limite otteniamo che
\begin{gather*}
 f(x)\geq u(x)
\end{gather*}
Questo implica che la funzione $f$ non può essere identicamente nulla e che:
\begin{gather*}
 \liminf_{x\to \bar y} f(x)\geq \lim_{x\to \bar y} u(x)=1
\end{gather*}
L'arbitrarietà di $\bar y\in D$ garantisce che
\begin{gather*}
 \liminf_{x\to \bar y} f(x)\geq 1
\end{gather*}
per ogni $\bar y\in D$.\\
Dato che per costruzione $\limsup_{x\to \bar y} f(x)\leq 1$, otteniamo che per ogni punto $\bar y\in D$, $f$ è continua e $f|_{D}=1$. Rimangono da considerare i punti sull'insieme $E\equiv\overline D\setminus D$. Sia a questo scopo $\tilde y\in E$. Per costruzione di $f$:
\begin{gather*}
 0\leq \liminf_{x\to \tilde y} f(x)\leq \limsup_{x\to \tilde y} f(x)=1
\end{gather*}
Consideriamo un numero $0<\lambda<1$. Ricordiamo che dato un insieme $I$, si dice
\begin{gather*}
 \lambda I=\{x\in \R^n \ t.c. \ x/\lambda \in I\}=\{\lambda x \ t.c. \ x\in I\}
\end{gather*}
Grazie all'invarianza per traslazioni di $\R^n$ e dell'operatore laplaciano standard, possiamo ridefinire per comodità il sistema di riferimento in modo che $\tilde y=0$.\\
Scegliamo $\lambda$ in modo che $\partial (\lambda B) = \lambda(\partial B) \subset B\setminus \Omega$ (vedi figura \ref{fig_2}).
\begin{figure}[ht!]
\begin{center}
     \includegraphics[width=60mm]{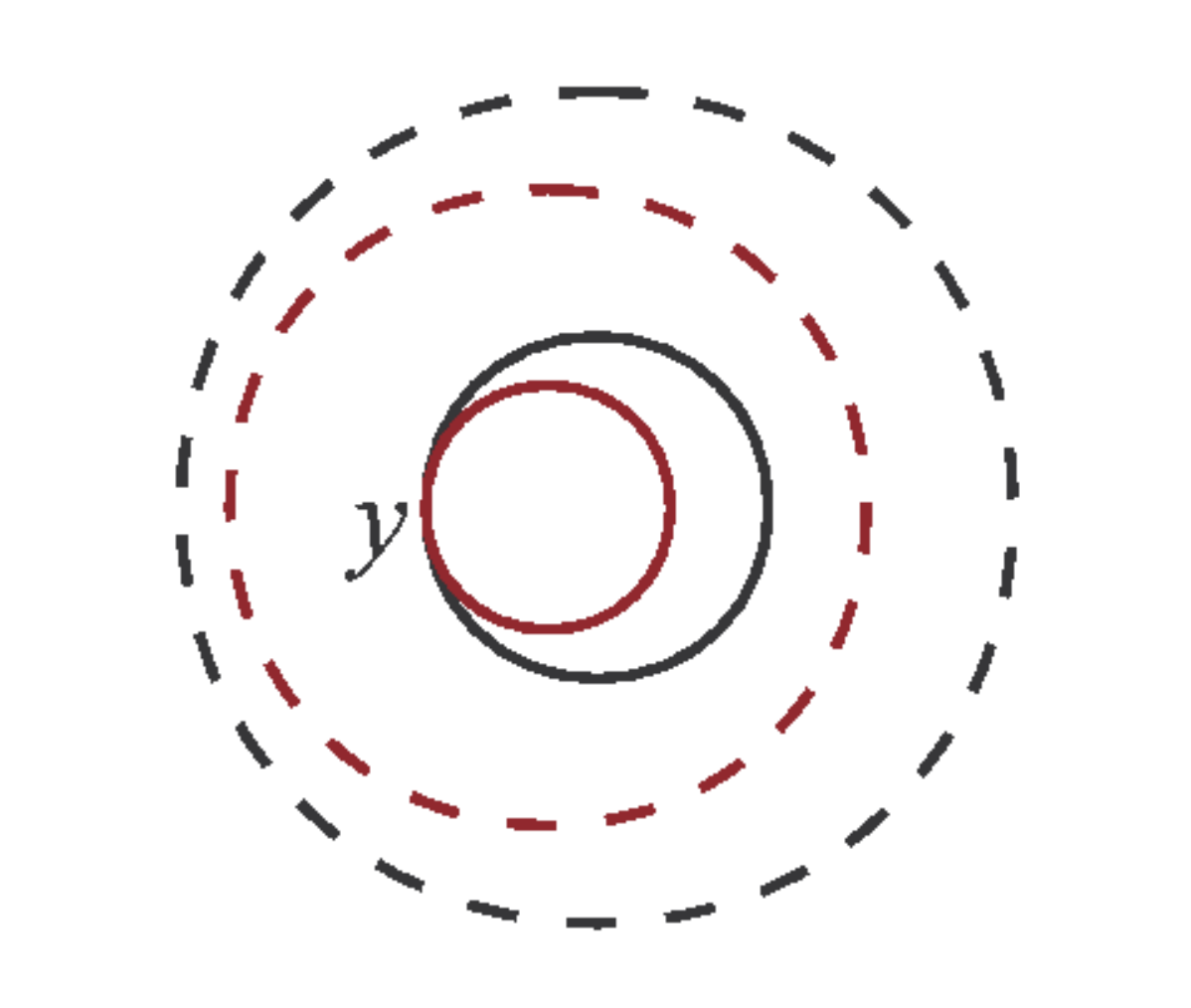}
\caption{Le linee piene rappresentano il bordo degli insiemi $m-1$ dimensionali, le linee tratteggiate il bordo degli insiemi $m$ dimensionali. Il colore rosso indica gli insiemi moltiplicati per $\lambda$.}\label{fig_2}
\end{center}
\end{figure}
La funzione $f$ è armonica su $\Omega$, quindi assume un valore minimo su $\lambda (\partial B)$. Grazie al principio del massimo \ref{prop_max2}, questo valore è strettamente positivo, diciamo $\epsilon>0$. Consideriamo la funzione
\begin{gather*}
 f_\lambda:\lambda B \to \R \ \ \ f_{\lambda}(x)=f(\lambda^{-1} x)
\end{gather*}
Osserviamo subito che $\liminf_{x\to \tilde y} f(x)=\liminf_{x\to \tilde y} f_\lambda (x)$ (ricordiamo che abbiamo scelto un riferimento dove $\tilde y=0$). Definiamo inoltre
\begin{gather*}
 f_\epsilon(x)=\frac{f_{\lambda}(x)+\epsilon}{1+\epsilon}
\end{gather*}
Questa funzione è armonica su $\lambda \Omega$, quindi essendo $\lambda<1$ è armonica su $\lambda B \setminus \overline D$. Inoltre vale che $f_\epsilon|_{\partial (\lambda B)}=\frac{\epsilon}{1+\epsilon}<\epsilon$, e $f_\epsilon\leq 1$ ovunque definita. Dato che la successione $f_n$ converge a $f$ monotonamente dall'alto, per ogni $n$ si ha che:
\begin{gather*}
 f_\epsilon |_{\partial (\lambda B)}< f_n|_{\partial (\lambda B)}
\end{gather*}
Inoltre visto che $f_\epsilon \leq 1$ ovunque definita, grazie al principio del massimo otteniamo che per ogni $n$:
\begin{gather*}
 f_n > f_{\epsilon} \ \ \Rightarrow \ \ f\geq f_{\epsilon}
\end{gather*}
Questo significa che
\begin{gather*}
 \liminf_{x\to \tilde y=0} f(x)\geq \liminf_{x\to \tilde y=0} f|_{\epsilon}(x)=\frac{\epsilon+\liminf_{x\to \tilde y=0} f(x)}{1+\epsilon}
\end{gather*}
se definiamo $\liminf_{x\to \tilde y=0}f(x)\equiv L$, otteniamo che:
\begin{gather*}
 L+\epsilon L \geq L+\epsilon \ \ \Rightarrow \ \ L\geq 1
\end{gather*}
Quindi:
\begin{gather*}
1\leq  \liminf_{x\to \tilde y} f(x)\leq \limsup_{x\to \tilde y} f(y)=1
\end{gather*}
ripetendo la costruzione per ogni $\tilde y\in E$, otteniamo la tesi.
\end{proof}
\end{prop}
Osserviamo che a patto di riscalare la funzione $f$, è possibile ottenere una funzione tale che:
\begin{gather*}
 f\in H(\Omega) \cap C(\overline \Omega)=C(\overline {B(\bar x,r)}) \ \ \ f|_{\partial B(\bar x,r)}=a, \ \ \ f|_{\overline{D(\bar x,r)}}=b
\end{gather*}
per qualsiasi valore fissato di $a$ e $b\in \R$.\\
Passiamo ora a dimostrare che ogni dominio $\Omega$ della forma descritta qui sopra è regolare rispetto al laplaciano standard.
\begin{prop}
 L'insieme $\Omega\equiv B(\bar x, r_2 )\setminus \overline{D(\bar x, r_1)}$ con $r_2>r_1$ è regolare per il laplaciano standard, cioè per ogni funzione $\phi:\partial \Omega\to \R$ continua, esiste una soluzione del relativo problema di Dirichlet.
\begin{proof}
Come sopra, possiamo assumere $\bar x=0$ e $r\equiv r_1<r_2=1$. Iniziamo con il considerare solo funzioni $\phi$ identicamente nulle su $\partial B$ e $0\leq \phi\leq 1$ ovunque su $\partial \Omega$, e dimostriamo che $P[\phi]$ risolve il relativo problema di Dirichlet.\\
Sappiamo che $0\leq P[\phi]\leq 1$ su tutto $\overline \Omega$. Inoltre grazie all'esistenza di barriere per ogni punto di $\partial B$, sappiamo che $P[\phi]$ è continua su un intorno di $\partial B$ (a priori disgiunto da $\overline D$) e $P[\phi]|_{\partial B}=0$. Sia $y\in \overline D$. Per continuità, dato $\epsilon>0$, esiste $\delta>0$ tale che
\begin{gather*}
 \phi(y)-\epsilon<\phi|_{B(y,\delta)\cap \overline D}<\phi(y)+\epsilon
\end{gather*}
Consideriamo ora un insieme $\Omega_y=B(x,\rho_2)\setminus \overline{D(x,\rho_1)}$, dove $x$, $\rho_1$, $\rho_2$ sono scelti in modo che $B(x,\rho_2)\subset B(y,\delta)$, $\overline D(x,\rho_1)\subset D$, $y\in \overline D(x,\rho_1)$. Ad esempio, se $y\in D$ è sufficiente scegliere $x=y$ e $\rho_1<\rho_2<\delta$, mentre se $y\in \overline D \setminus D$ (cioè se $y$ sta sul bordo del disco $D$), allora si può scegliere $D(x,\rho_1)$ in modo che sia un disco interno e tangente a $D$ nel punto $y$, sufficientemente piccolo in modo che $\overline {D(x,\rho_1)}\subset B(y,\delta)$.\\
Consideriamo ora due funzioni, $f_1$ e $f_2$, definite in questo modo:
\begin{gather*}
 f_1|_{\partial B(x,\rho_2)}=0 \ \ \ f_1|_{\overline D(x,\rho_1)}=\phi(y)-\epsilon \ \ \ f\in H(B(x,\rho_2)\setminus \overline{D(x,\rho_1}))\\
 f_2|_{\partial B(x,\rho_2)}=1 \ \ \ f_1|_{\overline D(x,\rho_1)}=\phi(y)+\epsilon \ \ \ f\in H(B(x,\rho_2)\setminus \overline{D(x,\rho_1}))
\end{gather*}
estendiamo le due funzioni a costanti su $\overline B\setminus B(x,\rho_2)$, in particolare $f_1=0$ e $f_2=1$ fuori dall'insieme $B(x,\rho_2)$. Entrambe le funzioni risultano continue su $\overline B$, $f_1$ è subarmonica su $\Omega$ (in quanto massimo tra due funzioni armoniche), mentre $f_2$ è superarmonica su $\Omega$ (in quanto minimo tra due funzioni armoniche). Inoltre per costruzione $f_1|_{\partial \Omega}\leq \phi$ e $f_2|_{\partial \Omega}\geq \phi$.\\
Per definizione di $P[\phi]$, otteniamo quindi che:
\begin{gather*}
 P[\phi]\geq f_1
\end{gather*}
quindi vale che:
\begin{gather*}
 \phi(y)-\epsilon =\lim_{x\to y} f_1(x)\leq \liminf_{x\to y}P[\phi](x)
\end{gather*}
Inoltre, qualunque funzione $u\in S_{\phi}$ \footnote{vedi \ref{deph_perron}} per il principio del massimo è minore della funzione $f_2$, quindi anche $P[\phi]$ conserva questa proprietà. Da questo ricaviamo che:
\begin{gather*}
 \phi(y)+\epsilon=\lim_{x\to y}f_2(x)\geq \limsup_{x\to y} P[\phi](x)
\end{gather*}
data l'arbitrarietà di $\epsilon$ e del punto $y$, possiamo concludere che $P[\phi]$ risolve il problema di Dirichlet considerato.\\
Osserviamo che con il procedimento appena descritto permette di costruire una barriera per ogni punto di $\overline D$. Indatti fissato $\bar y\in \overline D$, se scegliamo
\begin{gather*}
 \phi_{\bar y}(x)=\begin{cases}
                   1-d(x,\bar y) & se \ x\in \overline D \\
0 & se \ x\in \partial B
                  \end{cases}
\end{gather*}
grazie a quanto appena dimostrato la funzione $P[\phi]$ è armonica in $\Omega$, continua su $\overline B$, $P[\phi_y](y)=1$ e $P[\phi_y](x)<1$ per ogni $x\neq y$. Quindi la funzione $1-P[\phi_y]$ è una barriera per il punto $y$.\\
Dato che i punti di $\partial B$ sono regolari (hanno tutti una barriera locale), abbiamo dimostrato che il dominio $\Omega$ è regolare per l'operatore laplaciano standard.
\end{proof}
\end{prop}
Ricordiamo che grazie al lemma \ref{lemma_fiko}, i domini appena considerati sono regolari anche per l'operatore $\Delta_R$ (in carte locali).\\
Passiamo ora a dimostrare la proposizione principale di questo paragrafo:
\begin{prop}
Ogni dominio su una varietà riemanniana $R$ della forma $\Omega= A\setminus K$ dove $A$ è un dominio con bordo liscio e $K$ una sottovarietà regolare di codimensione $1$ con bordo regolare contenuta in $A$, è un dominio regolare per $\Delta_R$.
\begin{proof}
Dimostriamo che per ogni punto di $\partial \Omega$ è un punto regolare per il problema di Dirichlet.\\
Se $p\in \partial A$, $p$ è regolare grazie alla proposizione \ref{prop_dir_li}, quindi possiamo limitarci a considerare il caso $p\in K$. D'ora in avanti indicheremo con $\partial K$ il bordo di $K$ inteso come bordo della sottovarietà, non come bordo topologico, e $K^{\circ}$ sarà la parte interna sempre intesa nel senso di sottovarietà.\\
Se $p\in K^{\circ}$, consideriamo $(V(p),\psi)$ un aperto (nella topologia di $R$) coordinato intorno di $p$ tale che
\begin{gather*}
 V\cap K =\psi^{-1}\{(x_1,\cdots,x_{m-1},0)\}, \ \ \ \psi(p)=0
\end{gather*}
ad esempio possiamo considerare le coordinate di Fermi su un intorno di $p$ (vedi proposizione \ref{prop_coord_fermi}). Sia $r_2$ tale che $\overline {B(0,r_2)}\subset \psi(V)$ e sia $r_1<r_2$. Allora l'insieme
\begin{gather*}
 S=B(0,r_2)\setminus \overline {D(0,r_1)}
\end{gather*}
è un'insieme regolare per l'operatore $\Delta_R$, quindi esiste la soluzione del problema di Dirichlet
\begin{gather*}
 \Delta_R(\Phi)=0 \ su \ S, \ \ \ \Phi \in C(\overline S), \ \ \ \Phi(x)=d(x,0) \ \ \forall x\in \partial S
\end{gather*}
È facile verificare che la funzione $\Phi$ è una barriera locale per il punto $p$.\\
Se $p\in \partial K$, la dimostrazione della sua regolarità è molto simile. Consideriamo $(V,\psi)$ un'aperto con coordinate di Fermi tale che:
\begin{gather*}
 V\cap K=\psi^{-1}\{(x_1,\cdots,x_{m-1},0) \ t.c.  \ x_{m-1}\leq 0\} \ \ \ \psi(p)=0
\end{gather*}
cioè un sistema di coordinate che rappresenta $K$ come un semipiano in $\R^n$. Sia $r_2$ abbastanza piccolo da verificare $B(0,2r_2)\subset \psi(V)$. Se $r_1<r_2$, e $\bar x=(0,\cdots,0,-r_1,0)$, l'insieme
\begin{gather*}
 S= B(\bar x,r_2)\setminus D(\bar x,r_1)
\end{gather*}
è un'insieme regolare per l'operatore $\Delta_R$, quindi esiste la soluzione del problema di Dirichlet:
\begin{gather*}
 \Delta_R(\Phi)=0 \ su \ S, \ \ \ \Phi \in C(\overline S), \ \ \ \Phi(x)=d(x,0) \ \ \forall x\in \partial S
\end{gather*}
Come nel caso precedente, è facile verificare che la funzione $\Phi$ è una barriera locale per il punto $p$.
\end{proof}
\end{prop}

\subsection{Regolarità sul bordo}
In questo paragrafo riportiamo un lemma che garantisce sotto certe condizioni la regolarità delle soluzioni di particolari problemi di Dirichlet fino al bordo del dominio. Il lemma è basato sul teorema 17.3 pag 28 di \cite{32}, di cui riportiamo una versione semplificata senza dimostrazione.
\begin{teo}
 Sia $\Omega$ un dominio con bordo liscio, e $L$ un operatore uniformemente ellittico su $\overline \Omega$ con coefficienti $a, b,c$ lisci. Se $u$ è soluzione del problema di Dirichlet generalizzato
\begin{gather*}
 Lu=f \ in \ \Omega \ \ \ \ u|_{\partial \Omega}=0
\end{gather*}
con $f\in C^\infty(\overline \Omega)$, allora $u\in C^\infty(\overline \Omega)$.
\end{teo}
Ricordiamo che, grazie alla teoria sviluppata in \cite{32}, sotto le ipotesi citate il problema di Dirichlet generalizzato ha sempre un'unica soluzione. Osserviamo che considerando un dominio relativamente compatto con bordo liscio $\Omega\subset R$, il laplaciano su varietà soddisfa tutte le condizioni del teorema.\\
Grazie a questo teorema siamo in grado di dimostrare che:
\begin{lemma}\label{lemma_regbordo}
Siano $L$ e $\Omega$ come nel teorema precedente, sia $u\in C^2(\Omega)\cap C(\overline \Omega)$ con $Lu=0$. Se esiste una funzione $h\in C^\infty(\overline \Omega)$ tale che
$$
u|_{\partial \Omega}=h|_{\partial \Omega}
$$
allora $u\in C^\infty(\overline\Omega)$.
\begin{proof}
 La funzione $u-h$ soddisfa il problema di Dirichlet generalizzato:
\begin{gather*}
 L(u-h)=L(u)-L(h)=-L(h)\in C^\infty(\overline \Omega) \ \ \ \ \ (u-h)|_{\partial\Omega}=0
\end{gather*}
grazie al teorema precedente, la funzione $(u-h)\in C^\infty(\overline \Omega)$, quindi per differenza anche la funzione $u\in C^\infty(\overline\Omega)$.
\end{proof}
\end{lemma}
Nella tesi, spesso utilizzeremo questo lemma applicato a domini $\Omega$ con bordo costituito da 2 componenti connesse $\partial \Omega_0$ e $\partial \Omega_1$ \footnote{ad esempio un anello in $\R^2$ ha queste caratteristiche} e funzioni $u$ costanti su ogni componente \footnote{ad esempio $u|_{\partial \Omega_0}=0$ e $u|_{\partial \Omega_1}=1$}. In questo caso sappiamo che la funzione $h$ esiste, infatti:
\begin{prop}
 Sia $\Omega$ un dominio con bordo costituito da due componenti connesse $\partial \Omega_0$ e $\partial \Omega_1$. Allora esiste una funzione $h\in \overline \Omega$ tale che $h|_{\partial\Omega_0}=0$ e $h|_{\partial \Omega_1}=1$.
\begin{proof}
 Siano $V_0$ e $V_1$ due intorni disgiunti di $\partial \Omega_0$ e $\partial \Omega_1$ \footnote{gli insiemi $\partial \Omega_0$ e $\partial \Omega_1$ sono chiusi e disgiunti in $R$ spazio metrico}, e sia $\{\lambda_0,\lambda_1,\lambda\}$ uan partizione dell'unità subordinata al ricoprimento aperto $\{V_0,V_1, (\partial \Omega)^C\}$. Chiamiamo $f_0$ e $f_1$ due funzioni lisce con le caratteristiche descritte in \ref{lemma_T4liscio}, in particolare:
$$
\partial \Omega_0=f_0^{-1}(0) \ \ \ \partial \Omega_1=f_1^{-1}(0)
$$
È facile verificare che la funzione:
\begin{gather*}
 h\equiv \lambda_0\cdot f_0 + \lambda_1 \cdot (1-f_1)
\end{gather*}
è una funzione liscia con $h|_{\partial \Omega_0}=0 \ \ h|_{\partial \Omega_1}=1$.
\end{proof}
\end{prop}
Questo dimostra ad esempio che tutte le funzioni
$$
u\in H(\Omega), \ \ u|_{\partial\Omega_0}=0 \ \ u|_{\partial \Omega_1}=1
$$
sono funzioni $u\in C^\infty(\overline \Omega)$.

\chapter{Ultrafiltri e funzionali lineari moltiplicativi}\label{chap_filtri}
In questo capitolo introduciamo gli ultrafilti. L'introduzione è tratta da \cite{1}. Il risultato principale riguarda l'esistenza di funzionali lineari moltiplicativi sull'algebra di $l_{\infty}(\N)$, ed è tratto dall'articolo \cite{2}.\\
Lo scopo di questo capitolo è di descrivere un modo alternativo alle successioni per caratterizzare le topologie che non sono sequenziali, nelle quali cioè la chiusura per successioni non coincide con la chiusura topologica. Un esempio di questi spazi topologici è lo spazio $\Ro$ con la topologia $B$ descritto nella definizione \ref{deph_B}. Inoltre gli ultrafiltri possono essere utilizzati per definire dei caratteri sullo spazio $l_{\infty}(\N)$ (vedi paragrafo \ref{sec_charsucc}), e questo risultato può essere utilizzato per dare un esempio non banale di carattere sull'algebra di Royden. A questo scopo rimandiamo alla sezione \ref{subsec_charroy}.
\section{Filtri, ultrafiltri e proprietà}
Per prima cosa introduciamo la definizione di filtro.
\begin{deph}
 Sia $(X,\tau)$ uno spazio topologico. Un filtro in questo spazio è una collezione $\mathcal{A}=\{A_{\alpha} \ t.c. \ \alpha \in I\}$ di sottoinsiemi di $X$ tale che:
\begin{enumerate}
 \item $\forall \alpha\in I \ A_{\alpha}\neq \emptyset$
 \item $\forall \alpha, \ \beta \in I$ $\exists \gamma \in I \ t.c. \ A_{\gamma}\subset A_{\alpha}\cap A_{\beta}$
\end{enumerate}
 
\end{deph}
Osserviamo dalla definizione che ogni coppia di insiemi in $\mathcal{A}$ ha intersezione non vuota, e per induzione qualunque collezione finita di insiemi di $\A$ ha intersezione non vuota. Gli ultrafiltri generalizzano la nozione di convergenza di successioni e sono utili negli spazi non primo numerabili (come ad esempio l'algebra di Royden). A questo scopo definiamo:
\begin{deph}
 Dato un ultrafiltro $\A$ nello spazio $(X,\tau)$ diciamo che $\A$ converge a $x_0$ e scriviamo che $\A\to x_0$ se e solo se $\forall U(x_0)$ intorno aperto in $\tau$, esiste $\alpha$ tale che $\A_{\alpha}\subset U$.\\
Diciamo invece che $\A$ si accumula in $x_0$ e scriviamo $\A\triangleright x_0$ se e solo se $\forall U(x_0)$, $\forall \alpha$, $A_\alpha\cap U\neq \emptyset$
\end{deph}
Osserviamo subito che $\A\to x_0\Rightarrow \A\triangleright x_0$, e che negli spazi di Hausdorff $\A$ non si accumula in nessun altro punto \footnote{questo è il teorema 3.2 a pag. 214 di \cite{1}}. Questo è evidente dal fatto che ogni coppia di insiemi in $\A$ ha intersezione non vuota. Inoltre dalla definizione è evidente che $\A\triangleright x_0$ se e solo se $x_0\in \bigcap_{\alpha}\overline{A_{\alpha}}$.\\
Un esempio banale di filtro convergente è l'insieme degli intorni di un punto. Dato $x_0\in X$, indichiamo con $\U(x_0)$ il filtro dei suoi intorni. \`E evidente che $\U(x_0)\to x_0$.\\
\`E possibile definire in alcuni casi operazioni di unione e intersezione tra filtri:
\begin{deph}\label{deph_algb_filt}
 Dati due filtri $\A$ e $\B$ su $X$ definiamo i filtri:
\begin{enumerate}
 \item $\A\cup \B = \{A_{\alpha}\cup B_{\beta}\ t.c. \ \alpha\in I, \ \beta \in Y\}$
 \item se $A_{\alpha}\cap B_{\beta}\neq \emptyset$ per ogni scelta di $\alpha$ e $\beta$, allora definiamo $$\A\cap\B=\{A_{\alpha}\cap B_{\beta}\ t.c. \ \alpha\in I, \ \beta \in Y\}$$
\end{enumerate}
La facile verifica che questi insiemi sono filtri è lasciata al lettore.
\end{deph}

\begin{deph}
 Dati due filtri $\A=\{A_{\alpha} \ t.c. \ \alpha \in I\}$ e $\B=\{B_{\beta} \ t.c. \ \beta \in Y\}$, diciamo che $\B$ è subordinato ad $\A$ e scriviamo $\B\vdash\A$ se e solo se $\forall \alpha \ \exists \beta \ t.c. \ B_{\beta}\subset A_{\alpha}$
\end{deph}
\`E facile verificare che la relazione di subordinazione è una relazione di preordine, ma non di ordine parziale. Infatti è una relazione riflessiva e transitiva, ma non è vero che $\A\vdash\B \ \wedge \ \B\vdash\A \Rightarrow \A=\B$. Un controesempio si può trovare considerando $X=\N$, $\A=\{[2n,\infty) \ n\in \N\}$ e $\B=\{[n,\infty) \ n\in \N\}$.\\
Per la relazione di subordinazione valgono le seguenti proprietà:
\begin{prop}\label{prop_filt}
\begin{enumerate}
 \item $\A\subset \B\Rightarrow \B\vdash\A$
 \item Se $\B\vdash\A$, allora ogni elemento di $\B$ ha intersezione non vuota con ogni elemento di $\A$
 \item $\A\to x_0$ se e solo se $\A\vdash \U(x_0)$
 \item Se $\B\vdash\A$, allora $\A\to x_0 \Rightarrow \B\to x_0$ e anche $\B\triangleright x_0\Rightarrow \A\triangleright x_0$
\end{enumerate}
\begin{proof}
 La dimostrazione di questi punti è lasciata al lettore.
\end{proof}
\end{prop}
Per le successioni negli spazi topologici vale che una successione converge a un punto $p$ se e solo se da ogni sua sottosuccessione si può estrarre una sottosottosuccessione convergente a $p$, e una successione si accumula in un punto se e solo se esiste una sua sottosuccessione convergente a quel punto. Valgono risultati analogi per i filtri:
\begin{prop}
 Un filtro $\A$ converge a $x_0$ se e solo se per ogni $\B\vdash\A$, esiste $\mathcal{C}\vdash \B$ tale che $\mathcal C \to x_0$
\begin{proof}
 Dalla proposizione precedente risulta ovvia una delle due implicazioni. Per dimostrare l'altra, supponiamo che $\A$ non converga a $x_0$. Quindi esiste $U(x_0)$ intorno di $x_0$ tale che $\forall \alpha, \ A_\alpha \cap U^C\neq\emptyset$. La collezione $\B=\{A_\alpha\cap U^C \ t.c. \ \alpha\in I\}$ è un filtro subordinato a $\A$ e per il quale $x_0$ non è punto di accumulazione. Allora grazie alla proposizione precedente risulta che ogni filtro subordinato a $\B$ non può convergere a $x_0$.
\end{proof}
\end{prop}
\begin{prop}\label{prop_filt2}
 $\A\triangleright x_0$ se e solo se esiste $\B\vdash\A$ tale che $\B\to x_0$
\begin{proof}
 Se esiste $\B\vdash\A$ tale che $\B\to x_0$, è chiaro che $\A\triangleright x_0$. Supponiamo ora che $\A\triangleright x_0$. La collezione $\B=\A\cap \U(x_0)$ è un filtro dato che $\forall \alpha \ \forall U(x_0), \ A_{\alpha}\cap U \neq \emptyset$. Inoltre è chiaro che $\B\vdash\U(x_0)$. Allora per le proprietà della relazione di subordinazione $\B\to x_0$.
\end{proof}
\end{prop}
I filtri si possono definire anche tramite funzioni.
\begin{prop}\label{prop_filt5}
 Dato un filtro $\A$ su $X$ e una funzione $f:X\to Y$, la collezione $f(\A)=\{f(A_{\alpha} \ t.c. \ A_{\alpha}\in \A\}$ è un filtro su $Y$
\begin{proof}
 Il fatto che ogni elemento di $f(\A)$ sia non vuoto è ovvio. La seconda proprietà caratterizzante i filtri segue facilmente dalla considerazione che $f(A\cap B)\subset f(A)\cap f(B)$
\end{proof}

\end{prop}

Passiamo ora a definire gli ultrafiltri, o filtri massimali, e a dimostrarne l'esistenza.
\begin{deph}
 Un filtro $\M$ è detto massimale (o ultrafiltro) se non ha filtri propriamente subordinati, cioè se
\begin{gather*}
 \A\vdash\M \Rightarrow \M\vdash\A
\end{gather*}
\end{deph}
Non tutti gli autori condividono questa definizione. Su alcuni testi si definisce ultrafiltro un filtro che non è \textit{contenuto} propriamente in nessun altro filtro. Anche la definizione di subordinazione è diversa. Si definisce $\B$ subordinato a $\A$ se $\A\subset\B$. In questo modo la relazione di subordinazione diventa una relazione d'ordine parziale. Per gli scopi di questa tesi però non c'è differenza tra le due definizioni, quindi adotteremo quella meno restrittiva e più generale di \cite{1}. Passiamo ora a caratterizzare gli ultrafiltri e a dimostrarne l'esistenza.
\begin{prop}\label{prop_filt3}
Un filtro $\M=\{M_{\alpha} \ t.c. \ \alpha \in I\}$ è un ultrafiltro se e solo se per ogni $S\subset X$, uno dei due insiemi $S$ o $S^C$ contiene un elemento di $\M$.
\begin{proof}
 Dato un ultrafiltro $\M$ e un insieme $S$ è chiaro che non può accadere che sia $S$ che $S^C$ contengano un elemento di $\M$, altrimenti questi due elementi avrebbero necessariamente intersezione vuota. Assumiamo allora che $\forall \alpha, \ M_{\alpha}$ non sia contenuto in $S$, allora $\forall \alpha M_{\alpha}\cap S^C\neq\emptyset$. Questo implica che l'insieme $\M_1=\M\cap S^C=\{A_{\alpha}\cap S^C \ t.c. \ \alpha \in I\}$ sia un filtro tale che $\M_1\vdash \M$. Ma allora per ipotesi $\M\vdash\M_1$, quindi per ogni $\alpha$ esiste $\gamma$ tale che $A_{\gamma}\subset A_{\alpha}\cap S^C\subset S^C$.\\
L'altra implicazione si dimostra considerando un filtro $\M_1\vdash\M$. L'ipotesi assicura che per ogni $A^1_{\beta}\in \M_1$ esiste un $\A_{\alpha}\in \M$ tale che $\A_{\alpha}\subset A^1_{\beta}$ oppure $\A_{\alpha}\subset (A^1_{\beta})^C$. La seconda possibilità è esclusa dal fatto che $\M_1\vdash\M$, quindi $\M\vdash \M_1$, cioè $\M$ è massimale.
\end{proof}
\end{prop}
Grazie a questa caratterizzazione è facile dimostrare che dato $x_0\in X$, l'insieme $\M=\{A\subset X \ t.c. \ x_0\in A\}$ è un ultrafiltro, anche se è abbastanza banale e poco pratico. Per i nostri scopi è necessario dimostrare che ogni filtro è contenuto in un ultrafiltro, e per fare questo è necessario assumere il lemma di Zorn, che ricordiamo:
\begin{teo}[Lemma di Zorn]
 Ogni insieme preordinato in cui ogni catena contiene un limite superiore possiede almeno un elemento massimale
\end{teo}
Solitamente il lemma di Zorn è enunciato per insiemi parzialmente ordinati e non semplicemente preordinati. Si può dimostrare però che queste due versioni del teorema sono equivalenti (a questo scopo rimandiamo a \cite{3}, esercizio 1.16). Osserviamo che negli insiemi preordinati (con una relazione d'ordine $\geq$) un elemento $m$ è detto massimale se $a\geq m \Rightarrow m\geq a$.
\begin{teo}\label{teo_filt1}
 Dato un insieme $X$ e un filtro $\A$ su questo insieme, esiste un filtro massimale (ultrafiltro) $\M$ tale che $\M\vdash \A$.
\begin{proof}
 Chiamiamo $\hat B$ la famiglia dei filtri per i quali $\B\vdash \A$, e definiamo una relazione di preordine $\geq$ su questo insieme in questo modo:
\begin{gather*}
 \B_2\geq \B_1 \Longleftrightarrow \B_2\vdash\B_1
\end{gather*}
Ricordiamo che dalla proposizione \ref{prop_filt}, $\B_1\subset \B_2\Rightarrow \B_2\vdash\B_1$, quindi anche $\B_2\geq \B_1$. 
Consideriamo una catena $\B_{i}$ di elementi di $\hat B$ \footnote{ricordiamo che con catena in un insieme preordinato si intende un suo sottoinsieme tale che per ogni sua coppia di elementi $a$ e $b$, $a\geq b$ oppure $b\geq a$, dove è possibile che entrambe le affermazioni siano vere}. L'insieme $\tilde \B =\bigcup_i \B_i$ è un elemento massimale per questa catena. Per prima cosa dimostriamo che è un filtro. Ovviamente ogni suo elemento è non vuoto, inotlre dati due suoi elementi $A_1\in \tilde \B$ e $A_2\in \tilde \B$, esistono $\B_{i_1}$ e $\B_{i_2}$ filtri nella catena che li contengono. Supponiamo che $\B_{i_1}\vdash \B_{i_2}$, allora esiste $A_3\in \B_{i_1}$ tale che $A_3\subset A_2$, e quindi poiché $\B_{i_1}$ è un filtro, esiste anche $E\in \B_{i_1}$ tale che $E\subset A_1\cap A_3\subset A_1\cap A_2$. \`E poi evidente che $\tilde \B\geq \B_i \geq \A$ per ogni $i$, quindi $\tilde \B$ è un limite superiore per la catena. Applicando il lemma di Zorn otteniamo che esiste almeno un elemento massimale $\M$ in $\hat B$. $\M$ è un filtro massimale perché se $\mathcal{C}\vdash \M$, per transitività $\mathcal{C}\vdash\A$, quindi $\mathcal{C}\in \hat B$ e quindi per massimalità $\M\geq \mathcal{C}$, cioè $\M\vdash \mathcal{C}$.
\end{proof}
\end{teo}
Osserviamo ora due proprietà degli ultrafiltri che ci saranno utili nel seguito:
\begin{prop}\label{prop_filt4}
 Dato un ultrafiltro $\M$ su $X$, $\M\triangleright x_0$ se e solo se $\M\to x_0$
\begin{proof}
 Questa proposizione è un corollario della definizione di ultrafiltro e della proposizione \ref{prop_filt2}
\end{proof}
\end{prop}
\begin{prop}
 Data una qualunque funzione $f:X\to Y$ insiemi qualsiasi, se $\M$ è un ultrafiltro su $X$, allora $f(\M)=\{f(M_{\alpha}) \ t.c. \ M_{\alpha}\in \M\}$ è un ultrafiltro in $Y$.
\begin{proof}
 $f(\M)$ è un filtro grazie alla proposizione \ref{prop_filt5}, quindi rimane da dimostrare solo la sua massimalità. Utilizziamo la caratterizzazione \ref{prop_filt3} dei filtri massimali. Sia $S\in Y$. Visto che $f^{-1}(S^C)=(f^{-1}(S))^C$, esiste un elemento $M_{\alpha}$ contenuto o in $f^{-1}(S)$ o nel suo complementare. $M_{\alpha}\subset f^{-1}(S) \Rightarrow f(M_{\alpha})\subset S$ e lo stesso vale per $S^C$, quindi dato ogni sottoinsieme di $Y$, esiste un elemento di $f(\M)$ contenuto o in esso o nel suo complementare.
\end{proof}
\end{prop}
Il tipo di ultrafiltri massimali che ci interessano sono gli ultrafiltri \textit{non costanti}, ovvero gli ultrafiltri che non contengono nessun elemento finito
\begin{deph}
 Un ultrafiltro $\M$ è detto \textbf{non costante} se e solo se non contiene nessun insieme di cardinalità finita.
\end{deph}
Per dimostrare l'esistenza di ultrafiltri non costanti sull'insieme dei numeri naturali costruiamo un filtro con queste caratteristiche e lo estendiamo a massimale grazie al teorema \ref{teo_filt1} \footnote{dato che il teorema fa uso del lemma di Zorn, questa dimostrazione non è costruttiva. Non si possono cioè trovare esempi descrivibili facilmente di ultrafiltri massimali non costanti con questa tecnica.}.
\begin{oss}\label{oss_filt}
 La collezione $\F$ di sottoinsiemi a complementare finito è un filtro sull'insieme $\N$.
\begin{proof}
 La dimostrazione è immediata
\end{proof}
\end{oss}
Indichiamo con $\M_{\F}$ un ultrafiltro subordinato al filtro $\F$. Questo ultrafiltro necessariamente non conterrà nessun insieme a cardinalità finita. Sia per assurdo $S\in\M_{\F}$ insieme a cardinalità finita. $S^C\in \F$ per definizione. Ma visto che $\M_{\F}\vdash \F$, esiste un elemento di $\M_{\F}$ contenuto in $S^C$, ma allora questo elemento ha intersezione vuota con $S$, assurdo.\\
Prima di concludere questa sezione riportiamo la relazione tra compattezza di un insieme e filtri su quell'insieme.
\begin{prop}\label{prop_filt_comp}
 Per un insieme $K\subset (X,\tau)$ spazio topologico le seguenti affermazioni sono equivalenti:
\begin{enumerate}
 \item Ogni ricoprimento di aperti di $M$ ha un sottoricoprimento finito 
 \item Per ogni famiglia di chiusi $C_{i}$ tale che $\bigcap_i C_i=\emptyset$, esiste un numero finito di indici $i_1,\dots,i_n$ tali che $\bigcap_{k=1}^n C_{i_k}=\emptyset$
 \item Ogni filtro su $K$ ha un punto di accumulazione
 \item Ogni ultrafiltro in $K$ è convergente
\end{enumerate}
Se una qualsiasi delle proprietà precedenti è valida, l'insieme $K$ è detto compatto.
\begin{proof}
 (1) e (2) sono equivalenti grazie alle leggi di De Morgan, (3)$\Rightarrow$(4) grazie alla proposizione \ref{prop_filt4}, (4)$\Rightarrow$(3) grazie al teorema \ref{teo_filt1} e alla proposizione \ref{prop_filt}. (2)$\Rightarrow$(3) poiché dato un filtro $\A$, ogni collezione finita di suoi elementi $A_{\alpha_i}$ ha intersezione non vuota, quindi anche $\bigcap_{i=1}^{n} \overline{A_{\alpha_i}}\neq \emptyset $. Grazie a (2) allora $\bigcap_{\alpha} \overline{A_{\alpha}}\neq \emptyset $, quindi l'insieme dei punti di accumulazione di $\A$ è non vuoto. (3)$\Rightarrow$(2) perché data una qualsiasi collezione di chiusi $C_{\alpha}$ tale che ogni sua sottocollezione finita abbia intersezione non vuota, la collezione dei chiusi e di tutte le possibili intersezioni finite è un filtro, che per (2) ha un punto di accumulazione, quindi $\bigcap_{\alpha} \overline{C_{\alpha}}=\bigcap_{\alpha} C_{\alpha}\neq\emptyset$.
\end{proof}
\end{prop}
Vale anche la seguente caratterizzazione della continuità delle funzioni tramite ultrafiltri:
\begin{prop}
 Una funzione $f:X\to Y$ spazi topologici è continua in $x_0$ se e solo se per ogni filtro $\A\to x_0$, $f(\A)\to f(x_0)$
\begin{proof}
 Supponiamo che $f$ sia continua in $x_0$. Chiamando $\U(x)$ il filtro degli intorni aperti di $x_0$, dalla definizione di continuità si ha che per ogni intorno aperto $W$ di $f(x_0)$, esiste un aperto $U$ in $X$ tale che $f(U)\subset W$. Quindi per definizione $f(\U(x_0))\to f(x_0)$. Ora, se consideriamo un filtro $\A\to x_0$ qualsiasi, dalle proposizioni precedenti risulta che $\A\vdash \U(x_0)$, quindi anche $f(\A)\vdash f(\U(x_0))$, quindi necessariamente $f(\A)\to f(x_0)$.\\
Per dimostrare l'implicazione inversa, dato che $f(\U(x_0))\to f(x_0)$, per definizione per ogni intorno di $W(f(x_0))$, esiste $U(x_0)$ tale che $f(U(x_0))\subset W$, quindi $f$ è continua.
\end{proof}
\end{prop}
Per una successione in $\R$ possono accadere due cose: o la successione è limitata, quindi tutta contenuta in un compatto, oppure no. Se la successione è limitata e converge, allora converge a un punto al finito di $\R$, altrimenti una successione illimitata può ``convergere'' all'infinito. Per filtri e ultrafiltri /che giocano il ruolo delle successioni ``convergenti'' in senso generalizzato), vale una caratterizzazione simile.
\begin{deph}
 Un filtro $\A$ si dice avere supporto in un insieme $K\subset X$ se esiste $\alpha$ tale che $A_{\alpha}\subset K$
\end{deph}
Da questa definizione si ricava immediatamente che
\begin{prop}
 Se $\A$ ha supporto in $K$, allora $\A\triangleright x_0 \Rightarrow x_0\in \overline K$ (quindi anche $\A\to x_0 \Rightarrow x_0\in \overline K$)
\begin{proof}
 La dimostrazione è una conseguenza quasi immediata della definizione. Supponiamo per assurdo che $x_0 \not \in \overline K$. Allora esiste $U$ intorno di $x_0$ separato da $\overline K$. Ma allora $A_{\alpha}\cap U =\emptyset$ (dove $\alpha$ è l'indice tale che $A_{\alpha}\subset K$), quindi $x_0$ non può essere punto di accumulazione per $A$.
\end{proof}
\end{prop}
\begin{prop}
 Se $\A$ ha supporto in $K$, allora $\B\equiv \A\cap K$ è un filtro subordinato ad $\A$. Inoltre se $\A$ è un ultrafiltro, anche $\B$ è un ultrafiltro.
\begin{proof}
Visto che esiste $\alpha$ tale che $A_{\alpha}\subset K$, allora ogni insieme del filtro $\A$ ha intersezione non vuota con $K$, quindi $\B\equiv \A \cap K$ ben definisce un filtro, ovviamente subordinato ad $\A$. Il fatto che $\B$ sia un ultrafiltro se $\A$ lo è è una facile conseguenza della caratterizzazione \ref{prop_filt3}
\end{proof}
\end{prop}
Un filtro si dice avere \textit{supporto compatto} se esiste un insieme compatto che è un suo supporto. In uno spazio localmente compatto, ogni filtro convergente ha supporto compatto (basta considerare un intorno compatto del punto al quale il filtro converge). Se il filtro $\A$ non ha supporto compatto, vuol dire che per ogni compatto $K$ e per ogni indice $\alpha$, $A_{\alpha}\cap K^C\neq \emptyset$, quindi per ogni compatto $K$ il filtro $\B=\A\cap K^C$ è un filtro subordinato ad $\A$, e se $\A$ è un ultrafiltro, anche $\B$ lo è.
\section{Applicazioni: caratteri sulle successioni limitate}\label{sec_charsucc}
Come applicazione delle due sezioni precedenti, e soprattutto come esempio per il seguito, presentiamo l'esistenza di particolari caratteri sullo spazio $l_{\infty}(\N)$, cioè lo spazio di Banach delle successioni limitate a valori reali dotato della norma del sup. Questo spazio diventa un'algebra di Banach se definiamo la moltiplicazione tra successioni punto per punto, cioè date due successioni $x(n)$ e $y(n)$, definiamo $(x\cdot y)(n)\equiv x(n)\cdot y(n)$ \footnote{lasciamo le dovute verifiche al lettore}. \`E facile osservare che la successione costante uguale a 1 è l'unità di ques'algebra.\\
Tutti i funzionali $\phi_n:l_{\infty}(\N)\to\R$ definiti da $\phi_n(x)=x(n)$ sono caratteri. Cerchiamo però di definire altri funzionali, legati solo al comportamento di ogni successione all'infinito. Questa sezione è tratta dall'articolo \cite{2}.\\
Per prima cosa definiamo un'operazione di limite su $l_{\infty}(\N)$:
\begin{deph}
 Un'operazione di limite sullo spazio $l_{\infty}(\N)$ è un funzionale lineare $\phi:l_{\infty}(\N)\to \R$ tale che per ogni successione
\begin{gather*}
 \liminf_n x(n) \leq \phi(x) \leq \limsup_n x(n)
\end{gather*}
\end{deph}
Un'operazione di limite è quindi un funzionale lineare che vale $\lim_n x(n)$ se questo limite esiste. Grazie al teorema di Hahn-Banach è possibile dimostrare l'esistenza di alcune operazioni di limite (vedi ad esempio \cite{2}, oppure l'esercizio 4 cap.3 pag.85 di \cite{4}). \`E anche possibile dimostrare l'esistenza di operazioni di limite moltiplicative su $l_{\infty}(\N)$, e a questo scopo utilizzeremo la teoria sviluppata sugli ultrafiltri non costanti. Prima della proposizione ricordiamo che:
\begin{oss}\label{oss_sp}
 Ricordiamo la definizione di somma di insiemi. Dati due sottoinsiemi $A$ e $B$ di uno spazio vettoriale $V$, $A+B=\{a+b \ \ t.c. \ \ a\in A, \ b\in B\}$. Dalla definizione è chiaro che $(f+g)(A)\subset f(A)+g(A)$. In modo analogo alla somma di insiemi si può definire il loro prodotto, e vale ancora che $(f\cdot g )(A)\subset f(A)\cdot g(A)$.
\end{oss}

\begin{prop}\label{prop_car_N}
 Sull'insieme $l_{\infty}(\N)$ esistono operazioni di limite moltiplicative.
\begin{proof}
 Fissiamo un ultrafiltro non costante $\M$ sull'insieme $\N$ \footnote{l'esistenza di questo tipo di ultrafiltro è stata dimostrata nell'osservazione \ref{oss_filt}}. Una qualunque successione $x:\N\to \R$ permette di definire un ultrafiltro $x(\M)$ su $\R$, anzi visto che la successione è limitata, l'ultrafiltro sarà contenuto nell'insieme compatto $[-\norm x, \norm x]$, quindi grazie a \ref{prop_filt_comp} sarà convergente. Definiamo
\begin{gather*}
 \phi_{\M}(x)\equiv \lim_{\M} x \equiv \lim x(\M)
\end{gather*}
Per prima cosa osserviamo che se sostituiamo $\M$ con il filtro $\F$ degli insiemi a complementare finito, il limite non è sempre definito ma coincide con il limite di una successione in senso classico. Resta da dimostrare che $\phi_{\M}$ è lineare e moltiplicativo. A questo scopo consideriamo due successioni $x,\ y:\N\to \R$, e dimostriamo che $\phi_{\M}(x+y)=\phi_{\M}(x)+\phi_{\M}(y)$. Sia $\phi_{\M}(x)=s$ e $\phi_{\M}(y)=t$. Dalla definizione di limite di ultrafiltro otteniamo che per ogni $\epsilon>0$, esistono $\alpha$ e $\beta$ tali che
\begin{gather*}
 x(M_{\alpha})\subset B_\epsilon (s) \ \wedge \ y(M_{\beta})\subset B_\epsilon (t)
\end{gather*}
dove $M_{\alpha}, \ M_{\beta}\in \M$ e $B_{\epsilon}(s)$ indica l'insieme aperto $(s-\epsilon,s+\epsilon)$. Dato che $\M$ è un filtro, esiste $\gamma$ tale che $M_{\gamma}\subset M_{\alpha}\cap M_{\beta}$. Questo significa che per ogni $\epsilon>0$
\begin{gather*}
 (x+y)(M_{\gamma})\subset x(M_{\gamma})+y(M_{\gamma})\subset B_{\epsilon}(s)+B_{\epsilon}(t)=B_{2\epsilon}(s+t)
\end{gather*}
dove abbiamo usato l'osservazione \ref{oss_sp}. Questo dimostra che $$\lim_{\M}(x+y)=\lim_{\M}x+\lim_{\M}y$$
Consideriamo ora un qualunque numero $a\in \R$ e una successione $x\in l_{\infty}(\N)$. La dimostrazione del fatto che $\phi_{\M}(ax)=a\phi_{\M}(x)$ è del tutto analoga a quella appena mostrata, quindi lasciamo i dettagli al lettore \footnote{questo completa la dimostrazione che $\phi_{\M}$ è lineare.}.\\
La dimostrazione che $\phi_{\M}$ è moltiplicativo è anch'essa analoga a questa, l'unico punto un po' più difficile è capire come caratterizzare l'insieme $B_{\epsilon}(s)\cdot B_{\epsilon}(t)$. Osserviamo che
\begin{gather*}
 \abs{a-s}<\epsilon \ \wedge \ \abs{b-t}<\epsilon \ \Longrightarrow \\
\Longrightarrow \abs{ab-st}\leq \abs{ab-at}+\abs{at-st}\leq \epsilon(\abs{a}+t)\leq \epsilon(\abs{s}+\abs{t}+\epsilon)
\end{gather*}
quindi $B_{\epsilon}(s)\cdot B_{\epsilon}(t)\subset B_{[\epsilon(\abs{s}+\abs{t}+\epsilon)]}(st)$. Dato che $\epsilon(\abs{s}+\abs{t}+\epsilon)$ diventa piccolo a piacere al diminuire di $\epsilon$, seguendo un ragionamento del tutto analogo al precedente, si ottiene che per ogni $\epsilon>0$, esiste $\gamma$ tale che
\begin{gather*}
 (x\cdot y)(M_{\gamma})\subset x(M_{\gamma}) \cdot y(M_{\gamma})\subset B_{\epsilon}(s+t)
\end{gather*}
da cui per ogni ultrafiltro $\M$, $\phi_{\M}$ è lineare e moltiplicativo (quindi anche continuo) su $l_{\infty}(\N)$. Resta da dimostrare solo che è un'operazione di limite. A questo scopo procediamo per assurdo: sia $\phi_{\M}(x)>\limsup_n x(n)$, e scegliamo $p$ in modo che $\phi_{\M}(x)>p>\limsup_n x(n)$. L'insieme $(p,\infty)$ è un intorno di $\phi_{\M}(x)$, quindi contiene un elemento di $x(\M)$, quindi esiste un elemento di $\M$ contenuto in $x^{-1}(p,\infty)$. Dato che $p>\limsup_n(x)$, l'insieme $x^{-1}(p,\infty)$ ha cardinalità finita, quindi per definizione di ultrafiltro non costante è impossibile che $x^{-1}(p,\infty)$ contenga qualunque elemento di $\M$, da cui $\phi_{\M}(x)\leq\limsup_n x(n)$ $\forall x, \ \forall \M$ ultrafiltro non costante.\\
Osservando che $\liminf(x)=-\limsup(-x)$, si ottiene immediatamente anche l'altra disuguaglianza.
\end{proof}
\end{prop}
Quello che abbiamo ottenuto è un'operazione $\phi_{\M}$ che è lineare, moltilicativa, definita su ogni successione limitata e che coincide con il limite standard quando questo è definito. Il limite standard, oltre che ad essere lineare e moltiplicativo, è anche invariante per traslazioni, nel senso che $\lim_{n}x(n)=\lim_n x(n+1)$ quando è definito. Possiamo chiederci se vale una proprietà simile per $\phi_{\M}$.
\begin{oss}
 Nessuna operazione di limite lineare e moltiplicativa è anche invariante per traslazioni
\begin{proof}
 Dimostriamo questa affermazione per assurdo. Consideriamo la successione limitata $x:\N\to \R$, $x(n)\equiv(-1)^n$, e sia $\phi$ un'operazione di limite lineare e moltiplicativa. Dato che $\phi(x)^2=\phi(x\cdot x)=\phi(1)=1$, $\phi(x)=\pm1$. Per linearità $\phi(x(n))+\phi(x(n+1))=\phi(x(n)+x(n+1))=0$. Se $\phi(x(n))=\phi(x(n+1))$, allora dall'ultima uguaglianza si otterrebbe $\phi(x)=0$, assurdo.
\end{proof}
\end{oss}
Ricordiamo che esistono operazioni di limite lineari invarianti per traslazione (che evidentemente non possono essere moltiplicative). La loro esistenza è dimostrata ad esempio su \cite{4} nell'esercizio 4 cap.3 pag.85.

\chapter{Algebra e compattificazione di Royden}
Questo capitolo segue le orme del capitolo III di \cite{5} e dell'articolo \cite{6}. Lo scopo è quello di introdurre l'algebra di Royden su varietà Riemanniane (di dimesione finita qualsiasi) e quindi la compattificazione di Royden delle varietà.
\section{Funzioni di Tonelli}
\subsection{Definizioni e proprietà fondamentali}
Per prima cosa definiamo un'insieme particolare di funzioni, le funzioni di Tonelli. Lo scopo della sezione è dimostrare che queste funzioni su una varietà constituiscono un'algebra di Banach rispetto a una particolare norma. Iniziamo con il definire queste funzioni su rettangoli in $\R^m$
\begin{deph}
 Una funzione $f:\prod_{i=1}^m(a_i,b_i)\to \R$ si dice \textbf{di Tonelli} se soddisfa:
\begin{enumerate}
 \item $f$ è continua su $\prod_{i=1}^m(a_i,b_i)$
 \item per ogni $i$, $f^i(z)=f(\bar x_1,\dots,\bar x_{i-1},z,\bar x_{i+1},\dots,\bar x_m)$ è assolutamente continua rispetto a $z$ per quasi ogni valore di $(\bar x_1,\dots,\bar x_m)$ \footnote{rimandiamo alla sezione \ref{sec_ac} per la definizione e alcune proprietà delle funzioni assolutamente continue}
 \item Per ogni $i=1,\cdots,m$, $\frac{\partial f}{\partial x_i}$ è a quadrato integrabile su ogni sottoinsieme compatto di $\prod_{i=1}^m(a_i,b_i)$, cioè $\frac{\partial f}{\partial x_i}\in L ^2 _{loc} (\prod_{i=1}^m(a_i,b_i))$
\end{enumerate}
\end{deph}
Notiamo che tutte le funzioni $C^{\infty}(\prod_{i=1}^m(a_i,b_i),\R)$ sono ovviamente funzioni di Tonelli. L'insieme delle funzioni di Tonelli è ovviamente uno spazio vettoriale, e in un certo senso è l'insieme pià piccolo di funzioni continue per cui ha senso parlare di integrale di Dirichlet. A questo proposito dimostriamo che:
\begin{prop}
 Data una funzione $f:M\to \R$ che sia di Tonelli su un insieme compatto $K\Subset U$, dove $U=\phi^{-1}(\prod_{i=1}^m(a_i,b_i))$ è un insieme coordinato, allora detto
\begin{gather}
 D_K(f)=\int_K \abs{\nabla (f)}^2 dV = \int_{\phi(K)}  g^{ij} \frac{\partial f}{\partial x^i}\frac{\partial f }{\partial x^j} \sqrt{\abs g}dx^1\dots dx^n
\end{gather}
si ha che $D_K(f)$ è finito.
\begin{proof}
 Dalla definizione sappiamo che $\int_{\phi(K)}  \sum_{i=1}^m\abs{\frac{\partial f}{\partial x^i}}^2 dx^1\dots dx^n<\infty$. Iniziamo la dimostrazione col notare che per ogni punto $p$ e ogni m-upla di numeri $(y_1,\dots,y_m)$, esiste una costante $c_p$ (indipendente dai numeri $y_i$) tale che:
\begin{gather}\label{eq_g}
 c_p^{-1}\sum_{i=1}^m (y_i)^2\leq \sum_{i,j=1}^m g^{ij}y_i y_j  \leq c_p\sum_{i=1}^m (y_i)^2
\end{gather}
questo è vero grazie al fatto che $\sum_{i,j=1}^m g^{ij}y_i y_j$ definisce una norma su $\R^m$, e tutte le norme negli spazi finito-dimensionali sono equivalenti. Inoltre questa costante $c_p$ dipende con continuità da $p$ \footnote{le funzioni $g^{ij}$ dipendono con continuità da $p$}, quindi su ogni compatto $K\Subset M$, la funzione $c_p$ ha un massimo, che indicheremo semplicemente con $c$. Da queste considerazioni risulta che:
\begin{gather}\label{eq_g2}
 g^{ij} \frac{\partial f}{\partial x^i}\frac{\partial f }{\partial x^j} \leq \sum_{i=1}^m c \abs{\frac{\partial f}{\partial x^i}}^2
\end{gather}
dove gli indici ripetuti si intendono sommati. A questo punto, dato che su ogni compatto anche $\sqrt{\abs g}$ assume massimo finito (che indicheremo $G$), si ottiene facilmente che:
\begin{gather*}
 D_K(f)=\int_{\phi(K)}  g^{ij} \frac{\partial f}{\partial x^i}\frac{\partial f }{\partial x^j} \sqrt{\abs g}dx^1\dots dx^n\leq c G \sum_{i=1}^m\int_{\phi(K)}\abs{\frac{\partial f}{\partial x^i}}^2 dx^1\dots dx^m<\infty
\end{gather*}

\end{proof}

\end{prop}

La definzione di funzione di Tonelli può essere estesa in maniera naturale anche a funzioni $f:M\to \R$
\begin{deph}
 Una funzione $f:M\to \R$ si dice \textbf{di Tonelli} se è di Tonelli in ogni insieme parametrizzabile come rettangolo.
\end{deph}
La proposizione precedente assicura che su ogni compatto coordinato (contenuto in un rettangolo coordinato per la precisione), l'integrale di Dirichlet di una funzione di Tonelli è finito. La stessa cosa vale anche per un qualsiasi insieme compatto.
\begin{prop}
 Per ogni compatto $K\Subset M$, $D_K(f)<\infty$ se $f$ è una funzione di Tonelli
\begin{proof}
 La dimostrazione segue tecniche standard della geometria differenziale. Per ogni punto $p$ della varietà $M$ è possibile trovare un intorno compatto coordinato che descritto in carte locali abbia la forma $\prod_{i=1}^m[a_i,b_i]$. Su ciascuno di questi insiemi l'integrale di Dirichlet di $f$ è finito grazie alla proposizione precedente. Dato che l'insieme $K$ può essere ricoperto da un numero finito di questi insiemi, anche l'integrale di Dirichlet su $K$ sarà finito.
\end{proof}
\end{prop}
\subsection{Operazioni con le funzioni di Tonelli}
Oltre che somma e moltiplicazione per un numero reale, le funzioni di Tonelli ammettono anche altre operazioni tra loro. In particolare il valore assoluto di una funzione di Tonelli è di Tonelli, e quindi anche massimo e minimo tra due funzioni di Tonelli sono funzioni di Tonelli. In questo paragrafo ci occuperemo di queste operazioni.
\begin{prop}\label{prop_ton_1}
 Se $f$ è di Tonelli, anche $\abs f$ è di Tonelli. Inoltre su ogni insieme misurabile $S\subset M$, $D_S(\abs f)\leq D_S(f)$.
\begin{proof}
 Consideriamo una funzione $f:\prod_{i=1}^m (a_i,b_i)\to \R$ di Tonelli. Ovviamente $\abs f$ è una funzione continua. L'assoluta continuità delle funzioni $\abs{f}^i$ è garantita dal fatto che $\abs{ \abs{f}^i(z_1)-\abs f ^i (z_2)}\leq\abs{ f^i(z_1)- f ^i (z_2)}$ \footnote{vedi definizione \ref{deph_ac}}. Rimane da verificare solo la maggiorazione degli integrali di Dirichlet. Consideriamo come primo caso $K\Subset \prod_{i=1}^m (a_i,b_i)$.\\
 Sia $E$ un insieme di misura nulla tale che tutte le derivate $\partial \abs f /\partial x^i$ e $\partial f /\partial x^i$ esistano su $L\equiv E^C=\prod_{i=1}^m (a_i,b_i)\setminus E$. Siano ora $L_+=f^{-1}(0,\infty)\cap L$, $L_-=f^{-1}(-\infty,0)\cap L$, $L_0=f^{-1}(0)\cap L$. Questi insiemi sono ovviamente misurabili.\\
Sull'insieme $L_+$, $\partial \abs f /\partial x^i=\partial f/\partial x^i$, infatti per ogni punto $x\in L_+$, grazie al teorema di permanenza del segno esiste un intorno di $x$ sul quale $\abs f =f$. In modo analogo, su $L_-$, $\partial \abs f \partial x^i=-\partial f/\partial x^i$.\\
Su $L_0$ invece $\partial \abs f \partial x^i=0$. Infatti per definizione:
\begin{gather*}
 \frac{\partial \abs f}{\partial x^i}=\lim_{h\to 0} \frac{\abs {f(\bar x^1,\dots,\bar x^{i-1},x+h,\bar x^{i+1},\dots \bar x^m)}-0}{h}
\end{gather*}
dove abbiamo sfruttato il fatto che $x\in L_0$. Data la definizione di $L$, questo limite esiste necessariamente. Visto che scegliendo $h>0$ il limite è $\geq0$ e scegliendo $h<0$ succede il contrario, questo limite è necessariamente nullo.\\
Consideriamo ora un compatto $K\Subset \prod_{i=1}^m (a_i,b_i)$. Per additività degli integrali:
\begin{gather*}
 \int_K \abs{ \nabla \abs f}^2 dV=\int_{K\cap L_+}\abs{ \nabla \abs f}^2 dV+\int_{K\cap L_-}\abs{ \nabla \abs f}^2 dV+\int_{K\cap L_0}\abs{ \nabla \abs f}^2 dV
\end{gather*}
Visto che in coordinate $$\abs{\nabla f}^2 =g^{ij}\frac{\partial f}{\partial x^i}\frac{\partial f}{\partial x^j}$$ si ottiene che
\begin{gather*}
 \int_{K\cap L_+}\abs{ \nabla \abs f}^2 dV=\int_{K\cap L_+}\abs{ \nabla f}^2 dV \\
 \int_{K\cap L_-}\abs{ \nabla \abs f}^2 dV=\int_{K\cap L_-}\abs{ \nabla f}^2 dV \\
 \int_{K\cap L_0}\abs{ \nabla \abs f}^2 dV=0\\
\end{gather*}
da cui si ottiene $D_K(\abs f)\leq D_K(f)$.\\
Per un insieme $S\subset M$ qualsiasi, basta applicare la definizione di integrale con le partizioni dell'unità. A questo scopo sia $\{\lambda_n\}$ una partizione dell'unità di $M$ subordinata a un ricoprimento di aperti coordinati. Per definizione
\begin{gather*}
 \int_M \abs{\nabla f}^2 dV=\sum_{n=1}^{\infty} \int_{supp(\lambda_n)} \lambda_n\cdot\abs{\nabla f}^2 dV
\end{gather*}
con un ragionamento analogo a quello illustrato sopra, si ottiene $D_S(\abs f)\leq D_S(f)$ per ogni $S\subset M$ misurabile.
\end{proof}
\end{prop}
Grazie a questa proposizione, è immediato verificare che
\begin{prop}
 Date due funzioni di Tonelli $f$ e $g$, anche $\max\{f,g\}$ e $\min\{f,g\}$ sono funzioni di Tonelli
\begin{proof}
 La dimostrazione è immediata se si considera che
\begin{gather*}
\max\{f,g\}= \frac 1 2 (f+g)+\frac 1 2 \abs{f-g} \ \ \ \min\{f,g\}= \frac 1 2 (f+g)-\frac 1 2 \abs{f-g}
\end{gather*}

\end{proof}

\end{prop}

\section{Algebra di Royden}
In questa sezione ci occupiamo dell'algebra di Royden su una varietà Riemanniana $R$ e delle sue proprietà.
\subsection{Definizione}
\begin{deph}\label{deph_Ro}
 Data una varietà riemanniana $R$, definiamo l'algebra di Royden $\mathbb{M}(R)$ l'insieme delle funzioni $f:R\to \R$ tali che:
\begin{enumerate}
 \item $f$ è una funzione continua e limitata su $R$
 \item $f$ è una funzione di Tonelli
 \item $D_R(f)=\int_{R}\abs{\nabla f}^2 dV<\infty$
\end{enumerate}
\end{deph}
Notiamo subito che $\Ro$ è un'algebra commutativa di funzioni, cioè un insieme chiuso rispetto a somma, prodotto e prodotto per uno scalare.
\begin{prop}
 $\Ro$ è un'algebra commutativa di funzioni
\begin{proof}
 La verifica che $\Ro$ è chiusa rispetto a somma e prodotto per uno scalare è ovvia. Rimane solo da verificare che date due funzioni $f,\ g\in \Ro$, anche il loro prodotto appartiene all'algebra.\\
Per prima cosa, se $f$ e $h$ sono continue e limitate, $f\cdot h$ è continua e limitata da $\norm f _{\infty}\cdot \norm h _{\infty}$ \footnote{dove $\norm f _{\infty}$ è la norma del sup}. Se $f(\bar x^1,\dots,\bar x^{i-1},x,\bar x^{i+1},\dots,\bar x^m)$ e $h(\bar x^1,\dots,\bar x^{i-1},x,\bar x^{i+1},\dots,\bar x^m)$ sono entrambe assolutamente continue rispetto a $x$ (eventualità che si verifica quasi ovunque rispetto alle variabili barrate), allora la proposizione \ref{prop_ac_times} garantisce che anche il prodotto sia assolutamente continuo, quindi esiste quasi ovunque la derivata del prodotto e vale che:
\begin{gather*}
 \frac{\partial (fh)}{\partial x^i}=\frac{\partial f}{\partial x^i}h+f\frac{\partial h}{\partial x^i}
\end{gather*}
dove questa uguaglianza si intende quasi ovunque, quindi praticamente sempre dal punto di vista integrale.\\
Mancano da verificare le proprietà sull'integrale di Dirichlet. Iniziamo con il verificare che dato un compatto $K$ contenuto in un rettangolo coordinato, $D_K(fh)<\infty$. Questa verifica segue dalla catena di disuguaglianze:
\begin{gather*}
 D_K(fh)=\int_{\phi(K)} g^{ij}\ton{\frac{\partial f}{\partial x^i}h+f\frac{\partial h}{\partial x^i}}\ton{\frac{\partial f}{\partial x^j}h+f\frac{\partial h}{\partial x^j}}\sqrt{\abs g}dx^1\dots dx^m=\\
=\int_{\phi(K)} g^{ij}\frac{\partial f}{\partial x^i}\frac{\partial f}{\partial x^j}h^2 dV+
\int_{\phi(K)} g^{ij}\frac{\partial h}{\partial x^i}\frac{\partial h}{\partial x^j}f^2 dV+
2\int_{\phi(K)} g^{ij}\frac{\partial f}{\partial x^i}\frac{\partial h}{\partial x^j}fh dV
\end{gather*}
dove $dV=\sqrt{\abs g}dx^1\dots dx^m$ per semplicità di notazione. Considerando che la forma quadratica $g^{ij}$ è sempre definita positiva, l'argomento dei primi due integrali è sempre positivo, quindi il loro modulo può facilmente essere maggiorato da
\begin{gather*}
 \abs{\int_{\phi(K)} g^{ij}\frac{\partial f}{\partial x^i}\frac{\partial f}{\partial x^j}h^2 dV+
\int_{\phi(K)} g^{ij}\frac{\partial h}{\partial x^i}\frac{\partial h}{\partial x^j}f^2 dV}\leq\\
\leq \norm h_{\infty}^2 D_K(f)+\norm f_{\infty}^2 D_K(h)
\end{gather*}
Per l'ultimo integrale, applichiamo due volte la disuguaglianza di Schwartz (vedi ad esempio paragrafo 10.8 pag 210 di \cite{10}) in modo da ottenere:
\begin{gather*}
 \abs{\int_{\phi(K)} g^{ij}\frac{\partial f}{\partial x^i}\frac{\partial h}{\partial x^j}fh dV}\leq\norm f_{\infty}\norm h_{\infty}\int_{\phi(K)} \abs{g^{ij}\frac{\partial f}{\partial x^i}\frac{\partial h}{\partial x^j}} dV\leq\\
\leq\norm f_{\infty}\norm h_{\infty}\int_{\phi(K)} \abs{g^{ij}\frac{\partial f}{\partial x^i}\frac{\partial f}{\partial x^j}}^{1/2}\abs{g^{ij}\frac{\partial h}{\partial x^i}\frac{\partial h}{\partial x^j}}^{1/2} dV\leq\\
\leq\norm f_{\infty}\norm h_{\infty}\ton{\int_{\phi(K)} g^{ij}\frac{\partial f}{\partial x^i}\frac{\partial f}{\partial x^j}dV}^{1/2} \ton{\int_{\phi(K)}g^{ij}\frac{\partial h}{\partial x^i}\frac{\partial h}{\partial x^j}dV}^{1/2}=\\ \\
=\norm f_{\infty}\norm h_{\infty} D_K(f)^{1/2}D_K(g)^{1/2}
\end{gather*}
dove dalla prima alla seconda riga si sfrutta il fatto che $g^{ij}$ è definita positiva, quindi $g^{ij}x_iy_j$ è in prodotto scalare tra $x$ e $y$ \footnote{quindi possiamo applicare $\abs {\ps x y}\leq \norm x \norm y$}, e nella terza riga sfruttiamo la disuguaglianza di Schwartz per integrali.\\
Riassumendo, con queste disuguaglianze otteniamo che
\begin{gather}\label{DK_fg}
 D_K(fh)\leq\norm h_{\infty}^2 D_K(f) + 2\norm f_{\infty}\norm h _{\infty} (D_K(f)D_K(h))^{1/2} + \norm f_{\infty}^2 D_K(h)=\\
=\ton{\norm h_{\infty}D_K(f)^{1/2}+\norm f_{\infty}D_K(h)^{1/2}}^2\notag
\end{gather}
Sfruttanto le partizioni dell'unità, con un argomento del tutto analogo a quello utilizzato nella dimostrazione di \ref{prop_ton_1}, si ottiene la stessa disuguaglianza anche per gli integrali di Dirichlet estesi a tutta $R$, quindi:
\begin{gather}\label{D-fg}
  D_R(fh)\leq\ton{\norm h_{\infty}D_R(f)^{1/2}+\norm f_{\infty}D_R(h)^{1/2}}^2
\end{gather}

\end{proof}

\end{prop}
L'algebra $\Ro$ è quindi un'algebra dotata di unità (la funzione costante uguale a 1). Ha senso chiedersi quali siano i suoi elementi invertibili. \`E chiaro che l'inversa di una funzione $f$ è necessariamente la funzione $f^{-1}=1/f$, che esiste ed è continua solo se $f\neq 0$ ovunque. Però questa funzione è limitata solo se $\inf\abs f>0$. Questo suggerisce che
\begin{prop}\label{prop_inv}
 Data $f\in \Ro$, $f^{-1}\in \Ro$ se e solo se $\inf \abs f >0$.
\begin{proof}
 Supponiamo che $\inf \abs f>0$. Allora $f^{-1}$ esiste ed è continua e limitata. Inoltre se $f$ è derivabile in $x$, lo è anche $f^{-1}$, con
\begin{gather*}
 \left.\frac{\partial f^{-1}}{\partial x^i}\right|_{f(x)}=\left.-\frac{1}{f^2}\frac{\partial f}{\partial x^i}\right|_{x}
\end{gather*}
Dato che $\abs f\geq \inf \abs f>0$, si ottiene che
\begin{gather*}
 D_S(f^{-1})\leq \frac{1}{(\inf \abs f)^4} D_S(f)
\end{gather*}
per ogni sottoinsieme $S$ misurabile di $R$.\\
Supponiamo invece che $\inf \abs f =0$. Se esiste $x$ tale che $f(x)=0$, $f^{-1}$ non è definita, ma anche se la funzione non assume mai il valore $0$, $\inf \abs f =0 \Rightarrow f^{-1}$ non limitata.
\end{proof}
\end{prop}
\subsection{Topologie sull'algebra di Royden}
In questa sezione ci occupiamo di definire alcune topologie sull'algebra di Royden $\Ro$ e di caratterizzarne alcune proprietà, in particolare la completezza. Con una di queste topologie, $\Ro$ diventa un'algebra di Banach.\\
La prima topologia che introduciamo è $\tau_C$, la topologia della convergenza uniforme sui compatti.
\begin{deph}\label{deph_C}
 Data l'algebra $\Ro$, definiamo una base per una topologia $\tau_C$ gli insiemi della forma
\begin{gather*}
 V(f,\epsilon,K)=\{h\in\Ro \ t.c. \ \ \norm{f-h}_{\infty,K}\equiv \max_{p\in K} \{\abs{f(p)-h(p)}<\epsilon\}\}
\end{gather*}
dove $K\Subset R$ qualsiasi.
\end{deph}
Lasciamo al lettore la facile verifica che questa è una base per una topologia. Si nota subito che $\tau_C$ è primo numerabile, infatti fissata $f$, e fissata un'esaustione $K_n$ di compatti in $R$, una base di intorni è ad esempio
\begin{gather*}
 V_n(f)=V\ton{f,\frac 1 n, K_n}
\end{gather*}
Seguendo la traccia di \cite{4} (1.38 (c) a pagina 29), si ottiene che questa topologia è metrizzabile. Se $R$ è compatta, questa topologia è anche normabile, in caso contrario no poiché non è localmente limitata (vedi teorema 1.9 pag. 9 di \cite{4}). Una metrica per $\tau_C$ può essere ad esempio definita da
\begin{gather*}
 d(f,h)=\max_n \frac{1}{n}\frac{\norm{f-h}_n}{1+\norm{f-h}_n} \ \ \ \norm{f-h}_n=\max_{p\in K_n} \abs{f(p)-h(p)}
\end{gather*}
Sia nel caso $R$ compatta che nel caso $R$ non compatta questa metrica non è completa. L'idea è che l'algebra di Royden contiene solo funzioni derivabili in qualche senso, e una topologia che non chiede alcun controllo sulle derivate non può essere completa.
\begin{oss}\label{oss_cmpl_1}
 La topologia $\tau_C$ non è completa su $\Ro$.
\begin{proof}
 Le funzioni lisce limitate a supporto compatto appartengono all'algebra di Royden. Questo insieme però è denso nell'insieme delle funzioni continue su $R$ rispetto a $\tau_C$, quindi come sottospazio $\Ro$ non può essere chiuso, quindi $(\Ro,\tau_c)$ non può essere uno spazio completo.
\end{proof}
\end{oss}
La topologia $\tau_C$ essendo metrizzabile può essere caratterizzata dal comportamento delle sue successioni convergenti
\begin{prop}
 Rispetto alla topologia $\tau_C$, $f_n\to f$ se e solo se $f_n$ converge localmente uniformemente ad $f$, cioè se e solo se
\begin{gather*}
 \forall K\Subset R \ \ \lim_n \max_K \abs{f_n(p)-f(p)}=0
\end{gather*}
In questo caso scriviamo che
\begin{gather}
 f=C-\lim_n f_n
\end{gather}

\end{prop}
Un altro modo di definire una topologia vettoriale su $\Ro$ è il seguente:
\begin{deph}\label{deph_B}
 Fissata $K_n$ un'esaustione di $R$ \footnote{ricordiamo che esaustione significa che $K_n$ sono compatti, $K_n\Subset K_{n+1}^{\circ}$ e $\cup_n K_n=R$}, definiamo una base per la topologia $\tau_B$ di intorni di $0$ come:
\begin{gather}\label{eq_V1}
 V(0,E_n)=\{f:R\to \R \ t.c. \ \max_{p\in K_n} \abs{f(p)}< E_n\}
\end{gather}
dove $E_n$ è una qualunque successione di numeri strettamente positivi tali che $$\lim_{n\to \infty} E_n= \infty$$ La topologia su $R$ è univocamente determinata dalla sua invarianza per traslazioni.
\end{deph}
Visto che si verifica facilmente che
\begin{gather*}
 V(0,E_n)\cap V(0,'E_n) = V(0,\min\{E_n,E'_n\})
\end{gather*}
e che il minino di due successioni strettamente positive e tendenti a infinito mantiene queste proprietà, la definizione di base è ben posta.\\
Osserviamo subito che la definizione non dipende dalla scelta di $K_n$ (cioè al variare di $K_n$, tutte le topologie generate sono equivalenti tra loro), e nella definizione si può fare l'ulteriore richiesta che $E_n$ sia una successione crescente senza perdere di generalità.\\
\`E inoltre possibile descrivere in maniera alternativa questa topologia. A questo scopo, data una funzione $e:R\to \R$, diciamo che
\begin{gather*}
 \lim_{p\to \infty} e(p) =+\infty \ \ \Longleftrightarrow \ \ \forall N\in \N, \ \exists K\Subset R \ t.c. \ e(p)>N \ \forall p\not \in K
\end{gather*}
Con questa notazione la topologia $\tau_B$ può essere descritta da una base di intorni fatta così:
\begin{gather*}
 V(0,e)=\{f:R\to \R \ t.c. \ \abs{f(p)}< e(p)\}
\end{gather*}
Al variare di $e$ tra le funzioni continue positive che tendono a infinito, questa base di intorni definisce la stessa topologia definita dagli intorni in \ref{eq_V1}.\\
Un altro modo differente per definire questa topologia è il seguente. Dato $N\in \N$, definiamo l'insieme
\begin{gather*}
 X_N\equiv \{f\in \Ro \ t.c. \ \norm f _{\infty}\leq N\}
\end{gather*}
Dotiamo $X_N$ \footnote{che ovviamente non è uno spazio vettoriale} della topologia indotta da $(\Ro,\tau_C)$. Diciamo che un insieme $V$ è aperto in $\tau_B$ se e solo se $\forall N$ l'intersezione $V\cap X_N$ è aperto nella topologia $(X_N,\tau_C)$. Lasciamo al lettore la verifica che questa topologia è una topologia vettoriale.
\begin{prop}
 La topologia generata da questi aperti è la topologia $\tau_B$ definita sopra.
\begin{proof}
 Consideriamo $V$ aperto intorno di $0$ secondo la vecchia definizione. Allora esiste una successione $E_n\to \infty$ tale che
\begin{gather*}
 V(0,E_n)\subset V
\end{gather*}
L'intersezione $V(0,E_n)\cap X_N$ è l'insieme delle funzioni tali che:
\begin{gather*}
 V(0,E_n)\cap X_N =\{f\in X_N\ t.c. \ \norm{f}_{\infty,\ K_n}\leq E_n \ \forall n \}
\end{gather*}
Consideriamo $E_n$ successione crescente. Dato che $E_n\to \infty$, esiste $\bar n$ tale che $E_{n} \geq N$ per ogni $n\geq \bar n$, questo significa che
\begin{gather*}
 \norm f _{\infty}\leq N \ \Rightarrow \ \norm{f}_{\infty, \ K_{\bar n}^C}\leq E_n
\end{gather*}
per ogni $n\geq \bar n$. Quindi possiamo scrivere:
\begin{gather*}
 V(0,E_n)\cap X_N =\{f\in X_N \ t.c. \ \norm{f}_{\infty,\ K_n}\leq E_n \ \forall n \leq \bar n \}\supset \\
\supset\{f\in X_N\ t.c. \ \norm{f}_{\infty,\ K_{\bar n}}\leq E_1\}
\end{gather*}
e questo è per definizione un insieme della topologia $(X_N,\tau_C)$.\\
Sia ora un insieme aperto intorno di $0$ secondo la nuova definizione, cioè tale che per ogni $N$:
\begin{gather*}
 V\cap X_N \in \tau_C
\end{gather*}
Allora per $N=1$, esiste un insieme $K_1$ e un numero $E_1$ tale che
\begin{gather*}
 \{f \ t.c. \ \norm f_{\infty,\ K_1}\leq E_1\} \subset V\cap X_1\subset V
\end{gather*}
e lo stesso per ogni valore di $N$. Questo significa che esiste una successione $E_N$ tale che:
\begin{gather*}
 \{f\ t.c. \ \norm f_{\infty, \ K_N}\leq E_N \ \forall N\}\subset V
\end{gather*}
Supponiamo per assurdo che $E_N\not\to \infty$, quindi $E_N$ limitata da $M$. Questo è impossibile perchè altrimenti l'insieme:
\begin{gather*}
 \{f\ t.c. \ \norm f_{\infty, \ K_N}\leq E_N \ \forall N\}\cap X_M
\end{gather*}
non sarebbe aperto nella topologia di $(X_M,\tau_C)$.
\end{proof}

\end{prop}

La topologia $\tau_B$ ha lo svantaggio di non essere I numerabile (se $R$ non è compatta), quindi non è una topologia metrizzabile e soprattutto non può essere descritta in maniera completa dal comportamento delle sue successioni convergenti.
\begin{oss}
 $\tau_B$ non è I numerabile, cioè fissato un punto in $\Ro$, non esiste una base numerabile di intorni del punto.
\begin{proof}
 Supponiamo per assurdo che esista una base numerabile $V_k$ di intorni di $0$. Questi intorni avranno la forma
\begin{gather*}
 V_k\equiv V(0,E^{(k)}_n)
\end{gather*}
dove per ogni $k$ la successione $\{E^{(k)}_n\}$ è strettamente positiva e divergente. Costruiamo per induzione una successione strettamente crescente $k_m$ di interi tali che
\begin{gather*}
 \forall n\geq k_m, \ \ E^{(m)}_n \geq m
\end{gather*}
e definiamo la successione $B_m$ strettamente positiva e tendente a infinito come:
\begin{gather*}
 k_m\leq n < k_{m+1} \ \Rightarrow \ B_n=m/2
\end{gather*}
i termini $B_1,\cdots,B_{k_1 -1}$ possono essere assegnati casualmente.\\
Dalla definizione, otteniamo che per nessun valore di $k$ $B_n\leq E^{(k)}_n \ \forall n$, quindi non esiste $k$ per cui $V(0,E^{(k)}_n)\subset V(0,B_n)$, che contraddice il fatto che $\{V(0,E^{(k)}_n)\}$ sia una base di intorni di $0$.
\end{proof}
\end{oss}
Le successioni convergenti in questa topologia sono caratterizzate da
\begin{prop}
 Una successione $f_n$ converge a $f$ rispetto alla topologia $\tau_B$ se e solo se:
\begin{enumerate}
 \item esiste $M$ tale che $\norm{f_n}_{\infty}\leq M \ \forall n$
 \item $f=C-\lim_n f_n$
\end{enumerate}
\begin{proof}
 Dimostriamo prima che se valgono (1) e (2), allora la successione converge rispetto a $\tau_B$. Vista l'invarianza per traslazioni della topologia, è sufficiente verificare questa condizione con $f=0$. Consideriamo un qualsiasi intorno aperto di $0$ della forma $V(0,E_n)$. Dato che $E_n\to \infty$, esiste un numero $N_1$ tale che $$E_n>M \ \forall n\geq N_1$$ Inoltre poiché $f_n$ converge localmente uniformemente a $0$, esiste un numero $N_2$ tale che
\begin{gather*}
 \forall n \geq N_2, \ \max_{p\in K_{N-1}} \{\abs{f_n(p)}\} \leq \min_{1\leq i \leq N-1} \{E_i\}
\end{gather*}
quindi se $N= \max\{N_1,N_2\}$, $f_n\in V \ \forall n\geq N$, cioè $f_n$ converge a $f$ rispetto a $\tau_B$.\\
Per l'implicazione inversa, dato che $\tau_C\subset \tau_B$, la convergenza nella topologia $B$ implica la convergenza locale uniforme. Inoltre supponiamo per assurdo che $\norm{f_n}_\infty \to \infty$, quindi esiste $\{x_n\}\subset R$ tale che $\abs{f_n(x_n)}\to \infty$. Dato che $f_n$ converge localmente uniformemente a $0$, necessariamente $x_n\to \infty$ \footnote{cioè $x_n$ abbandona definitivamente ogni compatto}. Fissata un'esaustione $K_n$ di $R$, per ogni $n$ consideriamo
\begin{gather*}
 E_n =\min_{x_i\in K_n} \{\abs{f_i(x_i)}\}/2
\end{gather*}
Non è difficile verificare che $E_n\to \infty$, e che quindi a meno di un numero finito di termini \footnote{che in questo ragionamento sono ininfluenti}, $E_n>0$.\\
In questo modo la successione $\{f_n\}$ non è contenuta definitivamente nell'aperto $V(0,E_n)$, quindi non può convergere a $0$.
\end{proof}

\end{prop}

Indipendentemente dal fatto che $\tau_B$ non è primo numerabile, ha senso chiedersi se le successioni bastano a descrivere la topologia di questo spazio \footnote{spazi topologici di questo genere si chiamano spazi di \textit{Frechet-Urysohn}. Visti gli scopi della tesi, ci limitiamo a rimandare a \cite{17} capitolo I pag 13 per ulteriori approfondimenti}. La risposta a questa domanda è negativa, anzi qualunque topologia generi la convergenza $B$ è destinata a non essere una topologia ``strana''. Infatti in questi spazi la chiusura per successioni non è idempotente, cioè indicando $[A]_{seq}$ la chiusura per successioni di $A$, cioè l'insieme di tutti i possibili limiti di successioni in $A$, si ha che in generale $[[A]_{seq}]_{seq}\neq [A]_{seq}$, quindi necessariamente $[A]_{seq}\neq \overline A$.
\begin{prop}
 Data una successione $f_n\in \Ro$, se definiamo
\begin{gather*}
 f=B-\lim_n f_n \Longleftrightarrow f=C-\lim_n f_n \ \wedge \ \ton{\exists M \ t.c. \ \norm{f_n}_{\infty}\leq M}
\end{gather*}
allora se $R$ non è compatta, qualunque topologia generi questa convergenza non è di Frechet-Urysohn, quindi non è descrivibile con il comportamento delle sue successioni convergenti.
\begin{proof}
 La dimostrazione è ispirata all'esercizio 9 cap 3 pag 87 di \cite{4}. Consideriamo una successione discreta $\{x_n\}\subset \R$, $x_n\to \infty$, e sia $\{U_n\}$ una collezione di intorni disgiunti di $x_n$. Siano $g_n\in \Ro$ funzioni tali che $supp(g_n)\subset U_n$, $0\leq g_n\leq 1$ e $g_n(x_n)=1$. Definiamo $f_{n,m}\in \Ro$ come $f_{n,m}=g_n+ng_m$, e consideriamo
\begin{gather*}
 A=\{f_{n,m} \ n,m\in \N\}
\end{gather*}
Allora $0\not \in [A]_{seq}$, ma $0\in [[A]_{seq}]_{seq}$. Infatti qualunque successione $\{f_{n(k),m(k)}\}_{k=1}^{\infty}$ in $A$ che converge nel senso $B$ a $0$ necessariamente converge localmente uniformemente a $0$, quindi $n_k\to \infty$, ma allora anche $\norm{f_{n(k),m(k)}}_{\infty}=n_k\to \infty$, il che significa che questa successione non può convergere a $0$.\\
D'altro canto è facile verificare che tutte le funzioni $g_n\in [A]_{seq}$, basta considerare che la successione $f_{n(k),m(k)}\to g_n$ se $n(k)=n$, $m(k)\to \infty$. Ma dato che $g_n\to 0$ rispetto alla convergenza $B$, la tesi è dimostrata.
\end{proof}
\end{prop}
Questa proposizione impone di trattare la convergenza $B$ con particolare attenzione, perché non può essere descritta semplicemente dalle sue successioni (come accade per gli spazi metrici).

Un'altra topologia che possiamo definire su $R$ è la classica topologia della convergenza uniforme ovunque:
\begin{deph}\label{deph_U}
 La norma del sup è una norma sullo spazio $\Ro$. Definiamo la topologia descritta da questa norma $\tau_U$. Diciamo che 
\begin{gather}
f=U-\lim_n f_n 
\end{gather}
se $f_n$ converge ad $f$ rispetto a questa norma, quindi se
\begin{gather*}
 \lim_n \norm{f_n-f}_{\infty}=\lim_n \sup_R \abs{f_n(p)-f(p)}=0
\end{gather*}
\end{deph}
\begin{oss}
 \`E facile verificare che $\tau_C\subset \tau_B\subset \tau_U$. Le tre topologie coincidono nel caso $R$ compatta, mentre l'inclusione è stretta se $R$ non è compatta.
\end{oss}
Anche in questo caso vale un'osservazione molto simile a \ref{oss_cmpl_1} (e la dimostrazione è del tutto analoga, basta sostituire lo spazio delle funzioni continue con lo spazio delle funzioni continue limitate).
\begin{oss}
 La norma dell'estremo superiore non rende $\Ro$ uno spazio completo.
\end{oss}
Dalla definizione di $\Ro$, non è difficile immaginare che per rendere questa algebra un'algebra di Banach bisogna in qualche senso tenere in considerazione il comportamento delle derivate delle funzioni. A questo scopo introduciamo il concetto di $D-\lim$:
\begin{deph}
 Data una successione di funzioni $f_n\in \Ro$, diciamo che
\begin{gather}
 D-\lim_n f_n =f
\end{gather}
se e solo se $D_R(f_n-f)\to 0$
\end{deph}
Ricordiamo che l'integrale di Dirichlet è molto legato al concetto di norma nello spazio $\mathcal{L}^2(R)$, lo spazio di Hilbert delle 1-forme a quadrato integrabili su $R$ \footnote{a questo scopo rimandiamo alla sezione \ref{sec_L2}}.\\
Questi concetti di convergenza possono essere mischiati fra loro:
\begin{deph}
 Data una successione $f_n\in \Ro$, diciamo che
\begin{gather}
 f=QD-\lim_n f_n
\end{gather}
se e solo se $f=D-\lim_n f_n$ e anche $f=Q-\lim_n f_n$, dove $Q$ rappresenta una convergenza qualsiasi tra $C$, $B$, $U$.
\end{deph}
I due concetti di convergenza che saranno più usati in questo lavoro sono la $BD$-convergenza e la $UD$-convergenza. La seconda convergenza è generata da una norma, che indicheremo $\norm \cdot _R$, norma che rende $\Ro$ un'algebra di Banach.
\begin{teo}\label{teo_g1}
 Sull'algebra $\Ro$ definiamo la norma
\begin{gather*}
 \norm f _R \equiv \norm f _{\infty} + D_R(f)^{1/2}
\end{gather*}
questa norma genera la convergenza $UD$, e rispetto a questa norma $\Ro$ è un'algebra di Banach commutativa con unità.
\begin{proof}
Il fatto che $\norm \cdot _R$ sia una norma è facile conseguenza del fatto che $\norm \cdot _{\infty}$ è una norma, e $D_R(\cdot)^{1/2}$ è una seminorma.\\
Ovviamente $\norm 1 =1$, e grazie alla relazione (\ref{D-fg}) si ottiene:
\begin{gather*}
 \norm{fh}_R=\norm{fh}_{\infty}+D_R(fh)^{1/2}\leq \norm{f}_{\infty}\norm h _{\infty} + \norm h_{\infty}D_R(f)^{1/2}+\norm f_{\infty}D_R(h)^{1/2}\leq\\
\leq (\norm{f}_{\infty}+D_R(f)^{1/2})(\norm{h}_{\infty}+D_R(h)^{1/2})=\norm f _R \norm h_R
\end{gather*}
questo rende $(\Ro,\norm \cdot _R)$ un'algebra normata. Rimane da verificare la completezza.\\
Sia a questo scopo $\{f_n\}\subset \Ro$ una successione di Cauchy. Allora $\{f_n\}$ è di Cauchy uniforme, quindi esiste una funzione continua limitata sulla varietà  $R$ tale che $\norm{f_n-f}_{\infty}\to 0$. La parte più complicata è dimostrare che questa funzione è di Tonelli e che $D_R(f_n-f)\to 0$.\\
Ricordiamo che $\mathcal{L}^2(R)$, lo spazio delle 1-forme su $R$ normato con l'integrale del modulo quadro della forma, è uno spazio di Hilbert (vedi sezione \ref{sec_L2}). L'ipotesi che $f_n$ sia di Cauchy rispetto a $\norm \cdot _R$ implica che $df_n$ sia una successione di Cauchy nello spazio $\mathcal{L}^2(R)
$, quindi esiste una 1-forma $\alpha\in \mathcal{L}^2(r)$ tale che $df_n\to \alpha$. Se dimostriamo che $f$ è di Tonelli e $df=\alpha$, abbiamo la tesi. Osserviamo che queste due affermazioni hanno carattere locale, quindi fissiamo un qualunque aperto coordinato $(U,\phi)$ con l'accortezza che $\phi$ sia definita in un intorno di $\overline{U}$ e tale che $\phi(U)=\prod_{i=1}^m(a_i,b_i)$ e dimostriamo l'uguaglianza $df=\alpha$ in questa carta locale.\\
A questo scopo indichiamo con $\alpha_i(x)$ le componenti locali di $\alpha$, cioè:
\begin{gather*}
 \alpha=\alpha_i(x)dx^i
\end{gather*}
mentre per semplicità di notazione continuiamo a indicare con $f_n(x)$ e $f(x)$ le rappresentazioni locali di $f_n$ ed $f$ rispettivamente.\\
Dividiamo la dimostrazione in due parti: l'idea e i conti. L'idea è dimostrare $f=f^{(i)}$, dove $f^{(i)}$ sono funzioni per definizione assolutamente continue e la cui derivata parziale rispetto a $x^i$ è proprio $\alpha_i$. A questo scopo mostriamo che le funzioni $f_n$, che convergono uniformemente a $f$, convergono in norma $L^2(\phi(U))$ a $f^{(i)}$ (per tutti gli indici $1\leq i \leq m$). Quindi $f=f^{(i)}$ quasi ovunque (vedi proposizione \ref{prop_lp_1}). E questo conclude la dimostrazione.\\
Ora passiamo a conti. Definiamo le funzioni $f^{(i)}$ come:
\begin{gather}\label{eq_f_1}
 f^{(i)}(x^1,\cdots,x^i,\cdots,x^m)\equiv f(x^1,\cdots,c_i,\cdots,x^m)+\int_{c_i}^{x^i}\alpha_i(x^1,\cdots,t,\cdots,x^m) dt
\end{gather}
dove $c_i=(a_i+b_i)/2$ \footnote{in realtà per la dimostrazione va bene qualsiasi $c_i\in (a_i,b_i)$}.
Osserviamo che $\partial f_n/\partial x^i(x)$ converge in norma $L^2$ a $\alpha_i(x)$ quando $n$ tende a infinito. Infatti
\begin{gather*}
 \int_{\phi(U)} \ton{\frac{\partial f_n}{\partial x^i}(x)-\alpha_i(x)}^2 dx^1\cdots dx^m\leq
 \int_{\phi(U)} \sum_{i=1}^m \ton{\frac{\partial f_n}{\partial x^i}(x)-\alpha_i(x)}^2 dx^1\cdots dx^m\leq\\
\leq \int_{\phi(U)} C g^{ij}\ton{\frac{\partial f_n}{\partial x^i}(x)-\alpha_i(x)}\ton{\frac{\partial f_n}{\partial x^j}(x)-\alpha_j(x)} dx^1\cdots dx^m\leq\\
\leq \frac{C}{k}\int_{\phi(U)} g^{ij}\ton{\frac{\partial f_n}{\partial x^i}(x)-\alpha_i(x)}\ton{\frac{\partial f_n}{\partial x^j}(x)-\alpha_j(x)} \ \sqrt{\abs g} dx^1\cdots dx^m =\\
=\frac{C}{k} \norm{df_n-\alpha}_{\mathcal{L}^2(R)}\to 0
\end{gather*}
dove per passare dalla prima alla seconda riga abbiamo sfruttato la relazione \ref{eq_g}, scegliendo come $C$ il massimo dei $c_p$ al variare di $p\in \phi(\overline U)$, mentre $k$ è il minimo della funzione $\sqrt{ \abs g}$ sempre sullo stesso insieme compatto \footnote{per questo motivo è importante scegliere $\phi$ definita su un intorno di $\overline U$}.\\
Fissato $i$, tutte le funzioni $f_n(\bar x^i,\cdots,x^i,\cdots,\bar x^m)$ sono contemporaneamente assolutamente continue rispetto a $x^i$, quasi ovunque rispetto a $(\bar x^1,\cdots,\bar x^{i-1},\bar x^{i+1},\cdots,\bar x^m)$ \footnote{ogni funzione $f_n$ è assolutamente continua a meno di un insieme di misura nulla di $\bar x $, ma visto che l'unione numerabile di insiemi di misura nulla ha misura nulla, tutte le funzioni sono assolutamente continue contemporaneamente sullo stesso insieme con complementare di misura nulla}, quindi vale che
\begin{gather}\label{eq_f_2}
 f_n(\bar x^i,\cdots,x^i,\cdots,\bar x^m)=f_n(\bar x^i,\cdots,c_i,\cdots,\bar x^m)+\int_{c_i}^{x^i} \frac{\partial f_n(\bar x^i,\cdots,t,\cdots,\bar x^m)}{\partial t} dt
\end{gather}
Confrontando questa relazione con la relazione \ref{eq_f_1} otteniamo:
\begin{gather*}
\abs{f_n(\bar x^i,\cdots,x_i,\cdots,\bar x^m)-f^{(i)}(\bar x^i,\cdots,x_i,\cdots,\bar x^m)}^2=\\ \\
=\left|f_n(\bar x^i,\cdots,c_i,\cdots,\bar x^m)-f(\bar x^i,\cdots,c_i,\cdots,\bar x^m) +\right.\\ \\
\left. + \int_{c_i}^{x^i} \frac{\partial f_n(\bar x^i,\cdots,t,\cdots,\bar x^m)}{\partial t} dt-\int_{c_i}^{x^i} \alpha_i(\bar x^i,\cdots,t,\cdots,\bar x^m) dt\right|^2\leq\\ \\
\leq 2\abs{f_n(\bar x^i,\cdots,c_i,\cdots,\bar x^m)-f(\bar x^i,\cdots,c_i,\cdots,\bar x^m) }^2+\\ \\
+2\abs{\int_{c_i}^{x^i} \frac{\partial f_n(\bar x^i,\cdots,t,\cdots,\bar x^m)}{\partial t} dt-\int_{c_i}^{x^i} \alpha_i(\bar x^i,\cdots,t,\cdots,\bar x^m) dt}^2
\end{gather*}
per semplicità di notazione chiamiamo $2\abs{A_n}^2$ la prima riga dopo l'ultimo segno di disuguaglianza, e la seconda $2\abs{B_n}^2$. Grazie alle proprietà dell'integrale e alla disuguaglianza di Schwartz otteniamo:
\begin{gather*}
 \abs{B_n}^2\leq \ton{\int_{c_i}^{x^i} \abs{\frac{\partial f_n(\bar x^i,\cdots,t,\cdots,\bar x^m)}{\partial t} -\alpha_i(\bar x^i,\cdots,t,\cdots,\bar x^m) }dt}^2\leq\\
\leq \ton{b_i-a_i}\int_{a_i}^{b^i} \abs{\frac{\partial f_n(\bar x^i,\cdots,t,\cdots,\bar x^m)}{\partial t} -\alpha_i(\bar x^i,\cdots,t,\cdots,\bar x^m) }^2 dt
\end{gather*}
e questa relazione vale quasi ovunque rispetto a $(\bar x^1,\cdots,\bar x^{i-1},\bar x^{i+1},\cdots,\bar x^m)$.\\
Grazie a queste disuguaglianze ora siamo in grado di dimostrare che la successione $f_n$ (o meglio la successione delle rappresentazioni locali di $f_n$) converge in norma $L^2(\phi(U))$ alla funzione $f^{(i)}$ (con $1\leq i \leq m$ qualsiasi). Infatti:
\begin{gather*}
\int_{\phi(U)}\abs{f_m(x^1,\cdots,x^m)-f^{(i)}(x^1,\cdots,x^m)}^2 dx^1\cdots dx^m \leq\\
\leq 2\abs{A_n}^2 Vol(\phi(U))+\\
+2 (b_i-a_i)\int_{\phi(U)} \left\{\int_{a_i}^{b^i} \abs{\frac{\partial f_n(\bar x^i,\cdots,t,\cdots,\bar x^m)}{\partial t} -\alpha_i(\bar x^i,\cdots,t,\cdots,\bar x^m) }^2 dt        \right\}  dx =\\
=2\abs{A_n}^2 Vol(\phi(U))+2(b_i-a_i)^2 \norm{\frac{\partial f_n}{\partial x^i}-\alpha_i}^2_{L^2(\phi(U))}
\end{gather*}
Sappiamo che $A_n$ tende a 0 (poichè $f_n$ converge uniformemente a $f$), e anche $\norm{\frac{\partial f_n}{\partial x^i}-\alpha_i}^2_{L^2(\phi(U))}$ tende a zero come dimostrato prima.\\
\end{proof}
\end{teo}
Modificando leggermente questa dimostrazione, si può ottenere che:
\begin{prop}\label{prop_r1}
Se $f_n$ è una successione di funzioni in $\Ro$ tale che
\begin{gather*}
 f=C-\lim_n f_n
\end{gather*}
 dove $f$ è una funzione (continua) \textbf{limitata}, e se esiste $K<\infty$ tale che per ogni $n$:
\begin{gather*}
 D_R(f_n)\leq K <\infty
\end{gather*}
allora $f\in \Ro$, $D_R(f)\leq K$ e inoltre esiste una sottosuccessione $n_k$ tale che per ogni $g\in \Ro$:
\begin{gather*}
 D_R(f-f_{n_k},g)\to 0
\end{gather*}
\begin{proof}
 Grazie al teorema \ref{teo_BA}, esiste una sottosuccessione di $f_n$ che continueremo a indicare nello stesso modo tale che $df_n$ converge debolmente a $\alpha \in \mathcal L ^2 (R)$.\\
Consideriamo un rettangolo coordinato $T$ qualsiasi e sia $\phi \in C^{\infty}_0(T)$ lo spazio delle funzioni lisce a supporto compatto in $T$. Grazie a un'integrazione per parti otteniamo che:
\begin{gather*}
 \int_T f_n \frac{\partial \phi}{\partial x^i} dV = -\int_T \phi \frac{\partial f_n}{\partial x^i} dV
\end{gather*}
Ora poiché $df_n\to \alpha$ debolmente, si ha che:
\begin{gather*}
 \int_T \phi \frac{\partial f_n}{\partial x^i} dV \to \int_T \phi \alpha_i dV
\end{gather*}
e dato che $f=C-\lim_n f_n$ si ha che:
\begin{gather*}
 \int_T f_n \frac{\partial \phi}{\partial x^i} dV \to \int_T f \frac{\partial \phi}{\partial x^i} dV
\end{gather*}
Da queste consideraziono otteniamo che:
\begin{gather*}
 \int_T f \frac{\partial \phi}{\partial x^i} dV = \int_T \phi \alpha_i dV
\end{gather*}
Questo significa che le derivate distribuzionali di $f$ rispetto a $x^i$ coincidono con $\alpha_i$, e dato che $\alpha\in \mathcal{L}^2(R)$, grazie al lemma \ref{lemma_wder} (che riportiamo alla fine di questa proposizione), $f$ è di Tonelli con $df=\alpha$ nel senso standard e quindi $D_R(f)<\infty$, il che dimostra che $f\in \Ro$.\\
Il fatto che
\begin{gather*}
 D_R(f-f_{n},g)\to 0
\end{gather*}
è conseguenza diretta del fatto che $df_n$ converge debolmente nel senso di $\mathcal L^2(R)$ a $df$. Inoltre dalla teoria della convergenza debole, sappiamo che $$\norm{df}_{L_2}\leq \limsup_n\norm{df_n}_{L_2}$$ quindi ad esempio
\begin{gather*}
 \norm{df_n}_{L_2}\equiv D_R(f_n)\leq k \ \ \forall n \ \Rightarrow\ \norm{df}_{L_2}\equiv D_R(f)\leq k
\end{gather*}
\end{proof}
\end{prop}
\begin{lemma}\label{lemma_wder}
 Sia $f:T\to \R$ una funzione $f\in L ^2 (T)$ \footnote{osserviamo che tutte le funzioni continue sono a quadrato integrabile su insiemi relativamente compatti di misura finita, quindi sui rettangoli coordinati $T$ rispetto alla forma volume} tale che per tutti gli indici $i=1,\cdots,m$:
\begin{gather*}
 \frac{\partial f}{\partial x_i} \in L^2 (T)
\end{gather*}
dove le derivate sono intese in senso distrubuzionale. Allora $f$ è assolutamente continua rispetto a quasi tutti i segmenti paralleli agli assi coordinati in $T$, e le sue derivate standard coincidono quasi ovunque con le derivate distribuzionali. Quindi se $f$ è continua, è una funzione di Tonelli.
\begin{proof}
 Dato che questo teorema riguarda la teoria delle distribuzioni e gli spazi di Sobolev, che sono argomenti marginali in questa tesi, per la dimostrazione rimandiamo il lettore a \cite{33} (teorema 2.1.4 pag. 44), oppure a \cite{34} (teorema 1.41 pag. 22).
\end{proof}

\end{lemma}

\begin{oss}\label{oss_r1}
 Nella dimostrazione, il fatto che $f$ sia limitata è utile esclusivamente per dimostrare che $f\in \Ro$. Quindi se $f$ non è limitata, valgono tutte le conclusioni del teorema a meno dell'appartenenza all'algebra di Royden.
\end{oss}

\subsection{Densità di funzioni lisce}
In questo paragrafo dimostreremo la densità delle funzioni lisce nello spazio $\Ro$, e utilizzeremo questo risultato per dimostrare formule di Green generalizzate.
\begin{prop}\label{prop_dens}
 Sia $f$ una funzione di Tonelli su $R$. Per ogni $\epsilon>0$, esiste $f_{\epsilon}\in C^{\infty}(R,\R)$ tale che $\norm{f_{\epsilon} -f}_R<\epsilon$ \footnote{in questa proposizione non è richiesto che $f$ sia limitata, e che il suo integrale di Dirichlet $D_R(f)<\infty$, quindi potrebbe non avere senso $\norm{f}_{R}$}. Inoltre se $f$ ha supporto contenuto in un aperto $U\subset R$, anche $f_{\epsilon}$ può essere scelta con supporto contenuto in $U$.
\begin{proof}
Come spesso accade in questi casi, la dimostrazione si divide in due parti, una ``locale'' e una ``globale''. Nella prima parte dimostreremo il risultato per funzioni $f$ a supporto compatto in una carta locale, poi generalizzeremo il risultato utilizzando le partizioni dell'unità.\\
Sia quindi $f$ come nelle ipotesi e anche a supporto compatto in un intorno coordinato $(U,\phi)$ di $R$, e chiamiamo $\tilde f$ la rappresentazione di $f$ in questa carta locale. Consideriamo una successione $\Theta_n:\R^m\to \R$ di nuclei di convoluzione con supporto contenuto in $B_{1/n}(0)$. Grazie al lemma \ref{lemma_conv_2} la successione $\tilde f_n\equiv\Theta_n\ast \tilde f$ (che ha supporto definitivamente contenuto in $\phi(U)$) converge nella norma del sup a $\tilde f$. Inoltre grazie al lemma \ref{lemma_conv_1b}, si ha che:
\begin{gather*}
 \frac{\partial \tilde f_n}{\partial x^i}=\Theta_n \ast \frac{\partial \tilde f}{\partial x^i}
\end{gather*}
Passiamo ora a considerare $D_R(f_n-f)$. Grazie al fatto che esiste un insieme compatto $K$ tale che $supp(\tilde f_n)\subset K\Subset U$ definitivamente rispetto a $n$, e grazie alla relazione \ref{eq_g2}, otteniamo che:
\begin{gather*}
 g^{ij}\ton{\frac{\partial \tilde f_n}{\partial x^i}(x)-\frac{\partial \tilde f}{\partial x^i}(x)}\ton{\frac{\partial \tilde f_n}{\partial x^j}(x)-\frac{\partial \tilde f}{\partial x^j}(x)}\leq c \sum_{i=1}^m \ton{\frac{\partial \tilde f_n}{\partial x^i}(x)-\frac{\partial \tilde f}{\partial x^i}(x)}^2=\\
=c \sum_{i=1}^m \ton{\int_{B_n} \Theta_n(y)\ton{\frac{\partial \tilde f}{\partial x^i}(x-y)-\frac{\partial \tilde f}{\partial x^i}(x)}dV(y)}^2\leq\\
\leq c \sum_{i=1}^m \int_{B_n} \Theta_n(y)\ton{\frac{\partial \tilde f}{\partial x^i}(x-y)-\frac{\partial \tilde f}{\partial x^i}(x)}^2dV(y)
\end{gather*}
dove l'ultimo passaggio è giustificato dalla disuguaglianza di Jensen \footnote{infatti $\int \Theta_n(y) dV(y) =1$. Per riferimenti sulla disuguaglianza di Jensen, vedi teorema 3.3 pag 61 di \cite{12}}.\\
Queste considerazioni (assieme al teorema di Fubuni \footnote{vedi teorema 7.8 pag 140 di \cite{12}}) permettono di concludere:
\begin{gather*}
D_R(f_n-f)= \int_{K} g^{ij}\ton{\frac{\partial \tilde f_n}{\partial x^i}(x)-\frac{\partial \tilde f}{\partial x^i}(x)}\ton{\frac{\partial \tilde f_n}{\partial x^j}(x)-\frac{\partial \tilde f}{\partial x^j}(x)}dV(x)\leq\\
\leq c'\int_{K}dV(x)\sum_{i=1}^m \int_{B_n}dV(y) \Theta_n(y)\ton{\frac{\partial \tilde f}{\partial x^i}(x-y)-\frac{\partial \tilde f}{\partial x^i}(x)}^2=\\
=c' \sum_{i=1}^m\int_{B_n}dV(y) \Theta_n(y) \int_{K}dV(x) \ton{\frac{\partial \tilde f}{\partial x^i}(x-y)-\frac{\partial \tilde f}{\partial x^i}(x)}^2
\end{gather*}
dove $c'=\max_{x\in K} \{c \sqrt{\abs g(x)}\}$.\\
Chiamiamo
\begin{gather*}
 h(y)\equiv \int_{K}dV(x) \ton{\frac{\partial \tilde f}{\partial x^i}(x-y)-\frac{\partial \tilde f}{\partial x^i}(x)}^2
\end{gather*}
Grazie alla continuità dell'operatore traslazione in $L^2(R)$ \footnote{vedi teorema 9.5 pag 183 di \cite{4}}, e dato che $f$ è di Tonelli, possiamo concludere che:
\begin{gather*}
 \lim_{y\to 0} h(y) =0
\end{gather*}
Per definizione questo significa che
\begin{gather*}
 \lim_{n\to \infty} \max_{y\in B_n(0)} \abs{h(y)} = 0
\end{gather*}
E quindi anche:
\begin{gather*}
 \lim_{n\to \infty} \abs{\int_{B_n}dV(y) \Theta_n(y)h(y)}  \leq \lim_{n\to \infty} \max_{y\in B_n(0)} \abs{h(y)} \int_{B_n}dV(y) \Theta_n(y) = 0
\end{gather*}
Questo conclude la dimostrazione che $D_R(f_n-f)\to 0$.\\
Consideriamo ora una funzione $f$ di Tonelli qualsiasi, $\epsilon>0$, e sia $\lambda_n:R\to \R$ una partizione dell'unità liscia subordinata a aperti coordinati e a supporto compatto. Allora per ogni $n$, $f\cdot \lambda_n$ è una funzione del tipo descritto sopra. Scegliamo per ogni $n$ una funzione liscia $f_n$ \footnote{se $supp(f)\subset U$, allora $supp(f\cdot \lambda_n)\Subset U$, quindi possiamo scegliere $supp(f_n) \Subset U$} tale che $\norm{f_n-f\cdot \lambda_n}_R<\epsilon/2^n$. Allora la funzione $f_{\epsilon}\equiv \sum_{n=1}^{\infty} f_n$ soddisfa le richieste fatte.\\
$f_{\epsilon}$ è una funzione liscia per locale finitezza della partizione dell'unità $\lambda_n$, inoltre:
\begin{gather*}
 \norm{f_{\epsilon} -f}_R=\norm{\sum_{n=1}^{\infty} (f_n- f\cdot \lambda_n)}_R\leq \sum_{n=1}^{\infty}\norm{ f_n- f\cdot \lambda_n}_R < \epsilon
\end{gather*}
Inoltre se $supp(f_n)\Subset U$ per ogni $n$, allora anche $supp(f)\subset U$. Questo significa che se $f$ ha supporto compatto, visto che ogni compatto in $R$ è contenuto in un aperto relativamente compatto, allora $f_{\epsilon}$ può essere scelta a supporto compatto.
\end{proof}
\end{prop}
Vale anche una proposizione più forte rispetto a questa. Possiamo infatti chiedere che la funzione $f_{\epsilon}$ sia uguale ad $f$ su un insieme chiuso, ovviamente a patto di rilassare le ipotesi di regolarità sulla funzione $f_{\epsilon}$.
\begin{prop}\label{prop_dens2}
 Sia $f$ una funzione di Tonelli su $R$ e $A\subset R$ aperto con bordo $\partial A$ regolare. Per ogni $\epsilon>0$, esiste $f_{\epsilon}\in C^{\infty}(A,\R)$ di Tonelli su tutta la varietà tale che $\norm{f_{\epsilon} -f}_R<\epsilon$ e $f_{\epsilon}=f$ sull'insieme $A^C$.
\begin{proof}
 Riportiamo solo il filo conduttore della dimostrazione, lasciando alcuni dettagli al lettore. Questa dimostrazione è ispirata dal lemma 2.8 pagina 50 di \cite{18}.\\
Sia $K_n$ un ricoprimento di aperti relativamente compatti localmente finiti in $A$ \footnote{con localmente finito si intende un ricoprimento tale che per ogni $p\in A$ esiste un intorno che interseca solo un numero finito di insiemi $K_n$}, e sia $\epsilon_n$ una successione di numeri positivi che tende a 0 tale che $\epsilon_n\leq \epsilon$. Lo scopo è riuscire a creare una funzione $h_{\epsilon}\in C^\infty (A,\R)$ tale che
\begin{gather}\label{eq_hepsilon}
 \norm{h_{\epsilon}-f}_{K_n}\equiv \norm{h_{\epsilon}-f}_{\infty, \ K_n} + D_{K_n}(h_{\epsilon}-f)^{1/2}<\epsilon_n
\end{gather}
In questo modo se estendiamo la definizione di $h_{\epsilon}$ come
\begin{gather*}
 f_{\epsilon}(p)=\begin{cases}
               h_{\epsilon}(p)  & p\in A \\
	      f(p) & p\in A^C
              \end{cases}
\end{gather*}
questa funzione risulta di Tonelli su tutta $R$. Infatti $f_{\epsilon}$ è ovviamente continua e di Tonelli su tutti i punti di $R$ tranne che sul bordo $\partial A$. Consideriamo $p\in \partial A$, e consideriamo una successione $p_k\in A$, $p_k\to p$. Visto che $p_k$ converge al bordo di $A$, necessariamente $p_k\in K_{n(k)}$ dove $n(k)$ tende a infinito quando $k$ tende a infinito. Questo significa che
\begin{gather*}
 \abs{f_{\epsilon}(p_k)-f(p)}\leq \abs{f_{\epsilon}(p_k)-f(p_k)}+\abs{f(p_k)-f(p)}\leq \epsilon_{n(k)} +\abs{f(p_k)-f(p)}
\end{gather*}
 se $k$ tende a infinito, $\epsilon_{n(k)}$ tende a zero per ipotesi, mentre il secondo termine tende a zero per continuità di $f$. Questo prova che $f_{\epsilon}$ è continua.\\
Consideriamo un punto $p\in \partial A$ e una parametrizzazione locale $\phi(q)=(x^1,\cdots x^m)$ attorno a $p$ tale che $\partial A = \{x^m=0\}$. Allora la funzione $f_{\epsilon}$ è di certo assolutamente contiua quasi ovunque nell'insieme $x^m\neq 0$, quindi quasi ovunque. L'integrale di Dirichlet di $f_{\epsilon}$ è finito su ogni compatto che ha intersezione vuota con il bordo perché $f_{\epsilon}$ su questo insieme è una funzione liscia o una funzione di Tonelli per ipotesi. Se il compatto interseca il bordo di $A$, allora poiché $D_R(f-f_{\epsilon})$ è finito, necessariamente l'integrale di Dirichlet di $f_{\epsilon}$ su questo compatto è finita.\\
Resta da dimostrare che è possibile scegliere una funzione $h_{\epsilon}$ che soddisfi le condizioni \ref{eq_hepsilon}. La tecnica che utilizzeremo è un'adattamento della dimostrazione della proposizione \ref{prop_dens}.\\
Scegliamo una partizione dell'unità $\{\lambda_n\}$ di $A$ subordinata al ricoprimento $K_n$. Per ogni $n$, definiamo $$\alpha_m\equiv\frac{\epsilon_m}{\#n \ t.c. \ supp(\lambda_n)\cap K_m \neq \emptyset}$$ $$\delta_n\equiv \min\{\alpha_m \ t.c. \ K_m\cap supp(\lambda_n)\neq \emptyset\}$$
Per locale finitezza del ricoprimento $\{K_m\}$ tutte queste quantità sono ben definite e strettamente positive.\\
Ora grazie alle tecniche esposte nella dimostrazione della proposizione \ref{prop_dens}, scegliamo per ogni $n$ una funzione liscia $f_n$ con supporto in $K_n$ tale che:
\begin{gather*}
 \norm{f_n - f\lambda_n}_R < \delta_n
\end{gather*}
e definiamo $h_\epsilon \equiv \sum_n f_n$. Per locale finitezza della partizione dell'unità questa serie è ben definita e $h_\epsilon \in C^{\infty}(A,\R)$, inoltre:
\begin{gather*}
\norm{h_\epsilon - f}_{K_m} = \norm{\sum_n (f_n-f\lambda_n)}_{K_m} < \epsilon_m
\end{gather*}
da cui la tesi.
\end{proof}

\end{prop}

\subsection{Formule di Green e principio di Dirichlet}\label{sec_dirpri2}
Ricordando la notazione definita nella sezione \ref{sec_int}, grazie ai risultati di densità appena descritti possiamo dimostrare che:
\begin{prop}\label{prop_green2}
 Se $f$ è di Tonelli con $D_R(f)<\infty$, e $u\in H(\Omega)\cap C^\infty(\overline \Omega)$ \footnote{$H(\Omega)$ è l'insieme delle funzioni armoniche su $\Omega$. $C^\infty(\overline \Omega)$ è l'insieme delle funzioni che possono essere estese a funzioni lisce in un intorno di $\overline \Omega$}, dove $\Omega$ è un dominio regolare, abbiamo che:
\begin{gather*}
 D_{\Omega}(f,u)\equiv \int_{\Omega} \ps{\nabla f}{\nabla u} dV = \int_{\partial \Omega} f \ast du
\end{gather*}
\begin{proof}
 Consideriamo una successione $f_n\in C^{\infty}(R,\R)$ che converge in norma $\norm \cdot _R$ a $f$. Allora grazie alla proposizione \ref{prop_green1}, abbiamo che per ogni $m$:
\begin{gather*}
 D_{\Omega}(f_m,u)\equiv \int_{\Omega} \ps{\nabla f_m}{\nabla u} dV = \int_{\partial \Omega} f_m \ast du
\end{gather*}
Poiché $D_R(f_m-f)\to0$ e anche $\norm{f_m-f}_{\infty}\to0$, si ha che:
\begin{gather*}
 \abs{D_{\Omega}(f_m,u)-D_{\Omega} (f,u)}\equiv \abs{\int_{\Omega} \ps{\nabla (f_m-f)}{\nabla u} dV} \leq D_{\Omega}(u)D_{\Omega}(f_m-f)\to 0 
\end{gather*}
\begin{gather*}
 \abs{\int_{\partial \Omega} f_m \ast du - \int_{\partial \Omega} f_m \ast du}\leq\int_{\partial \Omega} \abs{f_m-f} \ast du \to 0
\end{gather*}
da cui la tesi.
\end{proof}
\end{prop}
Possiamo anche rilassare le ipotesi su $u$, e ottenere:
\begin{prop}\label{prop_green3}
 Sia $\Omega\subset R$ un dominio regolare, $f$ una funzione di Tonelli con $D_R(f)<\infty$, e $u\in H(\Omega)$ con $D_{\Omega}(u)<\infty$. Siano inoltre $\gamma_1$ e $\gamma_2$ due insiemi connessi disgiunti tali che $\gamma_1\cup \gamma_2 =\partial \Omega$. Se $f=0$ su $\gamma_1$, e $u\in C^\infty(\gamma_2)$, allora:
\begin{gather}\label{eq_green1}
 D_{\Omega}(f,u)\equiv \int_{\Omega}\ps{\nabla f}{\nabla u} dV = \int_{\gamma_2} f\ast du
\end{gather}
\begin{proof}
Per cominciare dimostriamo questa proposizione in un caso particolare. Sia $g$ una funzione di Tonelli con integrale di Dirichlet finito su $\Omega$, $g\geq 0$ su $\Omega$, $g|_{\gamma_1}=0$ e $g|_{\gamma_2}\geq \delta >0$.\\
Consideriamo ora una successione di insiemi $\Omega_n$ con bordo liscio tali che $\Omega_n\subset \Omega$, $\partial \Omega_n= \gamma_2\cup \beta_n$ e $\cup_n\Omega_n=\Omega$ \footnote{è possibile costruire una successione con tali caratteristiche ad esempio considerando una funzione liscia $h:\overline \Omega\to \R$ tale che $h(\gamma_1)=0$, $h(\gamma_2)=1$ e $0<h(x)<1$ per ogni $x\in \Omega$. Scegliendo una successione decrescente $r_n\searrow 0$ di valori regolari per $h$, gli insiemi $\Omega_n \equiv h^{-1}(r_n,1)$ hanno le caratteristiche cercate} e definiamo le funzioni
\begin{gather*}
 g_c(x)\equiv \max\{g(x)-c,0\}
\end{gather*}
per $0<c<\delta$. Per costruzione, queste funzioni coincidono con $g(x)-c$ su un intorno di $\gamma_2$ e sono tutte nulle in un intorno di $\gamma_1$, quindi fissato $c$, $g_c$ si annulla identicamente su $\beta_n$ definitivamente in $n$.\\
Grazie a queste considerazioni, e grazie alla proposizione \ref{prop_green2}, sappiamo che esiste $n$ per cui:
\begin{gather*}
D_\Omega(g_c,u)= D_{\Omega_n} (g_c,u)=\int_{\beta_n} g_c\ast du + \int_{\gamma_2} g_c \ast du= \int_{\gamma_2} g\ast du -c\int_{\gamma_2}\ast du
\end{gather*}
Facendo tendere $c$ a $0$, la seconda parte della disuguaglianza tende a $\int_{\gamma_2} g\ast du$, mentre la prima tende a $D_{\Omega}(g,u)$. Infatti se chiamiamo $A_c=\{x\in \Omega \ t.c. \ g(x)\leq~c\}$
\begin{gather*}
 \abs{D_{\Omega}(g_c,u)-D_{\Omega}(g,u)}=\abs{\int_{\Omega} \ps{\nabla(g_c-g)}{\nabla u} dV} =\abs{\int_{A_c} \ps{\nabla g}{\nabla u} dV}\leq\\
\leq \int_{A_c} \abs{\nabla g}^2 dV \int_{A_c} \abs{\nabla u}^2 dV \to 0
\end{gather*}
e dato che entrambi gli integrali di Dirichlet estesi a tutta $\Omega$ sono finiti, se l'insieme la misura dell'insieme di integrazione tende a zero, anche l'integrale tende a 0.\\
Questo dimostra la tesi su $g|_{\gamma_2}\geq\delta <0$. Se consideriamo una funzione $f\geq 0$ qualsiasi, la tesi si ottiene per linearità, infatti basta scegliere una funzione $g$ con le caratteristiche descritte sopra e applicare quanto appena per ottenere:
\begin{gather*}
 D_\Omega(f,u)=D_\Omega((f+g)-g,u)=D_\Omega(f+g,u)-D_\Omega (g,u)=\\
=\int_{\gamma_2} (f+g)\ast du - \int_{\gamma_2} g\ast du
\end{gather*}
Se la funzione $f$ non è positiva, la tesi si ottiene applicando il ragionamento appena esposto alle funzioni $f^+$ e $f^-$, in particolare:
\begin{gather*}
 D_\Omega(f,u)=D_\Omega(f^+-f^-,u)=\int_{\gamma_2} (f^+-f^-)\ast du =\int_{\gamma_2} f\ast du
\end{gather*}
\end{proof}
\end{prop}
Ora siamo in grado di dimostrare una generalizzazione del principio di Dirichlet esposto in \ref{prop_D1}.
\begin{prop}\label{prop_D2}
 Sia $f:R\to \R$ una funzione di Tonelli con $D_R(f)<\infty$, $u\in H(\Omega)\cap C(\overline \Omega)$, dove $\Omega$ è un dominio regolare in $R$ varietà riemanniana. Se $f|_{\partial \Omega}=u|_{\partial \Omega}$, allora $D_\Omega (u)<\infty$ e:
\begin{gather*}
 D_{\Omega}(f)=D_{\Omega}(u)+ D_{\Omega}(f-u)
\end{gather*}
\begin{proof}
 La tesi di questa proposizione è esattamente identica alla tesi della proposizione \ref{prop_D1}, la differenza è nelle richieste di regolarità su $f$, che in questo caso sono molto meno stringenti. Come negli esempi sopra, dimostreremo questo teorema approssimando $f$ con funzioni lisce.\\
Sia $f_m$ una successione di funzioni tali che $\norm{f_m-f}_R \to 0$, e siano $u_m$ le soluzioni del problema di Dirichlet su $\Omega$ con $f_m|_{\partial \Omega}$ come dato al bordo. La proposizione \ref{prop_D1} afferma che:
\begin{gather*}
  D_{\Omega}(f_m)=D_{\Omega}(u_m)+ D_{\Omega}(f_m-u_m)
\end{gather*}
Grazie al principio del massimo otteniamo anche che $$\norm{u_m-u}_{\infty}\leq \norm{f_m|_{\partial \Omega}-f|_{\partial \Omega}}_{\infty}\to 0$$ quindi $u=U-\lim_m u_m$. Inoltre:
\begin{gather*}
 D_{\Omega}(u_m-u_k)=D_{\Omega}(f_m-f_k)-D_{\Omega}(f_m-f_k+u_k-u_m)\leq D_{\Omega}(f_m-f_k)
\end{gather*}
e dato che $f_m\to f$ e che:
\begin{gather*}
 D_{\Omega}(f_m-f_k)^{1/2}\leq D_{\Omega}(f_m-f)^{1/2}+D_{\Omega}(f_k-f)^{1/2}
\end{gather*}
allora la successione $\{u_m\}$ è di Cauchy rispetto alla norma $\norm \cdot _R$, quindi per completezza di $\Ro$, $u_m$ converge e necessariamente converge a $u$. Questo implica che $D_{\Omega}(u)<\infty$. Osserviamo che:
\begin{gather*}
 D_{\Omega}(f)=\int_{\Omega} \abs{\nabla f}^2 dV=\int_{\Omega} \abs{\nabla (f-u)+u}^2 dV=\\
=\int_{\Omega} \abs{\nabla (f-u)}^2 dV +\int_{\Omega} \abs{\nabla u}^2 dV + 2 \int_{\Omega} \ps{\nabla (f-u)}{\nabla u}dV=\\ \\
=D_{\Omega}((f-u))+D_{\Omega}(u)+2D_{\Omega}((f-u),u)   
\end{gather*}
Applicando il lemma \ref{prop_green3} con $\gamma_1=\partial \Omega$ otteniamo che:
\begin{gather*}
 D_{\Omega}((f-u),u)= \int_{\Omega} \ps{\nabla (f-u)}{\nabla u}dV=0
\end{gather*}
dato che per ipotesi $(f-u)|_{\partial \Omega}=0$.
\end{proof}
\end{prop}
In realtà le ipotesi sulla regolarità del bordo di $\Omega$ non sono necessarie, infatti vale la proposizione
\begin{prop}\label{prop_D3}
  Sia $f:R\to \R$ una funzione di Tonelli con $D_R(f)<\infty$, $u\in~H(\Omega)\cap C(\overline \Omega)$, con $\Omega$ è un aperto relativamente compatto in $R$ varietà riemanniana. Se $f|_{\partial \Omega}=u|_{\partial \Omega}$, allora $D_\Omega (u)<\infty$ e:
\begin{gather*}
 D_{\Omega}(f)=D_{\Omega}(u)+ D_{\Omega}(f-u)
\end{gather*}
\begin{proof}
Per comodità, estendiamo l'insieme di definizione di $u$ a tutto $R$, ponendo $u|_{\Omega^C}=f|_{\Omega^C}$. In questo modo $u$ ha tutte le proprietà di una funzione in $\Ro$, tranne il fatto di avere integrale di Dirichlet finito \footnote{proprietà che seguirà dalla dimostrazione}.\\
Sia $\Omega_n$ un'esaustione regolare per $\Omega$, definiamo le funzioni:
\begin{gather*}
 u_n:R\to \R \ \ u_n\in H(\Omega_n)\cap C(R)\ \ \ \ u_n|_{\partial \Omega_n}=f|_{\partial \Omega_n} \ \ \ \ u_n|_{\Omega_n^C}=f|_{\Omega_n^C}
\end{gather*}
Dalla proposizione precedente è evidente $\forall n$ $D_\Omega(u_n)<\infty$, quindi $u_n\in \Ro$, inoltre per ogni $m>n$ vale che:
\begin{gather*}
 D_{\Omega_m}(u_n)=D_{\Omega_m}(u_m)+D_{\Omega_m}(u_n-u_m)
\end{gather*}
quindi considerando che $u_n|_{\Omega_m^C}=u_n|_{\Omega_m^C}=f|_{\Omega_m^C}$, si ha che $D_{\Omega_m^C}(u_n)=D_{\Omega_m^C}(u_m)$ e $D_{\Omega_m^C}(u_n-u_m)=0$, quindi:
\begin{gather*}
 \underbrace{D_{\Omega_m}(u_n)+D_{\Omega_m^C}(u_n)}_{D_R(u_n)}=\underbrace{D_{\Omega_m}(u_m)+D_{\Omega_m^C}(u_m)}_{D_R(u_m)}+\underbrace{D_{\Omega_m}(u_n-u_m)+D_{\Omega_m^C}(u_n-u_m)}_{D_R(u_n-u_m)}
\end{gather*}
Grazie a un ragionamento simile a quello riportato nella dimostrazione del teorema \ref{teo_1}, si ha che la successione $\{u_n\}$ è $D$-cauchy. Inoltre grazie al principio del massimo si può dimostrare che $u_n$ converge uniformemente a $u$ su $\Omega$, infatti%
\begin{gather*}
 \abs{u(x)-u_n(x)}\begin{cases}
        =0 & \text{se } x\in \Omega^C\\
=\abs{f(x)-u(x)} & \text{se } x\in \Omega\setminus \Omega_n\\
\leq \max_{x\in \partial \Omega_n}\{\abs{f(x)-u(x)}\} & \text{se } x\in \Omega_n
       \end{cases}
\end{gather*}
dove l'ultima riga segue dal principio del massimo. Da questa relazione ricaviamo la stima:
\begin{gather*}
 \norm{u-u_n}_\infty =\max_{x\in \Omega\setminus \Omega_n}\{\abs{u(x)-f(x)}\}
\end{gather*}
dato che $\Omega$ è compatto e $f|_{\partial \Omega}=u|_{\partial \Omega}$, la parte sinistra dell'ultima uguaglianza tende a $0$ quando $n$ tende a infinito, da cui $u$ è il limite uniforme di $u_n$.\\
Dato che la successione $\{u_n\}$ è di Cauchy rispetto alla metrica $CD$, converge in questa metrica, e per unicità del limite converge a $u$, in particolare:
\begin{gather*}
 \lim_{n\to \infty}D_R(u-u_n)=0
\end{gather*}
Per dimostrare la seconda parte della proposizione, applichiamo ancora una volta la proposizione precedente, ottenendo che:
\begin{gather*}
 D_{\Omega_n}(f)=D_{\Omega_n}(u_n)+D_{\Omega_n}(f-u_n)
\end{gather*}
e quindi anche
\begin{gather*}
 D_\Omega(f)=D_{\Omega}(u_n)+D_{\Omega}(f-u_n) 
\end{gather*}
La tesi si ottiene passando al limite per $n$ che tende a infinito.
\end{proof}

\end{prop}

\subsection{Ideali dell'algebra di Royden}\label{subsec_ideals}
In questo paragrafo parleremo di alcuni ideali dell'algebra di Royden.\\
Ricordiamo la definizione di ideale di un'algebra (commutativa)
\begin{deph}
 Data un'algebra $A$, un suo sottoinsieme $I\subset A$ è detto \textbf{ideale} se:
\begin{enumerate}
\item $I$ è un sottospazio vettoriale di $A$
\item per ogni $x\in A$ e $y\in I$, $xy\in I$
\end{enumerate}
\end{deph}
Un esempio di ideale su $\Ro$ è l'insieme delle funzioni che si annullano in un punto $p$. Questo ideale è per altro anche il nucleo del carattere $\tau(p)$. Gli ideali a cui saremo interessati in questa sezione però sono: $\Roo$ e $\Rod$, dove:
\begin{gather*}
 \Roo\equiv \{f\in \Ro \ t.c. \ supp(f) \ compatto\}
\end{gather*}
mentre $\Rod$ è l'insieme di tutte le funzioni che sono $BD-$limiti di successioni di funzioni in $\Roo$ \footnote{quindi $M_{\Delta}(R)$ è la \textbf{chiusura sequenziale} dell'insieme $M_0(R)$}. \`E abbastanza facile verificare che $\Roo$ è un ideale di $\Ro$, meno banale è la verifica che anche $\Rod$ è un ideale.
\begin{prop}
 $\Rod$ è un ideale di $\Ro$.
\begin{proof}
 Per dimostrare questa affermazione consideriamo una funzione $f\in \Rod$ e $h\in \Ro$. Allora esiste una successione $f_n\in \Roo$ tale che $$f=BD-\lim_n f_n$$ Ovviamente la successione $\{hf_n\}\subset \Roo$. Dimostriamo che la successione $hf_n$ converge nella topologia $BD$ a $hf$, in questo modo otteniamo la tesi.\\
Poiché $h$ è continua e limitata, la successione $hf_n$ converge localmente uniformemente a $hf$ ed è uniformemente limitata (quindi $hf=B-\lim_n hf_n$). Resta da dimostrare che $hf=D-\lim_n hf_n$. A questo scopo osserviamo che:
\begin{gather*}
 D_R(hf-hf_n)=\int_{R}g^{ij} \ton{\frac{\partial (hf-hf_n)}{\partial x^i}}\ton{\frac{\partial(hf- hf_n)}{\partial x^j}} dV =\\
=\int_{R}g^{ij} \ton{\frac{\partial h}{\partial x^i} f + h\frac{\partial f}{\partial x^i} - \frac{\partial h}{\partial x^i} f_n - h\frac{\partial f_n}{\partial x^i} }\ton{\frac{\partial h}{\partial x^j} f+ h \frac{\partial f}{\partial x^j}-\frac{\partial h}{\partial x^j} f_n- h \frac{\partial f_n}{\partial x^j}} dV  =\\
\int_{R}g^{ij} \ton{\frac{\partial h}{\partial x^i}( f-f_n) + h\ton{\frac{\partial f}{\partial x^i} -\frac{\partial f_n}{\partial x^i}} }\ton{\frac{\partial h}{\partial x^j}( f-f_n) + h\ton{\frac{\partial f}{\partial x^j} -\frac{\partial f_n}{\partial x^j}} } dV \leq \\
\leq 2 \int_R g^{ij} \frac{\partial h}{\partial x^i}\frac{\partial h}{\partial x^j} (f-f_m)^2 dV + 2 \int_R g^{ij} \ton{\frac{\partial f}{\partial x^i}-\frac{\partial f_n}{\partial x^i}}\ton{\frac{\partial f}{\partial x^j}-\frac{\partial f_n}{\partial x^j}}h^2 dV
\end{gather*}
dove l'ultima disuguaglianza segue dal fatto che per ogni norma $$\norm{A+B}^2\leq 2\ton{\norm A ^2 + \norm B ^2}$$ applicando questa disuguaglianza alla norma $\norm A = g^{ij}A_i A_j$ si ottiene il risultato.\\
Ora consideriamo che per ogni $K\Subset R$:
\begin{gather*}
 \int_R g^{ij} \frac{\partial h}{\partial x^i}\frac{\partial h}{\partial x^j} (f-f_m)^2 dV = \\
=\int_K g^{ij} \frac{\partial h}{\partial x^i}\frac{\partial h}{\partial x^j} (f-f_m)^2 dV + \int_{R\setminus K} g^{ij} \frac{\partial h}{\partial x^i}\frac{\partial h}{\partial x^j} (f-f_m)^2 dV
\end{gather*}
e poiché $f_n$ converge localmente uniformemente a $f$, il primo integrale converge a $0$ quando $m$ tende a infinito, mentre il secondo è limitato in modulo da:
\begin{gather*}
 \abs{\int_{R\setminus K} g^{ij} \frac{\partial h}{\partial x^i}\frac{\partial h}{\partial x^j} (f-f_m)^2 dV}\leq N D_{R\setminus K}(h)
\end{gather*}
dove $N$ è tale che $\abs{f_n(x)-f(x)}^2\leq N \ \forall n \ \forall x$ \footnote{poiché $f$ è limitata e $f_n$ è uniformemente limitata, questo numero esiste}.\\
Inoltre abbiamo che:
\begin{gather*}
\abs{ \int_R g^{ij} \ton{\frac{\partial f}{\partial x^i}-\frac{\partial f_n}{\partial x^i}}\ton{\frac{\partial f}{\partial x^j}-\frac{\partial f_n}{\partial x^j}}h^2 dV}\leq \\
\leq M\int_R g^{ij} \ton{\frac{\partial f}{\partial x^i}-\frac{\partial f_n}{\partial x^i}}\ton{\frac{\partial f}{\partial x^j}-\frac{\partial f_n}{\partial x^j}} dV=M D_R(f_n-f)
\end{gather*}
dove $M=\norm g _{\infty}^2$. Dato che $D_R(f_n-f)\to 0$, otteniamo che:
\begin{gather*}
 \limsup_{n\to \infty} D_R(hf-hf_n) \leq 2N D_{R\setminus K} (h)
\end{gather*}
 data l'arbitrarietà di $K$ e dato che $D_R(h)<\infty$, otteniamo che:
\begin{gather*}
 \lim_{n\to \infty} D_R(hf-hf_n) =0
\end{gather*}
Quindi $BD-\lim_n gf_n =gf$, il che dimostra che $gf\in \Rod$.
\end{proof}
\end{prop}
Vale una proposizione simile a \ref{prop_r1} per lo spazio $\Rod$:
\begin{prop}\label{prop_r2}
 Sia $f_n$ una successione in $\Rod$ tale che
\begin{gather*}
 f=C-\lim_n f_n
\end{gather*}
con $f$ funzione (continua) limitata e
\begin{gather*}
 D_R(f_n)\leq K
\end{gather*}
dove $K$ non dipende da $n$. Allora $f\in \Rod$.
\begin{proof}
 Grazie al teorema \ref{prop_r1}, $f\in \Ro$ e esiste una sottosuccessione $f_{n_k}$ (che per comodità continueremo a indicare con $f_n$) tale che per ogni $g\in \Ro$
\begin{gather*}
 D_R(f-f_n;g)\to 0
\end{gather*}
Sia $R_n$ un'esaustione regolare di $R$, e definiamo le funzioni $\phi_n$ in modo che:
\begin{enumerate}
 \item $\phi_n=0$ su $R\setminus \overline{R_2}$
 \item $\phi_n=f_n$ su $R_1$
 \item $\phi_n\in H(R_2\setminus \overline{R_1})$
\end{enumerate}
e sia $\phi$ definita in maniera analoga per $f$.\\
Dato che per ogni $n$, $\phi_n\in \Roo$, al posto che $f$ e $f_n$, nella dimostrazione possiamo considerare le funzioni $f-\phi$ e $f_n-\phi_n$, cioè possiamo assumere per ipotesi che $f_n=f=0$ sull'insieme $R_1$.\\
Sia $u_m$ una successione di funzioni definite da:
\begin{enumerate}
 \item $u_m=0$ su $R_1$
 \item $u_m=f$ su $R_m\setminus R_1$
 \item $u_m\in H(R_m\setminus \overline{R_1})$
\end{enumerate}
Chiaramente $u_m\in \Ro$, e per il principio del massimo $\norm{u_m}_{\infty}\leq \norm{f}_{\infty}$. Poiché $u_m=u_n=f$ su $\partial R_p$ e $u_p$ è armonica su $R_p$, grazie al principio di Dirichlet (vedi \ref{prop_D2}) otteniamo che se $p>m$:
\begin{gather*}
 D_R(u_p-u_m)=D_R(u_m)-D_R(u_p)
\end{gather*}
questo dimostra che la successione $D_R(u_m)$ è decrescente al crescere di $m$ \footnote{infatti se $p>m$, $D_R(u_m)-D_R(u_p)\geq0$}, e ovviamente limitata dal basso da $0$, quindi è convergente, e quindi se $m$ è sufficientemente grande e $p>m$, $D_R(u_m-u_p)$ è piccolo a piacere, cioè $\{u_m\}$ è $D$-Cauchy.\\
Grazie al principio \ref{prop_harnackpri}, esiste una sottosuccessione di $u_m$ (che continueremo a indicare nello stesso modo) tale che:
\begin{gather*}
 u=B-\lim_n u_n \ \ \Rightarrow \ \ u=BD-\lim_n u_n
\end{gather*}
allora $u\in \Ro$ con $u(R_1)=0$ e $u\in H(R\setminus R_1)$.\\
Dato che $f-u=BD-\lim_n (f-u_n)$, e per costruzione $(f-u_n)\in \Roo$, $f-u\in \Rod$.\\
Grazie alla formula di Green \ref{prop_green3}, si ha che:
\begin{gather*}
 D_R(f-u_n;u) =0
\end{gather*}
poiché $f-u_n=0$ sull'insieme $R\setminus R_n$, e anche su $R_1$ (quindi anche sui relativi bordi). Quindi anche $D_R(f-u;u)=\lim_n D_R(f-u_n;u)=0$. Quindi otteniamo che:
\begin{gather*}
 D_R(f-u;u)=0 \ \ \Rightarrow \ \ D_R(u)\equiv D_R(u;u)=D_R(f;u)=\lim_n D_R(f_n;u)
\end{gather*}
Dato che $f_n\in \Rod$ per ipotesi, per ogni $n$ esiste una successione $h_k$ di funzioni in $\Roo$ tali che $f_n=BD-\lim_k h_k$. Questo permette di osservare che
\begin{gather*}
 D_R(f_n;u)=\lim_k D_R(h_k;u) =\lim_k \int_{\partial R_1} h_k\ast du
\end{gather*}
dove l'ultima uguaglianza è conseguenza di \ref{prop_green3}. Dato che $\norm{h_k}_{\infty, \partial R_1}\to 0$ (poiché la funzione $f_n=0$ su $\partial R_1$), l'ultimo limite è $0$, quindi $D_R(u)=0$, il che significa che $u\equiv 0$, e quindi $f=f-u\in \Rod$, come volevasi dimostrare.
\end{proof}
\end{prop}
Come corollario a questa proposizione, osserviamo che:
\begin{oss}\label{oss_rodcmpl}
Se $\{f_n\}\subset \Rod$ e $f=BD-\lim_n f_n$, allora $f\in \Rod$.
\end{oss}

\section{Compattificazione di Royden}
Ora siamo pronti per introdurre il concetto di compattificazione di Royden. Data una varietà riemanniana $R$, la sua compattificazione di Royden $R^*$ è uno spazio compatto che in qualche senso contiene $R$. La compattificazione di Royden è l'iniseme dei caratteri sull'algebra $\Ro$, quindi l'integrale di Dirichlet delle funzioni gioca un ruolo fondamentale nella costruzione di questa compattificazione. Questo fa intuire che proprietà come parabolicità e iperbolicità vengano in qualche modo riflesse nelle proprietà di $R^*$.
\subsection{Definizione}
\begin{deph}\label{deph_R^*}
 Data una varietà riemanniana $(R,g)$, si definisce compattificazione di Royden uno spazio $R^*$ tale che:
\begin{enumerate}
 \item $R^*$ è uno spazio compatto di Hausdorff
 \item $R$ è un sottoinsieme aperto e denso di $R^*$
 \item ogni funzione in $\Ro$ può essere estesa per continuità a una funzione definita su $R^*$
 \item l'insieme $\overline \Ro$ delle funzioni che sono estensione delle funzioni in $\Ro$ separa i punti di $R^*$
\end{enumerate}
\end{deph}
Osserviamo subito che se $R$ è compatta, necessariamente $R^*=R$.
\begin{teo}\label{teo_Reu}
 Per ogni varietà riemanniana $R$ esiste $R^*$ ed è unica a meno di omeomorfismi che mantengaono fissi i punti di $R$.
\end{teo}
Spezziamo la dimostrazione di questo teorema in alcuni lemmi per facilitarne la lettura, e riportiamo nell proposizione \ref{prop_R_sumup} lo schema riassuntivo della dimostrazione.\\
Sia $R^*$ l'insime dei funzionali lineari moltiplicativi su $\Ro$ dotato della topologia debole-* rispetto a $\Ro$. Ovviamente questo spazio è uno spazio di Hausdorff \footnote{sia $p\neq q$, allora per definizione esiste $f\in \Ro$ tale che $\abs{p(f)- q(f)}=\delta>0$. Quindi gli intorni $V(p)=\{h\in R^* \ t.c. \ \abs{h(f)-p(f)}<\delta/3\}$ e $W(q)=\{h\in R^* \ t.c. \ \abs{h(f)-q(f)}<\delta/3\}$ sono intorni disgiunti dei due caratteri $p$ e $q$.}, e grazie alla proposizione \ref{prop_charK} è anche uno spazio compatto. $R^*$ è il candidato a compattificazione di Royden di $R$.\\
Per ogni punto $p\in R$, il funzionale $x_p(f)\equiv f(p)$ appartiene di certo all'algebra di Royden. Definiamo la funzione $\tau:R\to R^*$ come
\begin{gather*}
 \tau(p)=x_p \ \Longleftrightarrow \  \tau(p)(f)=f(p)
\end{gather*}
\begin{lemma}\label{lemma_R1}
 La funzione $\tau$ è un'omeomorfismo sulla sua immagine.
\begin{proof}
Infatti sia $p\neq q$, necessariamente $\tau(p)\neq \tau(q)$. Infatti poiché $R$ è una varietà riemanniana, per ogni coppia di punti disgiunti esiste una funzione $f$ liscia a supporto compatto (quindi $f\in\Ro$) tale che $f(p)\neq f(q)$ \footnote{questa funzione si può costruire ad esempio sfruttando le partizioni dell'unità}, quindi $\tau(p)(f)\neq \tau(q)(f)$.\\
Inoltre $\tau$ è continua. Consideriamo un qualunque aperto di $R^*$. Per definizione di topologia debole-*, questo aperto contiene un'intersezione finita di insiemi della forma
\begin{gather*}
 V(q,f,\epsilon)=\{p\in R^* \ t.c. \ \abs{p(f)-q(f)}<\epsilon\}
\end{gather*}
La controimmagine $\tau^{-1}(V)$ è aperta poiché $f$ è continua.\\
Inoltre la mappa $\tau$ è anche una mappa aperta sulla sua immagine. Consideriamo a questo scopo un punto $p\in R$ e un suo intorno aperto $V$ qualsiasi, e dimostriamo che $\tau(p)$ ha un intorno aperto contenuto in $\tau(V)\cap \tau(R)$. Sia a questo proposito $\psi$ una funzione liscia tale che $\psi(p)=0$ e $1-\psi$ abbia supporto compatto contenuto in $V$ (questo garantisce che $\psi \in \Ro$). L'insieme
\begin{gather*}
 A=\{q\in \tau(R) \ t.c. \ \abs{\psi(q)-\psi(p)}<1/2\}
\end{gather*}
è un aperto su $\tau(R)$ ed è ovviamente contenuto in $\tau(V)$.
\end{proof}
\end{lemma}
Per quanto riguarda l'estendibilità delle funzioni $f\in \Ro$, osserviamo che:
\begin{lemma}\label{lemma_R2}
Tutte le funzioni $f\in \Ro$ possono essere estese per continuità a $R^*$.
\begin{proof}
 Definiamo $\bar f$ l'estensione di $f$ a $R^*$, e definiamo in maniera naturale
\begin{gather*}
 \bar f (x) \equiv x(f) \ \ \forall x \in R^*
\end{gather*}
Se $p\in R$, allora $\bar f (p) = \tau(p)(f) = f(p)$, il che dimostra che $\bar f $ è un'estensione di $f$. La continuità di $\bar f$ segue dalla definizione di topologia debole-*, infatti:
\begin{gather*}
 \{y\in R^* \ t.c. \ \abs{\bar f (x)-\bar f(y)}<\epsilon\}= \{y\in R^* \ t.c. \ \abs{x(f)-y(f)}<\epsilon\}
\end{gather*}
dove l'ultimo insieme è aperto in $R^*$ per definizione.
\end{proof}
\end{lemma}
Per dimostrare la densità di $R$ in $R^*$, utilizzeremo quest'altro lemma:
\begin{lemma}\label{lemma_R3}
 Lo spazio $\overline{\Ro}=\{\bar f \ t.c. \ f\in \Ro\}$ è denso rispetto alla norma del sup nell'algebra delle funzioni continue da $R^*$ a $\R$, $C(R^*)$.
\begin{proof}
 Questa dimostrazione è una facile conseguenza del teorema di Stone-Weierstrass (vedi ad esempio teorema 7.31 pag 162 di \cite{15}). L'algebra $\overline \Ro$ infatti separa i punti su $R^*$ e contiene la fuzione costante uguale a 1 che non si annulla mai (quindi l'algebra non si annulla in nessun punto di $R^*$).
\end{proof}
\end{lemma}
\begin{lemma}\label{lemma_R4}
 L'insieme $\tau(R)$ è denso in $R^*$.
\begin{proof}
 Dimostriamo questa affermazione per assurdo: sia $\tilde x\in R^*\setminus \overline{\tau(R)}$. Allora per il lemma di Uryson (vedi teorema 2.12 pag 39 di \cite{12}) esiste una funzione continua $h:R^*\to \R$ tale che $h(\tilde x)=0$ e $h(\tau(R))=1$. Grazie al lemma precedente, esiste una funzione $\bar f\in \overline{\Ro}$ tale che $\bar f (\tilde x)=0$ e $\bar f(\tau(R))>1/3$ \footnote{infatti grazie al lemma precedente esiste una funzione $\bar h \in \overline\Ro$ tale che $\bar h(\tilde x)<1/3$ e $\bar h(\tau/R))>2/3$. La funzione $\bar f =\bar h -\bar h(\tilde x)$ ha le caratteristiche richieste}.Allora la sua restrizione $f:R\to \R$ ha estremo inferiore positivo, quindi è invertibile in $\Ro$. Questo significa che esiste $g$ tale che $(fg)(p)=1 \ \forall p\in R$. Ma per definizione di carattere $\overline{(fg)}(x)=x(fg)=1 \ \forall x\in R^*$, quindi anche $\overline{(fg)}(\tilde x )=\bar f(\tilde x) \bar g (\tilde x)=1$ \footnote{dato che $x(fg)=x(f)x(g)$, $\overline{(fg)}=\bar f \bar g$}, impossibile se $\bar f(\tilde x)=0$
\end{proof}
\end{lemma}

In seguito a questi lemmi siamo pronti a dimostrare che
\begin{prop}
 $R^*$ l'insieme dei caratteri su $\Ro$ con la topologia debole ereditata da questo spazio è una compattificazione di Royden per $R$.
\begin{proof}
 Il fatto che $R^*$ sia compatto e di Hausdorff è stato dimostrato a pagina \pageref{teo_Reu}. Che $R$ sia denso in $R^*$ è il contenuto dei lemmi \ref{lemma_R1} e \ref{lemma_R4}. Inoltre $R$ è aperto come sottoinsieme di $R^*$ grazie alla sua locale compattezza. Infatti sia $p\in R$ qualsiasi, e $U$ un suo intorno aperto relativamente compatto nella topologia di $R$, allora esiste un aperto $A\in R^*$ tale che $\tau(U)=A\cap R$. Ma allora $\tau(U)=A$, infatti se per assurdo non fosse così, $A\setminus \tau(\overline{U})$ sarebbe ancora un insieme aperto per compattezza di $\tau(\overline{U})$, e quindi se fosse non vuoto avrebbe intersezione non vuota con $R$ per densità, ma questo è impossibile, quindi $A\subset \tau(\overline{U})$, e cioè $A=\tau(U)$.\\
 Ogni funzione in $\Ro$ può essere estesa per continuità a una funzione definita su $R^*$ grazie al lemma \ref{lemma_R2}, e infine l'insieme $\overline \Ro$ delle funzioni che sono estensione delle funzioni in $\Ro$ separa i punti di $R^*$ per definizione di carattere, oppure come conseguenza del lemma \ref{lemma_R3}.
\end{proof}
\end{prop}
Resta da verificare l'unicità di $R^*$ a meno di omeomorfismi che tengano fissi gli elementi di $R$.
\begin{prop}\label{prop_R_sumup}
 Sia $X$ uno spazio con le proprietà (1),(2),(3),(4) di \ref{teo_Reu}. Allora la mappa $\sigma:X\to R^*$ definita da
\begin{gather*}
 \sigma(p) (f)=f(p)
\end{gather*}
dove $f\in \overline \Ro$, è un omeomorfismo che tiene fissi gli elementi di $R$.
\begin{proof}
 \`E ovvio che $\sigma$ tiene fissi gli elementi di $R$, e grazie a un ragionamento molto simile a quello del lemma \ref{lemma_R1}, $\sigma$ è un omeomorfismo sulla sua immagine. Essendo però $X$ compatto, $R$ denso in $R^*$, e anche $R^*$ compatto, allora necessariamente $\sigma(X)=R^*$, quindi $\sigma$ è anche suriettivo.
\end{proof}
\end{prop}
D'ora in avanti per comodità di notazione confonderemo la scrittura $f$ e $\bar f$, cioè indicheremo una funzione $f:R\to \R$ e la sua estensione a tutta la compattificazione di Royden nello stesso modo.
\subsection{Esempi non banali di caratteri}\label{subsec_charroy}
Come preannunciato nel capitolo \ref{chap_filtri}, utilizziamo gli ultrafiltri e in particolare i risultati ottenuti nella sezione \ref{sec_charsucc} per descrivere alcuni caratteri non banali sull'algebra di Royden (quindi punti di $R^*$).\\
Dalla definizione di carattere, ci si aspetta che sull'algebra di Royden esistano dei caratteri che in qualche senso rappresentino il limite della funzione in una certa direzione. Ad esempio se consideriamo una successione $x_n\in R$, l'operazione di limite lungo questa successione è un'operazione lineare e moltiplicativa su $\Ro$ quando definita, nel senso che date due funzioni $f$ e $g$ tali che esista il limite $\lim_n f(x_n)$ e $\lim_n g(x_n)$, allora vale che:
\begin{gather*}
 \lim_n (f\cdot g)(x_n)=\lim_n f(x_n) \cdot \lim_n g(x_n)
\end{gather*}
data una funzione in $\Ro$, non sempre è garantito che su ogni successione che tende a infinito esista il limite di $f(x_n)$. Per aggirare questo problema, utilizziamo i limiti lungo ultrafiltri.
\begin{oss}
 Con la tecnica mostrata nella sezione \ref{sec_charsucc} è possibile definire caratteri anche sulle algebre di Royden. Basta considerare una qualsiasi successione $\{x_n\}$ in $R$ e un qualsiasi ultrafiltro $\M$ su $\N$ e definire
\begin{gather*}
 \phi(f)=\lim_{\M} f(x_n)
\end{gather*}
Se l'ultrafiltro è non costante, allora $\phi$ è un'operazione di limite, quindi:
\begin{gather*}
 \liminf_n f(x_n)\leq \phi(f)\leq \limsup_n f(x_n)
\end{gather*}
Questo garantisce che l'operazione descritta sia un'estensione dell'operazione di limite standard, nel senso che se esiste $\lim_n f(x_n)$, allora $\phi(f)=\lim_n f(x_n)$, e anche che questo carattere non sia un carattere banale, non sia cioè un carattere immagine attraverso $\tau$ di un punto di $R$. Infatti è sempre possibile creare una funzione in $\Ro$ che valga $1$ su un qualsiasi punto fissato $\bar x$ e che valga $O$ definitivamente su una successione $x_n$ che tende a infinito \footnote{quindi $\phi(f)=0$, mentre $\bar x(f)\equiv f(\bar x) =1$}.\\
La compattificazione di Royden quindi contiene tutti i caratteri che rendono conto del comportamento di una funzione al ``limite'' lungo una successione.
\end{oss}
In realtà non è necessario passare attraverso le successioni, con un ragionamento molto simile a quello sviluppato nella sezione precendente, si osserva che
\begin{prop}
 Ogni ultrafiltro $\M$ su $R$ definisce un carattere sull'algebra di Royden, quindi un elemento di $R^*$
\begin{proof}
 La dimostrazione è la generalizzazione della dimostrazione di \ref{prop_car_N}.
Consideriamo una funzione continua limitata $f:R\to \R$. La collezione $f(\M)$ è un ultrafiltro in $\R$, anzi un ultrafiltro in $[\inf_{R}(f),\sup_R(f)]$, insieme compatto. Quindi per ogni funzione possiamo definire
\begin{gather*}
 \phi_{\M}(f)=\lim f(\M)
\end{gather*}
e con argomenti del tutto analoghi alla dimostrazione \ref{prop_car_N} ottenere che $\phi$ è un carattere su $M(R)$.
\end{proof}
\end{prop}
Quello che manca in questo caso rispetto a sopra è la proprietà che vale per gli ultrafiltri non costanti
\begin{gather*}
 \liminf_n x(n) \leq \phi(x) \leq \limsup_n x(n)
\end{gather*}
quindi una specie di controllo del carattere con il comportamento di $f$ all'infinito. Ovviamente non ha senso la definizione standard di $\liminf$ e $\limsup$ per una funzione che come dominio ha un'insieme con la potenza del continuo, ma comunque anche in questo caso esiste una proprietà simile, che verrà esplorata nella sezione successiva.
\subsection{Caratterizzazione del bordo}
In questa sezione ci occupiamo di due caratterizzazioni del bordo di $R^*$, una di natura funzionale e una di natura topologica.
\begin{deph}\label{deph_Gamma}
 Indichiamo con $\Gamma=R^*\setminus R$ il bordo di $R^*$
\end{deph}
Ovviamente $\Gamma$ è un'insieme compatto in $R^*$. Per prima cosa dimostreremo che per un carattere $p\in \Gamma$ $p(f)=p((1-\lambda) f)$ per ogni funzione $\lambda$ a supporto compatto in $R$, cioè il valore di un carattere in $\Gamma$ dipende solo dal comportamento ``all'infinito'' della funzione a cui è applicato.
\begin{prop}\label{prop_charG}
$p\in \Gamma$ se e solo se $p(f)=0$ $\forall f \in \Roo$, o equivalentemente se per ogni funzione $f\in \Ro$ e $\forall \lambda\in C^\infty(R,\R)$ a supporto compatto
\begin{gather*}
 p(\lambda f)=0 \ \ \Longleftrightarrow  \ \ p(f)=p((1-\lambda) f)
\end{gather*}
\begin{proof}
Questa dimostrazione è un'adattamento dell'esempio 11.13 (a) pag 283 di \cite{4}. Per prima cosa notiamo che vale una ovvia dicotomia per gli elementi di $R^*$. Consideriamo un'esaustione $K_n$ di $R$, e consideriamo una successione di funzioni di cut-off $\lambda_n$ tali che $\lambda_n(K_n)=1$ e $supp(\lambda_n)\subset K_{n+1}$. Allora $\forall p \in R^*$ vale che
\begin{gather*}
 \forall n \ p(\lambda_n)=0 \ \vee \ \exists n \ t.c. \ p(\lambda_n)\neq 0
\end{gather*}
Nel primo caso, consideriamo $\lambda$ una funzione a supporto compatto. Allora esiste $n$ tale che $supp(\lambda)\subset K_n$, quindi $\lambda \lambda_n f=\lambda f$ da cui:
\begin{gather*}
 p(\lambda f)=p(\lambda_n\lambda f)=p(\lambda_n)p(\lambda f)=0
\end{gather*}
Dimostriamo che nel secondo caso necessariamente $p=\tau(x)$ per qualche $x\in K_{n+1}$.\\
Per prima cosa osserviamo che nel secondo caso $\forall f\in \Ro$
\begin{gather}\label{eq_p1}
 p(f)=p(f\lambda_n)/p(\lambda_n)
\end{gather}
Supponiamo per assurdo che questo non sia vero, cioè che per ogni $x$ in $K_{n+1}$, $p\neq \tau(x)$. Questo significa che esiste $f\in\Ro$ tale che
\begin{gather*}
 x(f)=f(x)\neq p(f)\ \Rightarrow \ p(f-p(f)1)=0\neq x(f-p(f)1)=f(x)-p(f)
\end{gather*}
cioè per ogni punto $x$ in $K_{n+1}$, esiste $f_x$ tale che $p(f_x)=0$ e $f_x(x)\neq 0$. Visto che tutte le funzioni $f_x$ sono continue, per ogni $x$ esiste un intorno $U_x$ in cui $f_x\neq 0$. Quindi per compattezza di $K_{n+1}$ esiste un numero finito di punti $\{x_k\}\subset K_{n+1}$ tali che $\cup_k U_{k}\supset K_{n+1}$. Consideriamo la funzione
\begin{gather*}
 F=\sum_k f_k^2
\end{gather*}
questa funzione è strettamente positiva su tutto l'insieme $K_{n+1}$, quindi anche su un suo intorno compatto $A$. Sia $\psi$ una funzione di cut-off con $supp(\psi)\Subset A$ e $\psi(K_{n+1})=1$. La funzione
\begin{gather*}
 \tilde F = F\cdot \psi + (1-\psi)
\end{gather*}
è per costruzione una funzione strettamente positiva su tutta la varietà $R$ e $$\tilde F(x)=1 \ \forall x \not\in A$$ Quindi $\inf_R(\tilde F)>0$ e grazie alla proposizione \ref{prop_inv} è un'elemento invertibile di $\Ro$. Questo implica che necessariamente $p(\tilde F)\neq 0$, ma ciò nonostante:
\begin{gather*}
 p(\tilde F)=p(\lambda_n \tilde F)/p(\lambda_n)=p(\lambda_n F)/p(\lambda_n)=p(F)=\sum_{k} [p(f_k)]^2=0
\end{gather*}
dove abbiamo sfruttato la relazione \ref{eq_p1} e la definizione delle varie funzioni in gioco.
\end{proof}
\end{prop}
Passiamo ora alla caratterizzazione topologica del bordo $\Gamma$.
\begin{prop}
 $p\in \Gamma$ se e solo se $p$ non è un insieme $G_{\delta}$, cioè se e solo se $p$ non è intersezione numerabile di aperti.
\begin{proof}
 Supponiamo per assurdo che $p$ sia un $G_{\delta}$, quindi siano $V_n$ tali che $\overline{V_{n+1}}\subset V_n$ e sia $\{K_n\}$ un'esaustione di $R$. Dato che $p\in \Gamma$, $W_n\equiv V_n\cap K_n^C$ sono ancora intorni di $p$, e vale ancora che $\overline{W_{n+1}}\subset W_n$. Per ogni aperto $W_n\setminus \overline{W_{n+1}}$, scegliamo una funzione $w_n$ con supporto in questo aperto, con massimo $1$ e di integrale di Dirichlet $D_R(w_n)\leq 2^{-n}$ \footnote{possibile grazie alla proposizione \ref{prop_cap1}}. Allora la funzione
\begin{gather*}
 f=\sum_{n=1}^{\infty} w_n
\end{gather*}
converge localmente uniformemente a una funzione continua su $R$ con integrale di Dirichlet finito, quindi è estendibile con contintuità a $R^*$. Ma visto che per ogni intorno $W_n$ di $p$ la funzione oscilla tra $0$ e $1$, non è estendibile con continuità in $p$, contraddizione.
\end{proof}
\end{prop}
Questo prova che la topologia di $R^*$ non è I numerabile, quindi neanche metrizzabile.
\subsection{Bordo armonico e decomposizione}
In questa sezione definiamo il \textbf{bordo armonico} di $R^*$, concetto che sarà utile per dimostrare una versione del principio del massimo.
\begin{deph}\label{deph_ab}
 Definiamo $\Delta$ in bordo armonico di $R^*$ come:
\begin{gather*}
 \Delta= \{p\in R^* \ t.c. \ \forall f\in \Rod \ f(p)=0\} =\bigcap_{f\in \Rod} f^{-1}(0)
\end{gather*}
\end{deph}
Dalla definizione risulta evidente che $\Delta$ è un insieme chiuso, e che $\Delta \cap R=\emptyset$, cioè $\Delta \subset \Gamma$. Infatti per ogni punto di $R$ è sempre possibile trovare una funzione liscia a supporto compatto che vale 1 sul punto.\\
Il seguente lemma dimostra l'esistena di particolari funzioni in $\Rod$.
\begin{lemma}\label{lemma_frod}
 Sia $K$ un insieme compatto in $R^*$ tale che $K\cap \Delta=\emptyset$. Allora esiste una funzione $f\in \Rod$ identicamente uguale a $1$ su $K$.
\begin{proof}
 Dato che $K\cap \Delta=\emptyset$, per ogni punto $p\in K$ esiste una funzione $f_p\in \Rod$ tale che $f_p(p)>1$. Per compattezza di $K$, esiste un numero finito di punti $p_1,\cdots,p_n$ tali che
\begin{gather*}
 \bigcup_{i=1}^n f_{p_n}^{-1}(1,\infty) \supset K
\end{gather*}
Quindi la funzione
\begin{gather*}
 f_1(x)\equiv \sum_{i=1}^n f^2_{p_n}(x)
\end{gather*}
è una funzione appartenente a $\Rod$ (che ricordiamo essere un'ideale), e maggiore di $1$ sull'insieme $K$. La funzione $g(x)=\max\{1/2,f_1(x)\}$ è una funzione $g\in \Ro$ invertibile grazie al lemma \ref{prop_inv} e identicamente uguale a $f_1$ su $K$, quindi grazie al fatto che $\Rod$ è un'ideale di $\Ro$, si ha che la funzione
\begin{gather*}
f = f_1 \cdot g^{-1} 
\end{gather*}
soddisfa la tesi del lemma.
\end{proof}
\end{lemma}
Come corollario, osserviamo che:
\begin{oss}\label{oss_rodro}
 Se e solo se $\Delta=\emptyset$, allora $\Rod=\Ro$.
\begin{proof}
L'implicazione da destra a sinistra è evidente dalla definizione dell'insime $\Delta$. Per quanto riguarda l'altra implicazione, se $\Delta=\emptyset$, grazie alla proposizione precedente la funzione costante uguale a $1$ su tutto $R^*$ appartiene a $\Rod$, da cui segue la tesi.
\end{proof}

\end{oss}

Il seguente teorema permette di scrivere ogni funzione $f\in \Ro$ come la somma di due funzioni, una armonica e l'altra in $\Rod$. Questa decomposizione sarà utile sia in questa sezione per dimostrare una forma del principio del massimo che nel seguito.
\begin{teo}\label{teo_1}
 Sia $f\in \Ro$. Allora esistono una funzione $u\in HDB(R)$ \footnote{$HDB(R)$ è l'insime delle funzioni armoniche su $R$, limitate e con integrale di Dirichlet finito} che indicheremo per comodità $u=\pi(f)$ e una funzione $h\in \Rod$ tali che
\begin{enumerate}
 \item $f=\pi(f)+h$
 \item se $\Delta\neq \emptyset$, la decomposizione è unica, altrimenti unica a meno di costanti
 \item se $f\geq 0$, allora anche $\pi(f)\geq 0$
 \item $\norm{f}_{\infty}\geq \norm{\pi(f)}_{\infty}$
 \item per ogni funzione $\phi\in \Rod$, $D_R(\phi,\pi(f))=0$, quindi
\begin{gather*}
 D_R(f)=D_R(\pi(f))+D_R(h)
\end{gather*}
\end{enumerate}
\begin{proof}
 La dimostrazione di questo teorema è costruttiva. Sia $K_n$ un'esaustione di $R$ con domini regolari. Definiamo $u_n\in \Ro$ come $u_n=f$ sull'insieme $R\setminus K_n$, e $u_n\in H(K_n)$. Grazie al principio del massimo, notiamo subito che $\norm{u_n}_{\infty}\leq \norm{f}_{\infty}$ e che se $f\geq0$, anche $u_n\geq0$. Visto che la successione $u_n$ è uniformemente limitata, grazie al principio di Harnack \ref{prop_harnackpri} esiste una sua sottosuccessione (che continueremo a indicare con $u_n$) che converge localmente uniformemente a una funzione armonica $u$ definita su $R$. Inoltre, grazie all'identità di Green \ref{prop_green3}, se $k>n$:
\begin{gather*}
 D_{K_k} (u_k-u_n,u_k)=D_R(u_n-u_k,u_k)=0  
\end{gather*}
e quindi:
\begin{gather*}
 D_R(u_n)=D_R(u_k)+D_R(u_n-u_k)+2D_R(u_k,u_n-u_k)=D_R(u_k)+D_R(u_n-u_k)
\end{gather*}
in particolare, se $k>n$, $D_R(u_k)\leq D_R(u_n)$, quindi la successione $D_R(u_n)$ è decrescente e converge a un valore $\geq 0$.\\
Grazie all'ultima uguaglianza, questo significa anche che $D_R(u_n-u_k)$ è piccolo a piacere se $n$ e $k$ sono sufficientemente grandi, cioè la successione $\{u_n\}$ è $D-$Cauchy. Vista la completezza $BD$ dello spazio $\Ro$, la funzione $u\in \Ro$, quindi $$u\in HDB(R)$$ Consideriamo ora la successione $h_n=f-u_n$. Per le proprietà di $u_n$, $h_n\in \Roo$, e $BD-\lim_n h_n=f-u=h$, quindi chiaramente $h\in \Rod$.\\
Rimanne da dimostrare l'unicità della decomposizione. Supponiamo che esistano due decomposizioni distinte per $f$:
\begin{gather*}
 f=u+h=\bar u +\bar h
\end{gather*}
Definiamo $v=u-\bar u =\bar h -h$. Necessariamente $v\in HDB(R)\cap \Rod$. Quindi esiste una successione $v_n\in \Roo$ tale che $v=BD-\lim_n v_n$, quindi:
\begin{gather*}
 D_R(v)=\lim_n D_R(v_n,v)=\lim_n D_{K_k}(v_n,v)=\lim_n \int_{\partial K_k} v_n\ast dv =0
\end{gather*}
Dove $K_k$ è scelto in modo che $supp(v_n)\subset K_k$.\\
Questo dimostra che $v$ ha derivata nulla, quindi è costante su $R$, e se $\Delta\neq \emptyset$, allora $v=0$ dato che $v(\Delta)=0$.\\
Il punto (5) è una facile applicazione della formula di Green \ref{prop_green2}. Infatti se $\phi \in \Rod$, esiste una successione $\phi_n\in \Roo$ tale che $\phi=BD-\lim_n \phi_n$, quindi:
\begin{gather*}
 D_R(\phi,u)=\lim_n D_R(\phi_n,u)
\end{gather*}
Per ogni $n$, sia $K$ un compatto dal bordo liscio in $R$ tale che $supp(\phi_n)\subset K$, allora per \ref{prop_green2}.
\begin{gather*}
 D_R(\phi_n,u)=D_K(\phi_n,u)=0
\end{gather*}
poiché $\phi_n=0$ su $\partial K$.
\end{proof}
\end{teo}
Dalla dimostrazione è facile ricavare questo corollario:
\begin{prop}\label{prop_sub1}
 Data $f\in \Ro$, se esiste una funzione subarmonica $v$ definita su $R$ tale che $v\leq f$, allora
\begin{gather*}
 v\leq \pi(f)
\end{gather*}
allo stesso modo, se esiste $v$ superarmonica tale che $v\geq f$ allora
\begin{gather*}
 v\geq \pi(f)
\end{gather*}
\end{prop}
Questa proposizione può essere migliorata, infatti nella composizione si può chiedere che $\pi(f)=f$ su un predeterminato insieme compatto $K$ con bordo regolare a tratti.
\begin{teo}\label{teo_dec}
 Sia $f\in \Ro$ e sia $K$ un compatto non vuoto in $R$ con bordo regolare, allora esistono uniche una funzione $u\in HBD(R\setminus K)\cap \Ro$ e $g\in \Rod$, $g=0$ su $K$, tali che:
\begin{enumerate}
 \item $f=u+g$
 \item $D_R(u,\phi)=0$ per ogni $\phi\in \Rod$ e $\phi=0$ su $K$
 \item $D_R(f)=D_R(u)+D_R(g)$
 \item se $v$ è superarmonica su $R\setminus K$ e su questo insieme $v\geq f$, allora $v\geq u$ su $R\setminus K$
\item $\norm u_{\infty, R\setminus K} \leq \norm f _{\infty, R\setminus K}$
\end{enumerate}
\begin{proof}
 La dimostrazione è del tutto analoga alla dimostrazione precedente, quindi lasciamo i dettagli al lettore.
\end{proof}
\end{teo}
\begin{oss}\label{oss_decomp}
 Nelle ipotesi del teorema precedente, denotiamo la funzione $u$ con $\pi_K (f)$. Osserviamo che se $K\subset C$:
\begin{gather*}
 \pi_K(\pi_C(f))=\pi_K(f)
\end{gather*}
\begin{proof}
 La dimostrazione è una semplice applicazione del teorema precedente. Infatti sappiamo che:
\begin{gather*}
 f=\pi_K(f)+g_K\\
f=\pi_C(f)+g_C\\
\pi_C(f)=\pi_K(f)+g_K-g_C\equiv \pi_K(\pi_C(f))+g
\end{gather*}
dato che $g|_K\equiv (g_K-g_C)|_K=0$ e che $g\in \Rod$, per unicità della decomposizione abbiamo la tesi.
\end{proof}

\end{oss}

\subsection{Principio del massimo}\label{sec_max2}
Ora siamo pronti a presentare questa versione del principio del massimo.
\begin{teo}[Principio del massimo]\label{teo_max}
 Data una funzione $u\in HBD(R)$, si ha che:
\begin{gather*}
 \min_{p\in \Delta} u(p)=\inf_{p\in R^*} u(p) \ \ \ \max_{p\in \Delta} u(p)=\sup_{p\in R^*} u(p)
\end{gather*}
In particolare se $\Delta=\emptyset$, allora $u\in HBD(R)\Rightarrow u=cost.$
\begin{proof}
 Se $\Delta=\emptyset$ allora $\Ro=\Rod$ (vedi osservazione \ref{oss_rodro}), e si può dimostrare che $u=cost$ con un ragionamento del tutto simile a quello utilizzato nel teorema \ref{teo_1} per dimostrare che $v=cost$.\\
Supponiamo quindi che $\Delta\neq \emptyset$. Per dimostrare la tesi dimostriamo che $u(\Delta)\leq 0$ implica $u(R)\leq 0$. Grazie alla densità di $R$ in $R^*$, l'ultima affermazione implica che $u(R^*)\leq 0$, e se questo è vero, la tesi si ottiene considerando la funzione $u-~\max_{p\in \Delta} u(p)$, oppure la funzione
$-u+\min_{p\in \Delta} u(p)$.\\
Supponiamo quindi che $u(\Delta)\leq 0$. Fissato $\epsilon>0$, definiamo l'insieme
\begin{gather*}
 A=\{p\in R^* \ t.c. \ u(p)\geq \epsilon\}
\end{gather*}
Dato che ovviamente $A\cap \Delta=\emptyset$, per ogni $p$ \footnote{e per definizione di $\Delta$} esiste una funzione $f_p\in \Rod$ tale che $f_p(p)\geq 2$. Visto che l'algebra di Royden è chiusa rispetto all'operazione di massimo, possiamo fare in modo che $f_p\geq0$ su $R^*$ \footnote{basta considerare $\max\{f_p,0\}$}. Allora il ricoprimento
\begin{gather*}
 U_p=\{q\in R^* \ t.c. \ f_p(q)>1\}
\end{gather*}
è un ricoprimento aperto di $A$, e quindi per compattezza esiste un numero finito di punti $p_1,\cdots,p_n$ tali che $U_{p_i}$ ricopre $A$. Questo significa che la funzione
\begin{gather*}
 f(q)\equiv \sum_{i=1}^n f_{p_n}(q)
\end{gather*}
è una funzione $f\in \Rod$ tale che $f|_A \geq1$. Visto che per ipotesi la funzione $u$ è limitata da sopra, esiste un numero $M$ tale che $u-Mf \leq 0$ sull'insieme $A$, e quindi per definizione di $A$, $u-Mf-\epsilon \leq 0$ su tutta $R^*$.\\
Grazie al teorema di decomposizione appena dimostrato possiamo concludere che esistono (uniche visto che $\Delta\neq \emptyset$) due funzioni $v\in HBD(R)$ e $h\in \Rod$ tali che
\begin{gather*}
u-Mf-\epsilon = v+h 
\end{gather*}
e per la proposizione \ref{prop_sub1} $v\leq 0$ su $R^*$. Per unicità, $v=u-\epsilon$, e per arbitrarietà di $\epsilon$, $u\leq 0$ su $R^*$.
\end{proof}
\end{teo}
La dimostrazione di questo teorema può essere facilmente adattata per ottenere:
\begin{teo}\label{teo_max2}
  Data una funzione $u\in HBD(R\setminus K)$ con $K$ compatto in $R$, si ha che:
\begin{gather*}
 \min_{p\in (\Delta\cup \partial K)} u(p)=\inf_{p\in R^*\setminus K} u(p) \ \ \ \max_{p\in (\Delta\cup \partial K)} u(p)=\sup_{p\in R^*\setminus K} u(p)
\end{gather*}
\end{teo}

Come corollario a questo teorema possiamo facilmente dimostrare che
\begin{prop}
 Data una funzione subarmonica $u\in \Ro$ limitata, allora:
\begin{gather*}
 \max_{p\in \Delta} u(p)=\sup_{p\in R^*} u(p)
\end{gather*}
Data una funzione superarmonica $v\in \Ro$:
\begin{gather*}
\min_{p\in \Delta} v(p)=\inf_{p\in R^*} v(p)
\end{gather*}
\begin{proof}
 Dimostriamo solo il primo caso. La dimostrazione è identica al caso armonico, ma una volta trovata la decomposizione di $u-Mf-\epsilon=v+h$, non possiamo concludere che $v=u-\epsilon$, ma grazie alla proposizione \ref{prop_sub1}, possiamo concludere che $u\leq v$.
\end{proof}

\end{prop}
Un altro corollario di questo principio è il seguente:
\begin{prop}\label{prop_delta}
 Se $u\in H(R\setminus K)$ con $K$ insieme compatto \footnote{possibilmente vuoto} è una funzione limitata e $u\in \Rod$, allora $u=BD-\lim_n u_n$, dove $u_n$ sono le funzioni definite da:
\begin{gather*}
 u_n(x)=\begin{cases}
         u(x) & se \ x\in K\\
0 & se \ x\in K_n^C\\
u_n \in H(K_n\setminus K) 
        \end{cases}
\end{gather*}
dove $K_n$ è un'esaustione regolare di $R$ con $K\subset K_1$.\\
Questo implica che $\forall x \in R\setminus K$:
\begin{gather*}
 \min\{0,\min u|_{\partial K}\} \leq u(x) \leq \max\{0,\max u|_{\partial K}\}
\end{gather*}

\begin{proof}
 Grazie alle considerazioni fatte in precedenza, sappiamo che esiste il limite
\begin{gather*}
 v=BD-\lim_n u_n
\end{gather*}
la funzione $u-v$ è una funzione armonica su $R\setminus K$ con $(u-v)|_{K}=0$ e $(u-v)|_{\Delta}=0$, quindi grazie al principio del massimo appena dimostrato, $u=v$.\\
L'ultima considerazione segue dal principio del massimo applicato alle funzioni $u_n$ sull'insieme $K_n\setminus K$. Per le funzioni $u_n$ infatti massimo e minimo sono assunti su $\partial (K_n\setminus K)=\partial K \cup \partial K_n$. Passando al limite si ottiene la tesi.
\end{proof}

\end{prop}

La costruzione della compattificazione di Royden trova una buona giustificazione nel principio appena dimostrato. Su insiemi compatti con bordo, uno strumento fondamentale per lo studio delle funzioni armoniche è il principio del massimo che garantisce che una funzione armonica assume il suo massimo in un punto del bordo. Questo principio non è (ovviamente) applicabile se studiamo funzioni armoniche su varietà non compatte senza bordo. La compattificazione di Royden è il tentativo di ``rendere compatta'' una varietà in modo da poter applicare ancora i principi del massimo. Come ci si può aspettare, un insieme compatto con bordo e la compattificazione di Royden non si comportano esattamente nello stesso modo. Nelle ipotesi di quest'ultimo principio infatti compare l'ipotesi che $u\in HBD(R)$, cioè che l'integrale di Dirichlet della funzione sia finito, però questo principio dice anche dove esattamente cercare il massimo delle funzioni armoniche. Non in un punto qualsiasi del bordo come nel caso di insiemi compatti, ma è sufficiente controllare il comportamento di $f$ sul \textbf{bordo armonico} della varietà.\\
Concludiamo questo capitolo con una caratterizzazione di $\Rod$. Dalla definizione \ref{deph_ab}, il bordo armonico è l'insieme dove tutte le funzioni in $\Rod$ sono nulle. Vale anche una sorta di viceversa, nel senso che:
\begin{prop}
\begin{gather*}
 \Rod = \{f\in \Ro \ \ t.c. \ f(\Delta)=0\}
\end{gather*}
\begin{proof}
La dimostrazione segue facilmente dalla decomposizione descritta nella proposizione \ref{teo_1} e dal principio del massimo \ref{teo_max}.\\
Se $\Delta=\emptyset$, $\Rod=\Ro$ come dall'osservazione \ref{oss_rodro}. Supponiamo quindi che $\Delta\neq\ \emptyset$. Se $f\in \Rod$, è ovvio dalla definizione che $f(\Delta)=0$. Supponiamo quindi che $f(\Delta)=0$. Allora, grazie al teorema \ref{teo_1}, esistono due funzioni $u\in HBD(R)$ e $h\in \Rod$ tali che
\begin{gather*}
 f=u+h
\end{gather*}
Essendo $h(\Delta)=f(\Delta)=0$, si ha che $u(\Delta)=0$, quindi per il principio del massimo \ref{teo_max}, $u$ è identicamente nulla, da cui la tesi.
\end{proof}
\end{prop}

\chapter{Varietà paraboliche e iperboliche}\label{chap_parab}
In questo capitolo tratteremo proprietà e caratterizzazioni di varietà paraboliche e iperboliche. Il risultato principale è l'esistenza dei \textit{potenziali di Evans} sulle varietà paraboliche. In tutto il capitolo $R$ sarà una varietà riemanniana sensa bordo di dimesione $m$.
\section{Capacità}
In questa sezione introduciamo la capacità di una coppia di insiemi e definiamo di conseguenza le varietà paraboliche. I risultati principali sono tratti da \cite{20} e \cite{22} capitolo 7. 
\begin{deph}\label{deph_cap}
 Dati due insiemi $K$ compatto e $K\subset\Omega$ aperto paracompatto, definiamo la \textbf{capacità} di $K$ rispetto a $\Omega$:
\begin{gather*}
 \text{Cap}(K,\Omega)=\inf{\int_{R} (\nabla f)^2 dV}
\end{gather*}
dove l'inf è preso su tutte le funzioni di Tonelli \footnote{in realtà l'insieme di funzioni considerate può essere più generale senza cambiare il valore della capacità, vedi \cite{20}} tali che
\begin{gather*}
 f(K)=1 \ \ supp(f)\subset \overline{\Omega} \ \ 0\leq f \leq 1
\end{gather*}
\end{deph}
Data questa definizione, osserviamo subito che la capacità di una coppia di insiemi aumenta se rimpiccioliamo $\Omega$ e se ingrandiamo $K$, cioè se $K\subset K'$ e $\Omega\subset \Omega'$:
\begin{gather*}
 \text{Cap}(K,\Omega')\leq \text{Cap}(K,\Omega)\leq \text{Cap}(K',\Omega) \ ;  \ \text{Cap}(K,\Omega')\leq \text{Cap}(K',\Omega')\leq \text{Cap}(K',\Omega')
\end{gather*}
Visti gli scopi della tesi, d'ora in avanti considereremo solo insiemi $K$ e $\Omega$ con bordo regolare.\\
La prima domanda data questa definizione è se l'inf è raggiunto, e in tale caso da quale funzione. La risposta dipende dalla regolarità dei bordi di $K$ e $\Omega$, in particolare se entrambi i bordi sono lisci grazie al principio di Dirichlet possiamo dimostrare che:
\begin{prop}
 Il valore della capacità di una coppia di insiemi $K\Subset \Omega$ con bordi lisci è equivalente all'integrale di Dirichlet della funzione soluzione del problema
\begin{gather}\label{eq_dir_1}
 \ \begin{cases}
  \Delta u =0 \ \ in \ \Omega\setminus K \\
  u|_K=1 \\
  u|_{ \Omega^C}=0
 \end{cases}
\end{gather}
Chiamiamo la funzione $u$ \textbf{potenziale di capacità} della coppia $K$, $\Omega$.
\begin{proof}
Osserviamo che se i bordi degli insiemi sono lisci, è possibile risolvere il problema di Dirichlet che definisce $u$.\\
La dimostrazione è una diretta conseguenza del principio di Dirichlet riportato nella proposizione \ref{prop_D2}.
\end{proof}
\end{prop}
Nella seguente proposizione osserviamo che se è possibile risolvere il problema di Dirichlet \ref{eq_dir_1}, anche se $\Omega$ e $K$ non hanno bordo liscio $\text{Cap}(K,\Omega)=D_\Omega(u)$.
\begin{prop}
 Siano $K$ e $\Omega$ tali che $\Omega\setminus K$ sia regolare per il problema di Dirichlet \footnote{cioè tale che il problema di Dirichlet \ref{eq_dir_1} abbia soluzione, vedi sezione \ref{sec_dir} per condizioni sulla regolarità degli insiemi}. Se esiste una funzione di Tonelli $f$ tale che $f|_{K}=1$, $f|_{\Omega^C}=0$ e $D_R(f)<\infty$, allora $\text{Cap}(K,\Omega)=D_R(u)<\infty$.
\begin{proof}
 La dimostrazione è una diretta conseguenza del principio di Dirichlet riportato nella proposizione \ref{prop_D3}.
\end{proof}
\end{prop}

\begin{prop}
 Nella definizione di capacità, se $\Omega$ e $K$ hanno bordi regolari, allora possiamo sostituire l'insieme delle funzioni di Tonelli con l'insieme delle funzioni $f$ lisce in $\Omega\setminus K$ uguali a $1$ su $K$ e nulle su $\Omega^C$.
\begin{proof}
 Questo risultato segue dalla densità delle funzioni lisce in $\Omega\setminus K$ e uguali a $1$ su $K$ e nulle in $\Omega^C$ illustrato nella proposizione \ref{prop_dens2}.
\end{proof}
\end{prop}
Possiamo ulteriormente caratterizzare la capacità di una coppia di insiemi grazie alle formule di Green.
\begin{prop}\label{prop_cap2}
 Sia $u$ il potenziale di capacità di $(K,\Omega)$. Allora si ha che:
\begin{gather*}
 \text{Cap}(K,\Omega)=\int_{R} \abs{\nabla u}^2 dV =- \int_{\partial K} \ast du =- \int_{\partial \Omega} \ast du
\end{gather*}
\begin{proof}
 La dimostrazione della prima uguaglianza è una diretta conseguenza della formula di Green \ref{prop_green2}. La seconda segue dalla considerazione che $u$ è armonica in $\Omega\setminus K$, quindi grazie alla proposizione \ref{prop_green1}
\begin{gather*}
 0= \int_{\Omega\setminus K} \Delta u\ dV =\int_{\partial (\Omega\setminus K)} \ast du =\int_{\partial \Omega} \ast du + \int_{\partial K} \ast du
\end{gather*}
\end{proof}
\end{prop}

Oltre alla capacità di una coppia di insiemi, è possibile definire la capacità di un insieme compatto come:
\begin{deph}\label{deph_cap2}
 Dato un compatto $K\Subset R$ regolare per il problema di Dirichlet, definiamo
\begin{gather*}
 \text{Cap}(K)=\lim_n \text{Cap}(K,E_n)
\end{gather*}
dove $E_n$ è una qualsiasi esaustione regolare di $R$ \footnote{per l'esistenza di queste esaustioni, vedi \ref{prop_exreg}}. Questa definizione è ben posta, nel senso che non dipende dalla scelta dell'esaustione.
\begin{proof}
 Date due esaustioni regolari $E_n$ e $C_n$, definiamo
\begin{gather*}
 e_n\equiv \text{Cap}(K,E_n)\to e \ \ c_n\equiv \text{Cap}(K,C_n)\to c
\end{gather*}
le due successioni sono monotone decrescenti grazie alle considerazioni fatte prima. Inoltre per la compattezza di ogni $E_n$, esiste un intero $k$ tale che $E_n\Subset C_k$. Questo implica che $e_n\geq c_k\geq c$. Quindi anche $e\geq c$. Scambiando i ruoli di $E_n$ e $C_n$ si ottiene $e=c$.
\end{proof}
\end{deph}
Possiamo definire un potenziale di capacità anche per un singolo insieme compatto.
\begin{prop}
 Dato un insieme aperto relativamente compatto regolare per il problema di Dirichlet $K$ e una qualsiasi esaustione regolare $K_n$ con $K\subset K_1^{\circ}$ di $R$, sia $u_n$ il potenziale di capacità di $(K,K_n)$, allora la funzione
\begin{gather*}
 u=\lim_n u_n
\end{gather*}
è il potenziale di capacità per $K$, nel senso che $u$ è una funzione armonica su $R\setminus K$, $u=1$ su $K$ tale che
\begin{gather*}
 \text{Cap}(K)=D_R(u)
\end{gather*}
\begin{proof}
 La successione $u_n$ è una successione di funzioni armoniche su $K\setminus K_n$ strettamente positive su $K_n^{\circ}$ grazie alla proposizione \ref{prop_max3}. Dato che tutte queste funzioni sono uniformemente limitate da $1$, grazie alle proposizioni \ref{prop_harnackpri} e \ref{prop_conv_lu} $u_n$ converge localmente uniformemente a una funzione $u$ armonica in $R\setminus K$. Sempre grazie alla proposizione \ref{prop_conv_lu} si ha che
\begin{gather*}
 \lim_n D_R(u_n)=D_R(u)
\end{gather*}
come volevasi dimostrare.
\end{proof}

\end{prop}
Osserviamo che possiamo estendere la nozione di capacità e di potenziale di capacità anche a insiemi più generali di quelli utilizzati fino ad ora. In particolare, dato $C\subset R$ compatto, se è possibile risolvere il problema di Dirichlet
\begin{gather*}
 u|_K=1 \ \ \ u|_{\partial E_n}=0 \ \ \ u\in H(E_n\setminus K)
\end{gather*}
dove $E_n$ è un'esaustione regolare di $R$ con $K\subset E_1$, allora ha senso parlare di capacità di $K$ e di potenziale di capacità di $K$. In particolare:
\begin{oss}\label{oss_altri_dir}
 Se $K$ è una sottovarietà regolare compatta di $R$ di codimensione $1$ possibilmente con bordo, allora ha senso definire la sua capacità e il suo potenziale di capacità.
\begin{proof}
 Questo risultato segue dalle considerazioni della sezione \ref{sec_altri_dir}.
\end{proof}

\end{oss}

Prima di proseguire, riportiamo un'applicazione di quanto appena descritto. In particolare dimostriamo che in ogni insieme aperto, esiste sempre una coppia $(K,\Omega)$ contenuta in questo insieme con capacità grande a piacere e piccola a piacere (quindi una funzione armonica con integrale di Dirichlet grande a piacere e piccolo a piacere).
\begin{prop}
 Dato un qualsiasi aperto $A\subset R$, esistono $K\Subset \Omega \subset A$ tali che $\text{Cap}(K,\Omega)$ è grande a piacere.
\begin{proof}
 La dimostrazione di questa proposizione è costruttiva. Consideriamo un punto $p\in A$, allora esiste un intorno normale $B(p)\subset A$, intorno che possiamo dotare delle coordinate geodetiche $(r,\theta)$ (fuori da $p$) \footnote{vedi sezione \ref{sec_polar}}. Consideriamo
\begin{gather*}
 K=\{p\in B \ t.c. \ r(p)\leq R_1\} \ \ \Omega=\{p\in B \ t.c. \ r(p)< R_2\}
\end{gather*}
dove $R_1<R_2$ e $R_2$ è tale che $\Omega\subset B(p)$. Consideriamo una qualsiasi funzione di Tonelli $0\leq f \leq 1$ tale che $f(K)=1$ e $f=0$ fuori da $\Omega$. Per comodità di notazione, utiliziamo lo stesso simbolo $f$ con la sua rappresentazione in coordinate geodetiche. Abbiamo che
\begin{gather*}
 \int_{R} \abs{\nabla f}^2 \ dV= \int_{\Omega\setminus K} \abs{\nabla f}^2 \ dV=\int_{\Omega\setminus K} g^{ij}\frac{\partial f}{\partial x^i}\frac{\partial f}{\partial x^j} \ \sqrt{\abs g}dr d\theta^1 \cdots d\theta^{m-1}
\end{gather*}
dove $r=x^1$ e $\theta$ rappresenta tutte le altre coordinate. Grazie alla particolare forma che la metrica assume in coordinate polari geodetiche (vedi \ref{eq_gpol}) si ha che:
\begin{gather*}
 \int_{\Omega\setminus K} \abs{\nabla f}^2 \ dV\geq \int_{\Omega\setminus K} \ton{\frac{\partial f}{\partial r}}^2 \ \sqrt{\abs g} dr d\theta^1 \cdots d\theta^{m-1}\geq \\
\geq c \omega_{m-1} \int_{R_1}^{R_2} \ton{\frac{\partial f}{\partial r}}^2 dr
\end{gather*}
dove $c$ è un limite inferiore positivo per $\sqrt{\abs g}$ su $B(p)\setminus K$. Grazie alla disuguaglianza di Schwartz, otteniamo che:
\begin{gather*}
 \frac{1}{R_2-R_1}=\frac{1}{R_2-R_1} \ton{\int_{R_1}^{R_2}\frac{\partial f}{\partial r} dr }^2\leq \int_{R_1}^{R_2} \ton{\frac{\partial f}{\partial r}}^2 dr
\end{gather*}
il che dimostra che a patto di scegliere $R_2-R_1$ sufficientemente piccolo, la capacità dell'anello è abbastanza grande. Osserviamo che in questa dimostrazione è essenziale trovare il limite inferiore per $\sqrt{\abs{g}}$ \footnote{che ad esempio nel caso di $\R^m$ vale $\omega_{m-1}r^{m-1}$, dove $\omega_{m-1}$ è l'area della sfera $m-1$ dimensionale rispetto alla metrica euclidea standard}, quindi è essenziale che $R_1$ non tenda a $0$ quando scegliamo $R_2-R_1$ piccolo. In pratica per ottenere la tesi, fissato $R_1$ e trovata la costante $c$, scegliamo $R_2$ in modo che $R_2-R_1$ sia sufficientemente piccolo.
\end{proof}
\end{prop}
\begin{prop}\label{prop_cap1}
 Dato un qualsiasi aperto $A\subset R$, esistono $K\Subset \Omega \subset A$ tali che $\text{Cap}(K,\Omega)$ è piccola a piacere.
\begin{proof}
 Procediamo in maniera del tutto analoga a sopra, solo che in questo caso al posto di considerare una funzione di Tonelli qualsiasi, ne consideriamo una particolare. Infatti se vogliamo dimostrare che la capacità di un'insieme è piccola, basta trovare \textit{una} funzione che renda piccolo l'integrale di Dirichlet.\\
Sia come sopra $B(p)$ un intorno normale contenuto in $A$, e siano
\begin{gather*}
 K=\{p\in B \ t.c. \ r(p)\leq R_1\} \ \ \Omega=\{p\in B \ t.c. \ r(p)< R_2\}
\end{gather*}
dove $R_1<R_2$ e $R_2$ è tale che $\Omega\subset B(p)$. Consideriamo la funzione
\begin{gather*}
 f(p)=\begin{cases}
        1 & se \ p\in K\\
      0 & se \ p\in \Omega^C\\
      \frac{\log(r(p)/R_2)}{\log(R_1/R_2)} & se \ p\in \Omega\setminus K 
      \end{cases}
\end{gather*}
Dato che $f$ dipende solo da $r(p)$, il suo integrale di Dirichlet vale:
\begin{gather*}
 \int_{R}\abs{\nabla f}^2 dV= \int_{\Omega\setminus K} \ton{\frac{\partial f}{\partial r}}^2 \sqrt{\abs g}dr d\theta^1\cdots d\theta^{m-1}=\\
=\int_{\Omega\setminus K} \frac{1}{\log^2(R_1/R_2)}\frac{1}{r^2} \sqrt{\abs g}dr d\theta^1\cdots d\theta^{m-1}
\end{gather*}
Notiamo che l'integrale $\int_{r=k} \sqrt {\abs g } d\theta^1 \cdots d\theta^{m-1}$ è la superficie della sfera di raggio $k$, quindi per le proprietà della metrica, se $k$ è tale che $B_k\subset B(p)$, allora esiste una costante $C$ tale che:
\begin{gather*}
 \int_{r=k} \sqrt {\abs g } d\theta^1 \cdots d\theta^{m-1}\leq Ck^{m-1}
\end{gather*}
dove $m$ rappresenta la dimensione della varietà. Questa stima permette di concludere che:
\begin{gather*}
 \int_{R}\abs{\nabla f}^2 dV\leq \frac{C}{\log^2(R_1/R_2)}\int_{R_1}^{R_2}r^{m-3} dr
\end{gather*}
Se $m=2$, si ottiene:
\begin{gather*}
\int_{R}\abs{\nabla f}^2 dV\leq \frac{C}{\log(R_2/R_1)} 
\end{gather*}
quindi fissato $R_2$, è possibile scegliere $R_1>0$ in modo che questa quantità sia piccola a piacere.\\
In caso $m\geq 3$ invece si ottiene:
\begin{gather*}
 \int_{R}\abs{\nabla f}^2 dV\leq \frac{C}{\log^2(R_1/R_2)}(R_2^{m-2}-R_1^{m-2})
\end{gather*}
Anche in questo caso, fissato $R_2$ è sempre possibile scegliere $R_1$ in modo che $D_R(f)$ sia piccolo a piacere.
\end{proof}
\end{prop}
La capacità di una coppia di insiemi è legata al comportamento della funzione di Green su quell'insieme, in particolare si ha che \footnote{proposizione 4.1 pag 154 di \cite{20}}:
\begin{prop}\label{prop_cap3}
Siano $K, \ \Omega$ insiemi compatti dal bordo liscio in $R$ con $K\Subset \Omega^{\circ}$, e sia $p\in K$. Allora:
\begin{gather*}
 \min_{x\in \partial K} G_{\Omega}(p,x)\leq (\text{Cap}(K,\Omega))^{-1}\leq \max_{x\in \partial K} G_{\Omega}(p,x)
\end{gather*}
dove $G_{\Omega}$ indica la funzione di Green relativa al dominio $\Omega$ \footnote{vedi proposizione \ref{prop_fg1}}.
\begin{proof}
Siano
\begin{gather*}
 b\equiv \min_{x\in \partial K} G_{\Omega}(p,x) \ \ \ a\equiv\max_{x\in \partial K} G_{\Omega}(p,x)
\end{gather*}
e definiamo per ogni $c>0$ l'insieme compatto
\begin{gather*}
 F_c\equiv \{x\in \Omega \ t.c. \ G_{\Omega}(p,x)\geq c\}
\end{gather*}
Per prima cosa osserviamo che $F_a\subset K \subset F_b$. Questo è conseguenza del principio del massimo, infatti la funzione $G_{\Omega}(\cdot,p)$ è armonica su $\Omega\setminus K$, quindi assume il suo massimo su $\partial K$, cioè $G_{\Omega}(x,p)< a$ se $x\in \Omega \setminus K$. Inoltre dato che $G_{\Omega}$ è superarmonica su $K^{\circ}$, si ha che $G_{\Omega}(x,p)> b$ se $x\in K$.\\
Grazie alla ``monotonia'' della capacità, si ha che:
\begin{gather*}
 \text{Cap}(F_a,\Omega)\leq \text{Cap}(K,\Omega)\leq \text{Cap}(F_b,\Omega)
\end{gather*}
la tesi segue dalla considerazione che $\text{Cap}(F_c,\Omega)=1/c$.\\
Per il teorema si Sard, quasi ogni $c$ è un valore regolare per $G_{\Omega}(\cdot,p)$, quindi $F_c$ ha bordo regolare. Per questi valori di $c$, la funzione $G_{\Omega}/c$ (estesa costante uguale a $1$ dentro $F_c$ e uguale a $0$ fuori da $\Omega$) è il potenziale di capacità per la coppia $(F_c,\Omega)$, quindi grazie alla proposizione \ref{prop_cap2}:
\begin{gather*}
 \text{Cap}(F_c,\Omega)=-\frac{1}{c} \int_{\partial F_c} \ast dG_{\Omega}(\cdot,p) =\frac{1}{c}
\end{gather*}
Se $c$ non è un valore regolare di $G_{\Omega}(\cdot,p)$, allora esiste $c_n\nearrow c$ e $c_n'\searrow c$ successioni di valori regolari di $G_{\Omega}(\cdot,p)$. Dalla monotonia della capacità si ottiene che per ogni $n$:
\begin{gather*}
 \frac{1}{c_n}\leq \text{Cap}(F_c,\Omega)\leq \frac{1}{c_n'}
\end{gather*}
da cui $Cap(F_c,\Omega)=1/c$.
\end{proof}
\end{prop}
Grazie alla definizione di capacità appena data possiamo dividere in 2 categorie le varietà Riemanniane, quelle per cui ogni insieme aperto non vuoto relativamente compatto \footnote{noi considereremo solo compatti con bordo liscio per comodità} ha capacità positiva, e quelle per cui ogni insieme compatto ha capacità nulla. In questa sezione daremo delle caratterizzazioni equivalenti di queste proprietà.
\begin{deph}
 Una varietà Riemanniana $R$ si dice \textbf{parabolica} se ogni insieme aperto non vuoto relativamente compatto con bordo liscio ha capacità nulla, in caso contrario la varietà si definisce \textbf{iperbolica}
\end{deph}
Per prima cosa osserviamo che nella definizione non è necessario chiedere che \textit{ogni} insieme aperto relativamente compatto abbia capacità nulla, basta che un solo aperto relativamente compatto possieda questa proprietà e automaticamente tutti gli aperti non vuoti relativamente compatti hanno capacità nulla. Questo implica anche che se un solo aperto relativamente compatto ha capacità positiva, tuttigli aperti relativamente compatti hanno capacità postiva.
\begin{prop}
 Sia $R$ una varietà Riemanniana. Se un insieme aperto non vuoto relativamente compatto con bordo liscio $K\Subset R$ ha capacità nulla, allora tutti gli aperti non vuoti relativamente compatti in $R$ hanno capacità nulla.
\begin{proof}
Sia $K$ un aperto relativamente compatto con bordo liscio di capacità nulla. Allora tutti gli insiemi compatti contenuti in $K$ hanno capacità nulla grazie alla ``monotonia'' della capacità.\\
Consideriamo un aperto relativamente compatto con bordo liscio $C\supset K$ e chiamiamo $u$ il potenziale di capacità di $K$ e $u'$ il potenziale di capacità di $C$. Sia $K_n$ un'esaustione regolare di $R$ con $C\subset K_1^{\circ}$ e sia $u_n$ il potenziale di capacità di $(K,K_n)$, $u_n'$ il potenziale di capacità di $(C,K_n)$. Grazie al principio del massimo applicato all'insieme $K_n\setminus C$ (vedi \ref{prop_max1}), si ha che $u_n\leq u_n'$, quindi $u_n\leq u'$, disuguaglianza valida per ogni $n$. Passando al limite si ottiene che $u\leq u'$. Se $K$ ha capacità nulla, il suo potenziale di capacità è la funzione costante uguale a $1$, quindi necessariamente anche $u'$ è costante uguale a $1$, cioè la capacità di $C$ è nulla.\\
Se $C\cap K\neq \emptyset$, allora considerando una bolla $B\subset C\cap K$, $B$ ha capacità nulla poiché contenuta in $K$, e quindi $C$ ha capacità nulla poiché contiene $B$. Infine se $C\cap K=\emptyset$, allora $K\subset C\cup K$, quindi grazie a quanto appena dimostrato $C\cup K$ ha capacità e anche $C\subset C\cup K$ ha capacità nulla.
\end{proof}
\end{prop}
Osserviamo da questa definizione che modificando la metrica di $R$ su un insieme compatto, la varietà rimane parabolica o non parabolica come la varietà non modificata.
\begin{prop}
 Sia $(R,g)$ una varietà riemanniana, e sia $(R,g')$ un'altra varietà con $g=g'$ fuori da un compatto $K$. Allora $(R,g)$ è parabolica se e solo se $(R,g')$ lo è.
\begin{proof}
 Consideriamo un compatto $C$ che contenga $K$ nella sua parte interna. La capacità di questo compatto nelle due varietà è identica, perché su $C^C$ le due metriche coincidono, da cui la tesi.
\end{proof}
\end{prop}

\section{Bordo Armonico}
La prima caratterizzazione che diamo della parabolicità riguarda il bordo armonico di $R^*$.
\begin{prop}\label{prop_parrod}
 La varietà $R$ è parabolica se e solo se $1\in \Rod$.
\begin{proof}
 Una parte della dimostrazione è semplice. Se la capacità di un insieme compatto $K$ è nulla, questo implica che il suo potenziale armonico $u$ è costante uguale a $1$. Dato che $u$ per definizione è il limite locale uniforme del potenziale armonico $u_n$ della coppia $(K,K_n)$, e dato che $u_n$ è limitata da $1$, si ha che $u=B-\lim_n u_n$.\\
La successione $u_n$ inoltre è anche di Cauchy rispetto alla seminorma $D$, infatti se $m>n$:
\begin{gather*}
 D_R(u_n-u_m)=D_R(u_n)-2D_R(u_n;u_m)+D_R(u_m)
\end{gather*}
Poiché la funzione di Tonelli $u_m-u_n$ è nulla sul bordo di $K_m\setminus K$, grazie alla formula di Green \ref{prop_green3} osserviamo che
\begin{gather*}
 D_R(u_m-u_n;u_m)=D_{K_m\setminus K}(u_m-u_n;u_m)=0
\end{gather*}
Questo ci permette di concludere che:
\begin{gather*}
 2D_R(u_m) - 2D_R(u_m;u_n)=0
\end{gather*}
da cui
\begin{gather*}
  D_R(u_n-u_m)=D_R(u_n)-D_R(u_m)
\end{gather*}
poiché la successione $D_R(u_n)$ è decrescente e converge a $\text{Cap}(K)$, abbiamo dimostrato che $\{u_n\}$ è $D-$Cauchy. Quindi $1=BD-\lim_n u_n$, e poiché $u_n \in \Roo$, $1\in \Rod$.\\
Per dimostrare l'implicazione inversa, sia $u$ il potenziale di capacità di un compatto $K$ e sia $f_n$ una successione in $\Roo$ tale che:
\begin{gather*}
 1=BD-\lim_n f_n
\end{gather*}
Chiaramente:
\begin{gather*}
 u=B-\lim_n uf_n
\end{gather*}
Per dimostrare che $u=D-\lim_n uf_n$ consideriamo un compatto qualsiasi $C\Subset R$, e osserviamo che:
\begin{gather*}
 D_R(uf_n-u)=\int_R \abs{\nabla(u(f_n-1))}^2 dV \leq \\
\leq2\int_R \abs u ^2 \abs{\nabla (f_n-1)}^2 dV + 2 \int_R \abs{\nabla u}^2 \abs{f_n-1}^2 dV
\end{gather*}
il primo termine della somma converge a $0$ per ipotesi, mentre per il secondo termine osserviamo che:
\begin{gather*}
 \int_R \abs{\nabla u}^2 \abs{f_n-1}^2 dV= \int_{C} \abs{\nabla u}^2 \abs{f_n-1}^2 dV + \int_{R\setminus C} \abs{\nabla u}^2 \abs{f_n-1}^2 dV
\end{gather*}
dove $C$ è un compatto qualsiasi. Dato che $u\in \Ro$ (quindi ha integrale di Dirichlet finito) e $f_n$ converge localmente uniformemente a $1$, il primo addendo tende a zero indipendentemente dal compatto scelto. Detto $b=\limsup_n \abs{f_n-1}^2$ \ \footnote{sicuramente finito grazie all'uniforme limitatezza della successione $\{f_n\}$}, abbiamo che per ogni $C\Subset R$:
\begin{gather*}
 \limsup_n D_R(uf_n-u) \leq 2b \int_{R\setminus C} \abs{\nabla u}^2 dV
\end{gather*}
vista l'arbitrarietà di $C$, e dato che $\int_{R} \abs{\nabla u}^2 dV <\infty$, possiamo concludere che
\begin{gather*}
 \lim_n D_R(uf_n-u) =0
\end{gather*}
Sia $F_n$ un compatto dal bordo regolare che contiene $supp(f_n)$, allora grazie alla formula di Green \ref{prop_green3} applicata all'insieme $A_n=F_n\setminus K$ , osserviamo che
\begin{gather*}
 D_R((1-u)f_n,1-u)=D_{A_n}((1-u)f_n,1-u)=0
\end{gather*}
poiché la funzione di Tonelli $(1-u)f_n$ è nulla sul bordo di $A_n$.\\
Grazie alle ultime due considerazioni possiamo concludere che:
\begin{gather*}
 D_R(u)=D_R(1-u)\equiv D_R(1-u;1-u)=\lim_n D_R((1-u)f_n,1-u)=0
\end{gather*}
quindi la capacità di $K$ è nulla, il che dimostra la parabolicità di $R$.
\end{proof}
\end{prop}
Come corollario di questa proposizione osserviamo che:
\begin{prop}
 Una varietà $R$ è parabolica se e solo se
\begin{enumerate}
 \item $1\in \Rod$
 \item $\Rod=\Ro$
 \item $\Delta=\emptyset$
\end{enumerate}
\begin{proof}
 Il punto (1) è il contenuto della proposizione precedente. Dato che $\Rod$ è un'ideale di $\Ro$, se $1\in \Rod$ necessariamente $\Ro=\Rod$ e viceversa. L'equivalenza tra (2) e (3) è il contenuto dell'osservazione \ref{oss_rodro}.
\end{proof}
\end{prop}
Osserviamo che grazie ad una forma del principio del massimo (la forma riportata in \ref{teo_max}), una varietà parabolica non ammette funzioni armoniche limitate con integrale di Dirichlet finito che non siano costanti. Il viceversa non è vero, ad esempio $\R^n$ con la metrica euclidea standard non ammette funzioni armoniche limitate non costanti pur non essendo una varietà parabolica (per $n\geq3$) \footnote{vedi teorema 2.1 pag. 31 di \cite{24}}.\\
Possiamo migliorare questa osservazione, infatti non è necessario chiedere che la funzione armonica abbia integrale di Dirichlet finito, e neanche che sia limitata, basta una delle due condizioni per dimostrare che la funzione è costante. 

\begin{prop}
 Una varietà Riemanniana $R$ parabolica non ammette funzioni armoniche limitate non costanti.
\begin{proof}
 Sia $u$ una funzione subarmonica limitata dall'alto su $R$. A meno di una traslazione, possiamo supporre $u\geq 0$. Sia $a=\sup_{x\in R}
u(x)<\infty$. Consideriamo un'esaustione regolare $K_n$ $n=0,1,\cdots$, e siano $(1-\omega_n)$ i potenziali armonici della coppia $(K_0,K_n)$, e sia $b=\max_{x\in K_0} u(x)$. Allora grazie al principio del massimo applicato a $K_n\setminus K_0$ abbiamo che:
\begin{gather*}
 u(x)\leq b+ a \omega_n
\end{gather*}
per ogni $n$. Passando al limite, visto che $R$ è parabolica, otteniamo che:
\begin{gather*}
 u(x)\leq b+a0=b
\end{gather*}
Questo è valido per ogni $K_0$ dominio relativamente compatto con bordo liscio. Se consideriamo l'insieme delle bolle coordinate centrate in un punto qualsiasi $x_0$, otteniamo quindi che $u(x)\leq u(x_0)$. Ripetendo il ragionamento con $-u$, otteniamo che $u(x)=u(x_0)$ per ogni $x\in R$.
\end{proof}
\end{prop}
\begin{prop}
 Una varietà parabolica $R$ non ammette funzioni armoniche con integrale di Dirichlet finito che non siano costanti.
\begin{proof}
 Sia $u$ una funzione con integrale di Dirichlet $D_R(u)<\infty$. Dato che $R$ è parabolica, $1\in \Rod$ \footnote{vedi \ref{prop_parrod}}, quindi esiste una successione di funzioni $\phi_n\in\Roo$ tale che $1=BD-\lim_n \phi_n$. Consideriamo la funzione
\begin{gather*}
 u_m(z)\equiv \min\{\max\{u(z),-m\}m\}
\end{gather*}
cioè $u_m$ è la funzione $u$ troncata in $[-m,m]$. Grazie alla formula di Green \ref{prop_green2} se consideriamo un insieme compatto dal bordo liscio $K$ tale che $supp(\phi_n)\subset K$, otteniamo che:
\begin{gather*}
 D_R(\phi_n u_m, u) =D_K(\phi_n u_m, u)=0
\end{gather*}
poiché la funzione $\phi_n u_m |_{\partial K}=0$. Dato che $1=BD-\lim_n \phi_n$, $u_m =BD-\lim_n \phi_n u_m$. Infatti:
\begin{gather*}
 D_R(\phi_n u_m - u_m )=\int_{R} \abs{\nabla{(\phi_n-1)u_m}}^2 dV = \int_R \abs{\nabla(\phi_n-1) u_m + (\phi_n-1)\nabla u_m}^2 dV
\end{gather*}
dato che $u_m$ è limitata e ha integrale di Dirichlet finito, si ottiene facilmente che il limite per $n\to \infty$ di questa quantità è $0$. Questo dimostra che $u_m=D-\lim_n \phi_n u_m$, il fatto che $u_m=B-\lim_n \phi_n u_m$ è quasi scontato.\\
Dato che per ogni $n$, $D_R(\phi_n u_m, u)=0$, si ha che per ogni $m$:
\begin{gather*}
 D_R(u_m,u)=0
\end{gather*}
\`E facile verificare che $u=CD-\lim_m u_m$, quindi abbiamo che:
\begin{gather*}
 D_R(u)\equiv D_R(u,u)=\lim_m D_R(u_m,u)=0
\end{gather*}
da cui $u$ è costante.
\end{proof}

\end{prop}

\section{Funzioni di Green}\label{sec_green2}
Un'altra caratterizzazione delle varietà paraboliche riguarda l'esistenza di funzioni di Green definite su tutta la varietà.
\begin{prop}
 Consideriamo un'esaustione regolare $K_n$ di $R$, e siano $G_{n}\equiv G_{K_n}$ le funzioni di Green \footnote{estese a $0$ fuori da $K_n$} rispetto a questi compatti. La successione di funzioni
\begin{gather*}
 G_{n}(\cdot, p)
\end{gather*}
con $p\in K_1$ fissato converge a una funzione armonica positiva $G(\cdot,p)$ su $R\setminus \{p\}$ se e solo se la varietà $R$ non è parabolica.
\begin{proof}
Per $m>n$, sia
\begin{gather*}
 \delta_n^m(p) \equiv G_m(\cdot, p)-G_n(\cdot, p)
\end{gather*}
dove $\delta_n^m$ è definita sull'insieme $K_n$. Grazie al fatto che per $d(x,p)\to 0$ le due funzioni $G_n$ e $G_m$ hanno lo stesso comportamento asintotico e grazie alla proposizione \ref{prop_sing2} $\delta_n^m$ è una funzione armonica su tutto $K_n$ per ogni $m$.\\
Osserviamo che per ogni $n$, $\delta_n^{n+1}$ è una funzione strettamente positiva su $K_n$, infatti su $\partial K_n$ $G_n(\cdot,p)=0$, mentre per il principio del massimo $G_{n+1}(\cdot,p)>0$, quindi $\delta_n^{n+1}|_{\partial K_n}>0$, e sempre per il principio del massimo $\delta_n^{n+1}>0$.\\
Dato che
\begin{gather*}
 \delta_n^m=\sum_{i=n+1}^m \delta_n^i
\end{gather*}
otteniamo che al variare di $m$, $\delta_n^m$ è una successione di funzioni armoniche positive crescenti. Per il principio di Harnack (vedi \ref{prop_harnackpri}), la successione $\delta_n^m$ al variare di $m$ converge localmente uniformemente (in $K_n$) o diverge localmente uniformemente. \`E evidente che il comportamento delle successioni $\delta_n^m$ è indipendente dal parametro $n$, infatti sull'intersezione dei vari insiemi di definizione vale che
\begin{gather*}
 \delta_n^m=\delta_k^m-\delta_n^k
\end{gather*}
quindi se e solo se al variare di $m$ $\delta_n^m$ è limitata, anche $\delta_k^m$ lo è.\\
Consideriamo $K\subset K_1$ un compatto contenente $p$ come punto interno, allora dalla proposizione \ref{prop_cap3}, sappiamo che per ogni $m$:
\begin{gather*}
 \min_{x\in \partial K} G_{m}(x,p) \leq \text{Cap}(K,K_m)^{-1}\leq \max_{x\in \partial K} G_{m}(x,p)
\end{gather*}
facendo tendere $m$ a infinito, osserviamo che se $\text{Cap}(K)=\lim_m \text{Cap}(K,K_m)>0$, allora necessariamente $G_{m}|_{\partial K}(\cdot,p)$ è limitata, altrimenti tende a infinito. Dato che al variare di $m$ $\delta_n^m$ è limitata se e solo se $G_m(\cdot,p)|_{\partial K}$ lo è, allora indipendentemente da $n$, $\delta_n^m$ converge su $K_n$ a una funzione armonica se e solo se $R$ è non parabolica.\\
Quindi la successione $G_m(\cdot,p)$, che sul compatto $K_n$ è uguale a $G_n(\cdot,p)+\delta_n^m (\cdot)$ \footnote{se $m>n$}, converge localmente uniformemente in $R$ \footnote{anche se le funzioni $G_m$ non sono definite in $p$, la loro differenza può essere estesa in $p$, e in questo senso diciamo che la convergenza è uniforme anche su $p$} a una funzione armonica se e solo se $\delta_n^m$ converge, quindi se e solo se $R$ non è parabolica.
\end{proof}
\end{prop}
Da questa dimostrazione deduciamo che se $\Omega\Subset \Omega'$, allora $G_{\Omega}(\cdot,p)\leq G_{\Omega'}(\cdot,p)$, quindi che la funzione $G(\cdot,p)$ ottenuta come limite di $G_n$ è indipendente dalla scelta dell'esaustione $K_n$.
\begin{prop}
 La funzione $G(\cdot,p)=\lim_n G_n(\cdot,p)$ è indipendente dalla scelta dell'esaustione $K_n$ utilizzata per definire $G_n$.
\begin{proof}
 Siano $K_n$ e $K_m'$ due esaustioni regolari di $R$ con $p\in K_1\cap K_1'$, e siano $G_n(\cdot,p)$ e $G_m'(\cdot,p)$ i relativi nuclei di Green, dove
\begin{gather*}
 G(\cdot,p)\equiv \lim_n G_n(\cdot,p) \ \ \ G'(\cdot,p)\equiv \lim_m G_m'(\cdot,p)
\end{gather*}
Per ogni $n$, esiste $\bar m$ tale che $K_n\Subset K'_{\bar m}$, allora grazie all'osservazione precedente:
\begin{gather*}
 G_n(\cdot,p)\leq G_{\bar m}'(\cdot,p)\leq G'(\cdot,p)
\end{gather*}
quindi passando al limite su $n$, $G(\cdot,p)\leq G'(\cdot,p)$, e poiché i ruoli di $G_n$ e $G_m'$ sono simmetrici, vale anche il viceversa, quindi
\begin{gather*}
 G(\cdot,p)=G'(\cdot,p)
\end{gather*}
\end{proof}
\end{prop}
Osserviamo che anche per le funzioni di Green definite su tutta $R$ vale una proprietà analoga a \ref{prop_maxG}:
\begin{prop}\label{prop_maxG2}
Sia $G(\cdot,p)$ la funzione di Green ottenuta con il metodo di esaustione, allora se $K$ e $K'$ sono domini relativamente compatti con $p\in K\Subset K'$, vale che:
\begin{gather*}
G|_{K' \setminus \overline K}(\cdot,p) \leq \max_{x\in \partial K} G(x,p) \ \ \ \ \ \ G(\cdot,p)|_{R \setminus \overline K} \leq \max_{x\in\partial K}G(x,p)
\end{gather*}
\begin{proof}
Sia $K_n$ un'esaustione regolare di $R$ con $p\in K_1$ e siano $G_n(\cdot,p)$ le relative funzioni di Green. Se $n$ è abbastanza grande per cui $K'\Subset K_n$, allora vale che:
\begin{gather*}
 G_n|_{K' \setminus \overline K}(\cdot,p) < \max_{x\in\partial K}G_n(x,p)
\end{gather*}
poiché questa relazione vale definitivamente, passando al limite su $n$ otteniamo che:
\begin{gather*}
  G|_{K' \setminus \overline K}(\cdot,p) \leq \max_{x\in\partial K}G(x,p)
\end{gather*}
Data l'indipendenza da $K'$ di questa proprietà, possiamo concludere che:
\begin{gather*}
 G(\cdot,p)|_{R \setminus \overline K} \leq \max_{x\in\partial K}G(x,p)
\end{gather*}
come volevasi dimostrare.
\end{proof}

\end{prop}

Grazie alle tecniche usate, possiamo facilmente dimostrare anche alcune proprietà della funzione $G(x,p)=\lim_n G_n(x,p)$ che corrispondono alle proprietà delle funzioni $G_n$.
\begin{prop}
 Per la funzione $G(x,p)=\lim_n G_n(x,p)$ vale che:
\begin{enumerate}
 \item per ogni dominio relativamente compatto con bordo liscio $\Omega$ contenente $p$, vale che $G(\cdot,p)- G_{\Omega}(\cdot,p)$ è una funzione estendibile a una funzione armonica su tutto $\Omega$.
 \item $G$ è strettamente positiva su $\Omega$
 \item $G$ è simmetrica, cioè $G(p,q)=G(q,p)$
 \item Fissato $q\in \Omega$, la funzione $G(q,p)$ è armonica rispetto a $p$ sull'insieme $\Omega\setminus \{q\}$ e superarmonica su tutto $\Omega$.
 \item $G$ è soluzione fondamentale dell'operatore $\Delta$, cioé per ogni funzione liscia $f$ a supporto compatto in $R$:
\begin{gather*}
 \Delta_x \int_{R} G(x,y)f(y) dy =\int_{R} G(x,y)\Delta_y(f)(y) dy =-f(x)
\end{gather*}
questo significa che nel senso delle distribuzioni $\Delta_y G(x,y)=-\delta_x$
\item Il flusso di $G_{\Omega}(\ast,p)$ attraverso il bordo di un'insieme regolare $K\Subset \Omega$ con $p\not \in \partial K$ vale:
\begin{gather*}
 \int_{\partial K}\ast d G(\cdot,p)=\begin{cases}
                                   -1 & se \ p\in K\\
0 & se \ p\not \in K
                                  \end{cases}
\end{gather*}

\item La funzione $G$ ha un comportamento asintotico della forma:
\begin{gather*}
 G(x,y)\sim C(m) \begin{cases}
                  -log(d(x,y)) & m=2\\
		  d(x,y)^{m-2} & m\geq 3
                 \end{cases}
\end{gather*}
quando $d(x,y)\to 0$. La costante $C(m)$ dipende solo dalla dimensione della varietà e può essere determinata sfruttando la condizione (5).
\end{enumerate}
\begin{proof}
In tutta la dimostrazione supponiamo che $G$ esista, cioè che $R$ non sia parabolica.\\
 Il punto (1) è il punto chiave per dimostrare tutte le altre proprietà. Dalla costruzione nella proposizione precedente sappiamo che per ogni $n$ sull'insieme $K_n$ la funzione $$G(\cdot,p)-G_{n}(\cdot,p)=\lim_m\delta_n^m\equiv \delta_n$$ è una funzione armonica su $K_n$. Consideriamo una funzione di Green $G_{\Omega}(\cdot,p)$ con $\Omega$ qualsiasi, e sia $\Omega \Subset K_n$. Dalla proposizione \ref{prop_sing2} sappiamo che $G_n(\cdot,p)-G_{\Omega}(\cdot,p)$ è estendibile a una funzione armonica $\phi$ definita su tutto $\Omega$. Quindi:
\begin{gather*}
 G(\cdot,p)-G_{\Omega}(\cdot,p)=G(\cdot,p)-G_n(\cdot,p)+G_n(\cdot,p)-G_{\Omega}(\cdot,p)=\delta_n(\cdot) -\phi(\cdot)
\end{gather*}
è una funzione armonica su $\Omega$.\\
Questo implica in particolare che fissato $q\in \Omega$, $G(q,p)$ si può scrivere come la somma tra la funzione $G_\Omega(q,p)$ e una funzione armonica, il che dimostra anche il punto (4).\\ 
La positività e la simmetria di $G$ sono una ovvia conseguenza del fatto che queste proprietà valgono anche per tutte le $G_n$.\\
Per dimostrare il punto (5), consideriamo una funzione $f\in C^{\infty}_C(R)$. Sia $\Omega$ tale che $supp(f)\subset \Omega$, allora grazie al punto (1) sappiamo che $G(x,y)=G_{\Omega}(x,y) + \phi(x,y)$ con $\phi(\cdot,y)$ armonica in $\Omega$ per $y$ fissato \footnote{e viceversa per simmetria}, quindi:
\begin{gather*}
  \Delta_x \int_{R} G(x,y)f(y) dy =\Delta_x \int_{\Omega} G_{\Omega}(x,y)f(y) dy + \Delta_x \int_{\Omega} \phi(x,y)f(y) dy=-f(x)
\end{gather*}
Grazie al fatto che $\Delta_x \phi(x,y)=0$ e al teorema di derivazione sotto al segno d'integrale.\\
Allo stesso modo:
\begin{gather*}
  \int_{R} G(x,y)\Delta_y f(y) dy =\int_{\Omega} G_{\Omega}(x,y)\Delta_y f(y) dy + \int_{\Omega} \phi(x,y)\Delta_y f(y) dy
\end{gather*}
grazie a un'integrazione per parti su $\Omega$ (il termine al bordo è nullo grazie al fatto che $supp(f)\Subset \Omega$):
\begin{gather*}
\int_{\Omega} \phi(x,y)\Delta_y f(y) dy =-\int_{\Omega} \Delta_y \phi(x,y)\Delta_y f(y) dy=0
\end{gather*}
otteniamo che:
\begin{gather*}
  \int_{R} G(x,y)\Delta_y f(y) dy =\int_{\Omega} G_{\Omega}(x,y)\Delta_y f(y) dy =-f(x)
\end{gather*}
In maniera del tutto analoga, si dimostra il punto (6), basta considerare il fatto che $\ast dG =\ast dG_{\Omega}+ \ast d\phi$ e il flusso di una funzione armonica attraverso il bordo di un compatto regolare è nullo grazie alla formula di green \ref{prop_green1}.\\
Il punto (7) è conseguenza del fatto che la funzione armonica $\phi$ è limitata.
\end{proof}
\end{prop}
Osserviamo che anche nel caso di varietà paraboliche è possibile definire delle funzioni di Green, solo che in questo caso le funzioni non sono positive \footnote{e neanche limitate dal basso}. Dall'articolo \cite{25} riportiamo il teorema:
\begin{teo}\label{teo_gcont}
 Sia $R$ una varietà riemanniana completa non compatta e senza bordo. Allora esiste un nucleo di Green simmetrico $G(x,y)$ liscio sull'insieme $R\times~R\setminus~D$ dove $D=\{(x,x) \ t.c. \ x\in R\}$ è la diagonale di $R\times R$. In particolare $G$ soddisfa:
\begin{gather*}
 \Delta_x \int_{R} G(x,y)f(y) dy =\int_{R} G(x,y)\Delta_y(f)(y) dy =-f(x)
\end{gather*}
per ogni funzione $f$ liscia a supporto compatto in $R$.
\end{teo}
\subsection{Funzioni di Green sulla compattificazione di Royden}
In questa sezione dimostriamo che per una varietà iperbolica le funzioni di Green si possono estendere alla compattificazione di Royden $R^*$. Utilizzeremo i risultati ottenuti nella sezione \ref{sec_green2}.
\begin{prop}
 Per ogni $c>0$, la funzione
\begin{gather*}
 G_c(x,p)\equiv G(x,p)\curlywedge c \equiv\min\{G(x,p);c\}
\end{gather*}
appartiene a $\Rod$, inoltre
\begin{gather*}
  D_R(G(x,p)\curlywedge c)\leq c
\end{gather*}
\begin{proof}
Data un'esaustione regolare $K_n$ sia $G_n(x,p)$ la funzione di Green relativa a $K_n$ e $G(x,p)$ il suo limite.\\
L'obiettivo della dimostrazione è mostrare che per ogni $c>0$ \footnote{questo è vero anche per $c\leq0$, ma è privo di senso essendo $G(x,p)>0$ su $R$} la successione $G_n(\cdot,p)\curlywedge c $ converge localmente uniformemente a $G_c(\cdot,p)$ \footnote{questa considerazione è quasi ovvia} e ha integrale di Dirichlet limitato, quindi dato che ovviamente $G_n(x,p)\curlywedge c \in \Roo\subset \Rod$, grazie al teorema \ref{prop_r1} $G_c\in~\Rod$.\\
 Per dimostrare che $D_R(G_n(x,p)\curlywedge c)<M(c)$, per prima cosa osserviamo che se $c'\leq c $
\begin{gather*}
 D_R(G_n(x,p)\curlywedge c')=\int_{G_n(x,p)<c'} \abs{\nabla(G_n(x,p))}^2 dV \leq\\
\leq \int_{G_n(x,p)<c} \abs{\nabla(G_n(x,p))}^2 dV = D_R(G_n(x,p\curlywedge c)
\end{gather*}
quindi basta dimostrare che $D_R(G_n(x,p)\curlywedge c)<M(c)$ è valida per un insieme illimitato di $c$.\\
 Consideriamo un insieme aperto relativamente compatto $p\in \Omega$. Grazie a \ref{prop_maxG2}, $G(x,p)$ ristretta a $R\setminus \Omega$ assume il suo massimo $M$ in $\partial \Omega$. Allora l'insieme $U_c\equiv \{x\ t.c. \ G(x,p)>c\}$ è contenuto in $\Omega$ se $c>M$, e per monotonia anche gli insiemi
\begin{gather*}
 U_c^n\equiv \{x\ t.c. \ G_n(x,p)>c\}
\end{gather*}
sono contenuti in $\Omega$ se $c>M$, quindi sono relativamente compatti.\\
All'interno di $(M,\infty)$, consideriamo i valori di $c$ che sono valori regolari per tutte le funzioni $G_n(\cdot,p)$. Per il teorema di Sard \footnote{e per il fatto che l'unione numerabile di insiemi di misura nulla ha misura nulla} questi valori sono densi in $(M,\infty)$. D'ora in avanti consideriamo solo $c$ con queste caratteristiche.\\
Osserviamo ora che per ogni $n$, l'integrale:
\begin{gather*}
 D_R(G_n(x,p)\curlywedge c)= \int_{K_n\cap (U_c^n)^C} \abs{\nabla(G_n(x,p))}^2 dV = \\
\int_{\partial (K_n\cap (U_c^n)^C)} G_n\ast dG_n = - c \int_{\partial U_c^n} \ast dG_n=c
\end{gather*}
grazie alle proprietà di $G_n$ descritte in \ref{subsec_green1}. Questo e il teorema \ref{prop_r1} danno la tesi. Inoltre osserviamo che dall'ultima considerazione possiamo ottenere che:
\begin{gather*}
 D_R(G(x,p)\curlywedge c)\leq c
\end{gather*}
per ogni $c>0$.
\end{proof}

\end{prop}
Grazie al fatto che $G_c(x,p)\in \Rod$ per ogni $c$, possiamo estendere la definizione del nucleo di Green a $R\times R^*$. Sia a questo scopo $U_c\equiv \{x\in R \ t.c. \ G(x,z)>c\}$, e sia $c$ tale che $U_c$ è relativamente compatto in $R$. Allora grazie all'osservazione \ref{prop_charG}, per ogni funzione liscia a supporto compatto $\lambda$, se $p\in \Gamma$ e $z\in R$ fisso:
\begin{gather*}
 p(G_c(\cdot,z))=p(G_c(\cdot,z)(1-\lambda))
\end{gather*}
se scegliamo $\lambda$ in modo che $supp(\lambda)\supset U_c$, questa relazione mostra che per ogni $c'>c$ si ha che:
\begin{gather*}
 p(G_c(\cdot,z))=p(G_{c'}(\cdot,z))
\end{gather*}
quindi il valore che assume $G_c(\cdot,z)$ sul punto $p\in \Gamma$ è indipendente dalla scelta di $c$, e quindi ha senso definire:
\begin{gather*}
 G(p,z)=\lim_{q\to p} G(q,z)
\end{gather*}
Da questa definizione segue che:
\begin{prop}\label{prop_g'}
 La funzione $G'(z,p)$ con $(z,p)\in R\times R^*$ è una armonica in $R\setminus \{p\}$, continua su $R\times R^*\setminus D$, $G'_c(\cdot,p)\equiv\min\{G'(\cdot,p), c\}\in \Rod$ per ogni $c>0$, e $D_R(G_c'(\cdot,p))\leq c$.
\begin{proof}
\`E necessario dimostrare la continuità della funzione $G'(\cdot,\cdot)$ solo sui punti $R\times \Gamma$, dato che $G'$ è estensione di $G$ che è continua su $R\times R\setminus D$. Consideriamo quindi una coppia $(z_0,p_0)\in R\times \Gamma$, $\epsilon>0 $, e siano $U$ e $V$ due intorni di $z_0$ e $p_0$ aperti in $R^*$ e a chiusura disgiunta. Fissato $z_0$, la funzione $G'(z_0,\cdot)$ è continua su $R^*\setminus \{z_0\}$, questo significa che possiamo scegliere $V$ in modo che per ogni $p\in W$:
\begin{gather*}
 \abs{G'(z_0,p)-G'(z_0,p_0)}\leq \epsilon
\end{gather*}
Sia $s=\sup_{p\in R\setminus U}{G'(z_0,p)}=\sup_{p\in R^*\setminus U}{G'(z_0,p)}<\infty$ \footnote{l'ultima uguaglianza segue dal fatto che $R$ è denso in $R^*$}. Consideriamo la funzione di Harnack definita sull'aperto $U\times U$. Allora vale che:
\begin{gather*}
 k(z,z_0)^{-1}G'(z_0,p)\leq G'(z,p)\leq k(z_0,z)G'(z_0,p)
\end{gather*}
per ogni punto $p\in W\cap R$. Data la proposizione \ref{prop_harnackfun}, esiste un intorno $U'\subset U$ tale che $k(U',z_0)\subset[1,1+\epsilon/s)$. Allora per $z\in U'$ e per ogni $p\in W\cap R$, si ha che:
\begin{gather*}
 \abs{G'(z,p)-G'(z_0,p)}\leq \epsilon
\end{gather*}
Sempre per continuità su $R^*$ di $G'(z_0,\cdot)$, si ha che in realtà questa relazione vale per ogni $p\in W$. Riassumendo otteniamo che se $(z,p)\in U'\times W$, si ha che:
\begin{gather*}
 \abs{G'(z_0,p_0)-G'(z,p)}\leq\abs{G'(z_0,p_0)-G'(z_0,p)}+\abs{G'(z_0,p)-G'(z,p)}\leq 2\epsilon
\end{gather*}
data l'arbitrarietà di $\epsilon$ si ottiene la tesi.\\
La funzione $G'(\cdot,p)$ è il limite di funzioni armoniche positive, quindi grazie al principio di Harnack (vedi \ref{prop_harnackpri}) è una funzione armonica. Sappiamo che per ogni $p\in R$:
\begin{gather*}
  D_R(G(x,p)\curlywedge c)\leq c
\end{gather*}
Sia $p_0\in \Gamma$. Dato che $G'(\cdot,p_0)\curlywedge c =C-\lim_{p\in R, \ p\to p_0} G(\cdot,p)\curlywedge c$, al teorema \ref{prop_r1}, abbiamo che $G'(\cdot,p_0)\curlywedge c \in \Rod$ e anche $D_R(G'(\cdot,p_0)\curlywedge c)\leq c$.
\end{proof}
\end{prop}
Grazie alle proprietà fino a qui dimostrate, possiamo definire $G^*(\cdot,\cdot)$ l'estensione di $G'(\cdot,\cdot)$ a tutto $R^*\times R^*\setminus D$. Infatti, essendo $G'(\cdot,p)\curlywedge c \in \Rod\subset \Ro$ per ogni $c>0$, ha senso definire:
\begin{gather*}
 G^*(z_0,p_0)\curlywedge c =\lim_{z\in R, \ z\to z_0} G'(z,p_0)\curlywedge c \equiv \lim_{z\in R, \ z\to z_0} \lim_{p\in R, \ p\to p_0} G(z,p)\curlywedge c 
\end{gather*}
per ogni $c>0$, quindi:
\begin{deph}
 Definiamo $G^*(\cdot,\cdot):R^*\times R^*\setminus D \to [0,\infty)$ con il doppio limite:
\begin{gather*}
 G^*(z_0,p_0) =\lim_{z\in R, \ z\to z_0} G'(z,p_0) \equiv \lim_{z\in R, \ z\to z_0} \lim_{p\in R, \ p\to p_0} G(z,p)
\end{gather*}
\end{deph}
Definiamo inoltre il \textbf{bordo essenziale irregolare} di $R^*$ come:
\begin{deph}\label{deph_Xi}
 L'insieme dei punti in $\Gamma$
\begin{gather*}
 \Xi\equiv \{p\in \Gamma \ t.c. \ G(z,p)>0 \ z\in R\}
\end{gather*}
è detto bordo essenziale irregolare di $R$.
\end{deph}
Riassumendo le proprietà fin qui dimostrate, per la funzione $G^*$ valgono le seguenti proprietà:
\begin{prop}\label{prop_greenref}
Sia $R$ una varietà non parabolica senza bordo, allora per la funzione $G^*(\cdot,\cdot)$, l'estensione del nucleo di Green alla compattificazione $R^*$, valgono le proprietà:
\begin{enumerate}
\item $G^*(x,y)=G(x,y) \ \forall (x,y)\in R$
\item $G^*(\cdot,p)$ è continua sull'insieme $R^*\setminus \{p\}$
\item $G^*(\cdot,p)$ è una funzione armonica su $R\setminus \{p\}$
\item $D_R(G^*(\cdot,p)\curlywedge c)\leq c$ per ogni $c>0$
\item $G^*(\cdot,\cdot)|_{(\Delta\times R^*)}=0$
\item $G^*(\cdot,\cdot)|_{(R\cup \Xi)\times (R\cup \Xi)}>0$
\end{enumerate}
\begin{proof}
 (1) e (2) e (3) seguono dal fatto che $G^*$ è un'estensione di $G'$ che è un'estensione di $G$, e visto che sia $G$ che $G'$ sono in $\Ro$, l'estensione è continua su tutta $R^*$ per definizione della compattificazione di Royden.\\
Grazie alla \ref{prop_g'}, possiamo dedurre (4) e (5), infatti dato che $G'\in \Rod$, $D_R(G'\curlywedge c)\leq c$ e grazie alla proposizione \ref{prop_r2}, possiamo dedurre che $G^*\curlywedge c \in \Rod$, quindi, $G^*(\cdot,\cdot)|_{\Delta\times R^*}=0$.\\
Rimane da dimostrare (6). Come spesso accade, questa dimostrazione è un'applicazione del principio del massimo. Sappiamo che $G^*(x,y)>0$ se $(x,y)\in R\times R$, e anche se $(x,y)\in R\times \Xi$, resta il caso $(x,y)\in \Xi\times \Xi$. La funzione $G^*(\cdot,y)$ è armonica positiva su $R$ (come dimostrato sopra). Quindi se consideriamo un aperto $U\in R$ e un punto $p\in U$, esiste un numero positivo $a$ tale che $aG^*(\cdot,y)-G(\cdot,p)|_{\partial U}>0$. Consideriamo una successione $G_n(\cdot,p)$ di funzioni di Green relative a un'esaustione regolare $R_n$. Allora definitivamente $aG^*(\cdot,y)-G_n(\cdot,p)|_{\partial U}>0$ e anche $aG^*(\cdot,y)-G_n(\cdot,p)|_{\partial R_n}>0$. Quindi per il principio del massimo $aG^*(\cdot,y)-G_n(\cdot,p)|_{R_n\setminus U}\geq0$. Allora questa disuguaglianza vale anche per il limite su $n$, cioè
\begin{gather*}
 aG^*(\cdot,y)-G(\cdot,p)|_{\partial U}\geq 0
\end{gather*}
dato che $G(x,p)>0$ per definizione di $\Xi$, si ha la tesi.
\end{proof}
\end{prop}
D'ora in avanti confonderemo la notazione di $G, \ G', \ G^*$ quando non ci sia rischio di confusione.

\section{Potenziali di Evans}
In questa sezione diamo un'altra caratterizzazione della parabolicità di una varietà $R$ attraverso l'esistenza di particolari funzioni armoniche, i potenziali di Evans.
\begin{deph}
 Dato un compatto $K\Subset R$, un potenziale di Evans rispetto a questo compatto è una funzione armonica $f:(R\setminus K)\to \R$ tale che
\begin{gather*}
 \lim_{x\to \infty}f(x)=\infty \ \ \ f|_{\partial K}=0
\end{gather*}
L'ultimo limite può essere inteso equivalentemente in 2 sensi: per ogni successione $x_n\to \infty$, $f(x_n)\to \infty$, oppure per ogni $N>0$, esiste $K\subset K_N\Subset R$ tale che $f(K_N^C)\subset(N,+\infty)$.\\
Se la funzione $f$ è solo superarmonica in $R\setminus K$, viede definita potenziale di Evans superarmonico.
\end{deph}
I potenziali di Evans sono caratteristici delle varietà paraboliche, nel senso che una varietà è parabolica se e solo se per ogni compatto $K$ esiste un potenziale di Evans relativo a questo compatto. La dimostrazione di un'implicazione è quasi immediata, mentre l'altra implicazione è il contenuto fondamentale di questa tesi.\\
Osserviamo che per caratterizazione delle varietà paraboliche, è sufficiente chiedere che il potenziale di Evans sia superarmonico.
\begin{prop}
 Se esiste un compatto $K$ che ammette un potenziale di Evans superarmonico, allora la varietà è parabolica
\begin{proof}
 Sia $K$ un compatto e $f$ il relativo potenziale di Evans superarmonico. Consideriamo un compatto $C$ con $K\Subset C^{\circ}$ e bordo liscio. La funzione $f$ continua a essere superarmonica positiva su $C^C$ e tende a infinito. Inoltre grazie al principio del massimo \ref{prop_max2} $f|_{\partial C}>0$. Dimostriamo che il potenziale armonico $u$ di $C$ è costante uguale a $1$.\\
Sia $b=1-u$ funzione armonica definita su $C^C$ \footnote{estendibile per continuità a $0$ su $C$}, e sia $C_k$ un'esaustione regolare di $R$ con $C_1=C$. Per ogni $N>1$, esiste $\bar k$ tale che $f(C_{\bar k}^C)\subset(N,\infty)$. Dato che la funzione $b$ è necessariamente $\leq 1$ \footnote{indipendendemente dal fatto che $R$ sia parabolica}, e che $Nb|_{\partial C}<f|_{\partial C}$, per il principio del massimo abbiamo che per ogni $k>\bar k$, $Nb\leq f$ sull'insieme $C^k\setminus C$, quindi su tutto l'insieme $R\setminus C$. Data l'arbitrarietà di $N$, otteniamo che per ogni $N>1$
\begin{gather*}
 Nb\leq f \ \ \Rightarrow \ \ b\leq \frac{f}{N}
\end{gather*}
quindi necessariamente $b=0$, cioè $u=1$, come volevasi dimostrare.
\end{proof}
\end{prop}
I paragrafi seguenti si occupano di dimostrare l'implicazione inversa.

\subsection{Diametro transfinito}
In questo paragrafo definiamo due strumenti che saranno utili per costruire particolari funzioni armoniche, tra cui i potenziali di Evans. D'ora in avanti assumiamo che $R$ sia una varietà riemanniana non parabolica.
\begin{deph}
 Detto $G^*(\cdot,\cdot)$ il nucleo di Green su $R^*$, e dato $X\subset R^*$ definiamo:
\begin{gather}
 \binom{n}{2} \rho_n(X)\equiv \inf_{p_1,\cdots,p_n \in X} \sum_{i<j}^{1\cdots n} G(p_i,p_j)
\end{gather}
inoltre per convenzione $\rho(\emptyset)=\infty$
\end{deph}
\begin{prop}
 Fissato l'insieme $X$, $\rho_n(X)$ è una successione crescente al crescere di $n$.
\begin{proof}
 La dimostrazione è un semplice esercizio di algebra, non riguarda le caratteristiche di $R$.\\
Siano $p_1,\cdots,p_n$ punti qualsiasi in $X$, e consideriamo per ogni $1\leq k\leq n+1$:
\begin{gather*}
 \sum_{i<j }^{1\cdots n+1} G(p_i,p_j)= \sum_{i=1 }^{k-1} G(p_i,p_k)+ \sum_{j=k+1 }^{n+1} G(p_k,p_j)+ \sum_{i<j; \ i,j\neq k }^{1\cdots n} G(p_i,p_j)
\end{gather*}
e quindi dalla definizione di $\rho_n(X)$:
\begin{gather*}
 \sum_{i<j }^{1\cdots n+1} G(p_i,p_j)\geq \sum_{i=1 }^{k} G(p_i,p_k)+ \sum_{j=k+1 }^{n+1} G(p_k,p_j)+ \binom{n}{2}\rho_n(X)
\end{gather*}
Sommando tutte queste disuguaglianze al variare di $k$ tra $1$ e $n+1$, otteniamo che:
\begin{gather*}
(n+1) \sum_{i<j }^{1\cdots n+1} G(p_i,p_j)\geq 2\sum_{i<j }^{1\cdots n+1} G(p_i,p_j) + (n+1)\binom{n}{2}\rho_n(X)
\end{gather*}
cioè:
\begin{gather*}
 (n-1) \sum_{i<j }^{1\cdots n+1} G(p_i,p_j)\geq (n+1)\binom{n}{2}\rho_n(X)
\end{gather*}
vista l'arbitrarietà dei punti considerati, possiamo concludere che:
\begin{gather*}
 (n-1)\binom{n+1}{2}\rho_{n+1}(X)\geq(n+1)\binom{n}{2} \rho_n(X)
\end{gather*}
considerando che
\begin{gather*}
 (n-1)\binom{n+1}{2}=(n+1)\binom{n}{2}
\end{gather*}
otteniamo la tesi, cioè
\begin{gather*}
 \rho_{n+1}(X)\geq \rho_n(X)
\end{gather*}
\end{proof}
\end{prop}
Questa proposizione ci permette di definire il limite al variare di $n$ di $\rho_n(X)$, in particolare:
\begin{deph}
 Dato $X\in R^*$, definiamo il \textbf{diametro transfinito} di $X$:
\begin{gather}
 \rho(X)\equiv \lim_n \rho_n(X)
\end{gather}
\end{deph}
Assieme al diametro transfinito, definiamo la \textit{costante di Tchebycheff}, e esploriamo alcuni legami tra i due concetti.
\begin{deph}\label{deph_tau}
 Dato $X\subset R^*$, definiamo:
\begin{gather*}
 n\tau_n(X)=\sup_{p_1,\cdots,p_n\in X} \ton{\inf_{p\in X} \sum_{i=1}^n G(p,p_i)}
\end{gather*}
e per convenzione $\tau_n(\emptyset)=\infty$
\end{deph}
Anche in questo caso vale una relazione che ci consente di definire il limite di $\tau_n$.\\
Infatti osserviamo che:
\begin{gather*}
 \sum_{i=1}^{n+m} G(p,p_i)=\sum_{i=1}^{n} G(p,p_i)+\sum_{i=n+1}^{n+m} G(p,p_i)
\end{gather*}
quindi applicando l'$\inf$ a entrambi i membri otteniamo:
\begin{gather*}
 \inf_{p\in X} \sum_{i=1}^{n+m} G(p,p_i)\geq \inf_{p\in X}\sum_{i=1}^{n} G(p,p_i)+ \inf_{p\in X}\sum_{i=n+1}^{n+m} G(p,p_i)
\end{gather*}
cioè per definizione:
\begin{gather}\label{eq_tau1}
 (n+m)\tau_{n+m}(X)\geq n\tau_n(X)+m\tau_m(X)
\end{gather}
se consideriamo $n=m$ otteniamo il caso particolare:
\begin{gather}\label{eq_tau2}
 (qm)\tau_{qm}(X)\geq q\tau_{m}(X) \Rightarrow \tau_{qm}(X)\geq \tau_m(X)
\end{gather}
Grazie a questa relazione possiamo dimostrare che:
\begin{prop}
 Sia $\alpha\equiv \sup_n \tau_n(X)$, allora $\lim_n \tau_n(X)=\alpha$.
\begin{proof}
 La dimostrazione segue dalla relazione \ref{eq_tau1}. Infatti scegliamo un qualunque $\beta<\alpha$. Per definizione di $\sup$, esiste un $m$ tale che $\tau_m(X)>\beta$. \`E un fatto noto che ogni numero $n$ può essere scritto in maniera univoca come:
\begin{gather*}
 n=qm+r
\end{gather*}
dove $q,r\in \N$ e $0\leq r\leq m-1$. Quindi per $\tau_n(X)$ vale che:
\begin{gather*}
 n\tau_n(X)=(qm+r)\tau_{qm+r}(X)\geq qm \tau_{qm}(X) + r\tau_r(X) \geq qm\tau_m(X)
\end{gather*}
dove abbiamo utilizzato le relazioni \ref{eq_tau1} ed \ref{eq_tau2}. Quindi possiamo concludere che
\begin{gather*}
 \tau_n(X)\geq \frac{qm}{qm+r}\tau_m(X)>\frac{qm}{qm+r} \beta
\end{gather*}
Se $n$ tende a infinito, $q$ tente a infinito, mentre $0\leq r \leq n-1$ continua a valere, quindi:
\begin{gather*}
 \liminf_n \tau_n(X)\geq \liminf_{n} \frac{qm}{qm+r}\beta =\beta 
\end{gather*}
Cioè per ogni $\beta<\alpha=\sup_n \tau_n(X)$ abbiamo che:
\begin{gather*}
 \alpha\geq \limsup_n \tau_n(X) \geq \liminf_n \tau_n(X) \geq \beta
\end{gather*}
data l'arbitrarietà di $\beta$, otteniamo la tesi.
\end{proof}
\end{prop}
\begin{deph}
 Dato $X\subset R$ definiamo la sua \textbf{costante di \Tc} il limite:
\begin{gather*}
 \tau(X)\equiv \lim_{n\to \infty} \tau_n(X)
\end{gather*}
\end{deph}
Il diametro transfinito di un insieme e la costante di \Tc\ misurano in qualche senso la grandezza di un insieme. Più il nucleo di Green $G^*$ ha un valore alto sui punti di $X$, più queste due costanti hanno valore alto, ed è facile osservare che per insiemi di cardinalità finita $\tau(X)=\rho(X)=\infty$, e che se $X'\subset X$, allora $\tau(X')\geq \tau(X)$. Di seguito riportiamo un'altra proprietà di queste due costanti.
\begin{prop}\label{prop_cfrdiam}
 Per ogni $X\subset R^*$ si ha che:
\begin{gather}\label{eq_tau>rho}
 \tau(X)\geq \rho(X)
\end{gather}
\begin{proof}
Assumiamo che $X$ abbia cardinalità infinita \footnote{altrimenti abbiamo visto che $\tau(X)=\rho(X)=\infty$}. Per $n>1$ sia $r\equiv \frac{1}{n-1}$. \`E possibile scegliere $n$ punti $p_1,\cdots,p_n$ tali che per ogni $i=1,\cdots,n-1$ abbiamo che:
\begin{gather}\label{eq_rhotau1}
 \sum_{j=n+1-i}^n G(p_{n-i};p_j)\leq \inf_{p\in X} \sum_{j=n+1-i}^n G(p,p_j) + r
\end{gather}
Dimostriamo questa affermazione per induzione su $i$. Per $i=1$ l'affermazione segue direttamente dalla definizione di $\inf$. Infatti se scegliamo $p_n$ arbitrariamente e $p_{n-1}$ in modo che
\begin{gather*}
 G(p_{n-1},p_n)\leq \inf_{p\in X} G(p,p_n) +r
\end{gather*}
Supponiamo che la tesi sia vera per $1\leq i<n-1$ \footnote{quindi supponiamo di aver determinato $p_n,\cdots, p_{n-i+1}$} e consideriamo la funzione
\begin{gather*}
 f(p)=\sum_{j=n-i+1}^n G(p,p_i)
\end{gather*}
osserviamo che se $p=p_k$ per qualche $k=n-i+1,\cdots n$, allora $f(p)$ vale infinito. La funzione $f$ comunque è positiva, quindi possiamo trovare un punto $p_{n-i}$ tale che:
\begin{gather*}
 f(p_{n-i})\leq \inf_{p\in X} f(p) +r
\end{gather*}
quindi otteniamo che esistono punti $p_1,\cdots,p_n\in X$ tali che
\begin{gather*}
  \sum_{j=n+1-i}^n G(p_{n-i};p_j)\leq \inf_{p\in X}\sum_{j=n+1-i}^n G(p;p_j)\leq i\tau_i(X) +r
\end{gather*}
sommando queste disuguaglianze per $i=1,\cdots,n-1$ otteniamo che
\begin{gather*}
 \sum_{i<j}^{1\cdots n} G(p_i;p_j)\leq \sum_{i=1}^{n-1} i\tau_i(X) + 1
\end{gather*}
e dalla definizione di $\rho_n$ abbiamo:
\begin{gather*}
 \binom{n}{2}\rho_n(X) \leq \sum_{i=1}^{n-1} i\tau_i(X) + 1
\end{gather*}
quindi:
\begin{gather*}
 \rho_n(X) \leq \binom{n}{2}^{-1}\ton{\sum_{i=1}^{n-1} i\tau_i(X) + 1}
\end{gather*}
Passando al limite otteniamo la tesi, cioè:
\begin{gather*}
 \rho(X) = \lim_n \rho_n(X) \leq \lim_n\binom{n}{2}^{-1}\ton{\sum_{i=1}^{n-1} i\tau_i(X) + 1}=\tau(X)
\end{gather*}
dove l'ultimo passaggio è giustificato nel seguente lemma
\end{proof}
\end{prop}
\begin{lemma}
 Sia $a_n$ tale che $\lim_n a_n =a$, allora:
\begin{gather*}
 \lim_n \sum_{k=1}^{n-1} \frac{ka_k}{\binom{n}{2}}= a
\end{gather*}
\begin{proof}
Per definizione di limite, per ogni $\epsilon>0$, esiste $N$ tale che $a_k>a-\epsilon$ per ogni $k>N$, quindi per ogni $\epsilon>0$ e per ogni $n>N+1$ (quindi definitivamente):
\begin{gather*}
 \binom{n}{2}^{-1} \sum_{k=1}^{n-1} ka_k= \binom{n}{2}^{-1}\ton{\sum_{k=1}^N ka_k +\sum_{k=N+1}^{n-1} ka_k}>\\
>\binom{n}{2}^{-1}\ton{\sum_{k=1}^N ka_k} + (a-\epsilon) \binom{n}{2}^{-1}\ton{\sum_{k=N+1}^{n-1} k}
\end{gather*}
applicando il $\liminf_n$, e osservando che $\binom{n}{2}\to \infty$, $ \sum_{k=1}^{N} ka_k$ è costante al variare di $n$ e $\sum_{k=N+1}^{n-1} k= \binom{n}{2}-\binom{N+1}{2}$ otteniamo:
\begin{gather*}
 \liminf_n \binom{n}{2}^{-1} \sum_{k=1}^{n-1} ka_k \geq a-\epsilon
\end{gather*}
con un ragionamento del tutto analogo si ottiene anche:
\begin{gather*}
 \limsup_n \binom{n}{2}^{-1} \sum_{k=1}^{n-1} ka_k < a+\epsilon
\end{gather*}
e per l'arbitrarietà di $\epsilon$, otteniamo la tesi.
\end{proof}
\end{lemma}
\subsection{Stime per il diametro transfinito}
Lo scopo di questo paragrafo è ottenere la stima riportata nella proposizione \ref{prop_stima}, una stima tecnica che servirà nei paragrafi successivi a dimostrare che per ogni $\Xi'\Subset \Xi$, $\rho(\Xi')=\tau(\Xi')=\infty$. L'affermazione è vuotamente vera se $R$ è una varietà iperbolica regolare (cioè se $\Xi=\emptyset$), quindi in tutta la sezione assumeremo che $R$ sia una varietà iperbolica irregolare.\\
A questo scopo, introduciamo un insieme di funzioni ausiliarie che serviranno a stimare il diametro di $\Xi_n$.\\
Fissato un punto $z_0\in R$, sia $r_n$ una successione di numeri reali positivi tali che
\begin{gather*}
 r_n>r_{n+1} \ \ \  \lim_n r_n=0
\end{gather*}
inoltre chiediamo che $U_n\equiv\{z\in R \ t.c. \ G(z,z_0)>r_n\}$ non sia relativamente compatto e che $r_n$ sia un valore regolare di $G(\cdot,z_0)$ \footnote{in modo che $U_n$ sia un insieme con bordo regolare}. In questo modo possiamo descrivere l'insieme $\Xi$ come unione di $\Xi_n$, dove
\begin{gather*}
 \Xi_n=\overline{U_n}\cap \Gamma =\{z\in \Gamma \ t.c. \ G(z,z_0)\geq r_n\}
\end{gather*}
Osserviamo che l'insieme $U_n$ è necessariamente connesso, infatti se avesse una componente connessa non contenente $z_0$, per il principio del massimo e per il fatto che $G(\cdot,z_0)|_{\Delta}=0$, allora $G(\cdot,z_0)$ sarebbe necessariamente minore di $r_n$ su questa componente (assurdo per definizione di $U_n$.\\
Grazie alla proposizione \ref{prop_exreg3}, possiamo scegliere un'esaustione regolare $K_m$ di $R$ in modo che gli insieme 
\begin{gather*}
F_{nm}\equiv\overline{U_n}\cap \partial K_m 
\end{gather*}
siano per ogni $n$ e $m$ sottovarietà regolari di codimensione $1$ con bordo liscio.\\
Fissato un indice $m$, esiste una successione di insiemi $C_p$ aperti relativamente compatti in $R$ con bordo liscio tale che $\overline{C_{p+1}}\subset C_p$ e $K_m=\cap_p C_p$. Questa affermazione può essere dimostrata con argomentazioni simili a quelle riportate nella proposizione \ref{prop_exreg}. La successione $C_p$ è utile per dimostrare che:
\begin{lemma}
Con le notazioni introdotte qui sopra, esiste una funzione $w_{n, m, p}$ tale che:
\begin{enumerate}
 \item $ w_{n,m,p}\in H(U_{n+1}\setminus \overline{K_m})$
 \item $ w_{n,m,p}\geq 0$
 \item $w_{n,m,p}|_{\partial U_{n+1}\setminus C_p}=0$
 \item $w_{n,m,p}|_{\partial U_{n+1}\setminus K_m}\leq 1$
 \item $w_{n,m,p}|_{F_{n+1.m}}=1$
 \item $w_{n,m,p}|_{\Xi_n}\geq \sigma_n$
 \item $D_R(w)<\infty$
\end{enumerate}
dove $\sigma_n$ è un numero strettamente positivo indipendente da $m$ e $p$.
\begin{proof}
 La dimostrazione di questo lemma è molto tecnica. La sua utilità sarà illustrata nella proposizione seguente.\\
Costruiamo le funzioni $w_{n,m,p}$ su $U_{n+1}\setminus K_m$ per esaustione, dimostriamo che è possibile estendere queste funzioni a una funzioni in $\Ro$ (quindi ha senso parlare di $w_{n,m,p}|_{\Xi_n}$) e dimostriamo l'ultima disuguaglianza confrontando queste funzioni con funzioni di Green opportunamente modificate.\\
In tutta la dimostrazione considereremo fissati i valori di $n$, $m$, $p$, quindi $w_{n,m,p}\equiv w$, $U_{n+1}\equiv U$, $K_m\equiv K$ e $C_p\equiv C$.\\
Fissiamo una funzione liscia a supporto compatto $\lambda:R\to[0,1]$ tale che:
\begin{gather*}
 \lambda|_{(\partial K)\cap U}=1, \ \ \ \text{supp}(\lambda)\subset C
\end{gather*}
Per $k>m$, sia $u_k$ la soluzione del problema di Dirichlet
\begin{gather*}
 u_k\in H(U\cap (K_k\setminus \overline{K})), \ \ \ u_k|_{\partial K \cup [(\partial U)\cap K_k\setminus K]}=\lambda|_{\partial K \cup [(\partial U)\cap K_k\setminus K]}, \ \ \ u_k|_{\partial K_k}=0
\end{gather*}
Per il principio del massimo, la successione $u_k$ è una successione crescente e limitata da $1$, quindi grazie al principio di Harnack, $u_k$ ammette come limite una funzione armonica su $U\setminus \overline K$, in particolare:
\begin{gather*}
 \lim_k u_k =w
\end{gather*}
Osserviamo che grazie al principio di Dirichlet \ref{prop_D3}, $D_R(u_k)\leq D_R(\lambda)<\infty$ per ogni $k$, e inoltre se $i>k$
\begin{gather*}
 D_R(u_k)=D_R(u_i)+D_R(u_k-u_i)
\end{gather*}
da cui con un ragionamento simile a quello riportato nella dimostrazione del teorema \ref{teo_1}, otteniamo che la successione $\{u_k\}$ è $D$-Cauchy in $\Ro$, quindi per completezza $u_k$ ammette limite $CD$, e per unicità del limite
\begin{gather*}
 CD-\lim_k u_k =w
\end{gather*}
da cui $D_R(w)<\infty$, quindi (7) è dimostrata.\\
 Per quanto riguarda le altre proprietà, (2), (3), (4) e (5) seguono direttamente dal fatto che tutte le funzioni $u_k$ soddisfano queste proprietà. Per dare senso alla richiesta $w|_{\Xi_n}\geq \sigma_n$, estendiamo la funzione $w$ a una funzione $\tilde w$ definita su tutta $R$ in questo modo:
\begin{gather*}
 \tilde w(x)=\begin{cases}
              w(x) & \text{se } x\in U\setminus K\\
\lambda(x) & \text{se } x\in (U\setminus K)^C
             \end{cases}
\end{gather*}
Osserviamo che $\Xi_n\subset \overline U$, quindi il valore di $\tilde w|_{\Xi_n}$ è indipendente dall'estensione di $w$ che si sceglie.\\
È facile osservare che $\tilde w$ è continua e di Tonelli grazie a un ragionamento simile a quello riportato in \ref{prop_dens2}, e il suo integrale di Dirichlet è finito poiché gli integrali di Dirichlet di $\lambda$ e $w$ lo sono.\\
Questo dimostra che ha senso parlare di $w|_{\Xi_n}$ grazie al fatto che $w$ può essere estesa a tutto $R^*$ e che $w$ è definita su $U_{n+1}\setminus K_m$, quindi in un intorno di $\Xi_n$.\\
Per dimostrare (5), consideriamo la funzione
\begin{gather*}
h(x)=\frac{G(x,z_0)-r_{n+1}}{b-r_{n+1}}
\end{gather*}
dove $b$ è un numero positivo sufficientemente grande da rendere l'insieme
\begin{gather*}
 B\equiv \{ z\in R \ \ r.c. \ \ G(z,z_0)>b\}
\end{gather*}
contenuto nell'insieme $K_1\cap U_{n+1}$. In questo modo $h(x)\leq 1$ sull'insieme $U\setminus K$. Confrontiamo ora le funzioni $w$ e $h$ sull'insieme $U\setminus \overline K$. Osserviamo prima di tutto che queste sono entrambe funzioni armoniche su $U\setminus \overline K$ e continue fino al bordo di questo insieme. Per definizione di $U\equiv U_{n+1}$, la funzione $h$ è negativa su $\partial U$, e come osservato in precedenza $h\leq 1=w$ sull'insieme $\partial K\cap U$. Inoltre entrambe le funzioni sono nulle sull'insieme $\Delta \cap U$, quindi grazie al principio del massimo \ref{teo_max2}, $w\geq h$ su tutto l'insieme di definizione, quindi in particolare anche su $\Xi_n$.\\
Dato che $h|_{\Xi_n}\geq \frac{r_n-r_{n+1}}{b-r_{n+1}}\equiv \sigma_n$, si ha la tesi.
\end{proof}

\end{lemma}
Grazie alla funzione $w_{nm}$ possiamo ottenere una stima su $\rho(\Xi_n)$, infatti vale che:
\begin{prop}\label{prop_stima}
 Secondo le notazioni fino a qui introdotte:
\begin{gather}\label{eq_infty1}
 \rho(\Xi_n)\geq \sigma_n^2 \rho(F_{n+1,m})
\end{gather}
per ogni $m$.
\begin{proof}
 In tutta la dimostrazione, sottointendiamo che le sommatorie che scorrono su nessun indice sono nulle, nel senso che ad esempio:
\begin{gather*}
 \sum_{i=1}^0 a_i \equiv 0
\end{gather*}
Assumiamo inoltre che $\Xi_n$ abbia cardinalità infinita, in caso contrario la tesi è ovvia essendo $\rho(A)=\infty$ per ogni insieme $A$ di cardinalità finita.\\
Sia $k\geq 4$ un intero fissato e $p_1,\cdots,p_k$ punti arbitrari in $\Xi_n$. Ci prefiggiamo di trovare $k$ punti $z_1,\cdots,z_k\in F_{n+1,m}$ tali che
\begin{gather}\label{eq_rho1}
 \sigma^2_n \sum_{i<j}^{1\cdots t} G(z_i,z_j) + \sigma_n \sum_{i=1}^t \sum_{j=t+1}^k G(z_i,p_j) + \sum_{i<j}^{t+1,\cdots,k} G(p_i,p_j)\leq \sum_{i<j}^{1,\cdots,k} G(p_i,p_j)
\end{gather}
per ogni $t=1,\cdots,k$. Una volta dimostrato questo otteniamo in particolare per $t=k$ che:
\begin{gather*}
 \sigma^2_n \sum_{i<j}^{1\cdots k} G(z_i,z_j)\leq \sum_{i<j}^{1,\cdots,k} G(p_i,p_j)
\end{gather*}
e per definizione di $\rho(F_{n+1,m})$ questo implica che:
\begin{gather*}
 \sigma^2_n \binom{k}{2}\rho_k(F_{n+1,m})\leq \sum_{i<j}^{1,\cdots,k} G(p_i,p_j)
\end{gather*}
inoltre data l'arbitrarietà della scelta dei punti $p_1,\cdots,p_k$, si ha che:
\begin{gather*}
  \sigma^2_n \rho_k(F_{n+1,m})\leq  \binom{k}{2}\rho_k(\Xi_n)
\end{gather*}
passando al limite per $k$ che tende a infinito, si ottiene la tesi.\\
Resta da dimostrare la parte tecnica della prova, cioè la relazione \ref{eq_rho1}.\\
A questo scopo utilizzeremo l'induzione sull'indice $1\leq t\leq k$. Supponiamo di aver trovato dei punti $z_1,\cdots,z_{h-1}$ con $1\leq h \leq k-1$ per cui vale \ref{eq_rho1} \footnote{ovviamente l'ipotesi di induzione garantisce che questa relazione valga solo per $1\leq t \leq h-1$, perché per indici più grandi le equazioni non hanno senso non avendo determinato tutti punti $z_i$}. Definiamo
\begin{gather*}
 u_h(z)\equiv \sum_{j=h+1}^k G(z,p_j) + \sigma_n \sum_{i=1}^{h-1} G(z_i,z)
\end{gather*}
Dato che $u_h$ è continua e positiva su $R\setminus \{z_1,\cdots,z_{h-1}\}$, ammette un minimo positivo su $F_{n+1,m}$ assunto in $z_h$. Osserviamo che per ogni $\delta>0$:
\begin{gather*}
u_h(z)-(u_h(z_h)-\delta)>\delta
\end{gather*}
sull'insieme $\partial K_m \cap U_{n+1}$. Questo implica che esiste un intorno di questo insieme su cui questa funzione è strettamente positiva. Denotiamo questo intorno con il simbolo $A$. Per definizione della successione $C_p$, esiste un indice $p$ per il quale l'insieme $(C_p\setminus K_m)\cap U_{n+1}\subset A$. Consideriamo quindi la funzione $w_{n,m,p}$ definita nel lemma precedente. Sappiamo che 
\begin{gather*}
\phi_h^{\delta}(z)\equiv u_h(z)-(u_h(z_h)-\delta)w_{n,m,p}\geq 0
\end{gather*}
sull'insieme $\partial (U_{n+1}\setminus K_m)$, quindi per il principio del massimo questa funzione è positiva su tutto l'insieme $U_{n+1}\setminus K_m$, e in particolare sull'insieme $\Xi_n$, cioé:
\begin{gather*}
 u_h(z)\geq (u_h(z_h)-\delta)\sigma_n
\end{gather*}
per ogni $z\in \Xi_n$. Data l'arbitrarietà di $\delta$, otteniamo che
\begin{gather*}
 u_h(z)\geq u_h(z_h)\sigma_n
\end{gather*}
per ogni $z\in \Xi_n$.\\
Applicando le definizioni otteniamo che:
\begin{gather}\label{eq_rho2}
 \sum_{j=h+1}^k G(p_h,p_j) + \sigma_n \sum_{i=1}^{h-1} G(z_i,p_h) \geq \sigma_n \sum_{j=h+1}^k G(z_h,p_j) + \sigma_n^2 \sum_{i=1}^{h-1} G(z_i,z_h)
\end{gather}
questa relazione per $h=1$ è la dimostrazione di \ref{eq_rho1} per $t=1$. Per gli altri valori di $t$, dall'ipotesi induttiva sappiamo che per $t=h-1$:
\begin{gather*}
 \sum_{i<j}^{1\cdots,k} G(p_i,p_j)\geq \sigma^2_n \sum_{i<j}^{1\cdots h-1} G(z_i,z_j) + \sigma_n \sum_{i=1}^{h-1} \sum_{j=h}^k G(z_i,p_j) + \sum_{i<j}^{h,\cdots,k} G(p_i,p_j)=\\
=\sigma^2_n \sum_{i<j}^{1\cdots h-1} G(z_i,z_j) + \sigma_n \sum_{i=1}^{h-1} \sum_{j=h+1}^k G(z_i,p_j) +\\
+ \sum_{i<j}^{h+1,\cdots,k} G(p_i,p_j)+ \gr{\sigma_n \sum_{i=1}^{h-1} G(z_i,p_h) + \sum_{j=h+1}^k G(p_h,p_j)}
\end{gather*}
applicando la relazione \ref{eq_rho2} nelle parentesi graffe, otteniamo la \ref{eq_rho1} per $t=h$, il che completa il ragionamento di induzione.
\end{proof}
\end{prop}

\subsection{Pricipio dell'energia}\label{subsec_energy}
Anche questo paragrafo, come il precedente, riporta alcuni risultati tecnici utili per stimare il diametro transfinito degli insiemi compatti contenuti in $\Xi$. In particolare introduciamo l'energia e il potenziale di Green relativo a una misura, e troveremo come queste nozioni sono legate alla capacità di un insieme.\\
Alcuni risultati della teoria del potenziale sono stati cortesemente segnalati dal professor. Wolfhard Hansen (University of Bielefeld), che ringraziamo.\\
In tutto il paragrafo, $K$ indicherà una sottovarietà regolare $m-1$ dimensionale di $R$ possibilmente con bordo liscio.
\begin{deph}
 Sull'insieme delle misure di Borel regolari con supporto in $K$ definiamo il prodotto \footnote{il cui risultato potrebbe anche essere $\infty$}:
\begin{gather*}
\psb{\mu}{\nu}\equiv \int_{K}d\mu(x) \ \int_K \ d\nu(y) G(x,y) 
\end{gather*}
la simmetria di questo prodotto segue dalla simmetria di $G(x,y)$ rispetto a $x$ e $y$.\\
Chiamiamo $\psb{\mu}{\nu}$ l'\textbf{energia relativa} tra le misure $\mu$ e $\nu$, e definiamo \textbf{energia} di $\mu$ la quantità $\psb{\mu}{\mu}$.
\end{deph}
\begin{prop}
$\psb \mu \mu=0$ se e solo se $\mu=0$·
\begin{proof}
Supponiamo per assurdo che $\mu\neq 0$. Allora la funzione
\begin{gather*}
 f(x)\equiv \int G(x,y)d\mu(y)
\end{gather*}
 è strettamente positiva per ogni $x$, infatti su $K$, $G(x,\cdot)$ assume un minimo positivo $\lambda>0$, e quindi:
\begin{gather*}
 f(x)\geq \lambda \mu(K)
\end{gather*}
con un ragionamento analogo, si ottiene che anche $$\int f(x)d\mu(x)=\int d\mu(x)\int d\mu(y) G(x.y)>0$$
\end{proof}
\end{prop}

\begin{deph}
Dato un insieme $K$ con le caratteristiche descritte all'inizio del paragrafo, definiamo l'\textbf{energia minima} di $K$ la quantità
\begin{gather*}
 \epsilon(K)\equiv \inf_{\mu\in m_K} \psb{\mu}{\mu}
\end{gather*}
dove $m_K$ è l'insieme delle misure di Borel regolari positive unitarie con supporto su $K$ (cioè $\mu(K)=1$).\\
Per convenzione $\epsilon(\emptyset)=\infty$, e osserviamo immediatamente che per definizione di $\inf$ se $K\subset K'$, $\epsilon(K)\geq\epsilon(K')$. 
\end{deph}
\begin{oss}\label{oss_notinfty}
 Per gli insiemi $K$ con le caratteristiche descritte all'inizio del paragrafo, $\epsilon(K)<\infty$.
\begin{proof}
 Per dimostrare questa affermazione, basta trovare una misura $\mu$ di Borel regolare finita \footnote{il fatto che $\mu(R)\neq 1$ non è importante, infatti è sufficiente riscalare la misura} per la quale $\psb\mu\mu<\infty$.\\
Sia $(U,\phi)$ un insieme coordinato che interseca $K$ per il quale $$\phi(K)=(x_1,\cdots,x_{m-1},0)$$ e sia $B$ una bolla di raggio $\bar r$ in $\R^{m-1}$ contenuta in $\phi(K)$. Consideriamo la misura di superficie data dalla metrica riemanniana su questo insieme, misura che ha la forma
\begin{gather*}
 dV=\sqrt{\abs g }d\lambda
\end{gather*}
dove $\lambda$ è la misura di Legesue standard. Grazie alle proprietà del nucleo di Green, sappiamo che sull'insieme $B$ $G$ in coordinate locali si può scrivere come
\begin{gather*}
 G(x,y)=f(x,y)+C(m)\gamma(x,y)
\end{gather*}
dove $f(x,y)$ è una funzione continua nelle due variabili $x$ e $y$, $C(m)$ una costante che dipende solo dalla dimensione $m$ e $\gamma(x,y)$ è una funzione che dipende dalla dimensione $m$ della varietà $R$, in particolare:
\begin{gather*}
 \gamma(x,y)=\begin{cases}
                  -\log(d(x,y)) & m=2\\
		  d(x,y)^{m-2} & m\geq 3
                 \end{cases}
\end{gather*}
Per limitatezza di $f$ sull'insieme $B\times B$, è ovvio che:
\begin{gather*}
 \int_B f(x,y)dV(x)dV(y)=\int_B f(x,y)\sqrt{\abs g }d\lambda(x)\sqrt{\abs g }d\lambda(y) <\infty
\end{gather*}
Consideriamo quindi l'integrale 
\begin{gather*}
 \int_B \gamma(x,y) dV(x)dV(y) = \int \gamma(x,y)\sqrt{\abs g }d\lambda(x)\sqrt{\abs g }d\lambda(y) \leq M^2 \int_B \gamma(x,y)d\lambda(x)d\lambda(y)
\end{gather*}
dove $M$ è un limite superiore per la funzione $\sqrt{\abs g}$ sull'insieme relativamente compatto $B$. La funzione
\begin{gather*}
 \bar \gamma(x)\equiv \int_B \gamma(x,y)d\lambda(y)
\end{gather*}
è una funzione limitata in $x$, infatti scelto $x\in B$, detto $2B(x)$ la bolla di raggio $2\bar r$ centrata in $x$, si ha che:
\begin{gather*}
 \bar \gamma(x)=\int_B \gamma(x,y)d\lambda(y)\leq \int_{2B(x)} \gamma(x,y)d\lambda (y)
\end{gather*}
per $m=2$, si ha che:
\begin{gather*}
  \int_{2B(x)} -\log(\abs{x-y})d\lambda(y)=\int_{-\bar r}^{\bar r} -\log(r) dr =2\bar r(\log(\bar r)-1)<\infty
\end{gather*}
mentre per $m\geq 3$, passando alle coordinate polari centrate in $x$, si ha che
\begin{gather*}
 \int_{2B(x)} \frac{1}{\abs{x-y}^{m-2}}d\lambda(y)=\int_{S_{m-2}}d\theta \int_{0}^{\bar r} dr \ r^{m-2}\frac{1}{r^{m-2}}= \omega_{m-2} \bar r
\end{gather*}
dove $\omega_{m}$ è la superficie della sfera $m$ dimensionale. Essendo la funzione $\bar \gamma$ limitata, abbiamo che:
\begin{gather*}
 \int_{B\times B} \gamma(x,y)d\lambda(x)d\lambda(y) =\int_B \bar \gamma(x) d\lambda(x)\leq \lambda(B) \sup_B (\bar \gamma(x))<\infty
\end{gather*}
\end{proof}
\end{oss}
\begin{deph}\label{deph_greenpot}
Per una misura $\mu\in m_K$, definiamo il relativo potenziale di Green la funzione:
\begin{gather*}
 G_{\mu}(x)\equiv \int_{K}G(x,y)d\mu(y)
\end{gather*}
\end{deph}
Per prima cosa esploriamo le proprietà di $G_{\mu}$ per una qualsiasi misura di Borel positiva regolare a supporto compatto $S\in R$. 
\begin{prop}
 $G_{\mu}\in H(R\setminus S)$, $G_\mu$ è positiva, semicontinua inferiormente, e può essere estesa a tutto $R^*$ con $G_{\mu}|_{\Delta}=0$. Inoltre $G_{\mu}$ è superarmonica su $R$.
\begin{proof}
 La positività di $G_{\mu}$ è un'ovvia conseguenza della positività della misura e del nucleo di Green. Per dimostrare che $G_{\mu}\in H(R\setminus S)$, sfruttiamo il teorema di derivazione sotto al segno di integrale \ref{diid}. Sia $x\in R\setminus S$. Allora esiste un intorno compatto $V$ di $x$ disgiunto da $S$, e per le proprietà di $G(\cdot,\cdot)$, sappiamo che
\begin{gather*}
 G(x,y)\leq \sup_{z\in \partial V} G(x,z)<\infty
\end{gather*}
per ogni $y\in S$. Grazie al fatto che $G\in C^\infty (R\times R \setminus D)$, sappiamo che tutte le derivate $\partial _i G(x,y)$ e $\partial _i \partial _j G(x,y)$ sono uniformemente limitate se $(x,y)\in V\times S$, quindi possiamo applicare il teorema \ref{diid} e ottenere che:
\begin{gather*}
 \Delta G_{\nu}(x)=g^{ij}(x)\partial _i \partial _j \int_{S} G(x,y) d\mu(y)=\\
=g^{ij}(x)\partial _i  \int_{S} \partial _jG(x,y) d\mu(y)= \int_{S} g^{ij}(x)\partial _i \partial _j G(x,y) d\mu(y)=0
\end{gather*}
Con un raginamento del tutto analogo possiamo dimostrare che la funzione $G_{\mu}$ è continua sull'insieme $R\setminus S$.\\
Il fatto che $G_{\mu}$ è semicontinua inferiormente su $R$ segue dalla considerazione che grazie al teorema di convergenza monotona:
\begin{gather*}
 G_{\mu}(x)=\int_S G(x,y)d\mu(y)=\lim_n \int_S (G(x,y)\curlywedge n) d\mu(y) \equiv \lim_n G_{\mu}^n (x)
\end{gather*}
La successione $G_{\mu}^n$ è una successione crescente di funzioni continue \footnote{la continuità segue dalla continuità della funzione $G(x,y)\curlywedge n$ e dal teorema di convergenza dominata}, quindi $G_{\mu}(x)$ è necessariamente una funzione semicontinua inferiormente.\\
Sempre grazie al teorema di convergenza monotona, sappiamo che:
\begin{gather*}
 G_{\mu}(x)\curlywedge c=\int_S (G(x,y)\curlywedge c) d\mu(y)
\end{gather*}
Queste funzioni appartengono tutte a $\Ro$, infatti sono continue, di Tonelli, e:
\begin{gather*}
 D_R(G_{\mu}(x)\curlywedge c)=\int_R dx \abs{\nabla \int_S (G(x,y)\curlywedge c) d\mu(y)}^2\leq\\
\leq \int_R dx  \int_S \abs{\nabla (G(x,y)\curlywedge c)}^2 d\mu(y)= \int_S d\mu(y) \int_R dx  \abs{\nabla (G(x,y)\curlywedge c)}^2\leq\\
 \leq \int_S d\mu(y) c=c
\end{gather*}
dove abbiamo sfruttato il fatto che $\mu(R)=1$ \footnote{quindi $\abs{\int_S f(y)d\mu(y)}^2\leq \int_S\abs{f(x)}^2 d\mu(y)$} e il teorema di Fubuni per lo scambio di integrali quando l'integrando è positivo.\\
Sempre per il teorema di convergenza monotona, detto $G_n(x,y)$ i nuclei di Green rispetto a un'esaustione regolare di $R$, sappiamo che:
\begin{gather*}
 G_{\mu}(x)\curlywedge c=\lim_n \int_S (G_n(x,y)\curlywedge c) d\mu(y)
\end{gather*}
Poiché tutte le funzioni $G_n(x,y)\curlywedge c$ appartengono a $\Rod$ e
\begin{gather}\label{eq_1}
 D_R(G_n(x,y)\curlywedge c)\leq c
\end{gather}
grazie al teorema \ref{prop_r2}, $G_{\mu}(x)\curlywedge c\in \Rod$, quindi $G_{\mu}|_{\Delta}=0$.\\
La superarmonicità di $G_{\mu}$ segue dalle proposizioni \ref{prop_sup1} e \ref{prop_sup2}. Infatti sappiamo che:
\begin{gather*}
 G_{\mu}(x)=\int_{K} G(x,y)d\mu(y)=\lim_n \int_K (G(x,y)\curlywedge n) d\mu(y)
\end{gather*}
Le funzioni $\int_K (G(x,y)\curlywedge n) d\mu(y)$ sono superarmoniche grazie alla proposizione \ref{prop_sub1}, e anche il loro limite è superarmonico grazie a \ref{prop_sup2}.
\end{proof}
\end{prop}

\begin{prop}\label{prop_endeph}
 Se $\epsilon(K)<\infty$, esiste una misura $\nu\in m_K$ tale che:
\begin{gather*}
 \epsilon(K)=\psb \nu \nu
\end{gather*}
Cioè esiste una misura che realizza il minimo dell'energia.
\begin{proof}
Dalla definizione di $\inf$, esiste una successione $\nu_n\in m_K$ tale che:
\begin{gather*}
\lim_n \psb {\nu_n}{\nu_n}=\epsilon(K)
\end{gather*}
Dalla teoria della misura e degli spazi di Banach, sappiamo che esiste una sottosuccessione di $\{\nu_n\}$ (che per comodità continueremo a indicare con lo stesso indice) tale che $\nu_n$ converge debolmente a $\nu\in m_K$, cioè per ogni funzione $f\in C(K)$:
\begin{gather*}
 \lim_n \int f d\nu_n =\int fd\nu
\end{gather*}
Ora consideriamo le funzioni
\begin{gather*}
 \phi_n^c (x)\equiv \int_K (G(x,y)\curlywedge c) d\nu_n(y) \ \ \ \ \ \  \phi^c (x)\equiv \int_K (G(x,y)\curlywedge c) d\nu(y)
\end{gather*}
Ricordiamo che per convergenza monotona:
\begin{gather*}
\int_K G(x,y)d\nu(y)=\lim_{c\to \infty} \phi^c(x)
\end{gather*}
Osserviamo che la successione $\{\phi_n^c\}$ converge uniformemente su $K$ a $\phi^c$, infatti grazie alla definizione di $\nu$, c'è convergenza puntuale. Inoltre:
\begin{gather*}
 \abs{\phi_n^c(x)}=\int_K (G(x,y)\curlywedge c) d\nu_n(y)\leq c \nu_n(K) \leq c M \\
\abs{\phi_n^c(x_1)-\phi_n^c(x_2)}\leq \int_K \abs{(G(x_1,y)\curlywedge c)-G(x_2,y)\curlywedge c} d\nu_n(y)
\end{gather*}
Data la continuità uniforme della funzione $G(\cdot,\cdot)\curlywedge c$ sull'insieme $K\times K$, si ha che per ogni $\epsilon>0$, esiste $\delta$ per cui
$$
d(x_1,x_2)<\delta \ \ \Rightarrow \ \ \abs{G(x_1,y)\curlywedge c - G(x_2,y)\curlywedge c }<\epsilon
$$
quindi se $d(x_1,x_2)<\delta$:
\begin{gather*}
\abs{\phi_n^c(x_1)-\phi_n^c(x_2)}\leq \epsilon\nu_n(K)\leq \epsilon M
\end{gather*}
Queste osservazioni dimostrano l'uniforme limitatezza e l'equicontinuità della successione $\phi_n^c$ sull'insieme $K$, quindi grazie al teorema di Ascoli-Arzelà (vedi ad esempio appendice A5 pag 394 di \cite{4}) per ogni sottosuccessione di $\{\phi_n^c\}$, esiste una sua sottosottosuccessione convergente uniformemente su $K$. Dato che $\{\phi_n^c\}$ converge puntualmente a $\phi^c$, allora la convergenza è uniforme su $K$.\\
Questo implica che:
\begin{gather}\label{eq_mis}
\abs{ \int_K \phi_n^c d\nu_n -\int_K \phi^c d\nu}\leq \abs{ \int_K (\phi_n^c-\phi^c) d\nu_n}+\abs{\int_K \phi^c d\nu_n -\int_K \phi^c d\nu} \leq\\
\leq \norm{\phi_n^c-\phi^c}_{\infty,K} \nu_n(K) + \abs{\int_K \phi^c d\nu_n -\int_K \phi^c d\nu} \to 0 \notag
\end{gather}
Da queste considerazioni otteniamo che:
\begin{gather*}
 \psb{\nu}{\nu}=\lim_{c\to \infty} \int_Kd\nu(x)\int_Kd\nu(y) (G(x,y)\curlywedge c)=\lim_{c\to \infty}\lim_{n\to \infty}\int_K \phi_n^c d\nu_n=\\
\lim_{c\to \infty}\liminf_{n\to \infty}\int_K \phi_n^c d\nu_n\leq\liminf_{n\to \infty}\lim_{c\to \infty}\int_K \phi_n^c d\nu_n=\liminf_{n\to \infty} \ps{\nu_n}{\nu_n}
\end{gather*}
dove il penultimo passaggio, lo scambio tra limite e liminf, è giustificato nel lemma seguente (lemma \ref{lemma_asdf}).\\
Questa disuguaglianza ci permette di concludere che:
\begin{gather*}
 \epsilon(K)\leq \psb \nu \nu \leq \liminf_n \psb{\nu_n}{\nu_n}=\epsilon(K)
\end{gather*}
da cui segue la tesi.
\end{proof}
\end{prop}

\begin{lemma}\label{lemma_asdf}
 Sia $a_{n,m}$ una successione con due indici crescente in $m$ per ogni $n$ fissato. Allora:
\begin{gather*}
 \lim_{m\to \infty}\liminf_{n\to \infty} a_{n,m} \leq \liminf_{n\to \infty} \lim_{m\to \infty} a_{n,m}
\end{gather*}
\begin{proof}
Sia per definizione
\begin{gather*}
L\equiv  \lim_{m\to \infty}\liminf_{n\to \infty} a_{n,m}\ \ \ \ \ \ \ b_m\equiv  \liminf_{n\to \infty} a_{n,m}
\end{gather*}
Dalla definizione di limite, sappiamo che per ogni $\epsilon>0$, esiste $M$ tale che per ogni $m\geq M$, $b_m>L-\epsilon$, quindi dalla definizione di $\liminf$ sappiamo che fissato $\bar m$, esiste $N(\bar m)$ tale che per ogni $n\geq N(\bar m)$:
\begin{gather*}
 a_{n,\bar m}>L-\epsilon
\end{gather*}
Data la monotonia di $a_{n,m}$ rispetto all'indice $m$, sappiamo che per ogni $m\geq \bar m$ e per ogni $n\geq N(\bar m)$, $a_{n,m}>L-\epsilon$, quindi passando al limite in $m$:
\begin{gather*}
 \lim_m a_{n,m} \geq L-\epsilon
\end{gather*}
quindi passando al $\liminf$:
\begin{gather*}
 \liminf_n \lim_m a_{n,m}\geq L-\epsilon
\end{gather*}
la tesi si ottiene per arbitrarietà di $\epsilon$.
\end{proof}

\end{lemma}

Ora introduciamo una particolare misura positiva di Borel regolare sull'insieme $K$, ragionando con le funzioni armoniche.\\
Sia $K$ una sottovarietà regolare compatta di codimensione $1$ (possibilmente con bordo liscio). $R\setminus K$ avrà $n$ componenti connesse $R_1,\cdots,R_n$ \footnote{ad esempio, se $K$ è una sfera, $R\setminus K$ avrà due componenti connesse, se invece $K$ è una parte di piano, $R\setminus K$ è connesso}, inoltre ovviamente $K=\cup \partial R_i$. Su ogni insieme $R_i$ possiamo risolvere un problema di Dirichlet, nel senso che data una funzione $f$ continua su $\partial R_i$, esiste una funzione armonica in $R_i$ che indicheremo $H^i(f)$ che ristretta al bordo di $R_i$ uguagli $f$.\\
Questa caratterizzazione però è sufficiente solo se $R_i$ è relativamente compatto. In caso contrario indichiamo con $H^i(f)$ la funzione armonica ottenuta per esaustione di $R_i$. In dettaglio, sia $C_n$ un'esaustione regolare di $R$, e sia $H^i_n(f)$ la funzione
\begin{gather*}
H^i_n(f)\in H(R_i\cap C_n), \ \ H^i_n(f)|_{\partial R_i\cap C_n}=f, \ \ H^i_n(f)|_{\partial C_n\cap R_i} =0
\end{gather*}
grazie al principio del massimo si ottiene che $H^i_n(f)$ è una successione crescente in $n$ di funzioni armoniche, e poiché sempre per il principio del massimo
\begin{gather*}
 \norm{H^i_n(f)}_{\infty}\leq \norm{f}_{\infty}
\end{gather*}
grazie al principio di Harnack $H^i_n(f)$ converge a una funzione $H^i(f)$ armonica su $R_i$ e continua su $\overline{R_i}$.\\
Ora se scegliamo un punto qualsiasi $z_0\in R^i$, per quanto appena visto possiamo definire un funzionale lineare continuo positivo \footnote{per definizione, $\phi$ è positivo se $\phi(f)\geq 0$ per ogni $f\geq 0$} sullo spazio delle funzioni continue in $K$ in questo modo:
\begin{gather*}
 \phi_{z_0}(f)\equiv H^i(f)(z_0)
\end{gather*}
dalla teoria \footnote{vedi teorema di rappresentazione di Riestz, ad esempio su \cite{12}, teorema 2.14 pag 40} sappiamo che esiste unica una misura di Borel positiva regolare $\xi_{z_0}$ con supporto su $\partial R^i$ \footnote{quindi con supporto in $K$ a patto di estenderla a nulla su $K\setminus \partial R^i$} tale che:
\begin{gather*}
 \phi_{z_0}(f)\equiv H^i(f)(z_0)=\int_{K} f d\xi_{z_0}
\end{gather*}
per ogni funzione $f$ continua su $K$.
\begin{deph}
 La misura $\xi_{z_0}$ appena caratterizzata è detta la \textbf{misura armonica} dell'insieme $K$ rispetto a $z_0$.
\end{deph}
Nel seguente lemma ci occupiamo di un aspetto tecnico legato a questa misura, in particolare a cosa succede su un insieme $F\subset K$ di misura nulla rispetto a una di queste $\xi_{z_0}$.
\begin{lemma}\label{prop_9f}
 Se $\xi_{z_0}(F)=0$, esiste una funzione armonica positiva $h:R^i\to [0,\infty)$ tale che per ogni $x_0\in F$:
\begin{gather*}
 \lim_{x\to x_0} h(x)=\infty
\end{gather*}
Come corollario di questo lemma, possiamo dimostrare che la proprietà $\xi_{z_0}(F)=0$ NON dipende dalla scelta di $z_0$ (se $z_0$ viene scelto tra gli elementi dello stesso insieme $R^i$).
\begin{proof}
 Dato che $\xi_{z_0}\equiv \xi_0$ è una misura regolare, esiste una successione di insiemi $\{F_n\}$ aperti nella topologia di $K$ tale che $$F\subset F_{n}\subset \overline{F_n}\subset F_{n+1} \ \ \ \ \bigcap_n F_n =F \ \ \ \ \xi_0 (F_n)\leq 2^{-n}$$
Costruiamo le funzioni $q_n:K\to \R$ con queste proprietà:
\begin{enumerate}
 \item $q_n\in C(K,[0,1])$ per ogni $n$
\item $supp(q_n)\subset F_n$, e $q_n|_{F_{n+1}}=1$
\end{enumerate}
Sia $s_n\equiv \sum_{i=1}^n q_i$. Visto quanto abbiamo appena osservato, possiamo costruire la successione di funzioni armoniche positive
\begin{gather*}
 h_n\equiv H^i(s_n)
\end{gather*}
Per definizione di $s_n$, sappiamo che $h_n|_{F}=s_n|_{F} = n$, inoltre la successione di funzioni armoniche positive crescenti $h_n$ è tale che:
\begin{gather*}
 h_n(z_0)=H^i(s_n)(z_0) =\int_K s_n d\xi_{z_0}=\sum_{j=1}^n \int_K q_j d\xi_{z_0}\leq \sum_{i=j}^n \int_K \chi_{F_j}d\xi_{z_0}\leq \sum_{j=1}^n 2^{-j}\leq 1
\end{gather*}
quindi per il principio di Harnack \footnote{vedi \ref{prop_harnackpri}} la successione $h_n$ converge a una funzione armonica $h$ su $R^i$ semicontinua inferiormente su $\overline{R^i}$ \footnote{dato che la successione $h_n$ è una successione crescente di funzioni continue su $\overline{R^i}$} tale che $h(z)=\infty$ per ogni $z\in F$.\\
Consideriamo ora un punto $z_1\in R^i$, e sia $\xi_{z_1}$ la relativa misura armonica. Per la successione $h_n$ vale che;
\begin{gather*}
 h_n(z_1)=\int_{K}s_n d\xi_{z_1}
\end{gather*}
quindi passando al limite su $n$ otteniamo che:
\begin{gather*}
 \lim_n \int_{K}s_n d\xi_{z_1} = h(z_1)<\infty
\end{gather*}
se per assurdo $\xi_{z_1}(F)>0$, allora necessariamente
\begin{gather*}
 \int_K s_n d\xi_{z_1}\geq n\xi_{z_1}(F)\to \infty
\end{gather*}
che contraddice quanto appena scoperto, quindi se $\xi_{z_0}(F)=0$, allora per ogni $z_1\in R^i$ si ha che $\xi_{z_1}(F)=0$. Scambiando i ruoli di $z_0$ e $z_1$ si ottiene che per un insieme $F\subset K$, il fatto di avere misura armonica $\xi$ nulla è indipendente dalla scelta del punto $z_0\in R^i$.
\end{proof}
\end{lemma}

La misura armonica è legata all'energia dal fatto che
\begin{prop}\label{prop_a}
 Sia $z_0\not \in K$ e $\xi_{z_0}=\xi$ la misura armonica di $K$ relativa a $z_0$. Allora se un insieme di Borel $F\subset K$ ha misura armonica non nulla ($\xi(F)\neq 0$), allora ha energia finita.
\begin{proof}
 Consideriamo la misura di Borel $\xi_1\equiv \xi|_F$. Il potenziale di Green di questa misura è caratterizzato dal fatto che per ogni $x\not \in K$:
\begin{gather*}
 G_{\xi_1}(x)=\int_K G(x,y)d\xi_1(y)\leq \int_K G(x,y)d\xi(y)=H(G(x,\cdot)|_{z_0}\leq G(x,z_0)
\end{gather*}
infatti per superarmonicità, la funzione $G(x,\cdot)$ è maggiore della funzione armonica che assume valore $G(x,\cdot)$ sull'insieme $K$. Data la continutà di $G(\cdot,z_0)$ su $R\setminus \{z_0\}$ e data la semicontinuità inferiore di $G_{\xi_1}$, la disuguaglianza continua a valere anche se $x\in K$, dimostrando che:
\begin{gather*}
 G_{\xi_1}(x)\leq G(z_0,x)\leq M
\end{gather*}
dove $M$ è un maggiorante per la funzione continua $G(z_0,\cdot)|_K$, quindi:
\begin{gather*}
 \psb {\xi_1}{\xi_1}\leq M\xi_1(K)<\infty
\end{gather*}
\end{proof}

\end{prop}

Con la seguente proposizione ci occupiamo di alcune proprietà fini del potenziale $G_\nu$, dove $\nu$ è la misura che minimizza l'energia su $K$. La proposizione e la relativa dimostrazione sono tratte dal paragrafo 9G di \cite{5} e dal teorema III.12 di \cite{36}.
\begin{prop}\label{prop_alpha}
 Sia $K$ un'insieme con le caratteristiche descritte all'inizio del paragrafo, sia $\xi$ la misura armonica di $K$ rispetto a un punto qualsiasi $z_0\not \in K$ e sia $\nu\in m_K$ una misura tale che
\begin{gather*}
 E\equiv \epsilon(K)=\psb \nu \nu
\end{gather*}
Allora:
\begin{enumerate}
 \item $G_\nu(x)\leq E$ sul supporto $S_\nu$ della misura $\nu$
 \item $G_\nu(x)\geq E$ a meno di un insieme di misura armonica $\xi$ nulla
 \end{enumerate}
\begin{proof}
Sia $A_\delta=\{x\in K\ t.c. \ G_\nu(x)>E+2\delta\}$, dove $\delta >0$. Data la semicontinuità inferiore di $G_\nu$, $A_\delta$ sono insiemi aperti per ogni $\delta$, quindi Borel misurabili. Supponiamo per assurdo che $\nu(A_\delta)>0$ per qualche $\delta>0$. Per la regolarità della misura $\nu$, esiste un insieme $C\subset A_\delta$ compatto tale che $\nu(C)\neq 0$. Dato che
\begin{gather*}
 E=\epsilon(K)=\psb \nu \nu =\int_{S_\nu} G_\nu (x) d\nu(x)
\end{gather*}
esiste un punto $x_0\in S_\nu$ tale $G_\nu (x_0)\leq E+\delta$, inoltre l'insieme $$B=\{x\in K \ t.c. \ G_\nu(x)\leq E+\delta\}$$ è chiuso e $\nu(B_\delta)\neq 0$. Consideriamo la misura di Borel
\begin{gather*}
 \sigma(\cdot)= -\frac{\nu(\cdot\cap C)}{\nu(C)} +\frac{\nu(\cdot\cap B)}{\nu(B)}
\end{gather*}
Dato che
\begin{gather*}
0\leq \psb{\nu(\cdot\cap C)}{\nu(\cdot \cap C)}=\int_K G_{\nu(\cdot \cap C)}d\nu{\cdot \cap C} \leq \int_K G_\nu d\nu(\cdot \cap C)=\\
=\int_K G_{\nu(\cdot \cap C)} d\nu \leq\int_K G_\nu d\nu =\psb \nu \nu <\infty
\end{gather*}
e con un ragionamento analogo anche $\psb{\nu(\cdot \cap B)}{\nu(\cdot\cap B)}<\infty$, e dato che gli insiemi $B$ e $C$ sono compatti disgiunti, allora
\begin{gather*}
 0\leq \psb{\nu(\cdot \cap B)}{\nu(\cdot \cap C)}=\int_{K}G_{\nu(\cdot \cap B)}d\nu(\cdot\cap C) \leq M\nu(C)
\end{gather*}
dove $M$ è un maggiorante per la funzione $G_{\nu(\cdot \cap B)}$ continua sul compatto $C$,. Inoltre notiamo che $\psb \sigma \sigma\neq \infty$, infatti:
\begin{gather*}
 \psb \sigma \sigma=\psb{-\frac{\nu(\cdot\cap C)}{\nu(C)} +\frac{\nu(\cdot\cap B)}{\nu(B)}}{-\frac{\nu(\cdot\cap C)}{\nu(C)} +\frac{\nu(\cdot\cap B)}{\nu(B)}}=\\
=\frac{1}{\nu(C)^2} \psb{\nu(\cdot\cap C)}{\nu(\cdot \cap C)} -\frac{2}{\nu(C)\nu(B)}\psb{\nu(\cdot \cap C)}{\nu(\cdot \cap B)} +\\
+ \frac{1}{\nu(B)^2}\psb{\nu(\cdot\cap B)}{\nu(\cdot \cap B)}
\end{gather*}
Consideriamo ora per $0\leq \eta <\nu(C)$ la misura $\nu_\eta=\nu + \eta \sigma$. \`E facile verificare che $\nu_\eta$ è una misura di Borel regolare positiva e $\nu_\eta(K)=1$, quindi $\nu_\eta \in m_K$. D'altronde osserviamo che:
\begin{gather*}
 \psb{\nu_\eta}{\nu_\eta}-\psb \nu \nu=\psb{\nu+\eta \sigma}{\nu+\eta \sigma}-\psb \nu \nu=2\eta \psb \nu \sigma + \eta^2 \psb \sigma \sigma
\end{gather*}
Inoltre:
\begin{gather*}
 \psb \nu \sigma= -\frac 1 {\nu(C)}\int_C G_\nu d\nu +\frac 1 {\nu(B)} \int_B G_\nu d\nu\leq -E-2\delta +E+\delta <-\delta <0
\end{gather*}
Quindi per $\eta$ sufficientemente piccolo, $\psb{\nu_\eta}{\nu_\eta}<\psb \nu \nu$, contraddicendo la definizione di $\nu$.\\
La seconda parte della proposizione si dimostra in maniera analoga alla prima. Supponiamo per assurdo che esista $\delta>0$ tale che l'insieme $$A_\delta =\{x\in K \ t.c. \ G_\nu(x)\leq E-2\delta\}$$ abbia misura armonica positiva. Grazie alla proposizione \ref{prop_a}, l'energia di $A$ è finita, quindi esiste una misura positiva unitaria $\sigma_1$ con supporto in $A$ tale che $\psb {\sigma_1}{ \sigma _1}<\infty$. Sia $B=\{x\in K t.c. \ G_\nu(x)> E-\delta\}$. Dato che
\begin{gather*}
 \int_K G_\nu(x)d\nu(x)=E
\end{gather*}
$\nu(B)\neq 0$. Se definiamo $\sigma(\cdot)=\sigma_1(\cdot)-\frac{\nu(\cdot\cap B)}{\nu(B)}$, analogamente a quanto dimostrato sopra, si dimostra che $\sigma(K)=\sigma(R)=0$ e che per ogni $0\leq\eta<\nu(B)$, la misura $\nu+\eta \sigma\in m_K$.\\
La contraddizione segue dal fatto che:
\begin{gather*}
 \psb{\nu_\eta}{\nu_\eta}-\psb \nu \nu= 2\eta \psb \sigma \nu + \eta^2 \psb \sigma \sigma\\
\psb \sigma \nu =\int_K G_\nu(x) d\sigma = \\
=\int_{A_\delta} G_\nu(x) d\sigma_1 -\frac{1}{\nu(B)}\int_B G_\nu(x)d\nu \leq E-2\delta -E+\delta \leq -\delta <0
\end{gather*}
\end{proof}
\end{prop}
Nella seguente proposizione dimostriamo una proprietà di continuità del potenziale di Green $G_\mu$ rispetto a una qualsiasi $\mu\in m_K$.
\begin{prop}\label{prop_cont}
 Sia $\mu$ una misura di Borel positiva finita e $G_\mu$ il relativo potenziale di Green. Se $G_\mu$ ristretto al supporto $S_\mu$ della misura $\mu$ è continuo, allora $G_\mu$ è una funzione continua su tutta la varietà $R$.
\begin{proof}
 Dato che $G_\mu\in H(R\setminus S_\mu)$, è sufficiente dimostrare che $G_\mu$ è continua sull'insieme $S_\mu$, quindi che per ogni $x_0\in S_\mu$:
\begin{gather*}
 \lim_{x\to x_0} G_\mu (x)=G_\mu (x_0)
\end{gather*}
Sappiamo dalle proposizioni precedenti che $G_\mu$ è semicontinua inferiormente, quindi basta dimostrare che:
\begin{gather}\label{eq_boh}
 \limsup_{x\to x_0} G_\mu(x)\leq G_\mu(x_0)
\end{gather}
inoltre possiamo assumere che $G_\mu(x_0)<\infty$, quindi $\mu(\{x_0\})=0$.\\
Consideriamo un intorno coordinato di $U$ di $x_0$. Su questo intorno possiamo scrivere:
\begin{gather*}
 G(x,y)=f(x,y)+C(m)h(x,y)
\end{gather*}
dove $f(x,y)$ è una funzione continua in entrambe le sue variabili, $C(m)$ una costante che dipende solo da $m=dim(R)$, mentre
\begin{gather*}
 h(x,y)=\begin{cases}
         -\log(d(x,y)) \ & \ se \ m=2\\
d(x,y)^{-m+2} \ & \ se \ m\geq 3
        \end{cases}
\end{gather*}
Decomponendo la misura $\mu$ in $\mu_0(\cdot)\equiv\mu(\cdot\cap U)$ e $\mu'\equiv \mu-\mu_0$ possiamo scrivere il potenziale di Green come:
\begin{gather*}
 G_\mu(x)=G_{\mu_0}(x)+G_{\mu'}(x)
\end{gather*}
e dato che $G_{\mu'}(x)$ è una funzione continua in $x_0$, basta dimostrare \ref{eq_boh} per il solo potenziale $G_{\mu_0}$. Inoltre
\begin{gather*}
 G_{\mu_0}(x)=\int_K f(x,y)d\mu_0(y)+ C(m)\int_K h(x,y)d\mu_0(y)
\end{gather*}
e dato che la misura $\mu_0$ è finita, il primo integrale è una funzione continua della variabile $x$, quindi ancora basta dimostrare \ref{eq_boh} per la funzione
$$H(x)\equiv \int_K h(x,y)d\mu_0(y)$$
Poiché $\mu(\{x_0\})=0$, e dato che $\mu$ è una misura regolare, per ogni $\epsilon>0$, esiste un aperto $U_\epsilon\subset U$ contenente $x_0$ tale che $\mu(U_\epsilon)=\mu_0(U_\epsilon)<\epsilon$, quindi vale che:
\begin{gather*}
 H(x)=\int_{U_\epsilon} h(x,y)d\mu_0(y)+\int_{U\setminus U_\epsilon} h(x,y)d\mu_0(y)
\end{gather*}
Per ogni $x\in U$, definiamo $\pi(x)\in S_\mu$ un punto tale che $d(x,\pi(x))=\min_{y\in S_\mu} d(x,y)$. Osserviamo che in generale il punto $\pi(x)$ non è unico, ma comunque $$\lim_{x\to x_0}\pi(x)=x_0$$ se $x_0\in S_\mu$, inoltre vale che per ogni $y\in S_\mu$:
\begin{gather*}
 d(y,\pi(x))\leq d(x,y)+d(x,\pi(x))\leq 2d(x,y)
\end{gather*}
quindi:
\begin{gather*}
 h(x,y)\leq \begin{cases}
             \log(2)+h(\pi(x),y) & se \ m=2\\
	      2^{m-2} h(\pi(x),y) & se \ m\geq 3
            \end{cases}
\end{gather*}
Dividiamo la dimostrazione in due casi: se $n=2$:
\begin{gather*}
 H(x)\leq \int_{U_\epsilon} h(\pi(x),y)d\mu_\epsilon +\epsilon \log(2) + \int_{U\setminus U_\epsilon} h(x,y)d\mu_0(y)=\\
=H(\pi(x)) + \epsilon \log(2) + \int_{U\setminus U_\epsilon} (h(x,y)-h(\pi(x),y)d\mu_0(y)
\end{gather*}
Applicando il $\limsup_{x\to x_0}$ a entrambi i membri e tenedo conto dell'ipotesi di continuità di $G_\mu$ (quindi di $H$) ristretta all'insieme $S_\mu$ otteniamo:
\begin{gather*}
 \limsup_{x\to x_0}H(x)\leq H(x_0)+\epsilon \log(2) +\limsup_{x\to x_0} \int_{U\setminus U_\epsilon} (h(x,y)-h(\pi(x),y)d\mu_0(y)
\end{gather*}
Per convergenza dominata, l'ultimo addendo è nullo, e data l'arbitrarietà di $\epsilon$, otteniamo la tesi.\\
In maniera analoga, se $m\geq 3$:
\begin{gather*}
H(x)\leq \int_{U_\epsilon} h(\pi(x),y)d\mu_\epsilon +\alpha_m \int_{U_\epsilon} h(\pi(x),y)d\mu(y) + \int_{U\setminus U_\epsilon} h(x,y)d\mu_0(y)
\end{gather*}
dove $\alpha_m=2^{m-2}-1$. Applicando il $\limsup_{x\to x_0}$ a entrambi i membri e tenendo conto della continuità della funzione $H$ ristretta all'insieme $S_\mu$, otteniamo:
\begin{gather*}
 \limsup_{x\to x_0} H(x) \leq H(x_0) + \alpha_m \int_{U_\epsilon} h(x_0,y)d\mu(y) +\\
+ \limsup_{x\to x_0}\int_{U\setminus U_\epsilon} (h(x,y)-h(\pi(x),y)d\mu_0(y)
\end{gather*}
Come nel caso bidimensionale, l'ultimo limite è nullo per convergenza dominata \footnote{infatti dato che $x\to x_0$, $x$ e $\pi(x)$ appartengono definitivamente a un intorno compatto di $x_0$ contenuto in $U_\epsilon$, quindi se $y\in U\setminus U_\epsilon$, sia $h(x,y)$ che $h(\pi(x),y)$ sono uniformemente limitate da una costante}, mentre dato che $\int_U h(x_0,y)d\mu(y)<\infty $, grazie a una nota proprietà degli integrali se $\epsilon$ è sufficientemente piccolo $\int_{U_\epsilon} h(x_0,y)d\mu(y)$ può essere reso piccolo a piacere. Quindi ancora una volta otteniamo che:
\begin{gather*}
 \limsup_{x\to x_0} H(x) \leq H(x_0)
\end{gather*}
da cui la tesi.
\end{proof}

\end{prop}
\begin{oss}
 Osserviamo che nel caso $m=2$, la relazione \ref{eq_boh} segue dal ``principio di Frostman'' (vedi ad esempio la sezione 9E pag 320 di \cite{5}, o il teorema III.1 di \cite{36}) e non è necessario assumere che $G_\mu$ sia continuo quando ristretto all'insieme $S_\mu$. In dimensione maggiore, però, la dimostrazione di questo principio non è facilmente estendibile, per questa ragione riportiamo solo la versione generale della proposizione (che comunque è sufficiente per i nostri scopi).
\end{oss}
Questa proposizione e l'armonicità di $G_\mu$ su $R\setminus S_\mu$ ci permettono di dimostrare che:
\begin{prop}\label{prop_m1}
 Se $G_\mu$ è continua quando ristretta all'insieme $S_\mu$, allora il suo massimo è raggiunto su $S_\mu$, quindi se $G_\mu$ è continua quando ristretta a $S_\mu$ e $G_\mu(x)\leq c$ per ogni $x\in S_\mu$, allora $G_\mu(x)\leq c$ per ogni $x\in R$.
\begin{proof}
 La dimostrazione è una semplice applicazione del teorema \ref{teo_max2} (ricordiamo che $G_\mu|_{\Delta}=0$).
\end{proof}
\end{prop}
Questa proposizione può essere migliorata, nel senso che non è necessario chiedere la continuità di $G_\mu$ sull'insieme $S_\mu$. A questo scopo riportiamo i seguenti lemmi, tratto dai teoremi 3.6.2 e 3.6.3 di \cite{37}, e cortesemente segnalati dal professor. Wolfhard Hansen (University of Bielefeld):
\begin{lemma}
 Data una misura positiva regolare finita $\mu$ con supporto compatto $S_\mu$ tale che $G_\mu<\infty$ su $S_\mu$, per ogni $\epsilon>0$, esiste un insieme compatto $C\subset S_\mu$ tale che, detta $\mu|(\cdot)C\equiv \mu(\cdot\cap C)$, $G_{\mu|C}$ è continua su $C$ e $\mu(S_\mu \setminus C)<\epsilon$.
\begin{proof}
 Grazie al teorema di Lusin (vedi ad esempio teorema 2.23 pag 53 di \cite{12}), per ogni $\epsilon>0$ esiste un compatto $C\subset S_\mu$ tale che $G_\mu |_C$ è una funzione continua. Dato che:
\begin{gather*}
 G_{\mu|C}=G_\mu-G_{\mu-\mu|C}
\end{gather*}
e dato che su $C$ $G_\mu$ è continua, $G_{\mu|C}$ è semicontinua superiormente su $C$, e quindi continua su $C$ e automaticamente continua su tutto $R$ grazie alla proposizione \ref{prop_cont}.
\end{proof}
\end{lemma}
\begin{lemma}
 Se $G_\mu<\infty$ su $S_\mu$, allora esiste una successione di misure di Borel regolari $\mu_n$ tali che $\{G_{\mu_n}\}$ è una successione crescente di funzioni continue su $R$ e:
\begin{gather*}
 G_{\mu}(x)=\lim_n G_{\mu_n}(x)
\end{gather*}
\begin{proof}
 Scegliamo per induzione una successione di insiemi compatti $C_n$ tali che $C_n\subset C_{n+1}$, $\mu(S_\mu\setminus C_n)\leq 2^{-n}$ e $G_\mu$ continua sull'insieme $C_n$, e consideriamo $\mu_n =\mu|C_n$.\\
Dal lemma precedente, sappiamo che $G_{\mu_n}$ è una funzione continua, e per costruzione delle misure $\mu_n$, ovviamente $\{G_{\mu_n}(x)\}$ è una successione crescente $\forall x\in R$. Inoltre:
\begin{gather*}
 G_{\mu_n}(x)=\int_K G(x,y) \chi_{C_n}(y) d\mu(y)
\end{gather*}
e quindi per convergenza monotona, $\lim_nG_{\mu_n}(x)=G_\mu (x)$.
\end{proof}
\end{lemma}
Ora siamo pronti per dimostrare che:
\begin{prop}\label{prop_m2}
 Se $G_\mu(x)\leq c$ per ogni $x\in S_\mu$, allora $G_\mu(x)\leq c$ per ogni $x\in R$.
\begin{proof}
 Sia $\{\mu_n\}$ una successione di misure con le caratteristiche descritte nell'ultimo lemma, allora $G_{\mu_n}(x)\leq G_\mu (x)\leq c$ per ogni $x\in S_\mu$, e dato che $G_{\mu_n}$ sono funzioni continue, grazie alla proposizione \ref{prop_m1}, sappiamo che per ogni $x$ per ogni $n$:
$$G_{\mu_n}(x)\leq c$$
Passando al limite su $n$ otteniamo la tesi.
\end{proof}

\end{prop}

Grazie a questa proposizione, possiamo dimostrare questo teorema che garantisce lega l'energia di un insieme alla sua capacità:
\begin{teo}\label{teo_en2}
 Sia $K$ una sottovarietà regolare di codimensione $1$ di $R$ possibilmente con bordo liscio. Esiste una misura di Borel regolare unitaria \footnote{cioè $\nu(R)=1$} $\nu$ con supporto in $K$ tale che:
\begin{gather*}
 \epsilon(K)= \psb{\nu}{\nu}<\infty
\end{gather*}
Inoltre detto $u$ il potenziale armonico dell'insieme $K$, per il potenziale di Green relativo a $\nu$ vale che:
\begin{enumerate}
\item $G_{\nu}\in H(R\setminus K)$
\item $G_{\nu}(x)=\epsilon(K)u(x)$ per ogni $x\in R$
\item $G_{\nu}\in \Rod$
\item $D_R(G_{\nu})=\epsilon(K)$
\item $S_\nu=K$
\end{enumerate}
Inoltre la misura $\nu$ è unica.
\begin{proof}
Sia $u$ il potenziale di capacità dell'insieme $K$ (vedi osservazione \ref{oss_altri_dir}). Sia $\nu$ la misura che risolve il problema \ref{prop_endeph}. Sappiamo che il potenziale $G_\nu$ soddisfa $G_\nu|_{S_\nu}\leq \epsilon(K)u|_{S_\nu}=\epsilon(K)$, e grazie alla proposizione \ref{prop_m2}, sappiamo anche che $G|_\nu \leq \epsilon(K)$ su tutta la varietà $R$. Dato che entrambe le funzioni $G_\nu$ e $u$ sono armoniche in $R\setminus K$, su annullano su $\Delta$ e $G_\nu|_K\leq \epsilon(K)=\epsilon(K)u|_K$, grazie al principio \ref{teo_max2}, possiamo concludere che
\begin{gather*}
 G_\nu(x)\leq \epsilon(K)u(x) \ \ \ \forall x\in R
\end{gather*}
Per dimostrare questa disuguaglianza con il verso opposto, sappiamo che $G_\nu(x)\geq \epsilon(K)u(x)$ per $x\in K\setminus F$ dove $F$ è un'insieme di misura armonica $\xi$ nulla, quindi grazie al lemma \ref{prop_9f} esiste una funzione $w$ armonica su $R\setminus K$ che tende a infinito nei punti di $F$. Allora per ogni $\epsilon>0$, la funzione
\begin{gather*}
 G_\nu(x)+\epsilon w(x)\geq \epsilon(K)u(x)
\end{gather*}
per ogni $x\in R\setminus K$ grazie al teorema \ref{teo_max2}. Data l'arbitrarietà di $\epsilon$, possiamo concludere che
\begin{gather*}
 G_\nu(x)=\epsilon(K)u(x) \ \ \ \forall x\in R\setminus K
\end{gather*}
Dato che $u$ è continua su $R$ e $G_\nu$ è superarmonica, quindi semicontinua inferiormente, l'uguaglianza vale su tutto $R$. Sia infatti $x_0\in K$. Allora esiste una successione $\{x_n\}\subset R\setminus K$ che converge a $x_0$ e per la quale:
\begin{gather*}
 G_\nu(x_0)\geq \liminf_n G_\nu(x_n) \geq \lim_n \epsilon(K)u(x_n)=\epsilon(K)
\end{gather*}
dato che $G_\nu\leq \epsilon(K)$, necessariamente $G_\nu(x_0)=\epsilon(K)=\epsilon(K)u(x_0)$.\\
Per quanto riguarda la proprietà (4), osserviamo che per ogni aperto relativamente compatto con bordo liscio $C$ tale che $K\subset C$:
\begin{gather*}
 D_R(u)=-\int_{\partial C}\ast du
\end{gather*}
Consideriamo infatti una successione di aperti con bordo liscio $\{A_n\}$ tali che $A_n\subset A_{n-1}$ e $K=\cap_n A_n$ \footnote{un esempio di insiemi con queste caratteristiche è descritto nell'osservazione \ref{oss_an}}. Allora sappiamo che:
\begin{gather*}
 D_R(u)=\lim_n D_{R\setminus A_n} (u)= \lim_n \int_{\partial (R\setminus A_n)} u\ast du =-\lim_n \int_{\partial A_n} u\ast du
\end{gather*}
Grazie al fatto che $u$ è continua, è anche uniformemente continua su ogni insieme compatto, quindi ad esempio sull'insieme $\overline{A_1}$. Quindi per ogni $\epsilon>0$, esiste $\delta$ tale che
\begin{gather*}
 d(x,y)<\delta \ \ \Rightarrow \ \ \abs{u(x)-u(y)}<\epsilon
\end{gather*}
In particolare, se $d(x,K)<\delta$, otteniamo che $0\leq 1-u(x)\leq \epsilon$. Dato che $\cap_nA_n=K$, $A_n$ è contenuto definitivamente nell'aperto
\begin{gather*}
 K+\delta \equiv \{x\in R \ t.c. \ d(x,K)<\delta\}
\end{gather*}
quindi definitivamente in $n$ vale che:
\begin{gather*}
\int_{\partial A_n} (1-u)\ast du =\int_{\partial A_n} (1-u)(\ast du)^+ -\int_{\partial A_n} (1-u)(\ast du)^- \leq \\
\leq \epsilon \int_{\partial A_n}(\ast du)^+\leq \epsilon \int_{\partial A_n}\ast du\\
\int_{\partial A_n} (1-u)\ast du =\int_{\partial A_n} (1-u)(\ast du)^+ -\int_{\partial A_n} (1-u)(\ast du)^- \geq \\
\geq-\epsilon \int_{\partial A_n}(\ast du)^-\geq -\epsilon \int_{\partial A_n}\ast du
\end{gather*}
cioé:
\begin{gather*}
 -\epsilon\int_{\partial A_n} \ast du \leq \int_{\partial A_n} (1-u)\ast du \leq \epsilon \int_{\partial A_n}\ast du
\end{gather*}
dato che $u\in H(\overline{A_n}\setminus A_m)$: $$\int_{\partial A_n}\ast du =\int_{\partial A_m}\ast du$$ per ogni $n$ e $m$, quindi possiamo concludere che
\begin{gather*}
 \lim_n \int_{\partial A_n}(1-u)\ast du =0 \  \ \Rightarrow \\
\Rightarrow\ \ -D_R(u)= \lim_n\int_{\partial A_n}u\ast du=\lim_n \int_{\partial A_n} \ast du =\int_{\partial C} \ast du
\end{gather*}
Passiamo ora a considerare $D_R(G_\mu)$
\begin{gather*}
-D_R(G_\nu)=-D_R(\epsilon(K) u)=\epsilon(K)^2 \int_{\partial C} \ast du =\epsilon(K)\int_{\partial C} \ast d(G_{\nu})=\\
=\epsilon(K)\int_{\partial C}\ast d\ton{\int_K G(x,y)d\mu(y)}=\epsilon(K)\int_{\partial C} \int_{K} \ast d G(\cdot,y)d\mu(y)=\\
=\epsilon(K)\int_K d\mu(y)\int_C \ast d G(\cdot,y) = -\epsilon(K)
\end{gather*}
dove abbiamo sfruttato il fatto che $G(x,y)\in C^\infty (K\times C)$ per scambiare tra loro i segni di integrali e derivate, e la relazione (6) di \ref{deph_fgreen}.\\
Resta da dimostrare la proprietà (5). A questo scopo notiamo che per il principio del massimo forte, una funzione armonica su $R$ non può assumere il suo massimo in un punto interno del dominio, quindi se $x\in K\setminus S_\nu$, $G_\nu(x)<\max_{x\in R}G_\nu(x)=\epsilon(K)$. Ma dato che $G_\nu(x)=\epsilon(K)u(x)$, e dato che il potenziale di capacità $u$ è identicamente uguale a $1$ su $S_\nu$, necessariamente l'insieme $K\setminus S_\nu$ è vuoto.\\
Una dimostrazione alternativa di questa affermazione si può ottenere sfruttando una tecnica del tutto analoga a quella utilizzata nella dimostrazione della proposizione \ref{prop_alpha}. Supponiamo per assurdo che $S_\nu\not =K$. Allora, dato che $S_\nu$ è chiuso, sulla sottovarietà $K$ esiste una bolla chiusa $m-1$ dimensionale contenuta in $K$ ma disgiunta da $S_\nu$. Grazie all'osservazione \ref{oss_notinfty}, esiste una misura $\mu\in m_K$ con supporto contenuto nella bolla tale che $\psb \mu\mu <\infty$. Per ogni $0\leq \lambda \leq 1$, consideriamo la misura:
\begin{gather*}
 \nu_\lambda\equiv (1-\lambda)\nu +\lambda \mu
\end{gather*}
\`E facile verificare che $\nu_\lambda\in m_K$ per ogni $0\leq \lambda\leq 1$ inoltre:
\begin{gather*}
 \psb{\nu_\lambda}{\nu-\lambda}-\psb \nu \nu= \\
=[(1-\lambda)^2-1]\psb \nu \nu + \lambda^2 \psb \mu \mu +2\lambda (1-\lambda)\psb \nu \mu=\\
=\lambda^2(\psb \nu \nu +\psb \mu \mu-2\psb \nu \mu) + 2\lambda (\psb \nu \mu -\psb \nu \nu)
\end{gather*}
Dato che
\begin{gather*}
 \psb \nu \mu=\int_K G_\nu(x) d\mu \leq M < \epsilon(K)<\infty
\end{gather*}
dove $M$ è il massimo della funzione $G_\nu$ sull'insieme $S_\mu$, che è strettamente minore di $\epsilon(K)$ grazie al principio del massimo, $\psb \nu \mu -\psb \mu \mu<0$, quindi esiste un valore di $\lambda$ sufficientemente piccolo per il quale
\begin{gather*}
 \psb {\nu_\lambda}{\nu_\lambda}<\psb \nu \nu
\end{gather*}
contraddicendo l'ipotesi di minimalità di $\nu$.\\
Concludiamo con la dimostrazione dell'unicità della misura $\nu$. Supponiamo per assurdo che esista un'altra misura di minimo $\bar\nu$. Vale comunque che $G_\nu(x)=\epsilon(K)u(x)=G_{\bar\nu}(x)$. Grazie al punto (5) della proposizione \ref{prop_greenref}, per ogni funzione liscia a supporto compatto in $R$, abbiamo che:
\begin{gather*}
 \int_R f(y)d\nu(y)=\int_R\ton{\int_R G(x,y)\Delta f(x)d\lambda(x)}d\nu(y) 
\end{gather*}
dato che 
\begin{gather*}
\int_R \abs{\Delta f(x)}\int_R G(x,y)d\nu(y) d\lambda(x)\leq M \epsilon(K) \lambda(S)<\infty
\end{gather*}
dove $M$ è un maggiorante di $\abs{\Delta f(x)}$ e $S$ il supporto della funzione $f$, grazie al teorema di Fubini possiamo scambiare l'ordine di integrazione e ottenere:
\begin{gather*}
 \int_R f(y)d\nu(y)=\int_R \Delta f(x) G_\nu(x) d\lambda(x)=\int_R f(y)d\bar \nu(y)
\end{gather*}
Per densità delle funzioni lisce a supporto compatto nelle funzioni continue a supporto compatto, necessariamente $\nu=\bar \nu$.
\end{proof}
\end{teo}
\begin{oss}\label{oss_infty}
 Osserviamo subito che nelle ipotesi del teorema precedente, otteniamo che la capacità di $K$ è precisamente l'inverso di $\epsilon(K)$, infatti:
\begin{gather*}
 Cap(K)=D_R(u)=\frac{D_R(G_\mu)}{\epsilon(K)^2}=\epsilon(K)^{-1}
\end{gather*}
\end{oss}
Con la seguente stima leghiamo l'energia al diametro transfinito di un insieme.
\begin{prop}\label{prop_cfrrhoe}
Sia $K\subset R$. Allora $\rho(K)\geq \epsilon(K)$.
\begin{proof}
 Grazie alla definizione di $\rho_n(K)$, per ogni $n$ possiamo scegliere $n$ punti $p_1\cdots,p_n$ in $K$ tali che
\begin{gather}\label{eq_bi}
 \binom n 2 \rho_n(K) \geq \sum_{i<j}^{1,\cdots,n} G(p_i,p_j) -\frac{1}{n}
\end{gather}
Sia $\mu_n$ la misura data da
\begin{gather*}
 \mu_n=\frac1 n \sum_{i=1}^n \delta_{p_i}
\end{gather*}
dove $\delta_x$ è la misura che associa $1$ agli insiemi contenenti $x$ e $0$ agli altri. Quindi $\mu_n\in m_K$. Allora esiste una sottosuccessione di $\{\mu_n\}$ (che per comodità continueremo a indicare con lo stesso simbolo) che converge debolmente nel senso della misura a una misura $\mu$, cioè esiste una misura $\mu$ tale che per ogni funzione $\phi$ continua su $K$:
\begin{gather*}
 \lim_n \int_K \phi d\mu_n=\int_K \phi d\mu
\end{gather*}
quindi la misura $\mu$ è di Borel, positiva e unitaria. Definiamo le funzioni:
\begin{gather*}
 \phi_n^c (x)\equiv \int_K (G(x,y)\curlywedge c) d\mu_n(y) \ \ \ \ \ \  \phi^c (x)\equiv \int_K (G(x,y)\curlywedge c) d\mu(y)
\end{gather*}
e seguendo lo stesso ragionamento della dimostrazione della proposizione \ref{prop_endeph}, otteniamo la validità dell'equazione \ref{eq_mis}, cioè otteniamo:
\begin{gather*}
\lim_n \int_Kd\mu_n(x)\int_Kd\mu_n(y) (G(x,y)\curlywedge c) =  \int_Kd\mu(x)\int_Kd\mu(y) (G(x,y)\curlywedge c)
\end{gather*}
Moltiplicando ambo i membri di \ref{eq_bi} per $\binom n 2 ^{-1}$ otteniamo:
\begin{gather*}
 \rho_n(K)+\frac{1}{n}\binom n 2 ^{-1} \geq \binom n 2 ^{-1} \sum_{i<j}^{1,\cdots,n} G(p_i,p_j)\geq \frac 2{n^2} \sum_{i<j}^{1,\cdots,n} (G(p_i,p_j)\curlywedge c)
\end{gather*}
dove abbiamo sfruttato il fatto che $\binom n 2 \leq \frac{n^2}2$. Ora dalla definizione di $\mu_n$ otteniamo che:
\begin{gather*}
\sum_{i,j}^{1,\cdots,n} (G(p_i,p_j)\curlywedge c)=n^2 \int_K d\mu_n \int_K d\mu_n (G(x,y)\curlywedge c)
\end{gather*}
considerando che la funzione $G$ è simmetrica nei suoi due argomenti e che:
\begin{gather*}
 \sum_{i=j=1}^n (G(p_i,p_j)\curlywedge c)=\sum_{i=1}^n (G(p_i,p_i)\curlywedge c)=nc
\end{gather*}
otteniamo:
\begin{gather*}
 \sum_{i<j}^{1,\cdots,n} (G(p_i,p_j)\curlywedge c)=\frac{n^2}{2}\int_K d\mu_n \int_K d\mu_n (G(x,y)\curlywedge c) - \frac{nc}{2} 
\end{gather*}
quindi:
\begin{gather*}
 \rho_n(K)+\frac{1}{n}\binom n 2 ^{-1} \geq\int_K d\mu_n \int_K d\mu_n (G(x,y)\curlywedge c) - \frac{c}{2n} 
\end{gather*}
facendo tendere $n$ a infinito, otteniamo:
\begin{gather*}
 \rho(K)\geq \lim_n \int_K d\mu_n \int_K d\mu_n (G(x,y)\curlywedge c) = \int_K d\mu \int_K d\mu (G(x,y)\curlywedge c)
\end{gather*}
grazie all'arbitrarietà del parametro $c$, concludiamo:
\begin{gather*}
 \rho(K) \geq \int_K d\mu \int_K d\mu G(x,y)\geq \epsilon(K)
\end{gather*}
come volevasi dimostrare.
\end{proof}
\end{prop}

L'ultima proposizione ``tecnica'' della dimostrazione è:
\begin{prop}\label{prop_capinfty}
 La capacità di $F_{n+1,m}$ tende a $0$ se $m$ tende a infinito.
\begin{proof}
 Ricordiamo che la capacità di $F_{n+1,m}$ (insieme che in questa dimostrazione indicheremo per comodità con $F_m$) equivale per definizione all'integrale di Dirichlet della funzione $u_{n+1,m}$ (che indicheremo per semplicità $u_m$), funzione armonica su $R\setminus F_m$, $u_m\in \Rod$ e $u_m|_{\partial F_m}=1$.\\
Grazie al lemma \ref{lemma_frod}, sappiamo che esiste una funzione $f_m\in \Rod$, $f_m=1$ su 
\begin{gather*}
 E_{n+1,m}\equiv \overline{U_{n+1}\cap K_m^C}
\end{gather*}
Infatti per definizione:
\begin{gather*}
 U_n\equiv\{z\in R \ t.c. \ G(z,z_0)>r_n\}
\end{gather*}
e per \ref{prop_greenref}, $G(\cdot,z_0)|_{\Delta}=0$, cioè $U_n\cap \Delta=\emptyset \ \forall n$. Possiamo scegliere questa funzione armonica su $R\setminus E_{n+1,m}$, infatti definiamo $W_{n+1,m,p}$ con $p>m$ l'insieme
\begin{gather*}
 W_{n+1,m,p}\equiv U_{n+1}\cap (\overline{K_p}\setminus K_m)
\end{gather*}
e chiamiamo $\pi_{m,p}$ l'operatore definito in \ref{teo_dec} rispetto all'insieme $ W_{n+1,m,p}$. La successione $\{w_{m,p}\equiv \pi_{m,p}(f_m)\}$ è una successione di funzioni armoniche su $R\setminus E_{n+1,m}$, appartenenti all'insieme $\Rod$.\\
Grazie al punto (2) del teorema \ref{teo_dec}, si ha che:
\begin{gather*}
 D(w_{m,p+q}-w_{m,p},w_{m,p})=0 \ \ \Rightarrow \ \ D(w_{m,p+q}-w_{m,p})=D(w_{m,p})-D(w_{m,p+q})
\end{gather*}
Con un ragionamento analogo a quello esposto in \ref{prop_r2}, otteniamo che la successione $\{w_{m.p}\}$ è D-cauchy rispetto all'indice $p$.\\
Osserviamo che la successione $\{w_{m,p}\}$ è una successione crescente in $p$, infatti dato che $w_{m,p}=1$ su $W_{n+1,m,p}$ e $w_{m,p}|_{\Delta}=0$, grazie all'osservazione \ref{prop_delta} $w_{m,p}\leq 1$ su $R$, quindi in particolare su $W_{n+1,m,p+1}$ si ha che $w_{m,p}\leq w_{m,p+1}$. Sempre grazie alla proposizione \ref{prop_delta}, otteniamo che questa disuguaglianza è valida su tutto $R$. Quindi la successione $w_{m,p}$ converge monotonamente rispetto a $p$ a una funzione $w_m$ armonica su $R\setminus E_{n+1,m}$ per il principio di Harnack e identicamente uguale a $1$ su $E_{n+1,m}$ per costruzione. È facile verificare che la convergenza è locale uniforme su $R$.\\
Dimostriamo che:
\begin{gather*}
 0=BD-\lim_m w_m
\end{gather*}
Il punto (2) di \ref{teo_dec} garantisce che per ogni $p$:
\begin{gather*}
D_R(w_{m+q,\ p}-w_m,w_{m+q,\ p})=0
\end{gather*}
infatti $w_{m+q,\ p}-w_m=0$ su $W_{n+1,m,p}$ e questa funzione è in $\Rod$. Passando al limite su $p$, otteniamo che:
\begin{gather*}
 0=D_R(w_{m+q}-w_m,w_{m+q})=D_R(w_{m+q})-D_R(w_{m+q},w_m)
\end{gather*}
da cui:
\begin{gather*}
 D_R(w_{m+q}-w_m)=D_R(w_{m})-D_R(w_{m+q})
\end{gather*}
quindi grazie a un ragionamento simile a quello riportato in \ref{prop_r2}, otteniamo che la successione $\{w_{m}\}$ è D-cauchy.\\
La successione $w_n$ inoltre è una successione decrescente di funzioni. Infatti per ogni $p$ abbiamo che $w_{m+1,p}\leq 1$ sull'insieme $E_{n+1,m}$ grazie alla proposizione \ref{prop_delta}, e quindi $w_{m+1,p}\leq w_m$ sull'insieme $E_{n+1,m}$. La disuguaglianza vale su tutta la varietà $R$ (quindi per continuità anche su tutta $R^*$) grazie al teorema \ref{teo_max2}. Questo garantisce che esiste il limite
\begin{gather*}
 w=BD-\lim_m w_m
\end{gather*}
Grazie al principio di Harnack, questa funzione è armonica su tutta $R$, e per l'osservazione \ref{oss_rodcmpl} $w\in \Rod$. Grazie al principio del massimo \ref{teo_max}, otteniamo che $w=0$, quindi:
\begin{gather*}
 \lim_m D_R(w_m)=0
\end{gather*}
Questo risultato è utile in quanto
\begin{gather*}
 D_R(u_m)\leq D_R(w_m)
\end{gather*}
dove questa osservazione segue dal teorema \ref{teo_dec}. Consideriamo infatti l'insieme $F_{n+1,m}$ come l'insieme $K$ del teorema. La proiezione della funzione $w_m$ su questo insieme è la funzione $u_m$ dato che $(w_m-u_m)|_{F_{n+1,m}}=0$ e questa funzione è in $\Rod$. Quindi grazie al punto (3) del teorema abbiamo che:
\begin{gather*}
 D_R(w_m)=D_R(u_m)+D_R(w_m-u_m)
\end{gather*}
da cui la tesi.
\end{proof}

\end{prop}

\subsection[Il diametro transfinito è infinito]{Il diametro transfinito $\rho(\Xi_n)=\infty$}
Mettendo assieme le varie proposizioni e lemmi visti fino ad ora, siamo pronti per dimostrare che:
\begin{prop}
Per ogni $n$, $\tau(\Xi_n)=\rho(\Xi_n)=\infty$. Questo implica che ogni compatto contenuto nel bordo irregolare $\Xi$ ha diametro transfinito infinito.
\begin{proof}
 Grazie alla proposizione \ref{prop_cfrdiam} sappiamo che
\begin{gather*}
 \tau(X)\geq \rho(X)
\end{gather*}
per ogni insieme $X\subset R^*$. Quindi basta dimostrare che $\rho(\Xi_n)=\infty$. La relazione \ref{eq_infty1} combinata con la proposizione \ref{prop_cfrrhoe} ci permette di scrivere:
\begin{gather}
 \rho(\Xi_n)\geq \sigma_n^2 \rho(F_{n+1,m})\geq \sigma_n^2 \epsilon(F_{n+1,m})
\end{gather}
Questo assicura che se $\lim_m \epsilon(F_{n+1,m})=\infty$ otteniamo la tesi.\\
Per dimostrare questa uguaglianza, ricordiamo che $F_{n+1,m}$ è una sottovarietà di codimensione $1$ di $R$ con bordo liscio, quindi grazie a \ref{oss_infty}:
\begin{gather*}
 \epsilon(F_{n+1,m})=Cap(F_{n+1,m})^{-1}
\end{gather*}
e grazie alla proposizione \ref{prop_capinfty}, che assicura che:
\begin{gather*}
 \lim_m Cap(F_{n+1,m})=0
\end{gather*}
otteniamo la tesi.
\end{proof}
\end{prop}

\subsection{Funzioni armoniche che tendono a infinito sul bordo di R}
Utilizzando in fatto che $\tau(\Xi_n)=\infty$ per ogni $n$, in questa sezione costruiremo una funzione armonica su $R$ che converge a infinito su $\Xi$ (in un senso che analizzeremo meglio in seguito).
\begin{prop}
 Esiste una funzione armonica $E:R\to \R$ tale che $E(p)=\infty$ per ogni $p\in \Xi$.
\begin{proof}
 La dimostrazione è costruttiva. Dato che $\tau(\Xi_n)=\infty$ per ogni $n$, per definizione di $\tau$ \footnote{vedi \ref{deph_tau}} esiste una successione di interi $n_k$ tale che
\begin{gather*}
 \tau_{n_k}(\Xi_n)\geq 2^k
\end{gather*}
cioè sempre per definizione di $\tau_n$, esistono $n_k$ punti $p_1,\cdots,p_{n_k}\in \Xi_n$ tali che
\begin{gather*}
 \inf_{p\in \Xi_n} \sum_{i=1}^{n_k}G(p,p_i) >2^k n_k
\end{gather*}
definiamo la funzione
\begin{gather*}
 E_{n,k}(x)\equiv \sum_{i=1}^{n_k} \frac{1}{2^k n_k}G(x,p_{n_k})
\end{gather*}
Osserviamo che $E_{n,k}(p)>1$ per ogni $p\in \Xi_n$. Se definiamo
\begin{gather*}
 E_n (x) \equiv \sum_{k=1}^{\infty} E_{n,k} (x)
\end{gather*}
osserviamo che in ogni punto di $R$ la somma converge. Infatti per come è definita $E_{n,k}$, $E_n$ è una combinazione convessa delle funzioni $G(x,p_{n,k})$, quindi una serie di funzioni armoniche positive. Per il principio di Harnack, questa successione o converge localmente uniformemente in $R$ o diverge in tutti i punti di $R$.\\
Consideriamo $z_0\in R$ un punto qualsiasi. Sia $V$ un suo intorno relativamente compatto. Sappiamo che
\begin{gather*}
 \sup_{x\in R^*\setminus V} G(x,z_0)\leq \sup_{x\in \partial V}G(x,z_0)<\infty
\end{gather*}
questo vuol dire che la successione $G(z_0,p_k)$ è uniformemente limitata, quindi qualunque sua combinazione convessa converge \footnote{e dal principio di Harnack converge localmente uniformemente a una funzione armonica}. Un modo meno banale ma più veloce di dimostrare la stessa cosa, è osservare che $G(z_0,\cdot)$ è una funzione continua su $\Gamma$ insieme compatto, quindi assume massimo finito su questo insieme.\\
Inoltre osserviamo che per ogni $p\in \Xi_n$, $E_n(p)=\infty$. Per definire la funzione $E$ ripetiamo un ragionamento simile a quello fin qui esposto, in particolare definiamo:
\begin{gather*}
 E(x)\equiv \sum_{n=1}^{\infty} \frac{1}{2^n} E_n(x)
\end{gather*}
anche in questo caso $E$ è una combinazione convessa di funzioni della forma $G(\cdot,p)$ con $p\in \Xi$ \footnote{$p$ può stare in un $\Xi_n$ qualsiasi, quindi in generale $p\in \Xi$}, cioè esiste una successione di numeri positivi $t_k$ con somma $\sum_{k=1}^{\infty} t_k=1$ e una successione di punti $p_k$ in $\Xi$ tale che:
\begin{gather}\label{eq_E1}
 E(x)= \sum_{k=1}^{\infty} t_k G(x,p_k)
\end{gather}
Grazie al principio di Harnack, $E(x)$ è una funzione armonica su $R$ e per ogni $p\in \Xi$, $E(p)=\infty$.
\end{proof}
\end{prop}
Il fatto che $E(p)=\infty$ per ogni $p\in \Xi$ è assolutamente inutile se non consideriamo che:
\begin{oss}
 La funzione $E$ è semicontinua inferiormente su $R^*$, nel senso che
\begin{gather*}
 \liminf_{p\to p_0} E(p)\geq E(p_0)
\end{gather*}
per ogni $p$ in $R^*$.\\
Ricordiamo la definizione di $\liminf$ in una topologia non I numerabile:
\begin{gather*}
 \liminf_{p\to p_0} f(p)=L \ \Longleftrightarrow \ \forall \epsilon>0, \ \exists V(p_0) \ t.c. \ E(p)>L-\epsilon \ \forall p\in V(p_0) \ \wedge \\\wedge \ \forall \epsilon>0, \ \forall V(p_0), \ \exists p\in V(p_0) \ t.c. \ E(p)<L+\epsilon
\end{gather*}

\begin{proof}
Se $p\in R\cup \Delta$, la dimostrazione è ovvia essendo $E$ continua su $R$, positiva ovunque e uguale a $0$ su $\Delta$ \footnote{infatti tutte le funzioni $G(\cdot,p)$ si annullano sul bordo armonico $\Delta$, vedi \ref{prop_greenref}}.\\
Data la positività delle funzioni di Green, abbiamo anche che:
\begin{gather*}
 \inf_{q\in V\setminus\{q_0\}}\sum_{i=1}^{\infty}\alpha_i G(q,p_i)\geq  \inf_{q\in V\setminus\{q_0\}}\sum_{i=1}^{N}\alpha_i G(q,p_i)
\end{gather*}
per ogni scelta di $\alpha_i\geq 0$ e per ogni interno $N$.\\
Dato che le funzioni di Green sono semicontinue inferiormente \footnote{infatti se $q_0\neq p_i$, questo segue dalla continuità della funzione $G$, altrimenti dal fatto che $\lim_{V\to x} G(y,x)=\infty$, vedi sezione \ref{subsec_green1}}, cioè:
\begin{gather*}
 \lim_{V\to q_0}\inf_{V\setminus\{ q_0\}} G(q,p_i)= G(q_0,p_i)
\end{gather*}
possiamo osservare che:
\begin{gather*}
 \lim_{V\to q_0}\inf_{q\in V\setminus\{q_0\}}\sum_{i=1}^{\infty}\alpha_i G(q,p_i)\geq \lim_{V\to q_0} \inf_{q\in V\setminus\{q_0\}}\sum_{i=1}^{N}\alpha_i G(q,p_i)=\sum_{i=1}^{N}\alpha_i G(q,p_i)
\end{gather*}
visto che la relazione vale per ogni $n$, si ha la tesi.
\end{proof}
\end{oss}
Questa osservazione garantisce che per $p\to p_0, \ p\in R$, $E(p)\to \infty$. Riassumendo, abbiamo dimostrato l'esistenza di una funzione $E:R^*\to \R\cup \{\infty\}$ continua su $R$ e semicontinua inferiormente su $R^*$, tale che $E(p)=\infty$ per ogni $p\in \Xi$. Ora ci concentriamo sul dimostrare altre proprietà di questa funzione, come ad esempio caratterizzare l'integrale di Dirichlet $D_R(E\curlywedge c)$ per $c>0$ \footnote{come avevamo fatto per le funzion di Green $G^*(\cdot,p)$}.
\begin{prop}
 Per ogni $c>0$ vale che
\begin{gather*}
 D_R(E(\cdot)\curlywedge c)\leq c
\end{gather*}
\begin{proof}
 Siano $t_k$ e $p_k$ tali che valga la relazione \ref{eq_E1}. Definiamo la successione di funzioni
\begin{gather*}
 \psi_n(p)\equiv \sum_{k=1}^n t_k G(p,p_k)
\end{gather*}
cioè la successione delle somme parziali che definiscono $E(\cdot)$. Grazie alle considerazioni precedenti, sappiamo che $\psi_n$ converge localmente uniformemente a $E$, e se dimostriamo che:
\begin{gather*}
 D_R(\psi_n\curlywedge c)\leq c
\end{gather*}
allora grazie al teorema \ref{prop_r1} \footnote{e all'osservazione \ref{oss_r1}} abbiamo la tesi.\\
Resta da dimostrare la parte tecnica del teorema, cioè l'ultima relazione. Fissiamo $n$ e definiamo per comodità $\psi_n(p)\equiv \psi(p)$ \footnote{fino a quando non ci può essere confusione con l'indice $n$}. Sia $\alpha$ un valore regolare per la funzione $\psi$. 
Definiamo la quantità:
\begin{gather*}
 L(\alpha)\equiv \sum_{j=1}^n \int_{\psi=\alpha} \abs{\ast dG(\cdot,p_j)}=\sum_{j=1}^m\int_{\phi(\partial \Omega)}\abs{\frac{\partial G(\cdot,p_j)}{\partial x}(y)} \sqrt{\abs g} dy^1\cdots dy^{m-1}
\end{gather*}
dove $(x,y^1\cdots,y^{m-1})$ sono coordinate per la sottovarietà $\psi=\alpha$ dove $$\{\psi=\alpha\}=\{y^m=0\}$$ e dove la metrica assume la forma particolare \footnote{cioè per coordinatizzare un intorno di $\{\psi=\alpha\}$ scegliamo come coordinate $x$ e altre tutte ortogonali a $x$, possibile grazie al fatto che $\alpha$ è un valore regolare per la funzione $x$}
\begin{gather*}
 g=\begin{bmatrix} 1 & 0 \\ 0 & A(y) \end{bmatrix}
\end{gather*}
Ricordando che per ogni $n$-upla di numeri positivi
\begin{gather*}
\ton{ \sum_{i=1}^n a_i }^2\leq n\sum_{i=1}^n a_i^2
\end{gather*}
e la disuguaglianza di Schwartz, cioè:
\begin{gather*}
 \ton{\int \abs fg dx }^2\leq \int f^2 dx\int g^2 dx
\end{gather*}
abbiamo che:
\begin{gather*}
 L(\alpha)^2\leq\\
\leq n \sum_{j=1}^m\int_{\phi(\partial \Omega)} \ton{\frac{\partial G(\cdot,p_j)}{\partial x}(y)}^2\sqrt{\abs g} dy^1\cdots dy^{m-1} \int_{x=\alpha} \sqrt{\abs g} dy^1\cdots dy^{m-1}\leq\\
\leq n \sum_{j=1}^n\int_{\psi=\alpha} \abs{\nabla G(\cdot,p_j)}^2\sqrt{\abs g} dy^1\cdots dy^{m-1} \int_{\psi=\alpha} \ast dx
\end{gather*}
Dimostreremo nel lemma seguente (lemma \ref{lemma_w} e osservazione seguente) che vale la proprietà:
\begin{gather}\label{eq_w}
 D_R(\psi\curlywedge \alpha)=\alpha \int_{\psi=\alpha}\ast dx
\end{gather}
Consideriamo due numeri reali $c, \ c'$ con $0<c<\alpha<c'$ dove $c$ e $c'$ sono valori regolari della funzione $\psi$. Allora vale che:
\begin{gather*}
 \int_{\psi=\alpha}\ast dx = \alpha D_R(\psi\curlywedge \alpha) \leq c' D_R(\psi\curlywedge c')
\end{gather*}
quindi
\begin{gather*}
 L(\alpha) \leq n c' D_R(\psi\curlywedge c') \sum_{j=1}^n\int_{\psi=\alpha} \abs{\nabla G(\cdot,p_j)}^2\sqrt{\abs g} dy^1\cdots dy^{m-1}
\end{gather*}
Inoltre dato che l'insieme dei valori non regolari per $\psi$ ha misura nulla (vedi teorema di Sard), vale che:
\begin{gather*}
 \int_c^{c'} d\alpha \sum_{j=1}^n\int_{\psi=\alpha} \abs{\nabla G(\cdot,p_j)}^2\sqrt{\abs g} dy^1\cdots dy^{m-1}=\sum_{j=1}^n \int_{c<\psi<c'} \abs{\nabla G(\cdot,p_j)}^2 dV\leq\\
\leq \sum_{j=1}^n \int_{\psi<c'} \abs{\nabla G(\cdot,p_j)}^2 dV
\end{gather*}
Per definizione della funzione $\psi$, se $\psi(z)<c'$, allora necessariamente per ogni $j=1\cdots n$
\begin{gather*}
 G(z,p_j)<\frac{c'}{t_j}
\end{gather*}
quindi grazie alla proprietà (5) di \ref{prop_greenref}, si la stima sull'integrale di Dirichlet:
\begin{gather*}
 \int_{\psi<c'} \abs{\nabla G(\cdot,p_j)}^2 dV\leq D_R(G(\cdot,p_j)\curlywedge c't_j^{-1})\leq \frac{c'}{t_j}
\end{gather*}
Da questa stima, ricordando che $\psi\curlywedge c \in \Rod$, possiamo ottenere che:
\begin{gather*}
 \int_{c<\psi<c'} L(\alpha)^2 d\alpha \leq nc' D_R(\psi\curlywedge c') \sum_{j=1}^n \frac{c'}{t_j}<\infty
\end{gather*}
Il ragionamento fino a qui discusso ha l'utilità di dimostrare (grazie al teorema di Fubini) che a meno di un insieme di $\alpha$ di misura nulla rispetto alla misura di Lebesgue su $\R$, vale che:
\begin{gather*}
 \int_{\psi=\alpha} \abs{\ast dG(\cdot,p_j)} <\infty
\end{gather*}
D'ora in avanti considereremo valori di $\alpha$ regolari per la funzione $\psi$ e per i quali vale la relazione appena scritta per ogni $j=1\cdots,n$ \footnote{grazie alle considerazioni fatte fino ad ora, risulta evidente che questo insieme abbia complementare di misura nulla in $\R$, quindi sia denso in $\R$}.\\
Grazie alla relazione \ref{eq_w} e alla definizione della funzione $\psi$, sappiamo che:
\begin{gather}\label{eq_w2}
 D_R(\psi\curlywedge \alpha)=\alpha \int_{\psi=\alpha}\ast dx = \alpha \sum_{k=1}^n t_k\int_{\psi_\alpha} \ast dG(\cdot,p_k)
\end{gather}
Sia $c_k>\alpha/t_k$ un valore regolare per la funzione $G(\cdot,p_k)$. Grazie alla positività di tutte le funzioni $G(\cdot,p_k)$, abbiamo che:
\begin{gather*}
 K\equiv \{G(\cdot,p_k)\geq c_k\} \subset K'\equiv \{\psi(\cdot)\geq \alpha\}
\end{gather*}
applicando ancora il lemma seguente (lemma \ref{lemma_w}) alla funzione $(G(\cdot,p_k)\curlywedge c_k)/c_k)$, otteniamo con l'aiuto della proposizione \ref{prop_greenref}
\begin{gather*}
 \int_{\psi=\alpha} \ast d(c_k^{-1} G(\cdot,p_k)) =\int_{\partial K'}\ast d(c_k^{-1} G(\cdot,p_k)) = \int_{\partial K} \ast d(c_k^{-1} G(\cdot,p_k)) =\\
= \int_{\partial K} \ast d(c_k^{-1} (G(\cdot,p_k)\curlywedge c_k))= D_R (c_k^{-1} (G(\cdot,p_k)\curlywedge c_k))\leq c_k^{-2} c_k
\end{gather*}
e quindi:
\begin{gather*}
 \int_{\psi=\alpha} \ast dG(\cdot,p_k) \leq 1
\end{gather*}
da cui utilizzando \ref{eq_w2}:
\begin{gather*}
  D_R(\psi\curlywedge \alpha)\leq \alpha \sum_{i=k1}^n t_k\int_{\psi_\alpha} \leq \alpha
\end{gather*}
Se $\alpha$ non è un valore regolare per $\psi$, comunque esiste una successione $\alpha_n$ che converge dall'alto a $\alpha$ di valori regolari di $\psi$, quindi:
\begin{gather*}
 D_R(\psi\curlywedge \alpha)\leq D_R(\psi\curlywedge \alpha_n)\leq \alpha_n
\end{gather*}
dato che la relazione vale $\forall n$, abbiamo la tesi.
\end{proof}
\end{prop}

\begin{lemma}\label{lemma_w}
Sia $K$ un'insieme compatto in $R^*$ il cui bordo in $R$ sia liscio liscio e tale che $K\cap \Delta=\emptyset$, e sia $u$ una funzione limitata tale che $u\in \Rod$, $u|_K=1$, $u\in H(R\setminus K)$. Allora si ha che:
\begin{enumerate}
 \item $D_R(u)=-\int_{\partial K} \ast du$
 \item se $K'$ è un insieme compatto in $R^*$ con bordo liscio, disgiunto da $\Delta$ e la cui parte interna contiene $K$, allora se $\int_{\partial (K'\setminus K)} \abs{\ast du}<\infty$:
\begin{gather*}
 D_R(u)=-\int_{\partial K} \ast du=-\int_{\partial K'} \ast du
\end{gather*}
\end{enumerate}
\begin{proof}
 Sia $R_n$ un'esaustione regolare di $R$, definiamo $K_n\equiv K\cap \overline{R_n}$. Con le notazioni dell'osservazione \ref{oss_decomp}, abbiamo che:
\begin{gather*}
 \pi_{K_n}(\pi_{K_{n+p}}(u))=\pi_{K_n}(u)
\end{gather*}
per ogni $n$ e $p$ naturali. Grazie al punto (2) del teorema \ref{teo_dec}, si ha che:
\begin{gather*}
 D(\pi_{K_{n+p}}(u)-\pi_{K_{n}}(u),\pi_{K_{n}}(u))=0 \ \ \Rightarrow \\
\Rightarrow \ \ D(\pi_{K_{n+p}}(u)-\pi_{K_{n}}(u))=D(\pi_{K_{n+p}}(u))-D(\pi_{K_{n}}(u))
\end{gather*}
Con un ragionamento analogo a quello esposto in \ref{prop_r2}, otteniamo che la successione $\{\pi_{K_n}(u)\}$ è D-cauchy.\\
Osserviamo che la successione $\pi_{K_n}(u)$ è una successione crescente in $n$, infatti dato che $\pi_{K_n}(u)=1$ su $K_n$ e $\pi_{K_n}(u)|_{\Delta}=0$, grazie all'osservazione \ref{prop_delta} $\pi_{K_n}(u)\leq 1$ su $R$, quindi in particolare su $K_{n+1}$ si ha che $\pi_{K_n}(u)\leq \pi_{K_{n+1}}(u)$. Sempre grazie alla proposizione \ref{prop_delta}, otteniamo che questa disuguaglianza è valida su tutto $R$. Quindi la successione $\pi_{K_n}(u)$ converge monotonamente a una funzione $f$ armonica su $R\setminus K$ per il principio di Harnack e identicamente uguale a $1$ su $K$ per costruzione. È facile verificare che la convergenza è locale uniforme su $R$.\\
Dato che $(f-u)|_K=0$ e $(f-u)|_{\Delta}=0$, per la proposizione \ref{prop_delta}, $f=u$ su $R$, quindi
\begin{gather*}
 u=BD-\lim_n \pi_{K_n}(u)
\end{gather*}
Fissato $n$, per ogni $m>n$, definiamo la funzione $v^{(n)}_m$ tale che:
\begin{enumerate}
\item $v^{(n)}_m|_{K_n}=1 $
\item $v^{(n)}_m|_{R\setminus R_m}=0 $
\item $v^{(n)}_m\in H(R_m\setminus K_n) $
\end{enumerate}
dalla dimostrazione del teorema \ref{teo_1} \footnote{o meglio dal suo adattamento alla dimostrazione di \ref{teo_dec}}, si ha che
\begin{gather*}
 \pi_{K_n}(u)=BD-\lim_m v^{(n)}_m
\end{gather*}
Applicando la formula di Green \ref{prop_green3} all'insieme $R_m\setminus K $ otteniamo che:
\begin{gather*}
 D_R(v^{(n)}_m,u)= D_{R_m\setminus K}(v^{(n)}_m,u)= -\int_{\partial K} v_m \ast du
\end{gather*}
dove il segno $-$ viene dall'orientazione di $\partial K$, contraria a quella di $\partial (R_m\setminus K)$.\\
Grazie al fatto che $u|_K=1$ e $u|_{R\setminus K}\leq1$, possiamo dedurre che $\ast du \leq 0$ su $\partial K$, e dato che $v^{(n)}_m$ è una successione crescente in $m$ \footnote{semplice applicazione del principio del massimo}, grazie al teorema di convergenza monotona possiamo concludere che:
\begin{gather*}
 D_R(\pi_{K_n}(u),u)=\lim_n D_R(v^{(n)}_m,u)=\lim_n -\int_{\partial K}v^{(n)}_m \ast du =-\int_{\partial K} \pi_{K_n}(u)\ast du
\end{gather*}
La successione $\pi_{K_n}(u)$ è crescente sull'insieme $\partial K_n$ e in particolare tende a $1$ su questo insieme, quindi sempre per convergenza monotona:
\begin{gather*}
 D_R(u)\equiv D_R(u,u)=\lim_n D_R(\pi_{K_n}(u),u)=\lim_n -\int_{\partial K}\pi_{K_n}(u)\ast du =-\int_{\partial K} 1\ast du
\end{gather*}
Per dimostrare il punto (3), sfruttiamo il lemma \ref{lemma_frod}. Sia $f\in \Rod$ tale che $f_{K'}=1$. Allora per definizione di $\Rod$, esiste una successione di funzioni $f_n\in \Roo$ tale che $f=BD-\lim_n f_n$. Quindi vale che:
\begin{gather*}
 0=D_{K'\setminus K} (1,u)=D_{K'\setminus K} (f,u)=\lim_n D_{K'\setminus K} (f_n,u)
\end{gather*}
Sia $S_n\equiv supp(f_n)$ \footnote{compatto in $R$ per definizione di $\Roo$}, allora:
\begin{gather*}
 D_{K'\setminus K} (f_n,u)=D_{(K'\setminus K)\cap S_n} (f_n,u)
\end{gather*}
applicando la formula di Green \ref{prop_green3} otteniamo:
\begin{gather*}
 D_{(K'\setminus K)\cap S} (f_n,u)=\int_{\partial (K'\setminus K)} f_n\ast du
\end{gather*}
Per ipotesi, sappiamo che $\int_{\partial (K'\setminus K)} \abs{\ast du}<\infty$, quindi per uniforme limitatezza della successione $f_n$ vale anche che $\int_{\partial (K'\setminus K)}\abs{f_n} \abs{\ast du}<M<\infty$ dove $M$ è indipendente da $n$. Quindi possiamo applicare il teorema di convergenza dominata e ottenere che:
\begin{gather*}
 \lim_n \int_{\partial (K'\setminus K)} f_n\ast du = \int_{\partial (K'\setminus K)}\ast du
\end{gather*}
Riassumendo, abbiamo ottenuto:
\begin{gather*}
 0=\int_{\partial (K'\setminus K)}\ast du = \int_{\partial K'}\ast du - \int_{\partial K}\ast du
\end{gather*}
\end{proof}

\end{lemma}
Riassumendo abbiamo dimostrato l'esistenza di funzioni con queste caratteristiche:
\begin{prop}\label{prop_Eref}
 Data una varietà $R$ iperbolica irregolare ($\Xi(R)\neq \emptyset$), esiste una funzione armonica $E:R\to [0,\infty)$ tale che:
\begin{enumerate}
 \item $E\in H(R)$
 \item $E(\cdot)$ è la combinazione convessa di funzioni $G(p_i,\cdot)$ dove $p_i\in \Xi$
 \item $E(z)=\infty$ se $z\in \Xi$, $E(z)=0$ se $z\in \Delta$
 \item $E(z)$ è semicontinua inferiormente su $R^*$
 \item $D_R(E(z)\curlywedge c)\leq c$ per ogni $c>0$ 
\end{enumerate}

\end{prop}

\subsection{Varietà iperboliche irregolari}\label{subset_irreg}
Data una varietà iperbolica $R$, sappiamo che esiste una funzione di Green $$G(\cdot,\cdot):R^*\times R^* \to \R\cup\{\infty\}$$ con le proprietà descritte nella proposizione \ref{prop_greenref}. L'insieme $\Xi$ è definito come il sottoinsieme di $\Gamma$ tale che
\begin{gather*}
 \Xi\equiv \{p \ t.c. \ G(z,p)>0 \ z\in R\}
\end{gather*}
dove la definizione non dipende dalla scelta del punto $z$. In funzione del fatto che $\Xi$ sia vuoto o meno definiamo una varietà iperbolica \textit{regolare} o \textit{irregolare}
\begin{deph}
 Data una varietà iperbolica $R$, diciamo che è \textbf{regolare} se $\Xi=\emptyset$, altrimento diciamo che è \textbf{irregolare}.
\end{deph}
Una facile caratterizzazione delle varietà regolari è la seguente:
\begin{prop}
 Sia $G(z,\cdot)$ un nucleo di Green su $R$ iperbolica. $R$ è regolare se e solo se per ogni $c>0$ l'insieme
\begin{gather*}
 A_c\equiv \{p\in R \ t.c. \ G(z.p)\geq c\}
\end{gather*}
è compatto in $R$.
\begin{proof}
 La dimostrazione di questa affermazione è abbastanza immediata. Supponiamo che gli insiemi $A_c$ siano compatti. Allora per ogni $p\in \Gamma$, sappiamo che:
\begin{gather*}
 G(z,p)=p(G(z,\cdot))=p(G(z,\cdot) (1-\lambda_c))
\end{gather*}
dove $\lambda_c$ è una funzione liscia a supporto compatto tale che $\lambda_c=1$ sull'insieme $A_c$. Grazie all'osservazione \ref{prop_charG} vale la seconda uguaglianza scritta in formula, e inoltre grazie all'osservazione \ref{prop_maxG2} possiamo concludere che la funzione $G(z,\cdot) (1-\lambda_c)$ è compresa tra $0$ e $c$. Quindi per continuità, anche:
\begin{gather*}
 0\leq p(G(z,\cdot) (1-\lambda_c))=G(z,p)\leq c
\end{gather*}
data l'arbitrarietà di $c>0$, abbiamo che se $A_c$ è compatto per ogni $c$, allora fissato $z$, $G(z,p)=0$ per ogni $p\in \Gamma$, quindi $\Xi=\emptyset$, cioè la varietà iperbolica è regolare.\\
Supponiamo ora che esiste un insieme $A_c$ non compatto, allora esiste una successione $p_n\in R$, $p_n\to \infty$ tale che
\begin{gather*}
 G(z,p_n)\geq c \ \ \ \forall n
\end{gather*}
Consideriamo un qualunque carattere $p$ su $R$ descritto da un ultrafiltro su $z_n$. Dall'ultima equazione sappiamo che
\begin{gather*}
 G(z,p)\geq c
\end{gather*}
da cui $p\in \Xi\neq \emptyset$.
\end{proof}

\end{prop}
\begin{oss}
 Per le proprietà della funzione $G$, il fatto che per ogni $c>0$ $A_c$ sia compatto, è equivalente a chiedere che per ogni $c$ l'insieme $\partial A_c$ sia compatto.
\end{oss}
Per le varietà iperboliche irregolari, definiamo l'insieme delle successioni irregolari, cioé
\begin{deph}
 Data una varietà iperbolica $R$, definiamo:
\begin{gather*}
 \Sigma(R)=\{\{z_n\}\subset R \ t.c \ z_n \to \infty \ \ e \ \liminf_n G(z_0,z_n)>0\}
\end{gather*}
\end{deph}
Grazie a un ragionamento simile a quello riportato nella dimostrazione di \ref{prop_greenref}, si ottiene che questa condizione è indipendente dalla scelta di $z_0$, quindi $\Sigma(R)$ è ben definito. Inoltre dalla caratterizzazione appena dimostrata, è immediato verificare che $R$ è irregolare se e solo se $\Sigma(R)\neq \emptyset$.\\
Possiamo riformulare la proposizione \ref{prop_Eref} con questa nuova terminologia. In particolare:
\begin{prop}\label{prop_succirr}
Sia $R$ una varietà iperbolica irregolare, e sia $E$ la funzione descritta in \ref{prop_Eref}. Allora per ogni successione $\{z_n\}\in \Sigma(R)$:
\begin{gather*}
 \lim_{n\to \infty} E(z_n)=\infty
\end{gather*}
\begin{proof}
 Questa osservazione segue dalla semicontinuità di $E$ su $\Xi$.\\
Per prima cosa osserviamo che la proprietà appena enunciata è equivalente a:
\begin{gather}\label{eq_eq}
 \lim_{n\to \infty} \inf_{z\in V(z_0,a,n)} E(z) =\infty
\end{gather}
per ogni $a>0$, $a$ tale che $G(z_0,\cdot)\geq a$ non sia un insieme compatto in $R$, dove $V(z_0,a,n)\equiv \{z\in R \ t.c.\ G(z_0,z)> a\}\setminus K_n$ e $K_n$ è un'esaustione di $R$ \footnote{dalle considerazioni precedenti, sappiamo che questa proprietà è indipendente dalla scelta di $z_0\in R$}. Infatti, se questo è vero e consideriamo una successione $\{z_m\}\in \Sigma(R)$, dato che $$\liminf_m G(z_0,z_m) =\lambda>0$$ $z_m$ appartiene definitivamente a tutti gli insiemi $V(z_0,\lambda/2,n)$, e per $m$ che tende a infinito, esiste una successsione $n_m$ che tende a infinito tale che $z_n\in V(z_0,\lambda/2,n_m)$.\\
Per dimostrare l'altra implicazione, supponiamo che esista $a>0$ tale che:
\begin{gather*}
 \lim_{n\to \infty} \inf_{z\in V(z_0,a,n)} E(z) =M<\infty
\end{gather*}
questo significa che per ogni $n$ esiste un $z_n\in V(z_0,a,n)$ tale che $E(z_n)\leq M+1$. Per costruzione, la successione $z_n$ tende a infinito, e $G(z_0,z_n)\geq a$, quindi $\{z_n\}\in \Sigma(R)$. Questo completa la dimostrazione dell'equivalenza tra le due proprietà.\\
Dato che:
\begin{gather*}
 \bigcap_n V(z_0,a,n)=\Xi_a\equiv \{p\in \Gamma \ t.c. \ G(z_0,p)>a\}\subset \Xi
\end{gather*}
allora l'equazione \ref{eq_eq} è vera perché $E$ è semicontinua inferiormente sull'insieme $R\cup \Xi$ e $E(\Xi)=\infty$.
\end{proof}
\end{prop}
\subsection{Potenziali di Evans su varietà paraboliche}
Osserviamo che tutti gli spazi $R^n$ con $n\geq 3$ dotati della metrica euclidea standard sono iperbolici regolari, e gli stessi spazi privati di un punto qualsiasi sono iperbolici irregolari. Consideriamo ad esempio l'insieme $\R^3\setminus \{(0,0,1)\}$. Sappiamo che la funzione di Green su $\R^3$ è proporzionale all'inverso della distanza:
\begin{gather*}
 G(x,y)=C(3)\frac{1}{d(x,y)}
\end{gather*}
Si vede facilmente che la stessa funzione ristretta all'insieme $\R^3\setminus \{(0,0,1)\}$ è ancora un nucleo di Green. Fissiamo il punto $x=(0,0,0)$ e consideriamo la funzione $G(0,y)$. L'insieme $\{y\in \R^3\setminus \{(0,0,1) \ t.c. \ G(0,y) \geq C(3)/2\}\}$ NON è un insieme compatto in $\R^3\setminus \{(0,0,1)\}$, quindi questa varietà è iperbolica irregolare.\\
Altri esempi di varietà iperboliche irregolari possono essere costruiti togliendo un insieme compatto da una varietà riemanniana parabolica. Illustriamo questo risultato in due proposizioni.
\begin{prop}
 Sia $R$ una varietà Riemanniana e $K$ un compatto che sia chiusura della sua parte interna in $R$. Allora la varietà $R\setminus K$ dotata della metrica di sottospazio è una varietà iperbolica.
\begin{proof}
Consideriamo un compatto $C$ in $R'\equiv R\setminus K$ con bordo regolare, e dimostriamo che la capacità di questo compatto è necessariamente diversa da $0$. A questo scopo consideriamo un altro insieme compatto con bordo regolare $K'$ tale che:
\begin{gather*} 
 K\Subset K'^{\circ}\subset K' \ \ \ K'\cap C=\emptyset
\end{gather*}
Per ora trattiamo il caso di $K$ compatto con bordo regolare. Sia per definizione $w$ la funzione armonica in $K'\setminus K$ tale che:
\begin{gather*}
 w|_{\partial K}=0 \ \ \ w|_{\partial K'}=1
\end{gather*}
Consideriamo un'esaustione regolare $V_n$ di $R'$ tale che $C'\Subset V_1$, e sia $v_n$ il potenziale armonico della coppia $(C,V_n)$. Sappiamo che la capacità di $C$ è per definizione l'integrale di Dirichlet del limite di $v_n$, ed è nulla se e solo se $v=1$ costantemente.\\
Applicando il principio del massimo sull'insieme $V_n\cap K'$, otteniamo che per ogni $n$, $v_n\leq w$. Infatti sul bordo di $K'$, $w=1$ mentre $v_n\leq 1$, e sul bordo di $V_n$, $v_n=0$ mentre $w\geq 0$. Questo garantisce che
\begin{gather*}
 v\leq w \ \ \ su \ \ K'\setminus K
\end{gather*}
 quindi $v$ non può essere costante uguale a $1$.\\
Se $K$ non ha bordo liscio, è sufficiente ripetere la costruzione di $w$ con un insieme $K''\Subset K$ dal bordo liscio. Tutte le considerazioni fatte si applicano a questo caso senza complicazioni.
\end{proof}
\end{prop}
\begin{prop}\label{prop_last}
 Sia $R$ una varietà parabolica e sia $K$ un insieme compatto che sia chiusura della sua parte interna e abbia bordo regolare. Allora la varietà $R'\equiv R\setminus K$ è iperbolica irregolare. Inoltre tutte le successioni $\{z_n\}\subset R'\subset R$ che tendono a infinito in $R$ sono irregolari, appartengono all'insieme $\Sigma(R')$, e dato $z_0\in R'$, la funzione $G(z_0,\cdot)$ può essere estesa per continuità anche a $\partial K$ e su questo insieme è nulla.
\begin{proof}
 Dalla proposizione precedente sappiamo che $R'$ è una varietà iperbolica, quindi ammette nucleo di Green positivo. Sia $z_0\in R'$ qualsiasi, e consideriamo un compatto con bordo regolare $K'\Subset R$ \footnote{questo insieme NON è compatto in $R'$} tale che $K\subset K'^{\circ}$ e $z_0\in K'^{\circ}$. Sia inoltre $K_n$ un'esaustione regolare di $R$ con $K'\subset K_1$.\\
Osserviamo che il nucleo di Green $G(z_0,\cdot)$ su $R'$ è una funzione armonica sull'insieme $R\setminus K'$, e su $\partial K'$ assume minimo strettamente positivo $\lambda$ \footnote{$\lambda>0$ perchè il principio del massimo assicura che $G(z_0,\cdot)$ non può assumere il suo minimo in un punto interno all'insieme di definizione}. Allora grazie al principio del massimo possiamo confrontare $G(z_0,\cdot)$ con $v_n$, i potenziali armonici della coppia $(K',K_n)$. Per ogni $n$ si ha infatti che sull'insieme $R\setminus K'$:
\begin{gather*}
 G(z_0,\cdot)\geq \lambda v_n(\cdot)
\end{gather*}
passando al limite e ricordando che $R$ è parabolica (quindi $\lim_n v_n=1$) otteniamo che su $R\setminus K'$:
\begin{gather*}
 G(z_0,\cdot)\geq \lambda
\end{gather*}
Questo dimostra anche che qualunque successione tendente a infinito in $R$ è in $\Sigma (R')$.\\
Per dimostrare l'ultimo punto, consideriamo il massimo $\Lambda$ della funzione $G(\cdot,z_0)$ sul bordo di un dominio compatto \footnote{rispetto alla topologia di $R$} con bordo liscio $C$ tale che $K\subset C$ e $z_0\not\in C$. Sia $v$ il potenziale di capacità della coppia $(K,C)$. Grazie al principio del massimo, $G(z_0,\cdot)\leq \Lambda (1-v(\cdot))$, infatti $G(z_0,\cdot)$ è costruita come limite crescente di nuclei di Green $G_n$ su domini compatti, ed è facile verificare che
\begin{gather*}
 G_n(z_0,\cdot)\leq \Lambda (1-v(\cdot))
\end{gather*}
passando al limite su $n$ si ottiene la disuguaglianza desiderata. Questo dimostra che $G(z_0,\cdot)$ tende a zero se l'argomento tende a un punto qualsiasi di $\partial K$.\\
Osserviamo anche che grazie a questa proprietà, tutte le successioni che tendono a $\partial K$ \footnote{con questo si intendono tutte le successioni di elementi di $R'$ tali che per ogni intorno $U$ di $\partial K$, intorno rispetto alla topologia di $R$, esiste $N$ tale che $x_n\in U$ per ogni $n\geq N$. Ad esempio tutte le successioni che convergono rispetto alla topologia di $R$ a un punto di $\partial K$} non appartengono a $\Sigma(R')$.
\end{proof}
\end{prop}
Per dimostrare l'esistenza di potenziali di Evans relativi a un qualsiasi dominio compatto con bordo liscio $K$ in $R$ varietà parabolica, sfruttiamo quest'ultima proposizione, l'esistenza delle funzioni $E$ descritte in \ref{prop_Eref} e la proposizione \ref{prop_succirr}.
\begin{teo}\label{teo_evans}
 Data una varietà parabolica $R$ e un dominio compatto con bordo liscio $K$, esiste un potenziale di Evans su $R$ rispetto a $K$, cioè una funzione $E:R\setminus K^{\circ} \to[0,\infty)$ tale che:
\begin{enumerate}
\item $E\in H(R\setminus K)$
\item $E(z)=0$ se $z\in \partial K$
\item $E(z_n)\to \infty$ per ogni successione $z_n$ che tende a infinito in $R$.
\item $D_R(E\curlywedge c)\leq c$ per ogni $c>0$
\end{enumerate}
\begin{proof}
 La dimostrazione è semplicemente una raccolta di risultati precedentemente ottenuti.\\
Grazie a quanto appena dimostrato, $R\setminus K$ è una varietà iperbolica irregolare, quindi considerando la funzione $E$ costruita nella proposizione \ref{prop_Eref} e grazie alla proposizione \ref{prop_succirr}, otteniamo che $E$ soddisfa (1), (3) e (4).\\
Il punto (2) si ricava considerando che $E$ è una combinazione convessa di nuclei di Green sulla varietà $R\setminus K$, che grazie alla proposizione \ref{prop_last} si annullano sull'insieme $\partial K$.
\end{proof}
\end{teo}
\clearpage

\appendix
\chapter{Glossario}\label{chap_glossario}

\begin{tabular}{cl}
$A^C$ & complementare dell'insieme $A$\\
$\overline A$ & chiusura topologica dell'insieme $A$\\
d$\lambda$, d$\lambda^m$ & misura di Lebesgue su $\R$ o su $\R^m$\\
$(R,g)$ & varietà Riemanniana $R$ con tensore metrico $g$\\
$\Omega$ & dominio (insieme aperto connesso) in $R$\\
$\sqrt{\abs g}$ & radice quadrata del determinante di $g$\\
$dV$ &$dV=\sqrt{\abs g}dx^1\cdots dx^n$ elemento di volume su $(R,g)$\\
$L^2(\Omega)$ & spazio delle funzioni a quadrato integrabile su $\Omega$\\
$\mathcal L ^2 (\Omega)$ & spazio delle $1$-forme a quadrato integrabile su $\Omega$\\
& (vedi sezione \ref{sec_L2})\\
$\norm{f}_{\infty}$ & norma del sup della funzione $f$\\
$D_\Omega(f)$ & integrale di Dirichlet della funzione $f$ sull'insieme $\Omega$.\\ 
$\norm{f}_R$ & norma di $f$ nell'algebra di Royden $\Ro$\\
&(vedi teorema \ref{teo_g1})\\
$\Ro$ & Algebra di Royden su $R$ (vedi definizione \ref{deph_Ro})\\
$\Roo$ & Insieme delle funzioni a supporto compatto in $\Ro$\\
$\Rod$ & Completamento di $\Roo$ nella topologia $BD$\\
&(vedi paragrafo \ref{subsec_ideals})\\
$R^*$ & Compattificazione di Royden di $R$ (vedi definizione \ref{deph_R^*})\\
$\Gamma$ & $\Gamma=R^*\setminus R$ è il bordo di $R^*$ (vedi definizione \ref{deph_Gamma})\\
$\Delta$ & bordo armonico di $R$ (vedi definizione \ref{deph_ab})\\
$\Xi$ & bordo irregolare di $R$ (vedi definizione \ref{deph_Xi})\\
$f\ast h$& convoluzione della funzione $f$ con $h$ (vedi sezione \ref{sec_conv})\\
$\int f\ast du$ & vedi definizione \ref{deph_hodge}\\
$H(\Omega)$& spazio delle funzioni armoniche su $\Omega$\\
$HP(\Omega)$& spazio delle funzioni armoniche positive su $\Omega$\\
$HD(\Omega)$& spazio delle funzioni armoniche con $D_\Omega(u)<\infty$\\
%
$C(\Omega)=C(\Omega,\R)$& spazio delle funzioni continue a valori reali su $\Omega$\\
$C(\overline\Omega)=C(\overline\Omega,\R)$& spazio delle funzioni continue definite in un intorno di $\overline\Omega$\\
$\text{Cap}(K,\Omega)$ & capacità della coppia $(K,\Omega)$ (vedi definizione \ref{deph_cap})\\
$\text{Cap}(K)$ & capacità di $K$ (vedi definizione \ref{deph_cap2})
\end{tabular}

\section*{Ringraziamenti... e non solo}
\ 

I miei ringraziamenti vanno innanzitutto al prof. Alberto Giulio Setti e al prof. Stefano Pigola. In particolare ringrazio Alberto soprattutto per i consigli che mi ha dato e per l'estrema gentilezza, pazienza e diponibilità dimostrate, ben oltre quanto richiesto e non solo dal punto di vista accademico. Ringrazio anche il prof. Wolfhard Hansen (University of Bielefeld) per la disponibilità e per i suggerimenti dati riguardo alla teoria del potenziale.
\begin{flushright}
\textit{Vuoi studiare a Princeton? - ...}
\end{flushright}

\thispagestyle{empty} 

Ovviamente ringrazio tutti i miei amici dell'università, la compagnia del $4^\circ$ piano e della Norvegia, nanetti compresi. Non vedo l'ora del prossimo ciclo!
\begin{flushright}
\textit{C'è uno spiffero...}
\end{flushright}

Non posso dimenticare di ringraziare Madesimo, che mi ha fornito un posto idilliaco dove pensare e scrivere questa tesi, e non posso dimenticare tutti gli amici che mi hanno accompagnato in questa impresa.
\begin{flushright}
\textit{Ai larici a mezzogiorno?}
\end{flushright}

Un ringraziamento importante va a Internet e a tutti quelli che non credono nel diritto d'autore (senza i quali non avrei potuto studiare tutto quello che ho studiato), da \htmladdnormallink{Wikipedia}{http://www.wikipedia.org} a \htmladdnormallink{Ubuntu - Linux}{http://www.ubuntulinux.org} passando per l'indispensabile \htmladdnormallink{Gigapedia}{http://www.gigapedia.org}.

\begin{flushright}
\textit{A journey of a thousand sites begins with a single click.}
\end{flushright}


Ringrazio l'Unione Europea, per quello che è ma soprattutto per quello che potrebbe essere. Ringrazio in particolare tutti quelli che ci hanno lavorato seriamente, perché hanno aperto una strada che ora noi abbiamo la responsabilità di portare avanti e soprattutto perché ci hanno lasciato qualcosa che sprecare sarebbe una follia.
\begin{flushright}
\textit{European Union - United in diversity}
\end{flushright}

Un pensiero speciale va a Lucy, e a suo padre scomparso di recente.
\begin{flushright}
\textit{Cosa, cosa?}
\end{flushright}

Ringrazio Homer, Marge, Bart, Lisa e Maggie, che mi hanno divertito, intrattenuto, e dato da pensare.
\begin{flushright}
\textit{Farfalla vendetta!}\\
\end{flushright}

Ringrazio Kathryn, Seven, Benjamin, Jean-Luc e rispettive compagnie perché mi hanno fatto sognare.
\begin{flushright}
\textit{I look forward to it. Or should I say backward? - Don't get started!}
\end{flushright}

\ \newline
Un ringraziamento speciale a Oriana Fallaci, e alla maga.\\
\textit{\ [...] Ad esempio perché si trovasse qui, perché avesse scelto un mestiere che non si addiceva al suo carattere e alla sua struttura mentale cioè il mestiere di soldato, perché con quel mestiere avesse tradito la matematica. Quanto gli mancava la matematica, quanto la rimpiangeva! Massaggia le meningi come un allenatore massaggia i muscoli di un atleta, la matematica. Le irrora di pensiero puro, le lava dai sentimenti che corrompono l'intelligenza, le porta in serre dove crescono fiori stupendi. I fiori di un'astrazione composta di concretezza, di una fantasia composta di realtà [...] No, non è vero che è una scienza rigida, la matematica, una dottrina severa. \`E un'arte seducente, estrosa, una maga che può compiere mille incantesimi e mille prodigi. Può mettere ordine nel disordine, dare un senso alle cose prive di senso, rispondere ad ogni interrogativo. Può addirittuta fornire ciò che in sostanza cerchi: la formula della Vita. Doveva tornarci, ricominciare da capo con l'umiltà d'uno scolaro che nelle vacanze ha dimenticato la tavola pitagorica. Due per due fa quattro, quattro per quattro fa sedici, sedici per sedici fa duecentocinquantasei, e la derivata di una costante è uguale a zero. La derivata di una variabile è uguale a uno, la derivata di una potenza di una variabile... Non se ne ricordava? Sì che se ne ricordava! La derivata di una potenza di una variabile è uguale all'esponente della potenza moltiplicata per la variabile con lo stesso esponente diminuito di uno. E la derivata di una divisione? \`E uguale alla derivata del dividendo moltiplicato per il divisore meno la derivata del divisore moltiplicata per il dividendo, il tutto diviso il dividendo moltiplicato per sé stesso. Semplice! Bè, naturalmente trovare la formula della Vita non sarebbe stato così semplice. Trovare una formula significa risolvere un problema, e per risolvere un problema bisogna enunciarlo, e per enunciarlo bisogna partire da un presupposto... Ah perché aveva tradito la maga? Che cosa lo aveva indotto a tradirla? [...]}

\begin{flushright}
\thispagestyle{empty} 
 \textit{Oriana Fallaci: Insciallah. Atto primo, Capitolo primo}
\end{flushright}



\end{document}